\documentstyle[12pt,epsfig,amsfonts]{article}

\newtheorem{proposition}{Proposition}[section]
\newtheorem{definition}{Definition}[section]
\newtheorem{corollary}{Corollary}[section]
\newtheorem{coro}{Corollary}
\newtheorem{theorem}{Theorem}
\newtheorem{lemma}{Lemma}[section]

\newtheorem{claim}{Claim}[section]
\newtheorem{sublemma}{Sublemma}[section]
\begin{document}

\thispagestyle{empty}

\vspace*{0.7in}

\centerline{\Huge \bf Analysis of}

\bigskip

\centerline{\Huge \bf a Class of Strange Attractors}

\vspace{0.4in}

\centerline{\large Qiudong Wang\footnote{Department of Mathematics, UCLA, Los An
geles, CA 90095, email dwang@math.ucla.edu. This research is partially supported
 by NSF grant \#9970673 and an NSF Postdoctoral Research Fellowship} \ \ and \ \
 Lai-Sang Young\footnote{Courant Institute of Mathematical Sciences, 251 Mercer
St., New York, NY 10012, email lsy@cims.nyu.edu; and Department of Mathematics,
UCLA, Los Angeles, CA 90095, email lsy@math.ucla.edu. This research is partially
 supported by a grant from NSF and a Guggenheim Fellowship}}

\vspace*{0.5in}

\centerline{\bf \large TABLE OF CONTENT}
 
\vskip .3in
 
\ 1 \ \  Statements of Results

\ 2 \ \ Preliminaries
 
\medskip

\noindent PART I \ \ Controlling Nonhyperbolicity
 
\medskip
 
\ 3 \ \ The Critical Set
 
\ 4 \ \ Replication of Orbit Segments
 
\ 5 \ \ Pushing the Induction Forward
 
\ 6 \ \ Measure of Selected Parameters
 
\medskip

\noindent PART II \ \ Geometric and Statistical Properties 
 
\medskip
 
\ 7 \ \ Hyperbolic Behavior
 
\ 8 \ \ SRB Measures and Their Statistical Properties

\ 9 \ \ Global Geometry
 
\noindent \ \ \ \ 10 \ \  Symbolic Dynamics and Topological Entropy
 
\medskip
 
\noindent APPENDIX 

\medskip

\ A \ \ H\'enon Attractors and Homoclinic Bifurcations
 
\ B \ \ Computational Proofs

\pagestyle{myheadings}
\setcounter{page}{0}

\newpage

\newpage

The notion of a ``strange attractor" has been common knowledge
in dynamical systems for more than two decades and has captured the
attention of scientists in other disciplines. Rigorous mathematical
analysis, however, has not kept pace with these developments.
Among the examples that have been studied are
the Lorenz attractors (\cite{G}, \cite{Ro}, \cite{Ry}, \cite{W2}) and 
the H\'enon maps
(\cite{BC2}, \cite{BY1}, \cite{BY2}). In both of these
examples, the attractors are closely related to certain  
$1$-dimensional maps.
 
This paper is a general study of attractors that are derived,
in some fashion, from
$1$-dimensional maps. The unstable manifolds of the resulting attractors
are $1$-dimensional; the attractors themselves live in dimensions
$\geq 2$. We will limit ourselves to {\it discrete time} and 
{\it smooth} maps;
thus our study includes as a special case the H\'enon attractors
but not the Lorenz flows. Our investigation proceeds in several different
directions, ranging from local analysis to global geometry to statistical
properties. 
Although the attractors in our class have a number of features
in common with Axiom A
attractors and with piecewise monotonic maps in 1-dimension, the reader 
will find that 
these two theories together are far from adequate for
handling the new complexities that arise.

We now give a general description of the broad category of attractors
that are the objects of our study. Our results apply to a subset of
this class. Let $f:N \to N$ be a self-map of a circle or an interval,
and let $M=N \times D_n$ where $D_n$ is an $n$-dimensional disk.
Identifying $N$ with $N \times \{0\} \subset M$, we perturb $f$ into an
embedding of $N$ into $M$, and then extend it to an embedding $T$ of $M$
into itself. The attractor of interest to us is
given by $\Omega =\cap_{i \geq 0} T^i M$.
If $f(z)=z^2$ and $M$ is the solid torus $S^1 \times D_2$,
then $\Omega$ is the well known solenoid.
We propose to replace the map $f(z)=z^2$ in the standard solenoid
by an arbitrary
smooth  map. For lack of a better name, let us call these
{\bf (generalized) solenoidal attractors}.
 
In addition to the two examples  we have
encountered, namely the Axiom A solenoid \cite{Sm} (with $f(z)=z^2$) and the
H\'enon maps \cite{H} (with $f_a(x) = 1-ax^2, \ x \in [-1,1]$), other known examples
of solenoidal attractors
include dissipative twist maps \cite{Bi}, the most
standard of which can be realized as a suitable perturbation of
$f(x)=x+\frac{K}{2 \pi} \sin(2 \pi x)$, $x \in {\mathbb R}/ {\mathbb Z}$, and
certain periodically forced nonlinear oscillators (\cite{Lev}; see also 
\cite{GH}).

Since the picture is well understood when $f$ is uniformly expanding,
we are primarily interested in the case where $f$ has critical points.
When the critical orbits of $f$ tend to attractive cycles,
the dynamics of the $T$ (assuming the perturbations are small)
is also quite simple: the stable periodic orbits persist, and the complement
of their basins consists of ``horseshoes" and their stable manifolds.
We focus, therefore, on $1$-dimensional maps $f$ that are ``chaotic"
with no stable equilibria.
 We mention two important differences between
this situation and $f(z)=z^2$.
First, $T|\Omega$ in general cannot be realized as the inverse limit of $f$;
it is more complicated. Second, while $f(z)=z^2$ gives rise
to essentially one attractor -- in the sense that two different
perturbations $T$ and $T'$ can be conjugated by a homeomorphism
$C^0$-near the identity -- an arbitrary $f$  can (and does) give rise
to infinitely many ``different" attractors.
 
We now give a more precise description
of the setting to which our results apply.

\vskip .2in
 
\noindent {\large \bf Setting of this paper}

\bigskip

Our results are for attractors that arise from perturbations
of circle or interval maps. For definiteness, we assume $N=S^1$, 
$M$ is an annulus,
and impose the following conditions on $T$ to make the dynamics more tractable.
To ensure that $T$ is predominantly hyperbolic, it is necessary 
to start with a $1$-dimensional map with sufficiently strong
expanding properties. We assume $f$ satisfies 
the Misiurewicz condition, i.e. $f$ is an arbitrary piecewise monotonic map 
with the property that its forward critical orbits stay away from its critical points. We consider a 2-parameter 
family $\{T_{a,b}\}$ through $f$, using the parameter $a$ to control 
movements along the circle and $b$ to ``unfold" the $1$-dimensional
maps in the second direction. Mild transversality conditions are assumed
on the 2-parameter family, and the maps $T_{a,b}$ are required to be
diffeomorphisms for $b>0$.

This paper concerns the parameter range where $b$ is small, that is, 
where $T_{a,b}$ is {\bf strongly dissipative}. 
Detailed studies are  presented for maps corresponding to a positive measure 
set of parameters
in this range.

We mention some small generalizations. 
When $N=S^1$ and $|\deg(f)|>1$, at least three dimensions are needed
for $T_{a,b}$ to be globally injective. 
An extension of our techniques gives essentially the same results; 
details will appear elsewhere. Another possible extension, 
which we will not discuss,
is to replace $S^1$ with branched-1-manifolds (see \cite{W1}).

\vskip .2in

\noindent {\large \bf Overview of Results}

\bigskip
 
\noindent {\bf Selection of parameters and the critical set.} \ 
Given $\{T_{a,b}\}$, the goal of this step is to select a positive measure set 
of ``good" parameters corresponding to maps that one can control
in certain ways. Our criteria for parameter selection are 
similar to those 
of Benedicks and Carleson \cite{BC2}, which in turn draws its 
inspiration from previous work on $1$-dimensional maps, 
from \cite{BC1} and \cite{CE} in particular.
In this approach, one inductively identifies and controls an 
object called the {\it critical set}, which one hopes will play the role of 
critical points in $1$-dimension. Our inductive process gives 
new information not available in \cite{BC2}. 
We obtain a systematic description of the structure
of the map near the critical set, 
which we realize in a Cantor construction 
as the intersection of  a nested
sequence of sets each one of which is
a union of rectangles with known geometric properties. This detailed knowledge
of the critical set is crucial in many of our results.
Another departure from \cite{BC2} is that our analysis is based on   
simple geometric conditions, whereas the equations of the  
H\'enon maps are used in many computations there. We also give a more complete 
treatment of parameter-space issues than in previously published works.

\newpage

{\em The results below hold for maps corresponding to the parameters selected.} 
\bigskip

\noindent {\bf Hyperbolic behavior.} \  We prove that compact invariant 
sets disjoint from the critical set are uniformly hyperbolic, with
hyperbolicity getting weaker as one approaches the critical set. 
Our analysis also gives information on the nonuniform
character of hyperbolicity 
in the basin.
\bigskip

\noindent {\bf Statistical properties.} \
We construct Sinai-Ruelle-Bowen (SRB) measures on our attractors, bound the 
number
of ergodic SRB measures by the number of critical points of the generating $1
$-dimensional map, and show that with respect to Lebesgue measure, almost every
point in  the basin
is generic with respect to an ergodic SRB measure.
Appealing to the abstract results in \cite{Y3} and \cite{Y4}, we prove that
the attractors in our class enjoy a Central Limit Theorem and have exponential
decay of correlations on their mixing components.  The corresponding results for
Axiom A attractors have been known since the 1970s (\cite{S2}, \cite{R1},
\cite{R2}).
For the H\'enon family near $a=2, \ b=0$, SRB measures and their statistical
properties were studied in \cite{BY1} and \cite{BY2}, and the basin property
in \cite{BV}.

\bigskip
\noindent {\bf Global geometry.} \ 
The approximate shape and complexity of the Axiom A solenoid 
is given by a small tubular neighborhood of a simple closed 
curve winding around the solid torus $2^k$ times.
In analogy with piecewise monotonic maps in $1$-dimension, 
we introduce the notion of
``monotone branches" and show that our attractors have 
arbitrarily fine neighborhoods that are unions of finitely many of these
branches. From the way these branches fit together one obtains a certain 
insight into the differences between 
one and two-dimensional maps.

\bigskip

\noindent {\bf Symbolic dynamics and topological entropy.} \ 
The geometric considerations above make it possible to code
unambiguously all orbits on the attractor,
representing the dynamics of the map by a shift operator acting
on symbol sequences generated by a finite alphabet. 
In the non-Axiom A case this shift is not of finite type. 
The coding we give is an honest reflection of the true locations visited
by an orbit relative to the  ``components" of the critical set.
Kneading sequences for critical orbits are well defined, and
monotone branches correspond to cylinder sets.
This symbolic representation is nearly one-to-one, allowing us to deduce 
from it
the existence of equilibrium states
and various formulas for computing topological entropy.
For results in this direction for Axiom A attractors,
see \cite{Bo} and the references therein.

\vskip .2in

It is our hope that by formulating simple, checkable conditions as we have 
done in Sect. \ref{s1.1},
one can determine readily if the results of this paper
apply to a given situation. We illustrate this for generic {\bf homoclinic 
bifurcations} of 2-dimensional diffeomorphisms,
recovering the result in 
\cite{MV} (which extends \cite{BC2} to this setting) and 
obtaining immediately for the attractors in question the dynamical picture 
above.

\newpage

The results of this paper open the door to a host of questions for
the class of attractors being studied. For example, with the
information available, extensions of
the theory of equilibrium states to the present setting may be possible
(see e.g. \cite{Bo}, \cite{R2}, \cite{K}).
Our results on symbolic dynamics lead naturally to questions on the
zeta function (see e.g. \cite{Bal}, \cite{PP}, \cite{R3}).
With kneading sequences for critical orbits being well defined, it is
reasonable to consider the possibility of a kneading theory (see
\cite{C}, \cite{MT}). In a different direction, the notion of
monotone branches leads to questions about global topological structures and 
prime ends  (e.g. \cite{Bar}).
 
\medskip

We mention some other related works, omitting specific references to 
$1$-dimensional dynamics (see the reference in \cite{dMvS}).  
For results on piecewise uniformly
hyperbolic attractors, see e.g. \cite{CL}, \cite{I1}, \cite{I2}, 
\cite{M2} and \cite{Y1}. For the 
statistical properties of billiards, see e.g. \cite{S1}, \cite{BSC1}, 
\cite{BSC2} and \cite{Y3}. 
Closer to the setting of this paper 
are \cite{DRV} and \cite{MV}, in which Viana {\it et. al.} extended the 
analysis
in \cite{BC2} to ``H\'enon-like" maps and found applications for these 
extensions. See also \cite{V}.
Jakobson and Newhouse have announced that they have reproduced,
using different methods, the results in \cite{BC2} and \cite{BY1}.
We have been told that Luzzatto has done work in this direction, and
that Palis and Yoccoz have results for
certain non-attracting sets that arise in homoclinic bifurcations.

\bigskip

This paper is by and large self-contained --- with the exception of Section
\ref{s6}, where two results from 1-dimensional maps are quoted without
proof, and Section \ref{s10}, where previous work of the second-named 
author is used. 
Proofs that are computational in nature 
have been put in the Appendix so that they will not obstruct the main flow of 
ideas. In a paper as long as this one, it might be useful to indicate
the logical connections among the various sections. After Section \ref{s1},
we recommend at least looking through Section \ref{s2}, in which we 
introduce much of
the basic vocabulary for subsequent sections. The other sections are connected 
as indicated. (For example, the technical 
content of Section \ref{s6} is not needed for reading Sections 
\ref{s7}-\ref{s9}.)

\begin{picture}(13, 7)
\put(0,0){
\psfig{figure=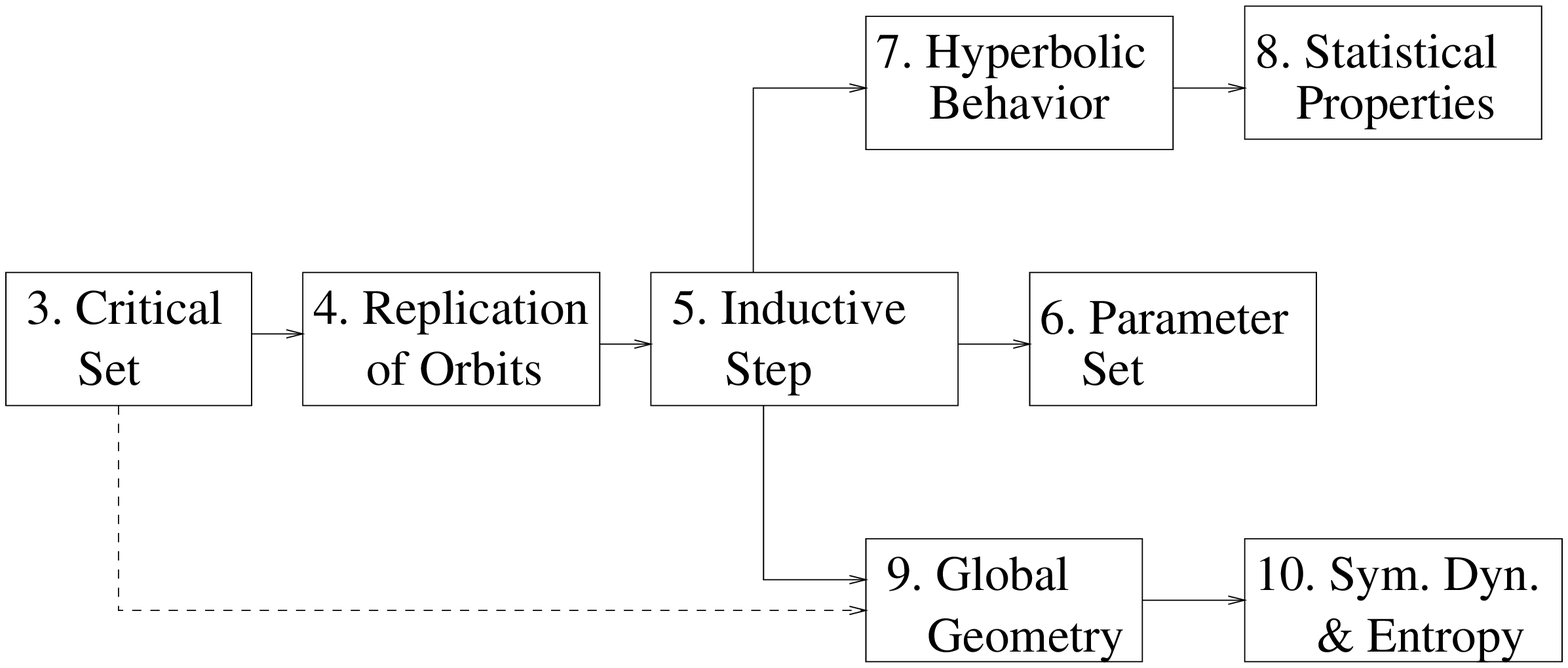,height = 6cm, width = 12cm}
}
\end{picture}
\newpage

%%%%%%%%%%%%%%%%%%%%Section 1
\section{Statements of Results}
\label{s1}
%%%%%%%%%%Section 1.1

\subsection{Setting}
\label{s1.1}

For definiteness, Theorems 1--7 are stated in the context of attractors 
that arise from perturbations of {\it circle} maps. For the interval case, 
see Sect. \ref{s1.5}.

Let $A=S^1 \times [-1,1]$. We consider 2-parameter families of maps
$\{T_{a, b}\}$ where for each $(a,b)$, $T_{a,b}:A \to A$ is a self-map of
$A$ and 
$(x,y,a,b) \mapsto T_{a,b}(x,y)$ is $C^3$. The class of 
2-parameter families $\{T_{a, b}\}$ to which our results apply
are constructed via the following four steps. 

\bigskip
\noindent {\bf Step I. } \ Let $f:S^1 \to S^1$ satisfy the following
{\bf Misiurewicz conditions}, i.e. letting $C=\{x:f'(x)=0\}$,
we assume:

1. $f^{\prime \prime}(x) \not = 0$ for all $x \in C$;

2. $f$ has negative Schwarzian derivative on $S^1 \setminus C$;

3. there is no $x \in S^1$ with 
$f^n(x)=x$ and $|(f^n)'(x)| \leq 1$;

4. for all $x\in C, \ \inf_{n>0}d(f^nx, C)>0$.

\bigskip

Observe that for $p \in S^1$ with $\inf_{n \geq 0}d(f^np, C)>0$,
if $g$ is sufficiently near $f$ in the $C^2$ sense, then there is a
unique point $p(g)$ having the same symbolic dynamics with respect to $g$ as
$p$ does with respect to $f$.
If $\{f_a\}$ is a 1-parameter family through
$f$, then for those $a$ for which it makes sense,
we will call $p(a)=p(f_a)$ the {\it continuation} of $p$.
For $x \in C$, we let $x(a)$ denote the corresponding 
critical point of $f_a$.
 
\bigskip
\noindent {\bf Step II.} \ Let $f$ be as in Step I, and let $\{f_a\}, \ a \in [a_0, a_1]$, be a 1-parameter family of maps
from $S^1$ to $S^1$ with $f=f_{a^*}$ for some $a^* \in [a_0, a_1]$.
We require that $\{f_a\}$ satisfy the following
{\bf transversality 
condition}\footnote{This 
transversality condition is used in \cite{TTY}.}: For every $x \in C$, if $p=f(x)$, then
%%%%Formula 1.1
\begin{equation}
\frac{d}{da} f_a(x(a)) \not = \frac{d}{da}p(a) \ \ \ \ \ \ {\rm at} \ \ 
a=a^*.
\label{formula1.1}
\end{equation}
 
\bigskip

\noindent {\bf Step III.} \ Let $\{f_a\}$ be as in Step II.
Identifying $S^1$ with $S^1 \times \{0\} \subset A$,
we extend $\{f_a\}$ to a $2$-parameter family $\{f_{a,b}\}, \
a \in [a_0, a_1], \ b \in [0,b_1]$, where $f_{a,b}:S^1 \to A$ is such that
$f_{a,0}=f_a$ and $f_{a,b}$ is an {\it embedding} for $b>0$.

\bigskip

\noindent {\bf Step IV.} \ Finally, we extend $f_{a,b}$ to 
$T_{a,b}:
A \to A$ in such a way that $T_{a,0}(A) \subset S^1 \times \{0\}$ and
for $b>0$, $T_{a,b}$ maps
$A$ diffeomorphically onto its image. We further impose the following
{\bf non-degeneracy condition}\footnote{This condition is not assumed in 
\cite{MV} or \cite{BV}. Their regularity condition on $|\det (DT)|$ and bound
on the perturbation term, however, imply a condition which is similar
(though not equivalent) to (\ref{non-dege}) and which serves a similar purpose.}
on the map $T_{a^*,0}$:
\begin{equation}
\partial_y T_{a^*,0}(x,0) \neq 0 \ \ \ {\rm whenever} \ \ \ f_{a^*}'(x)=0.
\label{non-dege}
\end{equation}

\bigskip
 
This completes our construction of admissible families $\{T_{a,b}\}$.
We remark that the transversality and non-degeneracy conditions in Steps II
and IV are generic. Thinking in terms of normal neighborhoods, one
constructs easily for a given $f_{a,b}$ extensions of the type
in Step IV;  the signs of the $\partial y$-derivatives at the
critical points of $f_{a^*}$ are determined by the orientations of the
turns of $f_{a,b}$ at the corresponding points. 
Step III is feasible if and only
if the degree of $f$ is $0, 1$ or $-1$. If $|\deg(f)|>1$, an extra
dimension is needed; this will be treated
in a separate paper.
 
Observe that for $b>0$, $T_{a,b}$ has the general form
\begin{equation}
\label{themap}
T_{a, b}: \left( \begin{array}{c} x \\ y \end{array} \right)
 \ \mapsto \ \left( \begin{array}{c} F(x, y, a) + b \ u(x, y, a, b)
\\ b \ v(x, y, a, b) \end{array} \right)
\end{equation}
where $F(x,y,a)=T_{a,0}^1(x,y)$, the first component of $T_{a,0}(x,y)$,
and the $C^2$ norms of
$(x,y,a) \mapsto u(x, y, a, b)$ and $ v(x, y, a, b)$
are uniformly bounded for all $b \in (0,b_1]$.\footnote{This is a 
calculus exercise: Observe that
$bu$ extends to a $C^3$ function $g$ on $\{b \geq 0\}$ with $g|\{b=0\}=0$.
Writing $\partial^2=\frac{\partial ^2}{\partial z_1 \partial z_2}$ where 
$z_i=x,y$ or $a$, we then check that $\partial^2 u$ extends to a continuous 
function $h$ on $\{b \geq 0\}$
with $h=\frac{\partial}{\partial b} \partial^2 g$ on $\{b=0\}$.} 

In terms of differentiability assumptions, the following slightly more 
technical formulation corresponds exactly to what is used:

\smallskip

(i) $(x, y, a) \mapsto T_{a, b}(x, y)$ has uniformly bounded $C^3$-norms in $b$,

(ii) $T_{a, b}$ is of the form (\ref{themap}) with uniformly bounded $C^2$
norms for $u$ and $v$.

\bigskip
{\bf Notation.} \ Given $\{T_{a,b}\}$, constants that are determined
entirely by the family $\{T_{a,b}\}$ will be referred to as
{\bf system constants}. Except where declared otherwise, the letter {\bf $K$}
is reserved throughout this article for use as a {\bf generic system constant}, meaning a system constant that
is allowed to change from statement to statement (the other system
constants are fixed). We will use $K_1, K_2$ etc. where $K$ appears in
more than one role in the same statement. 

\medskip
Let $K$ be such that $T_{a,b}(A) \subset R_0:= S^1 \times [-Kb, Kb]$
for all $(a,b)$. It is convenient
for us to work with $R_0$ instead of $A$.
For $T=T_{a,b}$, let $R_n =T^n R_0$.
Then $\{R_n\}$ is a decreasing sequence of neighborhoods
of the {\bf attractor} $\Omega:= \cap_{n=0}^\infty R_n =
\cap_{n=0}^\infty T^n A$.

%%%%%%%%%%Section 1.2 
\subsection{Critical set and hyperbolic behavior}
\label{s1.2}

Our first theorem identifies, for each map $T$ corresponding to a selected set 
of parameters, a fractal set ${\cal C}$ chosen to play the role of the 
critical set in $1$-dimension. This set will be called the {\bf critical set} 
of $T$. Our parameter selection imposes strong hyperbolic properties on orbits 
starting from ${\cal C}$ in the hope that these properties
will be passed on to the rest of the system. The geometric structure near
${\cal C}$, which is described in some detail in Theorem 1, is crucial for
many of our later results.

For $z_0 \in R_0$, let $z_i=T^i z_0$. If $w_0$ is a tangent vector at $z_0$, let $w_i=DT^i(z_0)w_0$. 
A curve in $R_0$
is called a {\bf $C^2(b)$-curve} if the slopes of its tangent vectors
are ${\cal O} (b)$ and its curvature is everywhere ${\cal O}(b)$.
%%%%%%%%%%%%%%%%%%%%%%%%%Theorem 1
\begin{theorem} 
{\bf (Parameter selection and the critical set)}  \
Given $\{T_{a,b}\}$ as in Sect. 1.1, there is a positive measure set
$\Delta \subset [a_0, a_1]\times (0,b_1]$ such that (1) and (2) below
hold for $T=T_{a,b}$ for all $(a,b) \in \Delta$. The set $\Delta$
is located near $a=a^*$ and $b=0$; it has the property that
for all sufficiently small $b,
\ \Delta_b:=\{a:(a,b) \in \Delta\}$ has positive 1-dimensional
Lebesgue measure. The constants $\alpha, \delta, c>0$ and $0<\rho<1$
below are system constants, and $b<<\alpha, \delta, 
\rho, e^{-c}$ for all $(a,b) \in \Delta$.
\begin{itemize}
\item[(1)] {\bf Critical regions and critical set.} There is a 
Cantor set ${\cal C} \subset \Omega$ called the critical set given by
${\cal C}=\cap_{k=0}^\infty {\cal C}^{(k)}$ where the ${\cal C}^{(k)}$ are
a decreasing sequence of neighborhoods of ${\cal C}$ called 
critical regions.
 
{\bf Geometrically}:

(i) ${\cal C}^{(0)}=\{(x,y) \in R_0: d(x, C) <\delta\}$ where $C$ is the set
of critical points of $f$.

(ii) ${\cal C}^{(k)}$ has a finite number of components called $Q^{(k)}$ each 
one of which is diffeomorphic to a rectangle. The boundary of $Q^{(k)}$ 
is made up of two $C^2(b)$ segments of $\partial R_k$ connected by two vertical lines: the horizontal boundaries are  
$\approx \min{(2\delta, \rho^k)}$ in length, and the Hausdorff distance between them is 
${\cal O}(b^{\frac{k}{2}})$.
 
(iii) ${\cal C}^{(k)}$ is related to ${\cal C}^{(k-1)}$ as follows:
$Q^{(k-1)} \cap R_k$ has at most finitely many components,
each one of which
lies between two $C^2(b)$ subsegments of $\partial R_k$ 
that stretch across $Q^{(k-1)}$ as shown.
Each component of $Q^{(k-1)} \cap R_k$ contains exactly one component of 
${\cal C}^{(k)}$.

\begin{picture}(12, 4.5)
\put(1,0.5){
\psfig{figure=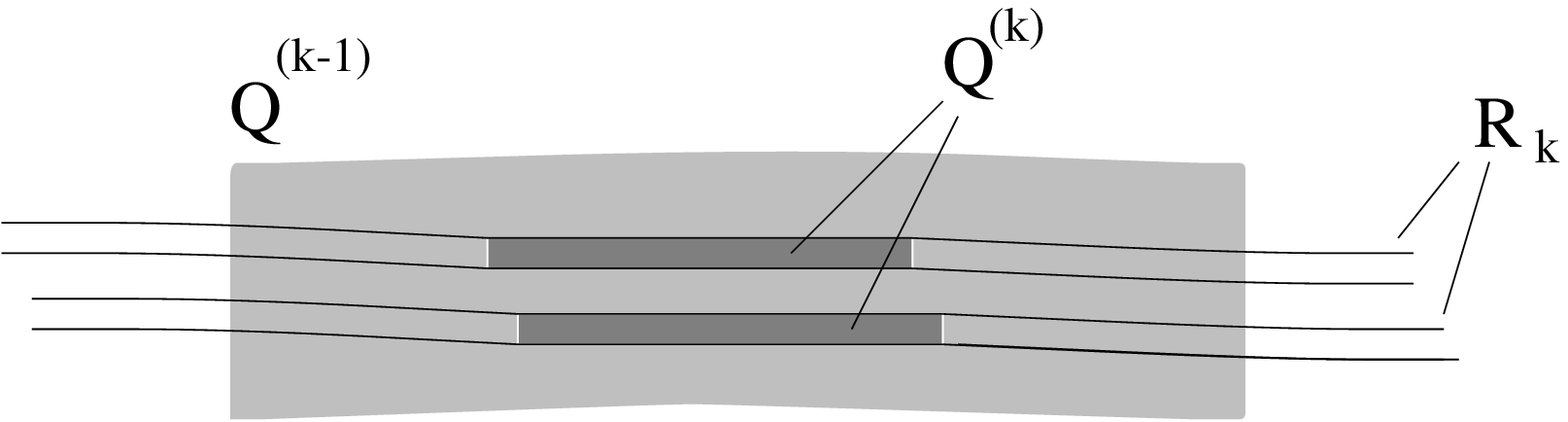,height = 3cm, width = 11cm}
}
\end{picture}

\medskip

\centerline{\rm Figure 1 \ \ Critical regions}

\medskip

{\bf Dynamically}: On each horizontal boundary $\gamma$ of $Q^{(k)}$  
there is a unique point $z$ located within ${\cal O}(b^{\frac{k}{4}})$ 
of the midpoint of $\gamma$ with the property that if $\tau$ is the unit 
tangent vector to
$\gamma$ at $z$, then $DT^n(z)\tau$ decreases in length
exponentially as $n$ tends to $\infty$.
 
\item[(2)] {\bf Properties of critical orbits.} For $z \in R_0$, 
let $d_{\cal C}(z)$ denote the following notion of 
``distance to the critical set":
If $z \not \in {\cal C}^{(0)}$,  
let $d_{\cal C}(z)=\delta$; if $z \in {\cal C}^{(0)} \setminus {\cal C}$,
let $k$ be the largest number with $z \in {\cal C}^{(k)}$,
and define $d_{\cal C}(z)$ to be the horizontal distance between $z$
and the midpoint of the component of ${\cal C}^{(k)}$ containing $z$.
Then  for all $z_0 \in {\cal C}$:
 
(i) \ $d_{\cal C}(z_j) \geq e^{-\alpha j}$ for all $j>0$;
 
(ii) $\|DT^j(z_0) {\tiny (\!\!\begin{array}{l} 0 \\ 1 \end{array}\!\!)}\| 
\geq K^{-1}e^{cj}$ for all $j>0$.
\end{itemize}
\label{theorem1}
\end{theorem}

\bigskip

{\bf Theorems \ref{theorem2}-\ref{theorem5} apply to 
$T=T_{a,b}$, $(a,b) \in \Delta$, where $\Delta$ is as in Theorem \ref{theorem1}.}

\bigskip

Our next theorem is about the abundance of hyperbolic behavior on the
attractor and in the basin. A compact $T$-invariant set $\Lambda$
is called {\bf uniformly hyperbolic} if there is a splitting of the
tangent bundle over $\Lambda$ into invariant subbundles $E^u \oplus E^s$
such that for some $C, \lambda>1$, we have, for all $n \geq 1$, $\|DT^nv\| \leq C\lambda^{-n}\|v\|$ for all $v \in E^s$ and $\|DT^{-n}v\| \leq C\lambda^{-n}\|v\|$
for all $v \in E^u$. 

%%%%%%%%%%%%%%%%%%%%%%%%%Theorem 2
\begin{theorem} 
{\bf (Hyperbolic behavior)} \
\begin{itemize}
\item[(1)] Let 
$$
\Omega_\varepsilon := \{z_0 \in \Omega: d_{\cal C}(z_n) \geq \varepsilon  \ 
\ \forall n \in {\mathbb Z} \}.
$$
\begin{itemize}
\item[(i)] For every $\varepsilon>0$, $\Omega_{\varepsilon}$ 
is uniformly hyperbolic. In fact, independent of $\varepsilon$, 
$\lambda$ in the 
definition of
hyperbolicity can be taken to be $\approx e^{\frac{c}{3}}$ where $c$ is as in
Theorem \ref{theorem1}. In particular, 
for evevry periodic point $z \in \Omega$ with $T^q z = z$, 
$\|DT^q | E^u(z) \| \geq K^{-1} e^{\frac{c}{3}q}$.

\item[(ii)] As $\varepsilon \to 0$, the hyperbolicity on $\Omega_{\varepsilon}$
deteriorates in the sense that 
$C \to \infty$ and the minimum angle between $E^u$ and $E^s$ tends to zero. 

\item[(iii)] $\Omega = \overline{\cup_{\varepsilon >0}
\Omega_{\varepsilon}}$ provided a surjective condition of the type
(*) below is assumed.
\end{itemize}

\item[(2)] Under the regularity conditions (**) below, we have
$$
\limsup_{n \to \infty} \frac{1}{n} \log \|DT^n(z_0)\| \geq \frac{c}{3}
$$
for Lebesgue-almost every $z_0 \in R_0$.
\end{itemize}
\label{theorem2}
\end{theorem}
The two technical conditions used in parts (1)(iii) and (2) of 
Theorem \ref{theorem2} are:

\bigskip

(*) \ Let $J_1, \cdots J_r$ be the intervals of monotonicity of $f$. Then for
each $i$, there

\ \ \ \ \ exists $j$ such that $f(J_j) \supset J_i$. 
 
\medskip
 
(**) \ There exist $K_1, K_2>0$ such that 
for all $z \in R_0$, 
$$ 
K^{-1}_1 b \leq | \det(DT_{a, b}(z))| \leq K_2 b
$$

\bigskip

We remark that Theorem \ref{theorem2}(1) confirms that ${\cal C}$ is the 
sole source of
nonhyperbolicity in the system. Part (2) expresses the fact that many orbits
experience at least some form of (nonuniform) hyperbolicity. 
A more detailed discussion is given in Section \ref{s7}.

%%%%%%%%%%%%%%%%%%%%%%%%%%%%%%%%%%%%%%%%%%%%%%%%%%%%%%%%%%%%%%%%
\subsection{SRB Measures and their Statistical Properties}
\label{s1.4}
\begin{definition} Let $g:M \to M$ be a diffeomorphism of a manifold.
A $g$-invariant Borel probability measure $\mu$ is called an
{\bf SRB measure} if $g$ has a positive Lyapunov exponents $\mu$-a.e.
and the conditional measures of $\mu$ on unstable manifolds are
absolutely continuous with respect to the Riemannian measure on
these manifolds.
\label{ndefinition1.1}
\end{definition}

In the absence of zero Lyapunov exponents, it follows from general
hyperbolic theory that an SRB measure has at most
a countable number of ergodic components, and that each ergodic
component has a positive measure set of generic points. A point
$z$ is said to be {\bf generic} with respect to $\mu$
if for every continuous function $\varphi$, $\frac{1}{n} \sum_{i=0}^n \varphi(g^
iz) \to \int \varphi d\mu$ as $n \to \infty$.
See  \cite{Led} and \cite{PS}.

\begin{theorem}{\bf (Existence and ergodic properties of SRB measures)}
\begin{itemize}
\item[(1)] $T$ admits an SRB measure. 
\end{itemize}
Assuming condition (**) above, we
have the following additional information:
\begin{itemize}
\item[(2)] $T$ admits at most $r$ ergodic SRB measures $\mu_i$ where $r$ is
the cardinality of the critical set of the 1-dimensional map $f$.
\item[(3)] Lebesgue-a.e. $z_0 \in R_0$ is generic with respect to some $\mu_i$; in fact, Lebesgue-a.e. $z_0 \in R_0$ lies in the stable manifold of a 
$\mu_i$-typical point in $\Omega$.
\end{itemize}
\label{theorem6}
\end{theorem}

We know from general hyperbolic theory that without zero
Lyapunov exponents, ergodic components of SRB measures are, up to finite
factors, mixing \cite{Led}.

\begin{theorem}
\label{theorem7}
{\bf (Decay of correlations and Central Limit Thoerem)}\
Let $\mu$ be an ergodic SRB measure, which, by taking a power of $T$
if necessary, we assume to be mixing. Then
\begin{itemize}
\item[(1)] for each  $\eta \in (0,1]$, there exists $\lambda=\lambda(\eta)<1$
such that if $\psi:A \to {\mathbb R}$ is H\"older continuous with exponent $\eta$ and 
$\varphi \in L^\infty(\mu)$, then there exists $K(\varphi, \psi)$ such that
$$
\left | \int (\varphi \circ T^n)
\psi d\mu - \int \varphi d\mu \int \psi d\mu \right | <K(\varphi, \psi) \lambda^
n
\ \ \ {\rm for \  all} \ n;
$$
\item[(2)] the Central Limit Theorem holds for all H\"older $\varphi$
with $\int \varphi d\mu =0$, i.e.
$$
\frac{1}{\sqrt{n}} \sum_{i=0}^{n-1} \varphi \circ T^i \ \ \to \ \
{\cal N}(0, \sigma)
$$
where ${\cal N}$$(0, \sigma)$ is the normal distribution with variance
$\sigma^2$; furthermore, $\sigma > 0$ if and
only if $\varphi \circ T \neq \psi \circ T - \psi$ for any $\psi$.
\end{itemize}
\end{theorem}

\newpage
We remark that the word ``attractor" has different meanings in the literature
(see \cite{Mil} for a discussion). In this article, it is convenient for us to
refer to $\Omega$ as ``the attractor". Theorem \ref{theorem6}
suggests, however, that from
a measure-theoretic point of view, it may be more appropriate to regard the
supports of the $\mu_i$ as attractors.

%%%%%%%%%%%%%%%%%%%%Section 1.3 
\subsection{Global geometry, symbolic dynamics and topological entropy}
\label{s1.3}
A {\bf monotone branch} of $R_n$ is a region diffeomorphic to a rectangle and 
bordered by two subsegments of $\partial R_n$.
Roughly speaking, it is the largest domain of this kind with the property that 
for
$0 \leq i \leq n$, the $x$-coordinates of its $T^{-i}$-image
stay inside some interval of monotonicity of $f$,
where $f$ is the initial 1-dimensional map from which $\{T_{a,b}\}$ is built. 
This notion is made precise in Section \ref{s8}, where a combinatorial tree
is introduced to describe the structure of a natural class of monotone branches.

 %%%%%%%%%%%%%%%%%%%%%%%%%Theorem 3 
\begin{theorem} 
{\bf (Coarse geometry of attractor)} \ \
There is a sequence
of neighborhoods $\tilde R_n$ of $\Omega$ with
$$
\tilde R_1 \supset \tilde R_2 \supset \tilde R_3 \supset \cdots
\ \ \ \ {\rm and}
\ \ \ \ \cap_i \tilde R_i = \Omega
$$
such that each $\tilde R_n$ is the union of a finite number of
monotone branches of $R_k$, $\ n \leq k \leq n(1+ K \theta)$,
where $\theta \sim \frac{-1}{\log b}$.
\label{theorem3}
\end{theorem}

Let $\{1,2, \cdots, k\}$ be a finite alphabet and let $\Sigma_k$
be the set of all bi-infinite sequences ${\bf s}=( \cdots, s_{-1}, s_0, s_1, 
\cdots)$ with $s_i \in \{1,2, \cdots, k\}$. The shift operator
$\sigma:\Sigma_k \to \Sigma_k$ is defined by $({\bf \sigma s})_i
=({\bf s})_{i+1}$. For $\Sigma \subset \Sigma_k$, we call
 $\sigma|\Sigma: \Sigma \to \Sigma$ a {\bf subshift} of the full shift
on $k$ symbols if $\Sigma$ is a closed $\sigma$-invariant subset of $\Sigma_k$.
 
Let $x_1 < x_2 < \cdots < x_r < x_{r+1}=x_1$ be the critical points
of $f$.
Let ${\cal C}_i^{(0)}$ be the component of ${\cal C}^{(0)}$ containing $x_i$
and let ${\cal C}_i={\cal C}\cap {\cal C}_i^{(0)}$. We remark that each 
${\cal C}_i$ is a fractal set
-- it is not contained in any smooth curve --
and that {\it a priori} there is no well defined notion of whether
a point lies to the left or to the right of ${\cal C}_i$.
%%%%%%%%%%%%%%%%%%%%%%%%%Theorem 4 
\begin{theorem} 
{\bf (Coding of orbits on attractor)} 
\begin{itemize}
\item[(1)] The critical set ${\cal C}$ partitions $\Omega \setminus {\cal C}$ 
into disjoint
sets $A_1, A_2, \cdots, A_r$ so that $z \in A_i$ has the interpretation
of being ``to the right" of ${\cal C}_i$ and ``to the left" 
of ${\cal C}_{i+1}$.
\item[(2)] There is a subshift $\sigma:\Sigma \to \Sigma$ of a full shift
on finitely many symbols and a continuous surjection $\pi: \Sigma
\to \Omega$ such that
$$
T \circ \pi \ = \ \pi \circ \sigma;
$$
$\pi$ is 1-1 except on $\cup_{i=-\infty}^\infty T^i {\cal C}$, 
where it is 2-1.
\item[(3)] Under the additional assumption that $f[x_j,x_{j+1}] \not \supset S^1$ for any $j$, the coding in (2) is given by (1), i.e. for all $z_0 \in \Omega \setminus 
\cup_{i=-\infty}^\infty T^i{\cal C}$, $\pi^{-1}(z_0)$ is the unique sequence $(s_i)_{i=-\infty}^{\infty}$ 
with $z_i \in A_{s_i}$.
\end{itemize}
\label{theorem4}
\end{theorem}

\begin{coro}
\label{ncoro1}
{\bf (Kneading sequences for critical points)}
\ For every $z_0 \in {\cal C}$, the itinerary of $\{z_1, z_2, \cdots \}$ is
uniquely represented by a sequence in $\Sigma$.
\end{coro}

Another consequence of Theorem \ref{theorem4} is the existence of 
equilibrium states. For a continuous map $g:X \to X$ of a compact metric space
and a continuous function $\varphi:X \to {\mathbb R}$, a $g$-invariant Borel 
probability measure $\mu$ on $X$ is called an {\bf equilibrium state}
for $g$ with respect to the potential $\varphi$ if $\mu$ maximizes the quantity
$$
\rm{sup} \ \{\ h_\nu(g)\ + \ \int \varphi d\nu \ \}
$$
where $h_\nu(g)$ denotes the metric entropy of $g$ with respect to
$\nu$ and the supremum is taken over all $g$-invariant Borel probability measures $\nu$. 

\begin{coro} 
\label{ncoro2}{\bf (Existence of equilibrium states)}\ 
T has an equilibrium state for every continuous 
$\varphi:\Omega \to {\mathbb R}$.
In particular, $T$ admits an invariant Borel probability measure 
maximizing entropy.
\end{coro}

The {\bf topological entropy} of $g$, written $h_{top}(g)$,
is usually defined in terms of open covers of arbitrarily small
diameters or in terms of $(n, \varepsilon)$-spanning or separated sets. 
For precise definitions, see \cite{Wa}. For the class of attractors studied
in this paper, $h_{top}(g)$ can be computated in more concrete ways.

In Theorem \ref{theorem4} 
we saw that every $z_0 \in \Omega$ can be unambiguously 
associated with one (and occasionally two) symbol sequences in $\Sigma$ 
determined by the locations of its iterates
with respect to the components of the critical set. 
We will show in Section \ref{s9} that in like manner
all the points in $R_0$ can be assigned symbol sequences --
except that this assignment is not unique. Let us temporarily refer to this as
the ``fuzzy" coding on $R_0$. Let

\begin{tabbing}
---\=$M^{\pm}_n$ \= = \= number of monotone segments in $\partial R^{\pm}_n$,
the two boundary components of \kill
     \>$N_n$  \> = \> number of distinct $n$-blocks in the coding of $\Omega$; \\
     \>${\tilde N}_n$ \> = \> number of distinct $n$-blocks in the ``fuzzy" coding of $R_0$; \\
     \>$P_n$ \> = \>number of fixed points of $T^n$; \\
     \>$M^{\pm}_n$ \> = \> number of monotone segments in $\partial R^{\pm}_n$, 
the two boundary components of \\ 
\> \> \> $R_n$ (see Sect. \ref{s8.1} for the precise definition).
\end{tabbing}

\newpage
\begin{theorem} {\bf (Formulas and inequalities for topological entropy)} \ 
\begin{itemize}
\item[(i)] 
$$
h_{\rm top}(T) \ = \ \lim_{n \to \infty} \frac{1}{n} \log N_n \ = \ 
\lim_{n \to \infty} \frac{1}{n} \log \tilde N_n \ = \ 
\lim_{n \to \infty} \frac{1}{n} \log P_n.
$$
\item[(ii)] 
$$
\limsup_{n \to \infty} \frac{1}{n} \log M_n^{\pm} \ \leq \ 
h_{\rm top}(T) \ \leq \ 
\left (\liminf_{n \to \infty}  \frac{1}{n} \log M_n^{\pm} \right )
(1 + \frac{K}{\log \frac{1}{b}}).
$$
\end{itemize}
\label{theorem5}
\end{theorem}

For a $1$-dimensional piecewise monotonic map $g$, it is a well known fact 
that $h_{\rm top}(g)$ is the growth rate of the number of 
intervals on which $g^n$ is monotonic \cite{MS}.
The factor $(1 + \frac{K}{\log \frac{1}{b}})$ gives, in a sense,
the potential {\it defect} in measuring the complexity of $T$ via 
the $1$-dimensional curves $\partial R_0^{\pm}$.

%%%%%%%%%%%%%%%%%
\subsection{H\'enon maps and homoclinic bifurcations}
\label{s1.5}

Theorems \ref{theorem1}--\ref{theorem5} are stated for attractors that arise 
from perturbations of circle maps. 
We state here, for the record, the corresponding results for 
{\it interval maps} and some of their applications. 
Reduction to the circle case is carried out in Appendix \ref{appendix-A.1}.

\begin{theorem} {\bf (Attractors arising from interval maps)}
\label{theorem8}
Let $I$ be a closed interval of finite length, and let $f:I \to I$ be 
a Misiurewicz map with $f(I) \subset$ int$(I)$.
Let $U$ be a neighborhood of $I \times \{0\}$ in ${\mathbb R}^2$, 
and let $\{T_{a,b}\}$ be a $2$-parameter family of maps with
$T_{a,b}:U \to {\mathbb R}^2$.
We identify $I$ with $I \times \{0\} \subset {\mathbb R}^2$, and 
assume that $\{T_{a,b}\}$ satisfies the conditions in Steps $II$, $III$ 
and $IV$ in Sect. \ref{s1.1} with $f_{a^*}=f$.
Then 
\begin{itemize}
\item[(i)] there exist $K>0$ and a rectangle $\hat \Delta=
[a_0, a_1] \times (0, b_1]$ arbitrarily near $(a^*, 0)$ such that 
for each $(a,b) \in \hat \Delta$, $T_{a,b}$ maps  $R:=I \times [-Kb, Kb]$ 
strictly into its interior, defining an attractor $\Omega:= 
\bigcap_{n \geq 0} T^n_{a,b} R$;
\item[(ii)] there is a positive measure 
set $\Delta \subset \hat \Delta$ such that the conclusions of Theorems 1--7
hold for $T=T_{a,b} \mid R$ for all $(a,b) \in \Delta$.
\end{itemize}
\end{theorem}

\begin{coro}  {\bf (The H\'enon family)} \ Let 
$$
T_{a,b}:(x,y) \mapsto (1-ax^2+y, \  bx),  \ \ \ \ (x,y) \in {\mathbb R}^2.
$$
Then for every $a^* \in [1.5, 2]$ for which $f_{a^*}:x \mapsto 1-a^* x^2$
is a Misiurewicz map, the conclusions of Theorem 8 hold.  
In particular, there is a 
positive measure set $\Delta$ near $(a^*, 0)$ such that the conclusions of
Theorems 1--7 hold for all $T=T_{a,b}$, $(a,b) \in \Delta$.
These results are valid for both $b > 0$ and $b < 0$.
\label{ncoro3}
\end{coro}

We remark that our method of proof does not distinguish between the orientation 
preserving and reversing cases of the H\'enon maps. 
When specialized to $a^*=2$ and $b>0$, the part of Corollary \ref{ncoro3} 
that corresponds to Theorem \ref{theorem1}, part (2), in this paper is 
a version of
the main result of \cite{BC2}.  
The parts of Corollary \ref{ncoro3} that correspond
to Theorem \ref{theorem6}(1),(2), Theorem \ref{theorem6}(3) 
and Theorem \ref{theorem7} are proved respectively in \cite{BY1}, 
\cite{BV} and \cite{BY2}.

\bigskip
Our last result concerns the application of Theorems 
\ref{theorem1}--\ref{theorem5} to 
homoclinic bifurcations.
Let $g_{\mu}, \ \mu \in [0, 1]$, be a $C^{\infty}$ one-parameter family of 
surface diffeomorphisms unfolding at $\mu=0$ a nondegenerate tangency of 
$W^u(p_0)$ and $W^s(p_0)$ where $p_0$ is a hyperbolic fixed point. We assume 
that the eigenvalues $\lambda$ and $\sigma$ of $Dg_0$ at $p_0$ satisfy 
$0 < \lambda < 1 < \sigma$ and $\lambda \sigma < 1$, and that they belong in 
the open and dense set of eigenvalue pairs that meet the hypotheses of 
Sternberg's linearization theorem. Under these conditions, it is well known
(see \cite{PT}) that for all sufficiently large $k$, 
there is an open 
set of parameters $\hat \Delta_k$ such that for all $\mu \in \hat \Delta_k$,
$g_{\mu}$ has a $k$-periodic attractor $\Omega_{\mu}$ all but finitely many
 of whose periodic components are located near the fixed point $p_{\mu}$. 

\begin{theorem}{\bf (Attractors arising from homoclinic bifurcations)} \ 
\label{theorem9}
Let $g_{\mu}$ be as above. Then for all sufficiently large $k$, there
is a positive measure set of parameters $\Delta_k \subset \hat \Delta_k$ 
for which the following hold: for all $\mu \in \Delta_k$,
there is a component $\Omega_{\mu}^0$ of $\Omega_\mu$ with the property 
that  if $T_{\mu}$ denotes the restriction of $g_{\mu}^k$ to a neighborhood of $\Omega_{\mu}^0$, then the conclusions of Theorems \ref{theorem1}--\ref{theorem5} hold for $T=T_\mu$.
\end{theorem}

Our proof of Theorem \ref{theorem9}, which is given in 
Appendix \ref{appendix-A.2}, consists of
observing that the maps $g_\mu^k$ meet the conditions of Theorem \ref{theorem8}. 
The part of Theorem \ref{theorem9} that corresponds to Theorem
\ref{theorem1}, part (2), in this paper
is the main result of \cite{MV}.

\newpage

%%%%%%%%%%%%%%%%%%%%Section 2
\section{Preliminaries}
\label{s2}
We gather in this section a collection of technical facts
used repeatedly in later sections. 
Most of the proofs are given in \ref{appendix-B}. 
Sects. \ref{s2.1}--\ref{s2.4} contain material not specific 
to the family $\{T_{a,b}\}$, and $K$ is not a ``system constant" in these 
subsections.

%%%%%%%%%%%%%%%%Section 2.1
\subsection{Linear algebra}
\label{s2.1}
Let $M$ be a $2 \times 2$ matrix.
Assuming that $M$ is not a scalar multiple of an orthogonal matrix,
we say that a unit vector $e$ defines {\bf the most contracted
direction} of $M$ if $\| M u \| \geq \| M e \|$
for all unit vectors $u$. 
For a sequence of matrices
 $M_1, M_2, \cdots$, we use $M^{(i)}$ to denote the matrix product $ M_i \cdots M_2 M_1$ and $e_i$ to denote the most contracted direction of $M^{(i)}$ when it makes sense. 

\bigskip
\noindent {\bf Hypotheses for Sect. \ref{s2.1}} \ The 
$M_i$ are $2 \times 2$ matrices; they satisfy $|\det(M_i)| \leq b$ and $\|M_i\| \leq K_0$ where $K_0$ and $b$ are fixed numbers with $K_0>1$ and $b<<1$. 

\medskip

%%%%Lemma 2.1
\begin{lemma}
\label{nlemma2.1}
There exists $K$ depending only on $K_0$ such that
if $\|M^{(i)}\| \geq \kappa^i$ and $\|M^{(i-1)}\| \geq \kappa^{i-1}$
for some $\kappa >> \sqrt b$, then $e_i$ and $e_{i-1}$ are well-defined, 
and
$$
\parallel e_i \times e_{i-1} \parallel
\leq (\frac{Kb}{\kappa^2})^{i-1}.
$$
\end{lemma}
%%%%Corollary 2.1 
\begin{corollary}
If for $1 \leq i \leq n$, $\|M^{(i)}\| \geq \kappa^i$ 
for some $\kappa >> \sqrt b$, then:

(a) $\parallel e_n - e_1 \parallel < \frac{Kb}{\kappa^2}$;

\medskip

(b) $\| M^{(i)} e_n \| \leq (\frac{Kb}{\kappa^2})^i$
for $1 \leq i \leq n$.
\label{ncorollary2.1}
\end{corollary}
\noindent {\bf Proof:} (a) follows immediately from Lemma 
\ref{nlemma2.1}.
For (b), since $\| e_n - e_i \| \leq (\frac{Kb}{\kappa^2})^i$,
we have $\| M^{(i)} e_n \| \leq \| M^{(i)}(e_n - e_i )\|
+ \| M^{(i)} e_i \|< K^i_0 \cdot (\frac{Kb}{\kappa^2})^i + 
(\frac{b}{\kappa})^i$. \hfill $\square$

\medskip

Next we consider for each $i$ a 2-parameter family
of matrices $M_i(s_1, s_2)$. For the purpose of the next corollary we 
make the
additional assumptions that for
$0 < j \leq 2$, $ \|\partial^j M_i(s_1, s_2)\| \leq K_0^i$ and
$|\partial^j \det(M_i(s_1, s_2))| <K_0^i b$ where $\partial^j$ represents
any one of the partial derivatives of order $j$ with respect to $s_1$ or $s_2$. 
Let $\theta_i(s_1, s_2)$ denote the angle $e_i(s_1, s_2)$ makes with
the
positive $x$-axis, assuming it makes sense.
%%%%Corollary 2.2
\begin{corollary}
Suppose that for some $\kappa >> \sqrt b$, 
$\|M^{(i)}(s_1, s_2)\| \geq \kappa^i$ for every $(s_1, s_2)$ and for
every 
$1 \leq i \leq n$.
Then for $j = 1, 2$, $ |\partial^j \theta_1| \leq K\kappa^{-(1+j)}$, and for 
$i \leq n$,
\begin{equation}
%%%%Formula 2.1
|\partial^j (\theta_i - \theta_{i-1}) | 
  <  (\frac{Kb}{\kappa^{(2+j)}})^{i-1},
\label{nformula2.1}
\end{equation}
%%%%Formula 2.2
\begin{equation}
\|\partial^j M^{(i)} e_n \| <
(\frac{Kb}{\kappa^{(2+j)}})^i.
\label{nformula2.2}
\end{equation}
\label{ncorollary2.2}
\end{corollary}

Our next lemma is a perturbation result.
Let $M_i, M'_i$ be two sequences of matrices,
let $w$ be a vector, and let $\theta_i$ and $\theta'_i$ denote the
angles $M^{(i)}w$ and $M'^{(i)}w$ make with the positive $x$-axis
respectively.
%%%%Lemma 2.2 
\begin{lemma} \ {\rm ([BC2], Lemma 5.5)} \ 
Let $\kappa, \lambda$ be such that $\frac{K b}{\kappa^2} < \lambda < K_0^{-12} 
\kappa^8$.
If for 
$1 \leq i \leq n, \ \|M_i-M'_i\| \leq \lambda^i$
and $\|M^{(i)}w\|\geq \kappa^i$, then

(a) $\| M'^{(n)}w \| \geq \frac{1}{2} \kappa^n$;
 
(b) $| \theta_n - \theta^{\prime}_n | < \lambda^{\frac{n}{4}}$.
\label{nlemma2.2}
\end{lemma}

Proofs of Lemmas \ref{nlemma2.1}, \ref{nlemma2.2} and 
Corollary \ref{ncorollary2.2} are given in Appendix \ref{app-B.1}.

\bigskip
\noindent {\bf Hypothesis for Sects. \ref{s2.2} and \ref{s2.3}:} \ 
$T:A \to A$ is an embedding of the form
$$
T(x,y)=(t_1(x,y), bt_2(x,y))
$$
where the $C^2$-norms of $t_1$ and $t_2$ are $\leq K_0$,
and $K_0>1$ and $b<<1$ are fixed numbers.

%%%%%%%%%%%%%%%Section 2.2 
\subsection{Stable curves}
\label{s2.2}

%%%%Lemma 2.3 
\begin{lemma}
Let $\kappa, \lambda$ be as in Lemma \ref{nlemma2.2} and $z_0 \in A$ be such 
that for 
$i=1, \cdots, n$, 
$\|DT^i(z_0)\| \geq \kappa^i$.
Then there is a $C^1$ curve
$\gamma_n$ passing through $z_0$ such that

(a) for all $z \in \gamma_n, \ d(T^iz_0, T^iz) \leq (\frac{Kb}{\kappa^2})^i$
for all $i \leq n$;

(b) $\gamma_n$ can be extended to a curve of length $\sim \lambda$  
or until it meets $\partial A$. 
\label{nlemma2.3}
\end{lemma}

A proof of this lemma is given in Appendix \ref{app-B.2}.

We call $\gamma_n$ a {\bf stable curve of order $n$}.
It will follow from this lemma that if $\|DT^i(z_0)\| \geq \kappa^i$
for all $i>0$, then there is a {\bf stable curve} $\gamma_\infty$
passing through $z_0$ obtained as a limit of the $\gamma_n$'s.
%%%%Section 2.3  
\subsection{Curvature estimates}
\label{s2.3}
Let $\gamma_0: [0, 1] \to A$ be a $C^2$ curve, and let
$\gamma_i(s) = T^i(\gamma_0(s))$. We denote the curvature of 
$\gamma_i$ at $\gamma_i(s)$ by $k_i(s)$.
%%%%Lemma 2.4
\begin{lemma}
\label{nlemma2.4}
Let $\kappa > b^{\frac{1}{3}}$. We assume that for every $s, \ k_0(s) \leq 1$ 
and 
$$\|DT^j(\gamma_{n-j}(s)) \gamma^{\prime}_{n-j}(s)\| 
\geq \kappa^j \|\gamma^{\prime}_{n-j}(s)\|
$$
for every $j < n$. Then 
$$
k_n(s) \leq \frac{Kb}{\kappa^3}.
$$
\end{lemma}

A proof is given in Appendix \ref{app-B.3}. 

%%%%Section 2.4 
\subsection{One-dimensional dynamics}
\label{s2.4}

We begin with some properties of maps satisfying the Misiurewicz condition.
Let $f$ be as in Sect. \ref{s1.1}, and let $C_\delta:=\{x \in S^1: 
d(x,C)<\delta\}$.

%%%%%%%%%%%Lemma 2.5
\begin{lemma} There exist $\hat c_0, \hat c_1>0$ such that the following
hold for all sufficiently small $\delta>0$: Let $x \in S^1$ be such that
$x, fx, \cdots, f^{n-1}x \not \in C_\delta$, any $n$. Then

(i) $| (f^n)'x | \geq \hat c_0 \delta e^{\hat c_1 n}$;

(ii) if, in addition, $f^nx \in C_\delta$, then 
$| (f^n)'x | \geq \hat c_0 e^{\hat c_1 n}$.
\label{nlemma2.5}
\end{lemma}

A proof is given in Appendix \ref{app-B.4}.

%%%%%%%%%Corollary 2.3
\begin{corollary} 
Let $c_0<\hat c_0$ and $c_1<\hat c_1$. Then for all sufficiently small $\delta$,
there exists $\varepsilon=\varepsilon(\delta)$ such that for all
$g$ with $\|g-f\|_{C^2}<\varepsilon$, (i) and (ii) above hold for $g$ with
$c_0$ and $c_1$ in the places of $\hat c_0$ and $\hat c_1$.
\label{ncorollary2.3}
\end{corollary}

\noindent {\bf Proof:} \ Let $N$ be such that $\delta e^{\hat c_1 N}
> e^{c_1 N}$, and choose $\varepsilon$ small enough so that for all 
$i \leq N$, if $x, gx, \cdots, g^{i-1}x \not \in C_\delta(g)$,
then $(g^i)'x \approx (f^i)'x$. \hfill $\square$

\medskip
The results in the rest of this subsection are not needed in this article.
We include them only as motivation for the corresponding results in
$2$-dimensions.

Temporarily write $C=C(g)$. To control $(g^n)'x$ when $g^ix \in C_\delta$ for some $i <n$,
we need to impose further conditions on $g$. Following \cite{BC1} and
\cite{BC2}, we assume there exist $\lambda>1$ and $0<\alpha <<1$ such that
for all $\hat x \in C$ and $n \geq 0$:

\medskip
(a) $d(g^n \hat x, C) \geq c_0 e^{-\alpha n}$  \ and

(b) $\mid (g^n)'(g \hat x)\mid \geq c_0 \lambda^n$.

\medskip
\noindent We define for each $x \in C_\delta$ a {\bf bound period}
$p(x)$ as follows. Fix $\beta>\alpha$. Let $\hat x \in C$ be such that
$|x-\hat x|<\delta$. Then $p(x)$ is the smallest $p$ such that
$$
|g^px-g^p \hat x|>c_0e^{-\beta p}.
$$

%%%%%%%%%Lemma 2.6
\begin{lemma}  {\bf (Derivative recovery)} \ There exists $K$ such that 
for $g$ satisfying the conditions above, if $|x-\hat x|=e^{-\mu}<\delta$ for some $\hat x \in C$, then

(i) \ $K^{-1}\mu \leq p(x) \leq K\mu$ ;

(ii) \ $K^{-1}(x-\hat x)^2 |(g^{i-1})'(g\hat x)|<|g^ix-g^i \hat x|<K(x-\hat x)^2 |(g^{i-1})'(g \hat x)|$;

(iii) \ $|(g^p)'x| \geq K^{-1}\lambda^{\frac{p}{2}}$  \ where $p=p(x)$.
\label{nlemma2.6}
\end{lemma}

\noindent {\bf Proof:} \ For this result there is no substantive difference
between the situation here and that of the
quadratic family $x \mapsto 1-ax^2$. See \cite{BC1} and \cite{BC2}, 
Section 2.
\hfill $\square$

\vskip .2in
\noindent {\bf Standing hypotheses for the rest of the paper:} \
$\{T_{a,b}\}$ is as in Sect. \ref{s1.1}. In particular, it has the form
$$
T_{a,b}(x,y) \ = \ (F_a(x,y)+bu_{a,b}(x,y), \ bv_{a,b}(x,y)).
$$
Where no ambiguity arises, we will write $T=T_{a,b}$. The phrase ``for $(a,b)$
sufficiently near $(a^*,0)$" will appear (finitely) many times in the next few sections. Each time it appears, the rectangle
in parameter space for which our results apply may have to be reduced.
From here on $K$ is the generic system constant as declared in Section \ref{s1}.

%%%%Section 2.5 
\subsection{Dynamics outside of ${\cal C}^{(0)}$}
\label{s2.5}
The first system constant to be chosen is $\delta$. A number of upper bounds
for $\delta$ will be specified as we go along. For now we think of it
as a very small positive number with $d(f^n \hat x, C)>>\delta$ for all $\hat x \in
C$ and $n > 0$. We assume also that $a$ is sufficiently near $a^*$
that the Hausdorff distances between the critical sets of $f_{a^*}$ and $f_a$ are $<<\delta$.

Recall that we will be working in $R_0=\{(x,y)\in A: |y| \leq Kb\}$. Our zeroth critical region ${\cal C}^{(0)}$ is defined to be
$$
{\cal C}^{(0)} =\{(x,y) \in R_0 :|x-\hat x|<\delta  \ \ \ {\rm for \  some} \  
{\hat x} \in C \}.
$$

Let $s(u)$ denote the slope of a vector $u$.
Assuming that $b^{\frac{1}{5}}<<\delta$, an easy calculation shows
that for $z \not \in
{\cal C}^{(0)}$, if $|s(u)|<\delta^2$, then 
$|s(DT(z)u)|={\cal O}(\frac{b}{\delta}) << \delta^2$.
Also,
if $\kappa_0 := \min{\|DT(z) u\|}$ where the minimum is taken over
all $z \not \in {\cal C}^{(0)}$
and unit vectors $u$ with $|s(u)|<\delta^2$, then  $\kappa_0>
K^{-1}\delta$.
Let $K(\delta):=\frac{K}{\kappa^3_0}$, so that $K(\delta)b$ is
the upper bound for $k_n$ in Lemma \ref{nlemma2.4}.
We call a vector $u$ a {\bf $b$-horizontal vector} if
$|s(u)|<K(\delta)b$. A curve $\gamma$  is called
a  {\bf $C^2(b)$-curve} if its tangent vectors are $b$-horizontal and its curvature is $\leq K(\delta)b$ at every point.  
%%%%Lemma 2.7  
\begin{lemma}
(a) For $z \not \in {\cal C}^{(0)}$, if $u$ is $b$-horizontal,
then so is $DT(z) u$.

(b) If $\gamma$ is a $C^2(b)$-curve outside of ${\cal C}^{(0)}$,
then $T(\gamma)$ is again a $C^2(b)$-curve.
\label{nlemma2.7}
\end{lemma}
\noindent {\bf Proof:} 
(a) has already been explained; (b) is an immediate consequence of (a) and 
Lemma \ref{nlemma2.4}.
\hfill $\square$

\medskip

Our next lemma describes the dynamics of $b$-horizontal vectors outside of 
${\cal C}^{(0)}$.
%%%%Lemma 2.8 
\begin{lemma}
There exist constants $c_0, c_1>0 $ independent of $\delta$ such that
the following holds for $T=T_{a,b}$ for all $(a,b)$ sufficiently near 
$(a^*,0)$. 
Let $z \in R_0$ be such that $z, Tz, \cdots, T^{n-1}z \not \in {\cal C}^{(0)}$, and let $u$ be a $b$-horizontal vector. Then

(i) $\| DT^n(z)u \| \geq c_0 \delta e^{ c_1 n}$;

(ii) if, in addition, $T^nz \in {\cal C}^{(0)}$, then 
$\| DT^n(z)u \| \geq c_0 e^{c_1 n}$.
\label{nlemma2.8}
\end{lemma}

\noindent {\bf Proof:} As with Corollary \ref{ncorollary2.3}, this follows from
Lemma \ref{nlemma2.5} by perturbation.
 \hfill $\square$
%%%%%%%%%%%%%%%Section2.6
\subsection{Critical points inside ${\cal C}^{(0)}$}
\label{s2.6} 
Wherever it makes sense, let $e_m$ denote the field of most contracted 
directions
of $DT^m$ and let $q_m$ be the slope of $e_m$.
When working with a curve $\gamma$ parameterized by arc length, 
we write $q_m(s)=q_m(\gamma(s))$.
We begin with some easy observations about $e_1$.

%%%%Lemma 2.9
\begin{lemma} For all $(a,b)$ sufficiently near $(a^*,0)$, $e_1$ is
defined everywhere on $R_0$, and there exists $K>0$ such that 

(a) $|q_1| > K^{-1}\delta$ outside of ${\cal C}^{(0)}$, and $q_1$ has opposite 
signs on adjacent components of $R_0 \setminus {\cal C}^{(0)}$;

(b) $$
 |\frac{dq_1}{ds}| > K^{-1}
$$
\label{nlemma2.9}
on every $C^2(b)$-curve $\gamma$ in ${\cal C}^{(0)}$.
\end{lemma}
\noindent {\bf Proof:}
The existence of $e_1$ follows from the fact that everywhere on $R_0$, 
$\|DT\|>K^{-1}$ (this uses the non-degeneracy condition in Step IV, Sect. 
\ref{s1.1}) while 
$|\det(DT)|={\cal O}(b)$. 
For  $a=a^*, \ b=0$ and $\{y=0\}$, the assertion in (a) is
obvious, and part (a) of Lemma \ref{nlemma2.9} follows by a perturbative 
argument.
The estimate for $|\frac{dq_1}{ds}|$ uses the non-degeneracy condition above 
and the fact that $f_{a^*}'' \neq 0$ on $C$. See Appendix \ref{app-B.5} 
for details.
 \hfill $\square$

%%%%Definition2.1 
\begin{definition}
Let $\gamma$ be a $C^2(b)$-curve in ${\cal C}^{(0)}$.
We say that $z_0$ is a {\bf critical point of order $m$ on $\gamma$} if

(a) $\|DT^i(z_0)\| \geq 1$ for $i=1,2, \cdots, m$;

(b) at $z_0$, $e_m$ coincides with the tangent vector to $\gamma$.
\label{ndefinition2.1}
\end{definition}

It follows from Lemma \ref{nlemma2.9} that on every $C^2(b)$-curve
that stretches across a component of
${\cal C}^{(0)}$, there is
a unique critical point of order $1$.
The next two lemmas are used in the ``updating" of existing critical
points and the creation of new ones. Their proofs are given in Appendix 
\ref{app-B.5}
%%%%Lemma2.10
\begin{lemma} \ {\rm ([BC2], p. 113)} \ 
Let $\gamma$ be a $C^2(b)$-curve in ${\cal C}^{(0)}$ where $\gamma(0)=z$
is a critical point of order $m$. We assume that

(a) $\|DT^i(z)\| \geq 1$ for $i=1,2, \cdots, 3m$;

(b) $\gamma(s)$ is defined for $s \in 
[-(Kb)^{\frac{m}{2}}, (Kb)^{\frac{m
}{2}}]$.

\noindent Then there exists
a unique critical point ${\hat z}$ of order $3m$ on $\gamma$, and 
$|{\hat z} - z| <(Kb)^m$.
\label{nlemma2.10}
\end{lemma} 
%%%%Lemma 2.11
\begin{lemma} \ {\rm ([BC2], Lemma 6.1)} \ 
For $\varepsilon > 0$, let $\gamma$ and $\hat{\gamma}$ be two 
disjoint $C^2(b)$-curves in
${\cal C}^{(0)}$ defined for 
$s \in [-4K_1\sqrt{\varepsilon}, 4K_1\sqrt{\varepsilon}]$ where $K_1$ is the 
constant $K$ in Lemma \ref{nlemma2.9}(b). We assume

(a) $\gamma(0)$ is a critical point of order $m$;

(b) the $x$-coordinates of $\gamma(0)$ and $\hat \gamma(0)$
coincide, and $\mid \gamma(0) - \hat{\gamma}(0) \mid < \varepsilon$.

\noindent Then there exists a critical point of order 
$\hat{m}$  at $\hat \gamma(\hat s)$ with $|\hat s| < 4K_1 \sqrt{\varepsilon}$ 
and $\hat{m} = min \{ m, K \log{\frac{1}{\varepsilon}} \}$.
\label{nlemma2.11}
\end{lemma}

%%%%%%%%%%%%%%%Section2.7
\subsection{Tracking $DT^n$: a splitting algorithm}
\label{s2.7}
The purpose of this section is to recall an algorithm
introduced in \cite{BC2} that gives, under suitable circumstances, a direct
relation between $DT^n$ and $1$-dimensional derivatives.
 
Let $z_0 \in R_0$, and let $w_0$ be a unit vector at $z_0$ that is 
$b$-horizontal. We write $z_n=T^nz_0$ and $w_n=DT^n(z_0)w_0$. 
In the case where
$z_i \not \in {\cal C}^{(0)}$ for all $i$, 
the resemblance to 1-d is made clear in Lemmas \ref{nlemma2.5} and 
\ref{nlemma2.8}.
Consider next an orbit $z_0, z_1, \cdots$ that visits
${\cal C}^{(0)}$ exactly once, say at time $t>0$. 
Assume:

\medskip

(a) there exists $\ell>0$ such that 
$\|DT^i(z_t){\tiny (\!\!\begin{array}{c} 0 \\ 1 \end{array}\!\!)}\| \geq 1$
for all $i<\ell$, so that in particular $e_\ell$, the most contracted
direction of $DT^\ell$, is defined at $z_t$, and

(b) $\theta(w_t, e_\ell)$, the angle between $w_t$ and
$e_\ell$, is $\geq b^{\frac{\ell}{2}}$.

\medskip

\noindent Then $DT^i(z_0)$ can be analyzed as follows. (Note that our
notation is different from that in \cite{BC2}.) We split $w_t$ into
$w_t=\hat w_t + \hat E$ where $\hat w_t$ is parallel to the vector 
${\tiny \left(\!\!\begin{array}{c} 0 \\ 1 \end{array}\!\!\right)}$
and $\hat E$ is parallel to $e_\ell$.
For $i \leq t$ and $i \geq t+\ell$, let $w_i^*=w_i$.
For $i$ with $t<i<t+ \ell$, let $w_i^*=DT^{i-t}(z_t)\hat w_t$.
We claim that all the $w_i^*$ are $b$-horizontal vectors, so
that $\{\|w_{i+1}^*\|/\|w_i^*\| \}_{i=0,1,2,\cdots}$ resemble a sequence
of 1-d derivatives. In particular, 
$\|w_{t+1}^*\|/\|w_t^*\| \sim \theta (w_t, e_\ell)$ simulates a drop in the 
derivative when an orbit comes near a critical point in 1-dimension.

To justify the statement about the slope of the $w_i^*$, we note that 
$DT(z_t){\tiny (\!\!\begin{array}{c} 0 \\ 1 \end{array}\!\!)}$ is
$b$-horizontal, so that
in view of lemma \ref{nlemma2.7} we need
only to consider $w_{t+\ell}^*$. We have
$$
\|DT^\ell(\hat{E})\| \leq b^\ell \frac{\|\hat w_t\|}{\theta (w_t, e_\ell)}
 \leq b^{\frac{\ell}{2}} \|\hat w_t\| \leq  b^{\frac{\ell}{2}}
\|DT^\ell(z_t)\hat w_t\|,
$$
the first and third inequalities following from (a) and the second from (b).
Since the slope of $DT^\ell(z_t)\hat w_t$ is smaller than $\frac{Kb}{2\delta}$,
it follows that $w_{t+\ell}^* 
=DT^\ell(z_t)\hat w_t+DT^\ell(z_t)\hat E$
remains $b$-horizontal.

The discussion above motivates the following
splitting algorithm introduced in
\cite{BC2}. Consider $\{z_i\}_{i=0}^{\infty}$, and let $t_1 < \cdots <
t_j < \cdots$ be the times when $z_i \in {\cal C}^{(0)}$. We let $w_0$
be a $b$-horizontal unit vector, and assume as before that
$e_{\ell_i}$ makes sense at $z_i$ for $i=t_j$.
Define $w_i^*$ as follows:
 
1. For $0 \leq i \leq t_1$, let $w_i^*=DT^i(z_0)w_0$.
 
2. At $i=t_j$, we split $w_i^*$ into
$$
w_i^*=\hat w_i+\hat E_i
$$
where
$\hat w_i$ is parallel to
${\tiny (\!\!\begin{array}{c} 0 \\ 1 \end{array}\!\!)}$ and 
$\hat E_i$ is parallel
to $e_{\ell_i}$.

3. For $i>t_1$, let
%%%%Formula2.6
\begin{equation}
w_i^*=DT(z_{i-1})\hat w_{i-1} \ + \ \sum_{j: \ t_j+\ell_{t_j}=i}
DT^{\ell_{t_j}}(z_{t_j})\hat E_{t_j}
\label{formula2.6}
\end{equation}
and let $\hat w_i=w_i^*$ if $i \neq t_j$ for any $j$.
 
This algorithm does not give anything meaningful in general.
It does, however, in the scenario of the next lemma.
%%%%Lemma 2.12
\begin{lemma}
Let $z_i, w_i$ and $w_i^*$ be as above. Assume
 
(a) for each $i=t_j$, $\theta(w_i^*, e_{\ell_i}) \geq b^{\frac{\ell_i}{2}
}$;
 
(b) the time intervals $I_j:=[t_j, t_j+\ell_{t_j}]$
are strictly nested, i.e. for $j \neq j'$, either $I_j \cap I_{j'} = \emptyset,
\ I_j \subset I_{j'}$, or
$I_{j'} \subset I_j$, and $t_j+\ell_{t_j} \neq t_{j'}+\ell_{t_{j'}}$.
 
Then $w_i=w_i^*$ for $i \not \in \cup_j I_j$, and
the $w_i^*$'s are all $b$-horizontal vectors.
The sequence $\{|w_i^*\|\}$ has the property that
\ $\|w_{i+1}^*\|/\|w_i^*\| \sim \theta (w_i^*, e_{\ell_i})$ for $i=t_j$,\
and \ $\|w_{i+1}^*\| \approx \|DT(z_i)w_i^*\|$
for $i \neq t_j$.
\label{nlemma2.12}
\end{lemma}

\noindent {\bf Proof}: The nested condition in (b) allows us to 
consider the $I_j$'s
one at a time beginning with the innermost time intervals. This reduces to the
case of a single visit to ${\cal C}^{(0)}$ treated earlier on. \hfill
$\square$

\vspace{.3in}
%%%%%%%%%%%%%%%%%%%%Section 3
\section{The Critical Set}
\label{s3} 

Many authors, including \cite{BC1}, \cite{CE}, \cite{J}, \cite{M1}, and 
\cite{NS}, have studied $1$-dimensional maps by controlling their 
critical orbits. These ideas were mimicked in \cite{BC2} 
where the authors developed techniques for identifying, for certain
H\'enon maps, a set they called the ``critical set". 
This is done via 
an inductive procedure involving parameter selection.
The first step in our analysis of the family $\{T_{a,b}\}$ is to
carry out a similar parameter selection, and the aim of
this section is to formulate suitable inductive hypotheses.

%%%%%%%%%%%%%%%Section3.1
\subsection{What is the critical set?}
\label{s3.1} 

In $1$-dimension, the critical set is where all previous expansion is
destroyed. Tangencies of stable and unstable manifolds play a similar
role in higher dimensions. Here is how we propose to capture the set
${\cal C}$ that we will prove in Section \ref{s7}
to be the {\it origin} of all {\it nonhyperbolic behavior}.

\medskip

Let ${\cal F}_0$ be the foliation on $R_0$ with leaves $\{y=$ constant$\}$,
and let ${\cal F}_k$ be its image under $T^k$.
In Sect. \ref{s2.5} we defined the $0$th critical region
${\cal C}^{(0)}$. Suppose that $T^i{\cal C}^{(0)} \cap {\cal C}^{(0)}
= \emptyset$ for all $i \leq n_0$. Then for $i \leq n_0$,
 ${\cal F}_i$ restricted to ${\cal C}^{(0)} \cap R_i$ consists of finitely
many bands of roughly horizontal leaves whose tangent vectors have been
expanded the previous $i$ iterates (Lemma \ref{nlemma2.8}).
From Corollaries \ref{ncorollary2.1}, \ref{ncorollary2.2}
and Lemma \ref{nlemma2.9}, we see also that in
${\cal C}^{(0)}$, $DT^i$ has a well-defined field of most contracted
directions, namely $e_i$, whose integral curves are roughly parabolas.
It is natural to take the set of tangencies in ${\cal C}^{(0)}$
between the leaves of ${\cal F}_i$ and the integral curves of $e_i$
to be our $i$th approximation of ${\cal C}$. Since these approximations
stabilize quickly with $i$, they would converge to ${\cal C}$ if this picture
could be maintained indefinitely, i.e. if the
``turns" of ${\cal F}_i$ could be prevented
from entering ${\cal C}^{(0)}$ for all $i$.

This, however, is impossible.
The ``turns" in ${\cal F}_1$ generated by what corresponds
to a single critical point of the $1$-dimensional map form a $1$-parameter
family of parabolas whose vertices lie on a roughly horizontal curve.
If this curve stays outside of ${\cal C}^{(0)}$,
it expands exponentially and therefore must intersect ${\cal C}^{(0)}$
after a finite number of iterates. What comes to our rescue is the observation
that the horizontal strips in ${\cal C}^{(0)} \cap R_i$ also become
exponentially thin with $i$, so that in all likelihood a roughly vertical curve
will intersect the attractor $\Omega$ in a very sparse Cantor set.
Since it is the ``turns"
{\it inside} $\Omega$ that count, it suffices to consider a Cantor set of
``turns", not the full $1$-parameter family.

These observations suggest that we modify our strategy as follows.
Since we do not know {\it a priori} the precise location of $\Omega$,
it is natural to consider a sequence of curves that limit on $\Omega$,
i.e. $\partial R_i$, $i=1, 2, \cdots$. We replace ${\cal F}_0$ by 
$\partial R_0$,
defining the $i$th approximation of the critical set for $i \leq n_0$
to be the set of tangencies between $\partial R_i$ and the integral curves
of $e_i$. 

Experience from $1$-dimension tells us that in order to retain a
positive measure set of parameters, we must allow
our ``turns" to approach the critical set slowly.
To maintain a picture similar to that for $i \leq n_0$, we shrink the 
critical regions {\it sideways} at a rate faster than this rate of approach.
As $i \to \infty$, the approximate critical sets converge to ${\cal C}$.

In order for the contractive fields above to be defined, it is necessary
that the derivative along orbits starting from $\cal C$ experience some
exponential growth. This growth, which is also useful for
controlling the movements of the ``turns",
is brought about in two ways: (i) by arranging for $n_0$ in the first
paragraph to be very large, growth is guaranteed for a long initial period;
(ii) when an orbit of ${\cal C}$ gets near a point $z \in {\cal C}$, it copies
the initial segment of the orbit of $z$, thereby {\it replicating} the growth
properties created in (i).

\medskip

A version of these ideas will be made precise in the inductive assumptions.
 
%%%%%%%%%%%%%%%%%%Section3.2
\subsection{Getting started}
\label{s3.2} 
The required initial growth in (i) above 
comes from the Misiurewicz property of $f$, the 1-dimensional map of
which $T$ is a perturbation. By choosing $(a,b)$ sufficiently near
$(a^*,0)$ and $\delta$ sufficiently small, $n_0$ can be arranged to be 
arbitrarily large.

Let $\Gamma_0$ be the set of all critical points of order $n_0$ on
$\partial R_0$. From Corollaries \ref{ncorollary2.1}, \ref{ncorollary2.2} 
and Lemma \ref{nlemma2.9}, we know that each connected segment of $\partial R_0 \cap {\cal C}^{(0)}$ contains exactly one point of $\Gamma_0$. These are our {\bf critical points of generation $0$}.

\medskip
In order to state properly our induction hypotheses, we introduce 
our main system constants. They are $\theta, \alpha, \beta,
\rho$, $c$, and $n_0$ and $\delta$ (which we have met): 

\smallskip
- There are two time scales, $N$ and a much slower one 
$\theta N$, where $\theta$ is chosen 

\ \ \ so that $b^\theta={\cal O}(1)$ and 
$<K^{-1}$ for some $K$ to be specified.

\smallskip
- $e^{-\alpha n}$ and $e^{-\beta n}$, with $\alpha<< \beta<<1$,
represent two small length scales.

\smallskip
- $c>0$ is our target Lyapunov exponent; it is $<c_1$ where 
$c_1$ is as in Lemma \ref{nlemma2.8}.

\smallskip
- Finally, $0<\rho<K^{-1}$ is an arbitrary number of order $1$. 
It determines the rate 

\ \ \ at which  our
critical regions decrease in size (see Sect. \ref{s3.1}).

\smallskip
The order in which these constants are chosen is as follows: $c, \rho, \alpha$
and $\beta$ are first fixed; $\delta$ is then taken as small as need to be.
The last constants to be determined are $n_0$ and $\theta$; observe that 
$n_0 \to \infty$ and $\theta \to 0$ corresponds essentially to $(a, b) \to 
(a^*, 0)$.

\bigskip

Parameters are deleted at each stage of our induction. {\bf Sections \ref{s3}
--\ref{s5}
are concerned with the dynamics of the maps corresponding to the
parameters retained. Issues pertaining to the measure of the set of
retained parameters} (including whether or not it is nonempty) {\bf are
postponed to Section \ref{s6}.}

%%%%%%%%%%%%%%%%%%%Section 3.3
\subsection{Inductive assumptions}
\label{s3.3}
Let $N \geq n_0$ be a large number, and let $\Delta_N$
be the set of parameters retained after $N$ iterates. We now
formulate a set of inductive assumptions
that describes the desired dynamical picture for $T=T_{a,b}, \
(a,b) \in \Delta_N$.
While we will continue to provide motivations and explanations,
(IA1)--(IA6) below are to be viewed as formal inductive hypotheses.
As before, let $z_i=T^iz_0$.
 
%%%%%%%%%%Section3.3.1 
\subsubsection{Critical points and critical regions}
\label{s3.3.1} 
 
\noindent {\bf (IA1)} \ {\bf (Structure of critical regions)} 
\  {\it For all $k \leq \theta N$,
the critical regions $ {\cal C}^{(k)}$ are defined
and have the geometric properties stated in {\rm (1)(i), (ii)} and 
{\rm (iii)} of
Theorem \ref{theorem1}. Moreover, on each horizontal boundary of each component of ${\cal C}^{(k)}$,
there is a critical point of order $N$ located within
${\cal O}(b^{\frac{k}{3}})$
of the midpoint of the segment.} 

\bigskip

Critical points on $\partial {\cal C}^{(k)}$ are called {\bf critical
points of generation $k$}. The set of critical points
of generation $\leq k$ is denoted by $\Gamma_k$.

%%%%%%%%%%%%%%%%%%Section 3.3.2
\subsubsection{Distance to critical set and loss of hyperbolicity}
\label{s3.3.2} 
If the critical set is where would-be stable and unstable directions
are interchanged, then distance to the
critical set might provide a measure of loss of hyperbolicity.
This is indeed the case under suitable circumstances and for a suitable notion of ``distance".

If $Q$ is a component of ${\cal C}^{(k)}$, we let $L_Q$ denote
the vertical line  midway between the two vertical boundaries of $Q$.

%%%%Definition 3.1 
\begin{definition} We say $z \in {\cal C}^{(0)}$ is {\bf horizontally
related} or simply {\bf h-related} to $\Gamma_{\theta N}$ 
if there exists a component $Q$ of 
${\cal C}^{(k)}$, $k \leq \theta N$, such that $z \in Q$ and $dist(z, L_Q)\geq b^{\frac{k}{20}}$. When this holds, we say $z$ is 
h-related to $z_0$ for all $z_0 \in \Gamma_{\theta N} \cap Q$.
\footnote{When studying the dynamics of $T$ on $\partial R_k$,
it will be convenient to include the following in the definition of
h-relatedness: Let $\gamma$ be a horizontal boundary of a component of
${\cal C}^{(k)}, \ k \leq \theta N$, and let $\hat z \in \gamma \cap
\Gamma_{\theta N}$. Then $z \in \gamma$ is also said to be h-related to 
$\hat z$.}
\label{ndefinition3.1}
\end{definition}

This is an attempt to describe the location
of a point relative to $\Gamma_{\theta N}$, which, as $N \to \infty$,
converges to a fractal set.
From Lemma \ref{nlemma4.1}, we see that $\Gamma_{\theta N} \cap Q$ is contained
in a region of width ${\cal O}(b^{\frac{k}{4}})$ in the middle of $Q$, so that $z$ and 
$\Gamma_{\theta N} \cap Q$ have a very obviously horizontal relationship.
We caution, however, that there may be points in $\Gamma_{\theta N}$
that are directly above or below $z$, and quite possibly both to its left and to
its right. Observe also that if $Q'$ is a component of ${\cal C}^{(k')}$ such that $z \in Q' \subset Q$, then
$dist(z, L_{Q'})\geq b^{\frac{k'}{20}}$.

\begin{definition}
For $z \in R_0$, we define its {\bf distance to the critical set},
denoted ${d_{\cal C}}(z)$, as follows: for $z \not \in {\cal C}^{(0)}$, 
let ${d_{\cal C}}(z)=\delta$;
for $z \in {\cal C}^{(0)}$, we let ${d_{\cal C}}(z)=dist(z,L_Q)$
where $Q$ is the component of ${\cal C}^{(k)}$ containing $z$
and $k$ is the largest number $\leq \theta N$ with $z \in  {\cal C}^{(k)}$.
We let $\phi(z)$ be one of the two points in $\partial Q \cap 
\Gamma_{\theta N}$ if $z$ is h-related to
$\Gamma_{\theta N}$.
\label{ndefinition3.2}
\end{definition}

For $z \in {\cal C}^{([\theta N])}$, the definitions of $d_{\cal C}(z)$ and $\phi(z)$ are temporary and will be modified as the
induction progresses.

 To secure growth properties for the orbits of $\Gamma_{\theta N}$,
we forbid them to approach the critical set too closely too soon.
(IA2) is a result of parameter selection.

\bigskip 
\noindent {\bf (IA2)} \ {\bf (Rate of approach to critical set)}
\ {\it For all $z_0 \in \Gamma_{\theta N}$
and all $i\leq N, \ d_{\cal C}(z_i) \geq min(\delta, e^{-\alpha i})$.}

\bigskip
(IA2) implies that for all $z_0 \in \Gamma_{\theta N}$ and $i \leq N$,
$z_i$ is h-related to $\Gamma_{\theta N}$ whenever it is in  ${\cal C}^{(0)}$. 
Intuitively, this is because $z_i$ is in a 
very ``deep" layer relative to its distance to $\Gamma_{\theta N}$.
Formally, let $z \in Q \subset {\cal C}^{(k)}$
where $Q$ and $k$ are as in Definition \ref{ndefinition3.2}. Then $k<<i$ since $\rho^k \geq e^{-\alpha i}$. Now $z_i \in R_i$. 
If $k< [\theta N]$, then $z_i \in Q \cap
R_{k+1}$, proving $d_{\cal C}(z_i) \geq \rho^{k+1}>>b^{\frac{k}{20}}$. 
If $k=[\theta N]$,
then $d_{\cal C}(z_i) \geq e^{-\alpha i}
\geq e^{-\alpha N} >>b^{\frac{1}{20}\theta N}$ provided that 
$b^\theta$ is chosen to be $<e^{-20\alpha}$.

%%%%Definition 3.3
\begin{definition} 
(a) For arbitrary $z \in {\cal C}^{(0)}$, we define its {\bf fold period} 
$\ell (z)$ to be the nonnegative integer $\ell \geq 1$ such that 
$b^{\frac{\ell}{2}}$ is closest to ${d_{\cal C}}(z)$.

(b) Given $z_0 \in R_0$ and unit vector $w_0$, 
we let $w^*_i, \ i=0,1, 2, \cdots,$ 
be given by the splitting algorithm in Sect. \ref{s2.7} with
$\ell_i=\ell(z_i)$ assuming 
$e_{\ell (z_i)}$ is defined at $z_i$.  
\label{ndefinition3.3}
\end{definition}

Recall that for $z_0 \in \Gamma_{\theta N}$, $\|DT^i(z_0)\| \geq 1$ for all 
$i \leq N$.
For $\ell \leq N$, Lemma \ref{nlemma2.2} gives an estimate on the size of the 
neighborhood of $\Gamma_{\theta N}$ on which $e_{\ell}$ 
is well defined. In particular, if $z$ is h-related to
$\Gamma_{\theta N}$, then
$e_{\ell (z)}$ is defined at $z$.

We fix $\varepsilon_0>0$ such that 
$\varepsilon_0<< |\frac{\partial q_1}{\partial x}|$ in ${\cal C}^{(0)}$ 
where $q_1$ is the slope of $e_1$.
For $z \in \partial R_k$, let $\tau(z)$ denote a tangent vector
to $\partial R_k$ at $z$. In the angle estimates below, $\tau$ and $e_\ell$ are
assumed to point in roughly the same direction as $w$. 

%%%%Definition 3.4
\begin{definition}
Let $z \in {\cal C}^{(0)}$ be h-related to
$\Gamma_{\theta N}$, 
and let $w$ be a vector at $z$. We say $w$ {\bf splits correctly} 
if $|\frac{w}{\|w\|}-\tau (\phi(z))| <\varepsilon_0 
d_{\cal C}(z)$.
\label{ndefinition3.4}
\end{definition}

\noindent {\bf (IA3)} \ {\bf (Correct splitting at returns)}
\ {\it For $z_0 \in \Gamma_{\theta N}, 
w_0 = {\tiny (\!\! \begin{array}{c} 0 \\ 1 \end{array}\!\!)}$ and $i \leq N$,
$w_i^*$ splits correctly whenever
$z_i \in {\cal C}^{(0)}$.}
\bigskip 

The sense in which this splitting is ``correct" is 
as follows. We wish to use Lemma \ref{nlemma2.12} to understand the evolution 
of $w_i$. First, (IA3) and Lemma \ref{nlemma2.9} together imply
condition (a) of this lemma. 
This is because
$|e_{\ell_i}(z_i)-\frac{w^*_i}{\|w^*_i\|}| \geq |e_{\ell_i}(z_i)-e_{\ell_i}(\phi(z_i))|
-|e_{\ell_i}(\phi(z_i))- \tau (\phi(z_i))| - |\tau (\phi(z_i))- \frac{w^*_i}{\|w^*_i\|}|
 \geq |\frac{\partial q_{\ell_i}}{\partial x}|d_{\cal C}(z_i)
- {\cal O}(b^{\ell_i}) -\varepsilon_0 d_{\cal C}(z_i)  \ \geq \ 
\frac{1}{2} |\frac{\partial q_1}{\partial x}|d_{\cal C}(z_i) 
\sim b^{\frac{\ell_i}{2}}$.
(For a comparison of $|\frac{\partial q_{\ell_i}}{\partial x}|$ and
$|\frac{\partial q_1}{\partial x}|$, see Corollary \ref{ncorollary2.2}.)
Condition (b) of Lemma \ref{nlemma2.12} is discussed in Sect. \ref{s4.1}.

\begin{picture}(13, 5)
\put(3, 0.5){
\psfig{figure=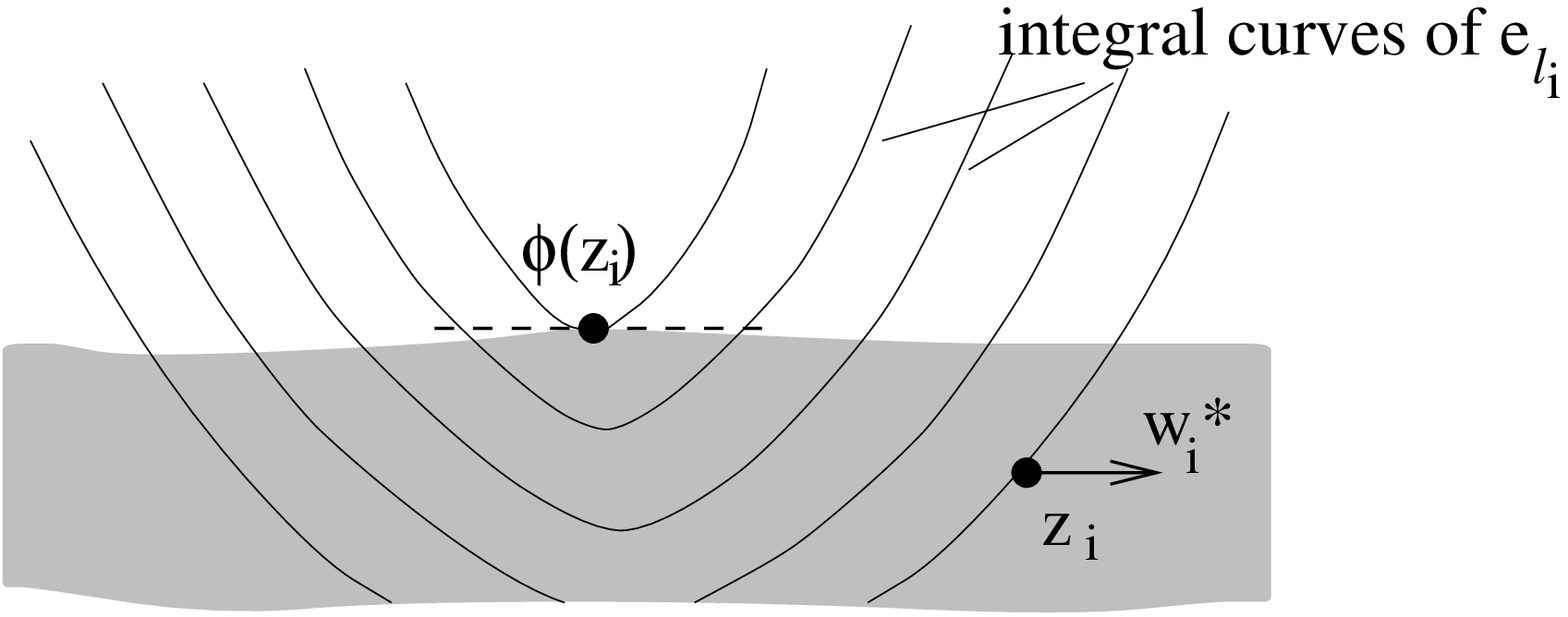,height=4.0cm,width=10cm}
}
\end{picture}

\centerline{Figure 2 \ Correct splitting of $w_i^*$}

%%%%%%%%%%%%%%%%Section 3.3.3
\subsubsection{Derivative along critical orbits}
\label{s3.3.3} 
We saw in the last paragraph that for $z_0 \in \Gamma_{\theta N}$,
as $z_i$ enters ${\cal C}^{(0)}$, $w_i^*$ suffers a loss of hyperbolicity
proportional to $d_{\cal C}(z_i)$. Combining this with (IA5)(c) below 
applied to an earlier step, we see
that this loss will be partially -- but not fully -- compensated for at 
the end of a certain period. To prevent a downward spiral in Lyapunov exponent, further parameter exclusion is needed. 

\bigskip 
\noindent {\bf (IA4)} \ {\bf (Derivative growth)} 
\ {\it For all $z_0 \in \Gamma_{\theta N}$ and 
$0\leq i \leq \frac{2}{3}N$, $\|w_i^*(z_0)\| \ > \ c_0e^{ci}$}.
\bigskip

In future steps of the induction, orbits of length $\frac{2}{3}N$
starting from $\Gamma_{\theta N}$ will be replicated; in other words,
they will serve as guides for other points that enter ${\cal C}^{(0)}$.

%%%%Definition 3.5
\begin{definition} 
\label{ndefinition3.5}
For arbitrary $\xi_0$ and $\xi_0^{\prime} \in {\cal C}^{(0)}$, 
we define their 
{\bf bound period} to be the largest integer $p$
such that for all $0<j \leq p$,
$$
|\xi_j- z_j| \ \leq \ e^{-\beta j}.
$$
\end{definition}

Observe that if $\xi_0^{\prime} = z_0 \in \Gamma_{\theta N}$, then for $j\leq p, \ |\xi_j- z_j| \ << d_{\cal C}( z_j)$.
We may assume $\delta$ is so small and $n_0$
so large that $d_{\cal C}(\xi_j)>\frac{\delta}{2}$ when $z_j$ is
outside of ${\cal C}^{(0)}$.
Our last two inductive assumptions deal with the
properties $z_0$ passes along to $\xi_0$.

\bigskip 
\noindent {\bf (IA5)} \ {\bf (Similarities with 1-dimensional maps)} 
\ {\it Let $z_0 \in \Gamma_{\theta N} \cap
\partial {\cal C}^{(k)}$, and let $\gamma:[0, \varepsilon] \to {\cal C}^{(0)}$
be a $C^2(b)$-curve with $\gamma(0)=z_0$
and $\gamma'(0)$ tangent to $\partial {\cal C}^{(k)}$. We regard all $\xi_0 \in \gamma$ as bound to $z_0$, and let $p(\xi_0)$ denote their bound periods. Then: 
 
(a) There exists $K$ such that for $\xi_0 \in \gamma$ with
$|\xi_0-z_0|=e^{-h}$,  
$$
\frac{1}{K} h \leq p(\xi_0) \leq K h \ \ \ {\it provided} \ \ \
Kh <\frac{2}{3}N;
$$
moreover, $p(\xi_0)$ increases monotonically
with the distance between $\xi_0$ and $z_0$;
  
(b) for $\ell \leq j\leq min(p, \frac{2}{3}N)$,
$|\xi_j-z_j| \approx |\xi_0-z_0|^2\|w_j(z_0)\|$ where ``$\approx$" means
up to a factor of $(1 \pm \varepsilon_1)$ for some $\varepsilon_1 > 0$;

(c) $\|w_p(\xi_0)\|\cdot |\xi_0-z_0| \geq e^{\frac{cp}{3}}$
provided $p<\frac{2}{3}N$.}

\bigskip
 
(IA5) describes the quadratic nature of the ``turn" as $\gamma$ is
mapped forward. For comparison with $1$-dimensional behavior, see 
Lemma \ref{nlemma2.6}.

\smallskip

The following distortion estimates are used in the proof of (IA5).
Let $w_0(\xi_0)=w_0(z_0)={\tiny (\!\!\begin{array}{c} 0 \\ 1 
\end{array}\!\!)}$, and let $\hat w^*_i(\xi_0)$ be 
given by Definition \ref{ndefinition3.3}(b) except that $e_{\ell (z_i)}$ 
(and not $e_{\ell (\xi_i)}$) is used for splitting at time $i$. 
(IA6) compares $w^*_i(z_0)$ and $\hat w^*_i(\xi_0)$. Let $M_i(\cdot)$ and  $\theta_i(\cdot)$ denote the
magnitude and argument of the vectors in question.  Define

%%%%Formula 3.1
\begin{equation}
\Delta_i(\xi_0, z_0) =
\sum_{s=0}^i (Kb)^{\frac{s}{4}} \mid \xi_{i-s} - z_{i-s} \mid.
\label{formula3.1}
\end{equation}

\medskip
\noindent {\bf (IA6)} \ {\bf (Distortion bounds)}
\ {\it Given $z_0 \in \Gamma_{\theta N}$
and any $\xi_0 \in {\cal C}^{(0)}$, we regard $\xi_0$ as bound to $z_0$
and let $p$ be the bound period. Then for $i \leq min \{p, N \}$, }
%%%%Formula 3.2
\begin{equation}
\frac{M_i(z_0)}{M_i(\xi_0)}, \ \ 
\frac{M_i(\xi_0)}{M_i(z_0)}  \ \leq \ 
\exp \{K \sum_{j=1}^{i-1} \frac{\Delta_j}{d_{\cal C}(z_j)} \}
\label{formula3.2}
\end{equation}
{\it and}
\begin{equation}
\mid \theta_i(\xi_0) - \theta_i(z_0) \mid \leq
(Kb)^{\frac{1}{2}}\Delta_{i-1} .
\label{formula3.3}
\end{equation}
 
\noindent {\it The estimates above also hold with $w^*_i(z_0)$ replaced by $\hat w^*_i(\xi_0')$
where $\xi_0'$ is another point in ${\cal C}^{(0)}$ also thought of as bound to
$z_0$, and $p$ is the minimum of the two bound periods.}

\bigskip
 
Let us return  for a moment to Definition \ref{ndefinition3.1}. 
From the geometry of
${\cal C}^{(k)}$ (see (IA1) and Lemma \ref{nlemma4.1}) it is an exercise in 
calculus
to show  that if $\xi_0$
is h-related to $z_0 \in \Gamma_{\theta N}$, then
it lies on a $C^2(b)$-curve through $z_0$ tangent to $\tau(z_0)$.
In particular, (IA5) applies.
 
Our rules of parameter exclusion, namely (IA2) and (IA4), are 
similar to those used in \cite{BC2}, but they are
applied to different orbits and with a different definition of 
``$d_{\cal C}(\cdot)$". The notions of bound and fold periods are borrowed 
from \cite{BC2}, as are (IA5) and (IA6).
Our construction of ${\cal C}$, however, has a distinctly different flavor.

%\vspace{.3in}
\newpage
%%%%%%%%%%%%%%%%%%%%%Section 4 
\section{Replication of Orbit Segments}
\label{s4} 
In Sect. \ref{s3.1} we outlined a scheme for obtaining derivative growth along
critical orbits, namely to choose a start-up geometry that guarantees 
some initial growth, and then to try to replicate this behavior.
Section \ref{s4} contains a detailed analysis of the replication process.
The main results are stated in Sect. \ref{s4.3}, after some
technical preparations in Sects. \ref{s4.1} and \ref{s4.2}, 
including amending slightly the definitions
of bound and fold periods. (IA1)--(IA6) are assumed
up to time $N$.

%%%%%%%%%%%%%%%Section 4.1
\subsection{Nested properties of bound and fold periods}
\label{s4.1} 

Consider $z_0 \in \Gamma_{\theta N}$. When $z_i$ enters
${\cal C}^{(0)}$, it is natural to assign to it a 
{\bf bound period} $p(z_i)$ defined using $\phi(z_i)$. An unsatisfactory aspect of this definition
is that two bound periods so defined may overlap without one being completely 
contained in the other. The purpose of this subsection is to adjust 
slightly the definition of
$p(z_i)$ to create a
simpler binding structure. A similar adjustment is made in \cite{BC2}.
 
First we fix some notation.
Let $Q^{(j)}$ denote  the components of ${\cal C}^{(j)}$, and 
let $\hat Q^{(j)}$ be the component of $R_j \cap {\cal C}^{(j-1)}$ containing $Q^{(j)}$.  For $z \in \partial R_j$,
let $\tau(z)$ be a unit vector at $z$ tangent to $\partial R_j$.

%%%%Lemma 4.1 
\begin{lemma}
\label{nlemma4.1}
For $z, z' \in \Gamma_{\theta N} \cap Q^{(k)}$, we have
$$
|z-z'|={\cal O}(b^{\frac{k}{4}})  \ \ \ \ {\it and} \ \ \ \ \| \tau(z) \times \tau(z^{\prime}) \| = {\cal O}(b^{\frac{k}{4}}).
$$
\end{lemma}

\noindent {\bf Proof:} \ Let $z^{(k)}$ be a critical point in $\partial Q^{(k)}$. For $k \leq i < [\theta N]$, 
let $z^{(i+1)}$ be a critical point of generation $i+1$ in 
$Q^{(i)}(z^{(i)})$, the component of $Q^{(i)}$ containing  $z^{(i)}$. 
From (IA1) we know that the Hausdorff distance between the two horizontal boundaries of $Q^{(i)}(z^{(i)})$ is 
${\cal O}(b^{\frac{i}{2}})$.  Lemma \ref{nlemma2.11} then tells us that $|z^{(i)}-z^{(i+1)}|={\cal O}(b^{\frac{i}{4}})$. The angle estimate
also follows from the proof of Lemma \ref{nlemma2.11} \hfill $\square$
%%%%Lemma 4.2
\begin{lemma}
Let $\xi_0$ be h-related to $z_0 \in \Gamma_{\theta N}$.
If during their bound period $z_i$ returns to ${\cal C}^{(k)}$, then
$\xi_i \in  \hat Q^{(k)}(z_i)$.
\label{nlemma4.2}
\end{lemma}
\noindent {\bf Proof:} \ Let $\gamma$ be a $C^2(b)$-curve joining $z_0$
and $\xi_0$. 
Then $T^i\gamma \subset R_i$. Since $e^{-\alpha i} \leq d_{\cal C}(z_i) \leq \rho^k$, we have $k < i$ and therefore $T^i\gamma \subset R_k$.
By the monotonicity of bound periods,
every point in $T^i\gamma$ is within a distance of $<e^{-\beta i}$
from $z_i$. This puts $\xi_i \in R_k \cap Q^{(k-1)}(z_i)$. \hfill $\square$
%%%%Lemma 4.3 
\begin{lemma}
Let $z_0 \in \Gamma_{\theta N}$ be such that $z_i \in {\cal C}^{(0)}$
at times $t_1<t_2< \cdots <t_r$, and that for each $j<r$
the bound period $p_j$ initiated at time $t_j$ extends beyond time $t_{j+1}$.
Then $p_j<(K\alpha)^{j-1}p_1$.
\label{nlemma4.3}
\end{lemma}
\noindent {\bf Proof:} \ Let $\tilde z_0=\phi(z_{t_1})$. We claim that
$|z_{t_2}-\phi(z_{t_2})| \approx |\tilde z_{t_2-t_1}-\phi(\tilde z_{t_2-t_1})|$, which is $>e^{-\alpha (t_2-t_1)}$. If true, this will
imply, by (IA5)(a), that $p_2<K\alpha (t_2-t_1)<K\alpha p_1$,
and the assertion in
the lemma will follow inductively. 
Since $|z_{t_2}-\tilde z_{t_2-t_1}|<e^{-\beta(t_2-t_1)}<<e^{-\alpha (t_2-t_1)}$, it suffices to show that $|\phi(\tilde z_{t_2-t_1})-
\phi (z_{t_2})| << |\tilde z_{t_2-t_1}-\phi(\tilde z_{t_2-t_1})|$.
Let $k$ be the largest number such that $\tilde z_{t_2-t_1}$
$\in {\cal C}^{(k)}$. By Lemma \ref{nlemma4.2}, 
$z_{t_2} \in Q^{(k-1)}(\tilde z_{t_2-t_1})$, 
so $\phi(\tilde z_{t_2-t_1})$ and $\phi(z_{t_2})$ must both be in
$Q^{(k-1)}(\tilde z_{t_2-t_1})$. By Lemma \ref{nlemma4.1}
they are $\leq b^{\frac{k-1}{4}}$ apart, and this is
$<<|\tilde z_{t_2-t_1}-\phi(\tilde z_{t_2-t_1})|$. \hfill $\square$

%%%%Definition 4.1 
\begin{definition} For $z_0 \in \Gamma_{\theta N}$ with $z_i \in
{\cal C}^{(0)}$, the {\bf adjusted bound period} $p^*(z_i)$
is defined to be the smallest number $p^* $ with the property
that for all $j$ with $i\leq j<i+p^*$, if $z_j \in {\cal C}^{(0)}$,
then $j+p(z_j) \leq i+p^*$.
\label{ndefinition4.1}
\end{definition}

Adjusted bound periods, therefore, have a nested structure by definition.

%%%%Corollary 4.1 
\begin{corollary}
(a) $p^* \leq p+K\alpha p$.

(b) For $z_i \in {\cal C}^{(0)}$ with $\phi(z_i)=\hat z_0$, we have
for all $j \leq p^*$,
$$
|z_{j+i}-\hat z_j|<e^{-\beta^*j}
$$
for some $\beta^*$ smaller than $\beta$ and $>>\alpha$.
\label{ncorollary4.1}
\end{corollary}
 
The proof is left as an exercise. We assume from here on that all bound periods for all critical orbits are adjusted, and write
$p$ and $\beta$ instead of $p^*$ and $\beta^*$.

This amended definition gives critical orbits the following simple structure of {\bf bound} and {\bf free states}. 
We call $z_i$ a {\bf return} if $z_i \in {\cal C}^{(0)}$.
Then $z_i$ is free for $i \leq n_1$ where $n_1>0$ is the time of the first return, and it is in bound state for $n_1<i \leq n_1+p_1$ where $p_1$ is the bound period initiated at time $n_1$. 
After time $n_1+p_1$, $z_i$ remains free
until its next return at time $n_2$, is bound for the
next $p_2$ iterates, and so on. The times $n_j$ are called {\bf free return} times. A {\bf primary bound period} begins at each $n_j$.
Inside the time interval $[n_j, n_j+p_j]$, there may be {\bf secondary
bound periods} which comprise disjoint time intervals, and so on.

Next we consider fold periods, which are denoted by $\ell$
and defined in Sect. \ref{s3.3.2}.
As with bound periods, if $z_i$ enters ${\cal C}^{(0)}$
at times $t_1$ and $t_2$ with $t_1<t_2 \leq N$,
and if the fold period begun at $t_1$ remains
in effect at $t_2$, then using Lemma \ref{nlemma4.2} we see that
$\ell_{t_2} < \frac{\alpha}{\log \frac{1}{b}} \ell_{t_1}$, so that 
{\bf adjusted fold periods}
can be defined similarly to give a nested structure.
This is condition (b) of Lemma \ref{nlemma2.12} .
 A further simplifying arrangement, which we
will also adopt, is that no fold periods expire at returns to ${\cal C}^{(0)}$
or at the step immediately after.
The proof of the following lemma is straightforward and will be omitted.
%%%%Lemma 4.4 
\begin{lemma} 
\ {\rm (cf. [BC2], Lemma 6.5)} \ 
Let $z_0 \in \Gamma_{\theta N}$. Then for every $i<N$,
there exist $i_1 \leq i \leq i_2$ with
$$
i_2 - i_1 < K \theta  \alpha i
$$
such that $i_1$ and $i_2$ are out of all fold periods.
\label{nlemma4.4}
\end{lemma}
%%%%%%%%%%%%%%%%%%%%%%Section 4.2
\subsection{Orbits controlled by $\Gamma_{\theta N}$}
\label{s4.2} 
In this subsection we consider $(z_0, w_0)$
where $z_0$ is an arbitrary point in $R_0$ and $w_0$ is a unit vector.
We write $z_i=T^iz_0$ and $w_i=DT^i(z_0)w_0$.

%%%%Definition 4.2
\begin{definition} We say $(z_0, w_0)$ is {\bf controlled} by
$\Gamma_{\theta N}$ up to time $m$ (with $m$ possibly $>N$)
if the following hold.
 
- Initial conditions: if $z_0 \not \in {\cal C}^{(0)}$, then $w_0$ is a
$b$-horizontal vector;  if $z_0 \in {\cal C}^{(0)}$, then
 
\ \ \ \ \ \ either $w_0={\tiny \left(\!\!\begin{array}{c} 0 \\ 1
\end{array}\!\!\right)}$, or $z_0$ 
is h-related to $\Gamma_{\theta N}$ and $w_0$ splits correctly.
 
- For $0<i \leq m$, if $z_i \in {\cal C}^{(0)}$, then $z_i$ is h-related
to $\Gamma_{\theta N}$ and  $w_i^*$ splits correctly.
\label{ndefinition4.2}
\end{definition}
 
No h-relatedness property is required for $z_0 \in {\cal C}^{(0)}$
when $w_0={\tiny \left(\!\!\begin{array}{c} 0 \\ 1 \end{array}\!\!\right)}$ 
because for practical purposes, one may think of
the sequence as starting with $(z_1, w_1)$.
 
Let $(z_0, w_0)$ be as above. Then the orbit of $z_0$ has a natural
bound/free structure defined as follows: If
$z_0 \in \Gamma_{\theta N}$,
then it is natural to regard $z_0, z_1, \cdots, z_i$ as free until $z_i$ 
returns to ${\cal C}^{(0)}$. For
$z_0 \in {\cal C}^{(0)} \setminus \Gamma_{\theta N}$, we may regard
$z_0$ as bound to any
$\hat z \in \Gamma_{\theta N}$ for a period $p$ provided that
$(\max \|DT\|)^p \ |z_0-\hat z|<e^{-\beta p}$. (This trivial bound
period is used to ensure that Lemma \ref{nlemma4.2} continues to work.)
When $z_i$ is h-related to $\Gamma_{\theta N}$, we take the bound
period to be that between $z_i$ and $\phi(z_i)$ (which is longer
than the trivial one).
Observe that Lemma \ref{nlemma4.3} is equally valid for controlled orbits
as for orbits starting from $\Gamma_{\theta N}$,
so that a nested structure can also be assumed for the bound
and fold periods
of controlled orbits.

In the language of Definition \ref{ndefinition4.2}, the situation can be summed up as follows. First, it follows from (IA2) and (IA3) that
for all $\hat z_0 \in \Gamma_{\theta N}$, $(\hat z_0, {\tiny \left(\!\!\begin{array}{c} 0 \\ 1 \end{array}\!\!\right)})$ is controlled by
$\Gamma_{\theta N}$ up to time $N$.
Second, for $(z_0, w_0)$ controlled by $\Gamma_{\theta N}$,
(IA5) and (IA6) apply to give information during its bound periods.
In particular, the orbit of $(z_0, w_0)$ has similar bound/free structures and ``derivative recovery" estimates as those of 
$(\hat z_0, {\tiny \left(\!\!\begin{array}{c} 0 \\ 1 \end{array}\!\!\right)})$, $\hat z_0 \in \Gamma_{\theta N}$, except that (IA2) and (IA4) need not hold.

\medskip
In the remainder of this subsection we record some technical facts 
on $w_i$ and $w^*_i$. Their proofs are given in Appendix \ref{app-B.6}. In 
Lemmas \ref{nlemma4.6}--\ref{nlemma4.8}, it is assumed that $(z_0, w_0)$ is controlled
by $\Gamma_{\theta N}$ up to time $m$, and  
all time indices are $\leq m$. 

%%%%Lemma 4.5
\begin{lemma}  Suppose $(z_0, w_0)$ satisfies the initial conditions
in Definition \ref{ndefinition4.2}, and for $0<i \leq m$, $z_i$ is h-related to
$\Gamma_{\theta N}$ at all returns. Then  $(z_0, w_0)$ is controlled
up to time $m$ if $w_i^*$ splits correctly at all free returns.
\label{nlemma4.5}
\end{lemma}

%%%%Lemma 4.6 
\begin{lemma}
\label{nlemma4.6}
Under the additional assumption that $d_{\cal C}(z_i)>
e^{-\alpha i}$ for all $i \leq m$, we have
$$
K^{-\varepsilon i}\|w_i^*\| \ \leq \ \|w_i\| \
\leq \ K^{\varepsilon i} e^{\alpha i} \|w_i^*\|, \ \ \ \ \ \varepsilon=K\alpha \theta.
$$
\end{lemma}

%%%%Lemma 4.7 
\begin{lemma}
There exists $c'>0$ such that for every $0\leq k <n$,
$$
\|w_n^*\| \geq K^{-1} d_{\cal C}(z_j)  e^{c'(n-k)} \|w_k^*\|
$$
where $j$ is the first time $\geq k$ when a bound period
extending beyond time $n$ is initiated. If no such $j$ exists,
set $d_{\cal C}(z_j)=1$.
\label{nlemma4.7}
\end{lemma}

%%%%Lemma 4.8 
\begin{lemma} 
\label{nlemma4.8}
Let $k<n$ and assume $z_n$ is free. Then
$$
\|w_n\|>  K^{-K\theta(n-k)} e^{c'(n-k)} \|w_k\|.
$$
\end{lemma}

%%%%%%%%%%%%%%%%%%%%%%Section 4.3
\subsection{Controlled orbits as ``guides" for other orbits}
\label{s4.3} 

(IA2)--(IA6) are about orbits starting from $\Gamma_{\theta N}$.
In Sect. \ref{s4.2} we introduced a class of orbits that successfully use
 orbits from $\Gamma_{\theta N}$ as their ``guides". We now let these
orbits serve as guides for other orbits 
and study the properties they pass along. This is the essence of the 
replication process.

Throughout Sect. \ref{s4.3} we assume that

\medskip
(1) $z_0 \in {\cal C}^{(0)}$, 
$w_0={\tiny \left(\!\!\begin{array}{c} 0 \\ 1 \end{array}\!\!\right)}$,
and $(z_0, w_0)$ is 
controlled by $\Gamma_{\theta N}$ up to time $m$;

(2) $d_{\cal C}(z_i)>e^{-\alpha i}$ for all $ 0 < i \leq m$.
\medskip
 
Our first order of business is to establish that for all $\xi_0$
bound to $z_0$, $\hat w_i^*(\xi_0)$ copies $w_i^*(z_0)$ faithfully.
A detailed proof of the following lemma is given in
Appendix \ref{app-B.7}. 

%%%%Lemma 4.9 
\begin{lemma} {\rm (cf. [BC2], Lemma 7.8)} \ 
Let $(z_0,w_0)$ be as above, and let $\xi_0 \in {\cal C}^{(0)}$ be an arbitrary
point which we think of as bound to $z_0$. Let $M_\mu( \cdot)$ and 
$\theta_\mu(\cdot)$ have the same meaning in (IA6). Then the estimates for
$$
\frac{M_\mu(\xi_0)}{M_\mu(z_0)}, \ \  \frac{M_\mu(z_0)}{M_\mu(\xi_0)}
\ \ \ {\it and} \ \ \ |\theta_\mu(\xi_0)-\theta_\mu(z_0)|
$$
as stated in (IA6) hold for all $\mu \leq min(p,m)$.
The corresponding distortion estimates for two points $\xi_0$ and $\xi_0'$ bound to $z_0$
apply as well.
\label{nlemma4.9}
\end{lemma}

In the rest of this subsection we consider the situation where $z_0$
is a critical point on a $C^2(b)$-curve in the sense of Sect. \ref{s2.6} and
study the quadratic behavior as this curve is iterated.
More precisely, let $e_m$ be the contractive field of order $m$, which we know 
from Lemmas \ref{nlemma4.6} and \ref{nlemma4.7}
is defined
at $z_0$. We assume 

\medskip
(3) $z_0$ lies on a $C^2(b)$-curve $\gamma \subset
{\cal C}^{(0)}$, and $e_m(z_0)$ is tangent to 
$\gamma$.
\medskip

For $\xi_0 \in \gamma$, let $p=p(\xi_0)$ denote the bound period between $z_0$ and $\xi_0$. We assume that during its bound period,
the orbit of $\xi_0$ inherits the secondary and higher order bound structures of the orbit of $z_0$.

%%%%Lemma 4.10 
\begin{lemma} In the part of $\gamma$ where $p<m$, $p$
increases monotinically with distance from $z_0$.
\label{nlemma4.10}
\end{lemma}
\noindent {\bf Proof:} \ Proceeding inductively, we assume that on
a connected subsegment $\gamma_k$ of $\gamma$ one of whose end points is
$z_0$, the minimum bound period
is $k$. It suffices to show that at time $k+1$, the part of $\gamma_k$ that remains bound to $z_0$ is connected. We may assume $T^k(\gamma_k)$
is not in a secondary fold period (otherwise all of $T^{k+1}(\gamma_k)$
will be in a bound period), and that $d_{\cal C}(\xi_0)>\frac{1}{2}\delta$ for all $ \xi_0 \in T^k(\gamma_k)$.
 
Let $T^k(\gamma_k)=\gamma^{(1)} \cup \gamma^{(2)}$ where 
$\gamma^{(1)}$ consists of points for which the
primary fold period remains in effect and $\gamma^{(2)}$ 
its complement. Then $\gamma^{(1)}$ is contained in a disk $B$ of radius $K^k b^{\frac{k}{2}}$ 
centered at $z_k$, and the bound period on no part of $B$ can expire
at time $k+1$. If the bound period of any part of $\gamma^{(2)}$ is to expire at time $k+1$, then the far end of $\gamma^{(2)}$ must be $>K^{-1}e^{-\beta (k+1)}$ from $z_k$. Also, 
its tangent vectors are $b$-horizontal.
One concludes that $T^k(\gamma) \setminus B$ is a $b$-horizontal
connected segment which will remain horizontal
in the next iterate, forcing the desired picture. \hfill $\square$

\bigskip

Let $s \to \xi_0(s)$ be the parametrization of $\gamma$ by 
arc length with $\xi_0(0)=z_0$. The following lemma, whose proof is 
given in  Appendix \ref{app-B.8}, contains a {\bf distance formula} for
$|\xi_\mu(s)-z_\mu|$. See Sect. \ref{s2.4} for comparison with $1$-d.

\begin{lemma} 
Let $\varepsilon_1>0$ be given. Then assuming $\delta$ is sufficiently small,
we have, for all $\mu \in {\mathbb Z}^+$ and $s>0$ satisfying $\mu \leq m$,
$(Kb)^{\frac{\mu}{2}}<s$ and $p(\xi_0(s)) \geq \mu$,
\begin{equation}
(1-\varepsilon_1)\ \|w_\mu(0)\| \ K_1 s^2 \ < \  |\xi_\mu(s)-z_\mu| 
\ < \ (1+\varepsilon_1)\ \|w_\mu(0)\| \ K_1 s^2
\label{nequation4.3(1)}
\end{equation}
where $K_1= \frac{1}{2}|\frac{dq_1}{dx}(z_0)|$.
\label{nlemma4.11}
\end{lemma}

%%%%Corollary 4.2
\begin{corollary}
Assume in addition to (1)--(3) above
that $\|w_j^*(z_0)\|>e^{cj}$ for all $j \leq m$. Let $\xi_0 \in \gamma$.
Suppose that $|\xi_0-z_0|=e^{-h}$ and $p(\xi_0) \leq m$. Then

(a) $\frac{h}{3K_2} \leq p \leq \frac{3h}{c}$ where $K_2=\log\|DT\|$;

(b) $\|w_p(\xi_0)\| \cdot |\xi_0-z_0| \geq e^{\frac{cp}{3}}$.
\label{ncorollary4.2}
\end{corollary}

\noindent {\bf Proof:}  (a) The lower bound for $p$ follows from the
fact that for all $j \leq \frac{h}{3K_2}$, 
$|\xi_j-z_j|< \|DT\|^j |\xi_0-z_0|<e^{-\frac{2h}{3}}
<<e^{-\beta \frac{h}{3K_2}}$. 
By Lemma \ref{nlemma4.11}, $p$ is the smallest $\mu$ such that
$\|w_\mu(0)\| \cdot |z_0 - \xi_0|^2 > K_1^{-1} e^{-\beta \mu}$. 
This must happen
for some $\mu \leq \frac{3h}{c}$ because 
$\|w_{\frac{3h}{c}}(z_0)\|\cdot  |z_0 - \xi_0|^2 
>K^{-\varepsilon \frac{3h}{c}} \|w^*_{\frac{3h}{c}}(z_0)\|\cdot  |z_0 - \xi_0|^2
> K^{-\varepsilon \frac{3h}{c}} e^{c\cdot\frac{3h}{c}} e^{-2h} > 1$.

(b) This follows from the fact that
$\|w_p(\xi_0)\| \approx \|w_p(z_0)\|$ (Lemma \ref{nlemma4.9}) and
$|z_0 - \xi_0| \cdot \parallel w_p(\xi_0) \parallel >
e^{-\frac{\beta}{2}p} \parallel w_p(\xi_0) \parallel^{\frac{1}{2}} >  
e^{-\frac{\beta}{2}p} e^{\frac{cp}{2}} >e^{\frac{cp}{3}}$.
 \hfill $\square$

\medskip
 
In analogy with Definition \ref{ndefinition3.3}, we define for
$\xi_0(s) \in \gamma$ the notion of a {\it fold period with
respect to $z_0$}.
This is the number $\ell$ such that $(Kb)^{\frac{\ell}{2}} \approx s$.
If $\tau_0(\xi_0)$, the unit tangent vector to $\gamma$ at $\xi_0$, is split
according to this definition, then the rejoining of the $E_i$-vector for
$\ell < i < p$ has negligible effect.
We may
assume also that as we iterate,
the sub-segment of $\gamma$ bound to $z_0$ acquires the same
fold periods as $z_i$, and think of these as {\it secondary
fold periods} for $\xi_i$ .

%%%%%%%Corollary 4.3 
\begin{corollary}
Let the assumptions and notation be as in Corollary \ref{ncorollary4.2}.
We let $p=p(\xi_0)$ where $|\xi_0-z_0|=e^{-h}$
and assume that $z_p$ is not in a fold period. Then
 
(a) the subsegment of $T^p\gamma$ between $\xi_p$ and $z_p$
contains a curve $\geq e^{-K\beta h}$ in length
 
\ \ \ \ \ and with $b$-horizontal tangent vectors;
 
(b) 
$$
\|\tau_p(\xi_0)\| \geq K^{-1}e^{h(1-\beta K)}.
$$
\label{ncorollary4.3}
\end{corollary}
 
\noindent {\bf Proof:} \ (a) By definition, $|\xi_p-z_p|>e^{-\beta p}$.
The part of $T^p\gamma$ in a fold period with respect to $z_0$
has length $\leq (Kb)^{\frac{p}{2}} \ \|DT\|^p$, and the rest have
$b$-horizontal tangent vectors.
To convert these estimates in $p$ into
bounds involving  $h$, use Corollary \ref{ncorollary4.2}(a).

(b) Splitting $\tau_0$ using $e_p$, we see that $\|w_p\|
\sim e^h\|\tau_p\|$. Combining this with Lemmas \ref{nlemma4.11}  and 
\ref{nlemma4.9}, we have
$e^h\|\tau_p(\xi_0)\| \sim \|w_p(\xi_0)\| \approx \|w_p(z_0)\|
>K^{-1}|\xi_p-z_p|e^{2h} \geq K^{-1} e^{-K\beta h}e^{2h}$.
 \hfill $\square$

\vspace{.3in}

%%%%%%%%%%%%%%%%%%%%Section 5 
\section{Pushing the Induction Forward}
\label{s5} 
The goal of this section is to define $\Delta_{3N}$ and to prove that 
(IA1)--(IA6) hold up to time $3N$ for parameters in $\Delta_{3N}$. 
The key to this inductive step
is the correct splitting of
the $w^*_i$-vectors at free returns (Proposition \ref{nproposition5.2}).
This is proved with the aid of another important fact, namely the control  
of points in 
$\partial R_k$ (Proposition \ref{nproposition5.1}). 
%%%%%%%%%%%%%Section 5.1
\subsection{Control of $\partial R_k, \ k\leq \theta N$}
\label{s5.1}
For $z \in \partial R_k$, let $\tau(z)$ denote a unit tangent vector to $\partial R_k$ at $z$.

%%%%Proposition 5.1 
\begin{proposition}
For every $\xi_0 \in \partial R_0$ and every $k \leq \theta N$, 
$(\xi_0, \tau_0)$ with $\tau_0=\tau(\xi_0)$ is controlled up to time $k$ 
by $\Gamma_k$.
\label{nproposition5.1}
\end{proposition}
\noindent {\bf Proof:} 
The proof proceeds by induction. The correctness of splitting of 
$\tau_0$ is evident. We assume all $(\xi_0, \tau_0)$
have been controlled up to time $k-1$, so that it makes sense to
speak of $\xi_k$ as being in a bound or free state. 
Suppose $\xi_k$ is bound to $z_i$ for some $z_0 \in \Gamma_{k-1}$.
Since $d_{\cal C}(z_i)>e^{-\alpha i}$, we have $z_i \in {\cal C}^{(j)} \setminus {\cal C}^{(j+1)}$ for some $j<<i \leq k$.
By Lemma \ref{nlemma4.2}, $\xi_k$ is h-related to $\Gamma_k$, and by 
Lemma \ref{nlemma4.5},
$\tau^*_k$ splits correctly, proving control at step $k$. 
Before proceeding to the free case, we state a
lemma of independent interest:

%%%%Lemma 5.1
\begin{lemma}
Let $\gamma$ be a subsegment of $\partial R_k$.
If all the points on $\gamma$ are free, then $\gamma$ is a $C^2(b)$-curve.
\label{nlemma5.1}
\end{lemma}

\noindent {\bf Proof:} \ 
That $\tau_k$ is a $b$-horizontal vector is an immediate
consequence of the splitting algorithm. As for curvature,
we appeal to  Lemma \ref{nlemma2.4} after using 
Lemma \ref{nlemma4.8} to establish 
that $\|\tau_k\|> K^{-K\theta(k-i)} \|\tau_i\|$ for all $i<k$. 
\hfill $\diamondsuit$

\bigskip

Returning to the proof of Proposition \ref{nproposition5.1}, let
$\xi_k$ be a free return, and let $\gamma$ be the maximal free subsegment 
of $\partial R_k$
containing $\xi_k$. Since the end points of $\gamma$ are in bound state, they cannot be in ${\cal C}^{(k-1)}$ as explained earlier.  This
leaves two possibilities for the relation between
$\gamma$ and ${\cal C}^{(k-1)}$.
 
\smallskip
{\it Case 1.} \ $\gamma$ passes through the entire length of
a component of ${\cal C}^{(k-1)}$. In this case we know from (IA1) that there is a critical point $z_0 \in \gamma$. 
To see that every $\xi' \in \gamma \cap {\cal C}^{(0)}$ is h-related
to $\Gamma_k$, start from $z_0$ and move away from it along $\gamma$.
Using the $C^2(b)$ property of $\gamma$, the structure of critical regions
(see (IA1)) and the fact that $\gamma \cap \partial R_i =
\emptyset \ \forall i<k$, we observe that after leaving $\partial Q^{(k)}(z_0)$ one gets into $Q^{(k-1)}(z_0)$, then $Q^{(k-2)}(z_0)$, 
and so on, with $d_{\cal C}(\xi') \geq \rho^i$ for $\xi' \in Q^{(i-1)}(z_0) \setminus Q^{(i)}(z_0)$. For the splitting
of $\tau(\xi')$, it follows from Lemma \ref{nlemma4.1} and the 
$C^2(b)$ property of $\gamma$ that for $\xi' \in \gamma \cap Q^{(i-1)}(z_0) \setminus Q^{(i)}(z_0)$,
$\angle(\tau(\xi'), \tau(\phi(\xi'))  \leq 
\angle(\tau(\xi'), \tau(z_0)) + \angle(\tau(z_0), \tau(\phi(\xi'))
<(Kb)|\xi'-z_0| +(Kb)^{\frac{i-1}{4}} <\varepsilon_0 d_{\cal C}(\xi')$.
 
\smallskip
{\it Case 2.} \ $\gamma$ does not intersect ${\cal C}^{(k-1)}$. Let $j<k$
be the largest integer such that $\gamma \cap {\cal C}^{(j-1)} \neq \emptyset$. Then there exists $z \in \gamma \cap (\hat Q^{(j)}\setminus
Q^{(j)})$ for some $Q^{(j)}$. Suppose for definiteness that $z$ lies
in the right component of $\hat Q^{(j)}\setminus Q^{(j)}$. Moving left
along $\gamma$ from $z$, we note that since $\gamma \cap Q^{(j)}=
\emptyset$, the left end point $\hat z$ of $\gamma$ must also be in the same 
component of $\hat Q^{(j)}\setminus Q^{(j)}$.  H-relatedness
and correct splitting are now proved 
as in Case 1 with $\hat z$ playing the role of
$z_0$. We know $\tau(\hat z)$ splits correctly because $\hat z$ is, by definition, in a bound state.
\hfill $\square$

%%%%%%%%%%%%%%%%%%%%%%Section 5.2
\subsection{Extending control of $\Gamma_{\theta N}$-orbits to time $3N$}
\label{s5.2} 
We continue to assume (IA1)--(IA6), which guarantee that
if $w_0={\tiny (\!\!\begin{array}{c} 0 \\ 1 \end{array}\!\!)}$, then for all 
$z_0 \in \Gamma_{\theta N}$,
$(z_0, w_0)$ is controlled
up to time $N$ by $\Gamma_{\theta N}$.
The next proposition plays a key role in the inductive process.
 
%%%%Proposition 5.2 
\begin{proposition} If $z_0 \in \Gamma_{\theta N}$ satisfies $d_{\cal C}(z_i)>e^{-\alpha i}$
for all $i \leq 3N$,
then $(z_0, w_0)$ is automatically controlled by $\Gamma_{\theta N}$
up to time $3N$.
\label{nproposition5.2}
\end{proposition}

\noindent {\bf Proof:} From the condition that 
$d_{\cal C}(z_i)>e^{-\alpha i}$, we have
that $z_0$ is h-related to 
$\Gamma_{\theta N}$ up to time $3N$ (see the remark following (IA2) 
in Sect. \ref{s3.3.2}), and that $p<K\alpha 3N << \frac{2}{3}N$.
It suffices therefore to prove the correct splitting property at free returns.
Proceeding inductively, we assume that $(z_0, w_0)$ is
controlled up to time $k-1$ for some $k$ with $N \leq k \leq 3N$,
and let $z_k$ be a free return.
Then either $z_k \in \hat Q^{(j)} \setminus Q^{(j)}$ for some
$j \leq \theta N$, or $z_k \in {\cal C}^{([\theta N])}$. In the latter case
we let $j=[\theta N]$ for purposes of the following arguments.
%%%%Claim 5.1 
\begin{claim}
\label{nclaim5.1}
There exists $j'$, $\frac{1}{3} j \leq j' < j$,
such that if
$$
\xi_0 = z_{k-j'} \ \ \ \ \ {\it and} \ \ \ \ \
u_0 = \frac{w_{k-j'}(z_0)}{\|w_{k-j'}(z_0)\|}\ ,
$$
then for $0 \leq s < j'$,
$$
\|DT^s(\xi_0) u_0\| \geq \|DT\|^{-s}.
$$
\end{claim} 

\noindent {\em Proof of Claim \ref{nclaim5.1}:}
We consider the graph 
${\cal G}$ of $i \mapsto \log\|w_i(z_0)\|$
for $k-j < i \leq k$. Let $L$ be the (infinite) line through
$(k,\ \log\|w_k\|)$ with slope $\log\|DT\|$. Then clearly,
all the points in ${\cal G}$ lie above $L$.
Let $P$ be the intersection of $L$ with the line $x=k - \frac{1}{3}j$.
We let $L$ be pivoted at $P$ and rotate it clockwise
until it hits some point in ${\cal G}$. (Draw a picture!)
Let $k-j'$ be the first coordinate of the first point hit.
Then $\frac{1}{3} j \leq j' < j$, and since all points in ${\cal G}$
lie above $L$, Claim \ref{nclaim5.1} is proved if we can show that
in its final position,
the slope of $L$ is $\geq -\log\|DT\|$.
This is true because $z_k$ being free, there must be some $j''$ with
$\frac{2}{3} j \leq j'' < j$ such that $z_{k-j''}$ is not in a
fold period, otherwise the bound period initiated at the same time
as this (very long) fold period would last beyond $z_k$.
By Lemma \ref{nlemma4.8}, $\|w_{k-j''}\| < \|w_k\|$.
Thus one cannot rotate $L$ to a slope  $< -\log\|DT\|$ without
first hitting the point $(k-j'', \ \log\|w_{k-j''}\|) \in {\cal G}$.
\hfill $\diamondsuit$

\medskip
 
Now by Lemma \ref{nlemma2.3}, there exists an integral curve $\gamma$ of the most contracted field of order $j'$ through $\xi_0$ having length ${\cal O}(1)$.
Since $\gamma$ follows roughly the direction of $e_1$, it 
has slope $>K^{-1}\delta$ outside of ${\cal C}^{(0)}$ and is roughly a
parabola inside ${\cal C}^{(0)}$ (Lemma \ref{nlemma2.9}). In both cases, 
$\gamma$
meets $\partial R_0$. 
Let $\xi_0^{\prime} \in \gamma \cap \partial R_0$. Then 
$$
| \xi_s - \xi_s^{\prime} | < (K^2b)^s
$$
for all $0 \leq s \leq j'$. Our next claim is made possible by
Proposition \ref{nproposition5.1}.
%%%%Claim 5.2 
\begin{claim}
$\xi_{j'}^{\prime}$ is a free return.
\label{nclaim5.2}
\end{claim} 
{\em Proof of Claim \ref{nclaim5.2}:} 
If not, then $\xi_{j'}^{\prime}$ would be 
bound to $\hat{z}$, a point on a critical orbit, 
and we would have $\xi_{j'}, \xi_{j'}^{\prime} \in \hat Q^{(i)}(\hat z)$
for some $i<<j' <j$ with $d_{\cal C}(\xi_{j'}) \approx d_{\cal C}(\xi_{j'}^{\prime}) \approx d_{\cal C}(\hat{z}) > e^{-\alpha j'}$.
This contradicts our assumption that $\xi_{j'} = z_k$ is in
$\hat Q^{(j)}$ or in ${\cal C}^{([\theta N])}$, for in either case,
$d_{\cal C}(z_k) < \rho^{j-1}$. \hfill $\diamondsuit$
%%%%Claim 5.3 
\begin{claim}
With $u_0$ as in Claim \ref{nclaim5.1}, let
$$
\tau_i = DT^i(\xi_0^{\prime}) \tau_0, \ \ \ \ \
u_i = DT^i(\xi_0) u_0,
$$
and let $\theta_i$ be the angle between $u_i$ and $\tau_i$. Then
$\theta_{j'} \leq b^{\frac{j'}{2}}$. 
\label{nclaim5.3}
\end{claim} 
{\em Proof of Claim \ref{nclaim5.3}:} 
Write $A = DT(\xi_{i-1})$ and
$A^{\prime} = DT(\xi_{i-1}^{\prime})$. Then
\begin{eqnarray*}
\theta_i =
\frac{\|\tau_i \times  u_i\|}{\|\tau_i\|\cdot \|u_i\|} 
&=& \frac{1}{\|\tau_i\|\cdot \|u_i\|}
\|A^{\prime} \tau_{i-1} \times A^{\prime} u_{i-1}
+ A^{\prime} \tau_{i-1} \times (A - A^{\prime}) u_{i-1}\| \\
&\leq & \frac{\| \tau_{i-1} \|}{\| \tau_i \|} \cdot
\frac{\|u_{i-1}\|}{\|u_i\|}
\cdot (|\det(A')| \theta_{i-1} + K |\xi_i-\xi_i'|) \\
&\leq & \frac{\| \tau_{i-1} \|}{\| \tau_i \|} \cdot
\frac{\|u_{i-1}\|}{\|u_i\|}
\cdot (b \theta_{i-1} + K (K^2b)^{i-1} ).
\end{eqnarray*}
Applying this relation for $\theta_i$ recursively, we obtain 
$$
\theta_{j'} < \left( \sum_{i=0}^{j'} \ \frac{\| \tau_i \|}{\| \tau_{j'} \|} 
\cdot
\frac{\|u_i\|}{\|u_{j'}\|} \right) (K^2 b)^{j'-1}.
$$
Since both $z_k$ and $\xi_{j'}$ are free returns, we may use Lemma 
\ref{nlemma4.8} 
to bound the sum in brackets by $\sum_i K^{4\theta (j'-i)}<2K^{4\theta j'}$, 
completing the proof of
Claim \ref{nclaim5.3}. \hfill $\diamondsuit$

\medskip
 
We are finally ready to prove that $w_k(z_0)$ splits correctly.
Recall that $\xi_{j'}=z_k \in \hat Q^{(j)}$ or $Q^{([\theta N])}$.
Since $| \xi_{j'} - \xi_{j'}^{\prime} | <(K^2b)^{j'}$,
$\xi_{j'}^{\prime} \in \partial R_{j'}$ and $j'<j$, we have $\xi_{j'}^{\prime}
\in \partial Q^{(j')}(z_k)$. By our inductive hypothesis, $\tau_{j'}(\xi'_0)$ splits correctly. Since 
$\angle(w_k(z_0), \tau(\xi'_{j'}))$ $\leq b^{\frac{j'}{2}}$ (Claim 
\ref{nclaim5.3}),
$\angle(\tau (\phi(\xi^{\prime}_{j'})), \tau (\phi(z_k)))= {\cal O}(b^{\frac{j'}{4}})$
and $|d_{\cal C}(\xi^{\prime}_{j'})- d_{\cal C}(z_k)|= {\cal O}(b^{\frac{j'}{4}})$  
(Lemma \ref{nlemma4.1}), it suffices to prove that $b^{\frac{j'}{4}}<<d_{\cal C}(z_k)^2$.
In the case where $z_k \in \hat Q^{(j)} \setminus Q^{(j)}$, this is trivial as
$d_{\cal C}(z_k) \sim \rho^j$. In the case where $z_k \in Q^{([\theta N])}$,
since $d_{\cal C}(z_k)>e^{-\alpha k}$, we have $d_{\cal C}(z_k)^2 > e^{-6\alpha N}$, which we may assume is $>> b^{\frac{1}{12}\theta N} \geq b^{\frac{1}{4}j'}$.
This completes the proof of Proposition \ref{nproposition5.2}. 
\hfill $\square$
%%%%%%%%%%%%%%%%%%%%%%%%%%%%Section 5.3
\subsection{Verification of (IA1)--(IA6) up to time $3N$}
\label{s5.3} 
 
\noindent {\bf Step 1} \ {\it Deletion of parameters.} \ 
We delete from $\Delta_N$ all $(a,b)$ for which there 
exists $z_0 \in \Gamma_{\theta N}$ and $i, \ N<i \leq 3N$, such
that 
$$
d_{\cal C}(z_i)<e^{-\alpha i} \ \ \ \ \ {\rm or} \ \ \ \ \  \|w_i^*(z_0)\| < e^{ci}.
$$
The set of remaining parameters is called $\Delta_{3N}$.
We do not claim in (IA1)--(IA6)
that $\Delta_{3N}$ has positive measure or even that
it is nonempty; this is discussed in Section \ref{s6}. 
Steps 2--5 below apply to $T=T_{a,b}$ for $(a,b) \in \Delta_{3N}$.

\bigskip
 
\noindent {\bf Step 2} \ {\it Updating of $\Gamma_{\theta N}$.} \ 
For each $z_0 \in \Gamma_{\theta N}$, since 
$\|w_i\|$ grows exponentially (Step 1 and Lemma \ref{nlemma4.6}), 
there exists a unique $z_0'$ on the component of $\partial {\cal C}^{(k)}$ containing $z_0$ that is
a critical point of order $3N$ (Lemma \ref{nlemma2.10}). 
Let $\Gamma'_{\theta N}$ be the set of these $z_0'$, 
i.e. $\Gamma^{\prime}_{\theta N}$ is a copy of 
$\Gamma_{\theta N}$ updated to order $3N$. 

\bigskip
 
\noindent {\bf Step 3} \ {\it Construction of $\Gamma_{3\theta N}$
and ${\cal C}^{(k)}$, $\theta N <k
\leq 3\theta N$.} \ 
We  establish control of $\partial R_k$ as in Sect. \ref{s5.1}, 
with one minor difference as explained in the next paragraph.
Assuming that all has been accomplished for $k-1$. Then $R_k$ meets each component $Q^{(k-1)}$ of ${\cal C}^{(k-1)}$ in at most a finite number of strips bounded by free, and hence $C^2(b)$, curves.
Let $\gamma$ be one of these curves.
By Lemma \ref{nlemma2.11}, there exists a critical point
$\hat{z}_0 \in \gamma$ of order $\hat{m} = min\{ 3N, -\log
d(z_0, \gamma)^{\frac{1}{2}} \}$ where
$z_0 \in \Gamma'_{\theta N}$ lies on the boundary of the
component $Q^{([\theta N])}$ containing $\gamma$. Since $d(z_0, \gamma)
= {\cal O}(b^{\frac{\theta N}{2}})$, we have, assuming $\theta$ is chosen with $e^{-3N} >K^{-N}>b^{\frac{\theta N}{4}}$, that $\hat z_0$ 
is of order $3N$.  

To continue, we need to set bindings for points
in $\partial R_k$. Technically, only $z_0 \in \Gamma_{\theta N}$ 
(and not the critical points on $\partial R_i, \ \theta N<i \leq k$)
can be used. This is of no concern to us for the following reason: for $k'$ with $k<k' \leq 3\theta N$, 
only those parts of $\partial R_{k'}$ that are free are involved
in the construction of ${\cal C}^{(k')}$; and for
$\xi_0 \in \partial R_k \cap {\cal C}^{([\theta n])}$, 
independent of which $z_0 \in Q^{([\theta n])}(\xi_0)$ 
we think of it as bound to, $\xi_i$ will remain in bound state
through time $3\theta N$ because $| \xi_i - z_i | \leq K^{3\theta N} \rho^{\theta N}
<< e^{-3\beta \theta N}$.
 
The newly constructed critical points in $\partial R_k, \ N < k
\leq 3N$, together with $\Gamma'_{\theta N}$ form $\Gamma_{3\theta N}$.
We have completed the verification of (IA1) up to time $3N$.
 
\bigskip
 
\noindent {\bf Step 4} \ {\it Updating the definitions of 
$d_{\cal C}(\cdot)$ and $\phi(\cdot)$}. \ 
Using $\Gamma_{3\theta N}$ and ${\cal C}^{(k)}$, $k \leq [3\theta N]$, 
we reset these definitions for $z \in {\cal C}^{([\theta N])}$
in accordance with Definition \ref{ndefinition3.2}. 
Since $|$old$\phi(z)-$new$\phi(z)|={\cal O}(b^{\frac{\theta N}{4}})$ 
and
$|\tau($old$\phi(z))-\tau($new$\phi(z))|={\cal O}(b^{\frac{\theta N}{4}})$
(Lemma \ref{nlemma4.1}), 
these changes have essentially
no effect on the correctness of splitting for points
with $d_{\cal C}(\cdot) >b^{\frac{3 \theta N}{20}}$. The 
relations in (IA5) are also not affected.

\bigskip
 
\noindent {\bf Step 5} \ {\it Verification of (IA2)--(IA6) for $i 
\leq 3N$}. \ This is carried out in 3 stages.
\begin{itemize}
\item[(1)] First we argue that for $z_0 \in \Gamma_{\theta N}$ 
(we really mean $\Gamma_{\theta N}$, not $\Gamma'_{\theta N}$),
(IA2)--(IA6) hold for $i \leq 3N$: (IA2) and (IA4) hold by design;
(IA3) is given by Proposition \ref{nproposition5.2}, and (IA5) and (IA6) 
are proved in Sect. \ref{s4.3} with $m=3N$.

\item[(2)] With the properties of $\Gamma_{\theta N}$ in (1) having
been established, we observe that continuing to use $\Gamma_{\theta N}$ as the 
source of control,  the material in Sects. \ref{s4.2} and \ref{s4.3} 
are now valid for times up to $min(m,3N)$.
 
\item[(3)] We are now ready to argue that (IA2)--(IA6) hold for all
$z_0' \in \Gamma_{3\theta N}$. For each $z_0' \in \Gamma_{3\theta N}$, 
whether it is in $\Gamma'_{\theta N}$ or of generation $> \theta N$, there exists $z_0 \in \Gamma_{\theta N}$ such that 
$|z_0'-z_0|={\cal O}(b^{\frac{\theta N}{4}})$. 
This implies, for $i \leq 3N$, that
$|z_i'-z_i|<b^{\frac{\theta N}{4}}\|DT\|^{3N}<<e^{-\beta 3N}$
provided $\theta$ is chosen so that $b^{\frac{\theta}{4}}\|DT\|^3< \frac{1}{2}e^{-\beta}$. 
(IA2) follows immediately from the corresponding condition for $z_0$.
Regarding $z'_0$ as bound to $z_0$ for at least $3N$ iterates,  
(IA3) and (IA4) follow from property (IA6) of $z_0$.
Finally, regarding $(z'_0, {\tiny (\!\!\begin{array}{c} 0 \\ 1 \end{array}\!\!)})$ 
as controlled by $\Gamma_{\theta N}$ up to time $3N$,
we obtain (IA5) and (IA6) from Lemmas \ref{nlemma4.9}-\ref{nlemma4.11}  and 
Corollary \ref{ncorollary4.2}.
\end{itemize} 

\bigskip
 
\noindent {\large \bf Conclusions from Sections \ref{s3}--\ref{s5}:} \ 
After letting $N$ go to infinity, we have defined for each $T=T_{a,b}$
with $(a, b) \in \Delta := \cap_N \Delta_N$ a set
${\cal C}$ given by ${\cal C}=\cap_{i \leq 0}{\cal C}^{(i)}$.
This is the {\it critical set} in Theorem \ref{theorem1}. 
Let $\Gamma$ be the set to which $\Gamma_{\theta N}$ converges
as $N \to \infty$.
An equivalent characterization of ${\cal C}$
is that it is the set of accumulation points of $\Gamma$.
Clearly, the properties that $d_{\cal C}(z_i) \geq e^{-\alpha i}$ and  $\|w_i\|$  grows exponentially are passed on to points 
in ${\cal C}$. We have thus completed the proof  of Theorem \ref{theorem1}
modulo the positivity of the measure of $\Delta$.

\vspace{.3in}
\section{Measure of Selected Parameters}
\label{s6}

In this section we fix $b>0$ and consider the $1$-parameter family
$a \mapsto T_{a,b}$. Let $\Delta_b=\{a: (a,b) \in \Delta \}$.
The Lebesgue measure of a set $A \subset {\mathbb R}$ is denoted by $|A|$.
More generally, we use $|\cdot|$ to denote the measure on curves induced
by arc length. The purpose of this section is to prove that
$|\Delta_b|>0$ for all sufficiently small $b>0$.

\subsection{Phase-space dynamics and curves of critical orbits}
\label{s6.1}

Assuming $\delta=e^{-\mu^*}$ for some $\mu^* \in {\mathbb Z}^+$, let ${\cal P}=\{I_{\mu j}\}$ be the partition of the interval $(-\delta, \delta)$ defined as follows: for $\mu\geq \mu^*$, let $I_\mu=(e^{-(\mu+1)}, e^{-\mu})$, 
and let each $I_\mu$ be further subdivided into $\mu^2$ 
subintervals of equal length called $I_{\mu j}, \ j=1,2, \cdots \mu^2$; 
for $\mu \leq -\mu^*$, let $I_{\mu j}=-I_{(-\mu) j}$. 

Next let $\gamma$ be a curve with nearly horizontal tangent vectors.
We assume for simplicity that $\gamma$ meets only one component $Q^{(0)}$ 
of ${\cal C}^{(0)}$, and let $\hat z =(\hat x, \hat y)$ be a point 
near the center of $Q^{(0)}$. The partition
 ${\cal P}_{\gamma, \hat z}$ on $\gamma$ is defined 
to be $(\psi^{-1}{\cal P})| \gamma \cup \{I^\pm\}$ where  
$\psi(x,y)= x-\hat x$ and $I^\pm$ are the two
components of $\gamma \setminus \psi^{-1}(-\delta, \delta)$.
An element of ${\cal P}_{\gamma, \hat z}$ is said to have ``full length"
if its image under $\psi$ is either equal to some $I_{\mu j}$ or 
longer than all the $I_{\mu j}$'s.
When $\gamma$ and  $\hat z$ are understood, we often refer to 
${\cal P}_{\gamma, \hat z}$ simply as ${\cal P}$ and
$(\psi^{-1}I_{\mu j}) \cap \gamma$ as  $I_{\mu j}$.

\medskip
Before proceeding to the estimation of $|\Delta_b|$, we consider first 
the following problem in phase-space dynamics. The estimation of 
$|\Delta_b|$ includes an argument parallel to and more
complicated than this.

\bigskip

\noindent {\bf A model problem in phase-space dynamics}

\medskip
Let $T=T_{a,b}$ with $(a,b) \in \Delta$. Recall from the proof of 
Proposition \ref{nproposition5.1} that if $\gamma \subset \partial R_k$ is a 
maximal free
segment meeting some $Q^{(0)}$, then either $\gamma \cap Q^{(0)}$
contains a critical point $\hat z \in \Gamma$ or the entire segment
$\gamma \cap Q^{(0)}$ is h-related to some $\hat z \in \Gamma$.
In both cases, ${\cal P}_{\gamma, \hat z}$ is the partition of choice
on $\gamma$. Note that for $z \in I_{\mu j}, \ d_{\cal C}(z)
\approx e^{-|\mu|}$.

Let $\omega_0$
be a subsegment of $\partial R_0$, and write
$\omega_i:=T^i\omega_0$. We assume that (i) for all $z_0 \in \omega_0,
\ d_{\cal C}(z_i)>e^{-\alpha i}$ for all $i \leq N$, 
(ii) each $\omega_i$, $i<N$, is contained in three consecutive $I_{\mu j}$,
and (iii) 
$\omega_N$ is free and is approximately equal to some 
$ I_{\mu_0 j_0}$. The problem is to find a lower estimate
for the measure of
$\{z_0 \in \omega_0: d_{\cal C}(z_i)>e^{-\alpha i}$ for all $i\}$. 

We may assume that all the points in $\omega_N$ have the same bound period, 
and let $i_1>N$ be the first moment in time after the expiration of this bound period
when $\omega_{i_1}\cap {\cal C}^{(0)}$ contains an $I_{\mu j}$ of
full length. This must happen at some point, for the length of $\omega_i$
grows by a factor $>K$ between successive free returns 
(Corollary \ref{ncorollary4.3}). 
 It is easy to check that  $d_{\cal C}>e^{-\alpha i}$ is not violated between 
times $N$ and $i_1$. Let $\{\omega\}$ be the partition 
${\cal P}$ on $\omega_{i_1}$ with end segments attached to
their neighbors if they are not of full length. We delete those 
$\omega$'s that contain some $z$ with $d_{\cal C}(z)<e^{-\alpha i_1}$. 
For each $\omega$ that is kept, 
we repeat the procedure above with $\omega$ in the place of
$\omega_N$, that is, we
iterate until $\omega$ makes a free return at time $i_2=i_2(\omega)$
with $T^{i_2-i_1}\omega$ containing an $I_{\mu j}$ of full length.
We then partition $T^{i_2-i_1}\omega$, discard subsegments that
violate $d_{\cal C}>e^{-\alpha i_2}$, and continue to iterate the rest. 

We estimate the fraction of $\omega_{i_1}$ deleted at time $i_1$ 
as follows.
Since $\omega_N \approx I_{\mu_0 j_0}$, the bound period $p$ is $\leq 
K |\mu_0|$. From Corollary \ref{ncorollary4.3}, we see that $\omega_{i_1}$ has length 
$>\frac{K^{-1}}{\mu_0^2}e^{-\beta K|\mu_0|} > e^{-2\beta |\mu_0| K}$.
Now $|\mu_0| \leq \alpha N$ and $i_1>N+p_0$ where $p_0>0$ is a lower
bound for all bound periods. Then
$$
\frac{|\{z \in \omega_{i_1}:d_{\cal C}(z)<e^{-\alpha i_1}\}|}
{|\omega_{i_1}|} \ < \ 
\frac{2e^{-\alpha (N+p_0)}}{e^{-2K\alpha \beta N}}
\ < \ e^{-\frac{1}{2}\alpha N}
$$
assuming $N$ is sufficiently large. Similarly, for each subsegment 
$\omega \approx I_{\mu j}$ of $\omega_{i_1}$ 
that is kept, the fraction of $T^{i_2-i_1}\omega$ deleted at time $i_2$ is 
$<e^{-\frac{1}{2}\alpha i_1}<
e^{-\frac{1}{2}\alpha (N+p_0)}$, and so on. 
To estimate the total measure of $\omega_0$ deleted, these fractions
have to be pulled back to $\omega_0$. This involves a {\it distortion estimate} for $DT^i$ along certain subsegments of $\partial R_k$.
Using the fact that this distortion is uniformly bounded 
(Lemma \ref{nlemma10.2}), we see that the 
fraction of $\omega_0$ deleted in this procedure is
$<K\sum_i e^{-\frac{1}{2}\alpha (N+ip_0)} < Ke^{-\frac{1}{2}\alpha N}$.

We remark that the scheme in this paragraph 
relies on the fact that $\omega_N$
has a certain minimum length depending on $N$, otherwise
the entire segment may be obliterated before time $i_1$ is reached.

\bigskip

\noindent {\bf Strategy for estimating $|\Delta_b|$}

\medskip
Since $b$ is fixed throughout this discussion, let us for notational 
simplicity omit mention of it and write $\Delta, \Delta_N$ and $T_a$ instead of $\Delta_b, \Delta_b \cap \Delta_N$ and $T_{a,b}$.
 Let $N$ be fixed. The problem is to estimate the measure of parameters
deleted between times $N$ and $3N$.
Our strategy is as follows: For $\hat a \in \Delta_N$ and $z_0 \in
\Gamma_{\theta N}(\hat a)$, let $a \mapsto z_0(a)$ be defined
on an interval containing ${\hat a}$. We consider
$$
\gamma_0 \ \to \ \gamma_1 \ \to \gamma_2 \ \to \ \cdots
\ \ \ \ \ \ \ \ {\rm where} \ \ \ \ \ \ \ \ 
\gamma_i(a) :=z_i(a)=T^i_a(z_0(a)),
$$
and estimate the measure of the set of $a$ for which
$z_i(a)$ violates (IA2) or (IA4).

\smallskip
The idea behind this line of proof is that {\it qualitatively, 
the evolution of $\gamma_0$ is similar to that of $\omega_0$} in the model
phase-space problem.
If this is true, then the measure deleted on account of (IA2) can be
estimated analogously.
To understand why the $\gamma_i$'s behave like phase curves, i.e.
curves that are obtained through the iteration of $T_a$, 
observe the way in which $\frac{d}{da}\gamma_i$, the tangent vector 
to the curve $a \mapsto \gamma_i(a)$, is transformed: 
if $\|\frac{d}{da}\gamma_i(a)\|>>1$, then $\frac{d}{da}\gamma_{i+1}(a)
\approx DT_a(\gamma_i(a)) \frac{d}{da}\gamma_i(a)$; that is to say,
$\gamma_{i+1} \approx T_a \circ \gamma_i$  near $\gamma_i(a)$. 

\bigskip
\noindent {\bf Issues to be addressed}

\medskip

\noindent 1. {\it Similarity of space- and $a$-derivatives.} \ 
This is the first and most important 
step in justifying the thinking in the last paragraph. 
Let $\gamma_0$ be as above. 
In Sect. \ref{s6.2}, we show that 
$\frac{d}{da} \gamma_i \sim DT^i
{\tiny \left(\!\!\begin{array}{c} 0 \\ 1 \end{array}\!\!\right)}$.
As we will see, this is made possible by our transversality condition 
on $\{f_a\}$ in Sect. \ref{s1.1}.
The only other prerequisite for this comparison is that
the slopes of $\gamma_0$ be suitably bounded.
This is verified in Sect. \ref{s6.3} 
for curves corresponding to critical points of all generations and
all orders.

\medskip

\noindent 2. {\it Dynamics of the curves $a \mapsto \gamma_i(a)$.} \ 
Our next step is to show that {\it as curves parametrized by $a$}, the $\gamma_i$ have properties similar to those of $\omega_i$.
For example, with $\Gamma_{\theta N}$ moving with $a$, 
how is $d_{\cal C}(z_i(a))$ affected? 
Other properties include the geometry of free segments,
quadratic behavior of the type in Sect. \ref{s4.3}, 
distortion estimates along $\gamma_i$ etc. 
These questions are discussed in Sect. \ref{s6.4}.

\medskip

\noindent 3. {\it Deletions of parameters in violation of (IA2) or (IA4).} \ 
We consider $z_0 \in \Gamma_{\theta N}$ one at a time, and let
$\gamma_0$ be the corresponding curve of critical points.
Assuming the success of the last step, deletions on $\gamma_0$
on account of (IA2) are estimated following the scheme
outlined in the model problem. Estimates for the measure of
parameters deleted on account of (IA4) are discussed in Sect. \ref{s6.5}.

\medskip
\noindent 4. {\it Combined effect of deletions corresponding to all
$z_0 \in \Gamma_{\theta N}$.} \ Obviously, we need to multiply
the measure of the parameters deleted on each $\gamma_0$
by the cardinality of $\Gamma_{\theta N}$, but there are technical 
considerations: As in our phase-space model, to get started
we need $\gamma_N$ to have a certain minimum length. This raises the 
question of the length of the parameter interval on which each 
$a \mapsto z_0(a)$ can be continued 
(this problem appears already in Sect. \ref{s6.3}). 
Also relevant is the combined effect of deletions on all
critical curves prior to time $N$. The final estimate is made 
in Sect. \ref{s6.6}.

\bigskip

The idea to relate parameter-space dynamics to phase-space dynamics is, of
course, not new. Two results on 1-dimensional maps are cited without proof and
used in this section: a transversality condition from \cite{TTY} is
used in Sect. \ref{s6.2} and a large deviation estimate from \cite{BC2} is
used in Sect. \ref{s6.5}.
%%%%%%%%%%%%%%%%%%%%%%%%%%%%%%%%%%%%%%%%%%%%%%%%%%%%%%%%%%%%

\subsection{Equivalence of space- and $a$-derivatives}
\label{s6.2}
The setting of this subsection is as follows: For fixed $b>0$, let 
${\hat a} \in \Delta_N$ for some $N$, and let  
$z_0=z_0(\hat a) \in \Gamma_{\theta N}({\hat a})$. Let $n \leq N$. Then $z_0$ 
obeys (IA2) and 
(IA4) and the conclusions of Lemmas \ref{nlemma4.6}--\ref{nlemma4.8}
up to time 
$n$. We assume also that
$z_0(\hat a)$ has a smooth continuation $a \mapsto z_0(a)$ 
to an $a$-interval containing $\hat a$.
Let $w_i=DT_{\hat a}^i(z_0(\hat a)){\tiny \left(\!\!\begin{array}{c} 0 \\ 1 \end{array}\!\!\right)}$ and
$\tau_i = \frac{dz_i}{da}(\hat a)$. The goal 
of this subsection is to compare $w_i$ and $\tau_i $.
Let $\tau_0=(\tau_{0,1}, \tau_{0,2})$.

%%%%%%%%%%%%%%%%%%%%%%Proposition 6.1
\begin{proposition}
\label{nproposition6.1}
Given $\bar \tau>0$, there exist constants 
$\lambda_2 > \lambda_1 > 0$ and a small $\varepsilon > 0$ 
such that the following holds:
If $(\hat a, b)$ is sufficiently
near $(a^*,0)$, $z_0(\hat a)$ is as above, 
$\|\tau_0\|<\bar \tau$ and $|\tau_{0,2}| < \varepsilon$,  
then for all $i \leq n$, 
$$
\lambda_1 \leq \frac{\| \tau_i\|}{\|w_i\|} \leq \lambda_2.
$$
\end{proposition}

We will show below that once we have $\| \tau_i\| \sim \|w_i\|$
for some $i$ with $\|w_i\|$ sufficiently large, then this relationship will hold
from there on. The estimate for the initial stretch is 
guaranteed by our transversality condition on the $1$-dimensional family
$\{f_a\}$. We recall a relevant result from $1$-dimension:

\medskip

Let $f$ and $\{f_a\}$ be as in Sect. \ref{s1.1}.
Let $x_0$ be a critical point of $f$, and let $p=f(x_0)$.
Since $f=f_{a^*}$, we write $x_0(a^*)=x_0, \ p(a^*)=p$, and
let $a \mapsto x_0(a)$ and $a \mapsto p(a)$ be the continuation of
$x_0$ and $p$
as defined in Sect. \ref{s1.1}. Let $x_k(a)=f_a^k(x_0(a))$.
We will use $(\cdot)'$ to denote differentiation with respect to $x$.

%%%%%%%%%%%%%%%%Lemma 6.1
\begin{lemma} {\rm ([TTY], Proposition VII.7)} 
\label{nlemma6.1}
As $k \to \infty$,
$$
Q_k(a^*) \ := \ \frac{\frac{dx_k}{da}(a^*)}{(f^{k-1})'(x_1(a^*))}
\ \ \to \ \ \lambda_0 := \ \frac{dx_1}{da}(a^*)-\frac{dp}{da}(a^*).
$$
\end{lemma}

The transversality condition in Sect. \ref{s1.1}, Step II, states that 
$\lambda_0
\neq 0$. We will also need the following technical lemma the proof of which 
is given in Appendix \ref{app-B.9}.

%%%%%%%%%%%%%%%%Lemma 6.2
\begin{lemma} 
\label{nlemma6.2}
There exist constants $K$ and $c^{\prime} > 0$ such that for every $0\leq s<i$, we have
$$
\| DT^{i-s}(z_s) \| \leq K e^{-c^{\prime} s} \|w_i\|.
$$
\end{lemma}

\medskip

\noindent {\bf Proof of Proposition \ref{nproposition6.1}:} \ Since
$$
\tau_i = DT(z_{i-1}) \tau_{i-1} + \psi(z_{i-1})
$$
where $\psi(z)=\frac{\partial (T_a z)}{\partial a}(\hat a)$, we have inductively
$$
\tau_i = DT^i(z_0) \tau_0 + \sum_{s=1}^i DT^{i-s}(z_s) \psi(z_{s-1}).
$$
The upper estimate for $\frac{\| \tau_i\|}{\|w_i\|}$ follows
from Lemma \ref{nlemma6.2} and the uniform boundedness of $\|\psi(\cdot)\|$:
$$
\frac{\| \tau_i\|}{\|w_i\|} \leq 
\frac{\|DT^i(z_0) \tau_0\|}{\|w_i\|} + \sum_{s=1}^i 
\frac{\|DT^{i-s}(z_s) \psi(z_{s-1})\|}{\|w_i\|}
$$
$$
< K\|\tau_0\| + K \sum_{s=1}^\infty e^{-c^{\prime} s} := \lambda_2.
$$
To obtain a lower bound for $\frac{\| \tau_i\|}{\|w_i\|}$, we pick $k_0$ 
large enough that $|Q_{k_0}(a^*)|>\frac{1}{2}|\lambda_0|$ where $Q_{k_0}$
and $\lambda_0$ are as in Lemma \ref{nlemma6.1}, and 
decompose $\tau_i$ into $\tau_i = I + II$ where
$$
I = DT^i(z_0) \tau_0 + \sum_{s=1}^{k_0}  DT^{i-s}(z_s) \psi(z_{s-1}),
$$
$$
II = \sum_{s=k_0+1}^i DT^{i-s}(z_s) \psi(z_{s-1}).
$$
Again by Lemma \ref{nlemma6.2}, we have
$$
\frac{\| II \|}{\|w_i\|} < \sum_{s=k_0+1}^{\infty} K e^{-c^{\prime} s}.
$$
We will show $\frac{\| I \|}{\|w_i\|}>K_0^{-1}|\lambda_0|$ for some $K_0$, 
and assume $k_0$ is chosen so that $\sum_{s>k_0} K e^{-c^{\prime} s}$
$<< K_0^{-1}|\lambda_0|$. Write
$$
I = DT^{i-k_0}(z_{k_0})V
$$
where
$$
V = DT^{k_0}(z_0) \tau_0 + \sum_{s=1}^{k_0} DT^{k_0-s}(z_s) \psi(z_{s-1}).
$$

%%%%%%%%%%%%%%%%%%Claim 6.1
\begin{claim} 
$$
\|V\|\ > \ \frac{1}{3} \ \frac{\|w_{k_0}\|}{\|w_1\|} \ |\lambda_0|,
$$ 
and the second component of $V$ tends to $0$ as $({\hat a},b) \to (a^*, 0)$.
\label{nclaim6.1}
\end{claim}

\noindent {\it Proof of Claim \ref{nclaim6.1}:} Let $z_0 \to (x_0,0)$ as  
$({\hat a},b) \to (a^*, 0)$. The two terms of $V$ are estimated as follows:

(i) $\|DT^{k_0}(z_0) \tau_0\|<K|\tau_{0,2}|$ for $(\hat a,b)$ sufficiently
near $(a^*, 0)$.
This is because $k_0$ is a system constant, and writing  
$T^{k_0}_{a^*, 0}=(T^1, T^2)$,
we have
$$
DT^{k_0}_{\hat a,b}(z_0)\tau_0 \ \to \ \left (\frac{\partial T^1}
{\partial x}(x_0,0)\tau_{0,1}
+\frac{\partial T^1}{\partial y}(x_0,0)\tau_{0,2}, \ 0 \right ) 
= \left (\frac{\partial T^1}{\partial y}(x_0,0)\tau_{0,2}, \ 0 
\right ).
$$

(ii) For $({\hat a},b)$ sufficiently near $(a^*, 0)$, $z_s$ stays out of ${\cal C}^{(0)}$ for $>k_0$ iterates, and
$$
\frac{\sum_{s=1}^{k_0} DT^{k_0-s}(z_s) \psi(z_{s-1})}
{\|w_{k_0}\|/\|w_1\|} \ \to \ 
\left (\frac{\sum_{s=1}^{k_0} (f^{k_0-s})'(x_s(a^*))
\frac{d}{da}(f_a(x_{s-1}))(a^*)}
{\pm (f^{k_0-1})'(x_1(a^*))},  \ 0 \right )  
$$
$$
= \ \left (\pm \sum_{s=1}^{k_0} \frac{\frac{d}{da}(f_a(x_{s-1}))(a^*)}
{(f^{s-1})'(x_1(a^*))}, \ 0 \right ),
$$
which by a simple computation is equal to $(\pm Q_{k_0}(a^*),\  0)$.
\hfill $\diamondsuit$
\medskip

Assuming further that the $x$-coordinate of $z_{k_0}$ is in a small neighborhood
of $x_{k_0}$ (which is bounded away from the critical set),
and that $z_s$ stays outside of ${\cal C}^{(0)}$ for $>>k_0$ iterates,
we have that the slope of $e_{i-s}(z_s)$ is bounded below by some
$K^{-1}$. This together with Claim \ref{nclaim6.1} gives
$$
\| DT^{i-k_0}(z_{k_0})V\|  
> K^{-1} \| DT^{i-k_0}(z_{k_0})w_{k_0}\|\ \frac{\|V\|}{\|w_{k_0}\|}
> K_0^{-1}\| w_i \| \ |\lambda_0|.
$$

\hfill $\square$

\bigskip
We will also need an estimate on the angle between $\tau_i$ and $w_i$, 
which we denote by $\theta_i$. The assumptions are as in 
Proposition \ref{nproposition6.1} .

%%%%%%%%%%%%%%Lemma 6.3
\begin{lemma}
If $z_i$ is a free return, then
$\theta_i < \frac{K}{\|\tau_i\|}.$
\label{nlemma6.3}
\end{lemma}

\noindent {\bf Proof:} 
$$
\mid \sin{\theta_i} \mid \ \leq \ \frac{1}{\|\tau_i\|}\left( 
\sum_{s=1}^i \frac{1}{\|w_i\|} \| w_i \times DT^{i-s}(z_s) \psi
(z_{s-1})\| + \frac{\|w_i \times DT^i(z_0)\tau_0\|}{\|w_i\|} \right)
$$
$$
\leq  \frac{1}{\|\tau_i\|}\left(
\sum_{s=1}^i \frac{\|w_s\|}{\|w_i\|} 
\ \left \|\frac{w_s}{\|w_s\|} \times \psi(z_{s-1}) \right \| \ b^{i-s} +
\frac{\|\tau_0\|}{\|w_i\|} \ b^i \right)
\leq \frac{K}{\|\tau_i\|}\sum_{s=0}^{\infty} b^s.
$$
The last inequality is valid if, for example, $\|w_s\| \leq \|w_i\|$
for all $s \leq i$, which is the case at free returns. \hfill $\square$

%%%%%%%%%%%%%%%%%%%%%%%%%%%%%%%%%%%%%%%%%%%%%%%%%%%%%%%%%%%%%%

%%%%%%%%%%%%%%%%%%%%%%%%%%%Section 6.3
\subsection{Initial data for critical curves}
\label{s6.3}
The goal of this subsection is to verify the conditions on $\tau_0$
in Proposition \ref{nproposition6.1} for critical curves of all generations 
and all orders. Our plan of proof is as follows:

\medskip

\noindent 1. We obtain information on the slopes of critical curves of
generation $i$ by comparing them to critical curves of generation $i-1$.
Following \cite{BC2}, this is done using a lemma of Hadamard,
which requires that the intervals of definition of the critical curves
be sufficiently long. We are thus led to the following question: on how
long of a parameter interval can one continue a critical curve with
reasonable properties?

\medskip
\noindent 2. As the order of a critical point tends to infinity,
the length of the parameter interval on which it is defined goes to zero.
This makes it necessary for us to prove our results in two steps,
to first work with critical points having orders commensurate with
their generations, and then to pass the bounds on to curves corresponding 
to higher orders. 

\subsubsection{Stability of critical regions}
\label{s6.3.1}
In Sections \ref{s3}--\ref{s5}, we construct for $N=N_0, 3N_0, 3^2N_0, \cdots$
a parameter set $\Delta_N$ such that for $a \in \Delta_N$, $\Gamma_
{\theta N}$ is well defined and consists of critical points of 
generation $\theta N$ and order $N$. 
Let us denote this set by $\Gamma_{\theta N, N}$. In the discussion to follow, 
it will be convenient to consider $\Gamma_{i, n}$ for arbitrary $i \leq n$. 
We define these sets
formally 
as follows:

\medskip

First we fix $a \in \Delta_N$, and define $\Gamma_{i,N}, \ \theta N <i
\leq N$, 
inductively by carrying out the steps in Section \ref{s5} in a slightly 
different order. Assuming that 
$\Gamma_{i-1,N}$ is defined and all the points in $\partial R_0$
are controlled for $i-1$ iterates, we define ${\cal C}^{(i)}$ and 
$\Gamma_{i,N}$. Immediately, we observe that
the newly constructed critical points are controlled by 
$\Gamma_{\theta N, N}$. In particular, they satisfy (IA2) and (IA4) (with possibly slightly weaker constants) and can be used for
binding. For free segments of $\partial R_i$ that lie in ${\cal C}^{(0)}$,
we may then set binding as in the proof of Proposition \ref{nproposition5.1}, 
and proceed to 
step $i+1$.

For $n$ with $N <n \leq 3N$, let $\Delta_n:=\{a \in \Delta_N $: (IA2) and
(IA4) are satisfied up to time $n$ for orbits from $\Gamma_{\theta N} \}$. 
A slight extension of the argument
above defines $\Gamma_{i,n}$ for all $a \in \Delta_n$ and
$i \leq n$. 

Finally, we introduce for each $n$ the parameter set $\tilde \Delta_n$,
which has the same definition as $\Delta_n$ except that in the definition of 
$\tilde \Delta_N$, $N = N_0, 3N_0, \cdots$, (IA2) and (IA4)
are replaced by $d_{\cal C}(z_j) >\frac{1}{2}e^{-\alpha j}$ and
$\|w^*_j\|>\frac{1}{2} c_0 e^{cj}$. One checks easily that all the
results in Sections \ref{s3}--\ref{s5} are valid under these slightly relaxed 
rules, 
as is the discussion in the last two paragraphs, so that 
$\Gamma_{i,n}$ is defined for all $a \in \tilde \Delta_n$ and
$i \leq n$.

\medskip

We remark before proceeding further that built into our definition of 
$\Gamma_{i,n}$ for $\frac{N}{3} < n < N$ is the 
property that $z_0 \in \Gamma_{i,n}$ 
has all the properties of $\tilde z_0 \in \Gamma_{\theta N, N}$ 
(except for the factor $\frac{1}{2}$) up to time $n$. In particular, 
Proposition \ref{nproposition6.1} applies to ${\hat a} \in \tilde \Delta_n$ and
$z_0 = z_0(\hat a) \in \Gamma_{i, n}$.

%%%%%%%%%%%%%%%%%%%%%%Definition 6.1 
\begin{definition} 
For $i \leq n$, an interval $J \subset \tilde \Delta_n$ and $\hat a
\in J$, we say $\Gamma_{i,n}(\hat a)$ has a {\bf smooth continuation} to
$J$ if there is a map $g:\Gamma_{i,n}(\hat a) \times J \to R_0$
such that 

- $g(\cdot, a)=\Gamma_{i,n}(a)$ for all $a$ and 

- for each $z \in \Gamma_{i,n}(\hat a)$, 
$a \mapsto g(z, a)$ is smooth.

\noindent Likewise one has the notion of the {\bf critical regions} 
${\cal C}^{(i)}$ {\bf deforming continuously} as $a$ ranges over $J$.
\label{ndefinition6.1}
\end{definition}

%%%%%%%%%%%%%%%%%%%%Lemma 6.4
\begin{lemma}
Let $\hat a \in \Delta_n$ and  
$J=[\hat a-\rho^{2n}, \hat a +\rho^{2n}]$. Then $J \subset \tilde 
\Delta_n$; moreover, 
$\Gamma_{n,n}(\hat a)$ has a smooth continuation to $J$, and ${\cal C}^{(i)}, 
\ i \leq n$, deform continuously on $J$.
\label{nlemma6.4}
\end{lemma}

The structual stability of the critical regions
comes from the fact that the components of ${\cal C}^{(i)}$ are stacked 
together in a very rigid
way, and their relations to the components of ${\cal C}^{(i-1)}$ are
equally rigid. As $a$ varies over $J$, the entire structure may move
up or down by amounts $>>b^{\frac{i}{2}}$, the maximum height of the 
components of ${\cal C}^{(i)}$, but it takes a relatively
large horizontal displacement to slide these components past each
other. A proof of Lemma \ref{nlemma6.4} is given in Appendix \ref{app-B.10}.

%%%%%%%%%%%%%%%%Section 6.3.2
\subsubsection{Comparing $\tau_0$-vectors for different critical curves}
\label{s6.3.2}

%%%%%%%%%%%%%%%%%%%%%lemma 6.5
\begin{lemma}
\label{nlemma6.5}
There exists $K$ such that the following holds 
for all $n$: Consider $\hat a \in \Delta_n$ and $J=[\hat a-\rho^{2n}, \hat a+\rho^{2n}]$. Let $z^{(n)}\in \Gamma_{n,n}(\hat a)$, $z^{(n-1)}\in \Gamma_{n-1,n-1}(\hat a) \cap Q^{(n-1)}(z^{(n)})$, and let
$z^{(n)}(a)$ and $z^{(n-1)}(a)$ be the continuations of $z^{(n)}$
and $z^{(n-1)}$ on $J$. Then 
$$
\|\frac{dz^{(n)}}{da}(a)-\frac{dz^{(n-1)}}{da}(a)\| \ < \ (Kb)^{\frac{n}{9}}.
$$
\end{lemma}

From this lemma it follows inductively that $\|\frac{dz^{(n)}}{da}-
\frac{dz^{(0)}}{da}\|<Kb^{\frac{1}{9}}$ where $z^{(0)}$ is a critical
point of generation $0$ and order $1$ lying 
in $Q^{(0)}(z^{(n)})$. Since there is only a finite number of
critical curves of generation $0$ and order $1$, and for them
$\tau_{0,2}=0$, Lemma \ref{nlemma6.5} proves that the hypotheses on
$\tau_0$ in Proposition \ref{nproposition6.1} are met for curves corresponding 
to all $z^{(n)} \in
\Gamma_{n,n}$. It remains to pass these properties to critical curves
of higher order.

%%%%%%%%%%%%%%%%%%%%%%lemma 6.6
\begin{lemma}
\label{nlemma6.6} 
Let $m>n$, $\hat a \in \Delta_m$, and let $z^m \in \Gamma_{n,m}(\hat a)$ be the updating of $z^n \in \Gamma_{n,n}
(\hat a)$ to order $m$. Then for all $a \in  
[\hat a-\rho^{2m}, \hat a+\rho^{2m}]$,
$$
\|\frac{dz^m}{da}(a) - \frac{dz^n}{da}(a)\| < (Kb)^{\frac {n}{4}}.
$$ 
\end{lemma}

Lemmas \ref{nlemma6.5} and \ref{nlemma6.6} are proved 
in Appendix \ref{app-B.10}.

%%%%%%%%%%%%%%%%%%%%%%%%%%%%%%%%%%%%%%%%%%%%%%%%%%%%%%%%%%%%%

%%%%%%%%%%%%%%%%%%%%%%Section 6.4
\subsection{Dynamics of critical curves}
\label{s6.4}

We fix a parameter interval $J$ and a critical point
$z_0$ which we assume can be smoothly continued to all of $J$.
As usual, let $\gamma_i(a)=z_i(a)$.
The purpose of this subsection is to make precise the parallel between
the dynamics of $\gamma_0 \to \gamma_1 \to \gamma_2 \to \cdots $ and the action of $T_a^i$ on $\partial R_0$. Let 
$\tau_i(a)=\frac{d \gamma_i}{da}(a)$.

%%%%Lemma 6.7
\begin{lemma}
There is a small number $k(\delta)>0$ such that for all $i >$ some
$i_0$, if $\gamma_i(a) \not \in {\cal C}^{(0)}$ and $|{\rm slope}
(\tau_i(a))|<k(\delta)$, then
(i)  $|{\rm slope}(\tau_{i+1}(a))|<k(\delta)$; 
 (ii) \ $\tau_{i+1}(a) \approx DT_a(\gamma_i(a))\tau_i(a)$.
Thus $\gamma_i$ with $|{\rm slope}(\tau_i)|<k(\delta)$ grows exponentially
in length as long as it stays outside of ${\cal C}^{(0)}$.
\label{nlemma6.7}
\end{lemma}

\noindent {\bf Proof:} (ii) is evident once $\|\tau_i\|$ is 
sufficiently large. By Proposition \ref{nproposition6.1}, this happens 
after some $i_0$. (i) is a consequence of (ii) and Lemma \ref{nlemma2.7}.
The exponential growth follows from Lemma \ref{nlemma2.8}. \hfill $\square$  

\medskip
We assume $(a,b)$ is sufficiently near $(a^*,0)$ that
$z_s$ remains outside of ${\cal C}^{(0)}$ for $>i_0$ iterates. 

Next we allow $\gamma_i$ to intersect ${\cal C}^{(0)}$.
For each fixed $a$, we have introduced in Sections \ref{s3}--\ref{s5} 
definitions of distance to the critical set, 
binding point, bound period, etc. 
To emphasize their
dependence on $a$, we write
$d_{{\cal C}(a)}(\cdot)$, $\phi_a(\cdot)$ and $p_a(\cdot)$
when referring to definitions that belong to the map $T_a$. 
Even for a fixed map, these quantities depend sensitively
on the location of the point in question; vertical displacements
of $z$, for example, may
dramatically change $\phi_a(z)$. In the ``dynamics" of critical curves,
the problem is all the more delicate, for not only does
 $z_i(a)$ move with $a$, the entire critical set moves as well.
The goal of the next few lemmas is to establish some viable notions of 
$d_{\cal C}(\cdot)$ and bound/free states that work in a coherent fashion
for all points in $\gamma_i$.

\smallskip
We assume for the rest of this subsection that 

\noindent (i) $J \subset \tilde \Delta_{K\alpha n}$, so that for each $a$ the binding structure is in place for points with 

\ \ $d_{{\cal C}(a)}(\cdot) >e^{-\alpha n}$;

\noindent (ii) $z_0$ obeys (IA2) and (IA4) up to time $n$, and 

\noindent (iii) all time indices are $\leq n$.
\smallskip

 In the next lemma, we let $| \cdot -\cdot|_h$ denote the horizontal distance between two points, 
and assume
for simplicity that $\gamma_i$ is contained in  one component 
of ${\cal C}^{(0)}$.
%%%%Lemma 6.8 
\begin{lemma}
\label{nlemma6.8}
Suppose $|{\rm slope}(\tau_i)| <k(\delta)$.
Then there exists $\bar z \in {\cal C}^{(0)}$ such that whenever
$d_{{\cal C}(a)}(\gamma_i(a))>\frac{1}{2}e^{-\alpha i}$,
$$
\left | |\gamma_i(a) - \bar z|_h - d_{{\cal C}(a)}(\gamma_i(a)) \right |
<Ke^{-ci}d_{{\cal C}(a)}(\gamma_i(a)).
$$
\end{lemma}

Thus we may put the partition ${\cal P}_{\gamma_i, \bar z}$ on $\gamma_i$
and define $d_{\cal C}(\cdot)=
|\cdot - \bar z|_h$ (the precise definition of $d_{\cal C}(\gamma_i(a))$
is irrelevant for $a$ with $d_{{\cal C}(a)}(\gamma_i(a))<\frac{1}{2}
e^{-\alpha i}$).

%%%%Lemma 6.9 
\begin{lemma} 
Let $\gamma_i$ be as above. We assume further that $z_i(a)$ is a free return for every $a$. Then for each $\omega_0=I_{\mu j} \in 
{\cal P}_{\gamma_i, \bar z}$ with $|\mu|<\alpha i$, there exists 
$\tilde p =\tilde p (\omega_0)< K|\mu|$
such that for all $a,\ a'$ with $z_i(a), z_i(a') \in \omega_0$,

(a) $|z_{i+j}(a)-z_{i+j}(a')| < e^{-\beta j}$ for $j \leq \tilde p$;

(b) $z_{i+\tilde p}$ is out of all fold periods,
$|{\rm slope} (\tau_{i+\tilde p})| <k(\delta)$ and 
$|\omega_{\tilde p}| \geq \frac{1}{\mu^2}e^{-\beta K|\mu|}$;

(c) $\|w_{i+\tilde p}\| > K^{-1}e^{\frac{\tilde p}{3}}\|w_i\|$, and
$\|\tau_{i+\tilde p}\| > K^{-1}e^{\frac{\tilde p}{3}}\|\tau_i\|$.
\label{nlemma6.9}
\end{lemma}

Lemma \ref{nlemma6.9} allows us to define a natural notion of bound/free states 
for the curves $\gamma_i$ that agrees essentially with the 
dynamical notion previously
defined for each $z_i(a)$.

%%%%Proposition 6.2 
\begin{proposition}
\label{nproposition6.2}
We assume the following hold for all $a \in J$
and $i \leq n$:

(i) for each $i$, the entire segment $\gamma_i$ is bound or free simultaneously, and $\gamma_i$ is  

\ \ \ \ contained in three contiguous $I_{\mu j}$'s at all free returns;

(ii) $\gamma_n$ is a free return.

\noindent Then there exists $K$ (independent of $\gamma_0$ or $n$) such that 
for all $a, a' \in J$,
$$
\frac{1}{K} \leq \frac{\| \tau_n(a) \|}{\| \tau_n(a^{\prime})\|} \leq K.
$$
\end{proposition}

Lemmas \ref{nlemma6.8} and \ref{nlemma6.9} are proved in 
Appendix \ref{app-B.11}.
Proposition \ref{nproposition6.2} is proved in Appendix \ref{app-B.12}.

%%%%%%%%%%%%%%%%%%%%%%%%%%%%%%%%%%%%%%%%%%%%%%%%%%%%%%%%%%%%%%%

%%%%%%%%%%%%%%%Section 6.5
\subsection{Deletions on account of (IA4)}
\label{s6.5}
Let $J \subset \tilde \Delta_{3K\alpha N}$, and let $z_0$ be a critical 
point with a smooth continuation on $J$.
We assume that for all $a \in J$, (IA2) and (IA4) hold up to time $N$,
and that $\gamma_N \approx I_{\mu j}$ 
for some $\mu$ with $|\mu| < \alpha N$. Let

\bigskip

\centerline{$E_{N, z_0} := \{a \in J: \exists n, \ N<n\leq 3N$, 
such that (IA2) is satisfied up to time $n$ }

\centerline{ and (IA4) is violated at time $n$\}.}

\bigskip
\noindent The set $E_{N, z_0}$ consists of parameters for which  $z_i$ has an 
abnormally high frequency of close returns between times $N$ and $3N$.

%%%%%%%%%%%%%%%%Proposition 6.3
\begin{proposition} 
\label{nproposition6.3}
Given $\varepsilon >0, \ \exists \delta_0
=\delta_0(\varepsilon)$ 
such that if $\delta<\delta_0$ and $(a,b)$ is sufficiently near $(a^*,0)$, then 
$$
|E_{N, z_0}|<e^{-\varepsilon n}|J|.
$$
\end{proposition}

A $1$-dimensional version of this result is proved in
\cite{BC2}, page 81-86. After the groundwork in Sect. \ref{s6.4}, 
the adaptation of
this result to our setting is straightforward.

%%%%%%%%%%%%%%%%%%%%%%%%%%%%%%%%%%%%%%%%%%%%%%%%%%%%%%%%%%%%%%

%%%%%%%%%%%%%%%%%%%%Section 6.6
\subsection{Estimating $|\Delta|$}
\label{s6.6}
The initial parameter set $\Delta_0$ is chosen as follows. 
Let $C=\{x_i\}$ be the critical set 
of $f$, and let $\delta_1$ be the minimum distance between $C$ and
$f^n x_i, \ n>0$. We assume that $\delta_1>>\delta$. 
Let $n_0$ be the number of iterates the critical orbits of $T_{a,b}$
are required to stay outside of ${\cal C}^{(0)}$.
We assume $n_0$ is as large as need be and prespecified.
Then there exists $\varepsilon>0$ such that for all $a \in
[a^*-\varepsilon, a^*+\varepsilon]$ and for all small enough
$b$, the first $n_0$ iterates of all the generation $0$ critical points 
stay $>\frac{\delta_1}{2}$ away from ${\cal C}^{(0)}$. 
We let $\Delta_0=[a^*-\varepsilon, a^*+\varepsilon]$, and let
 $b$ be fixed in the rest of the discussion.

We recall briefly the induction process: For $n=n_0, n_0+1, \cdots, 3N_0$
where $N_0=[\frac{1}{\theta}]$, we consider for each fixed $a \in \Delta_0$ critical orbits of generation $0$ and make deletions according to 
(IA2) and (IA4). By continuously updating $\Gamma_0$, we show in 
Sections \ref{s3}--\ref{s5} and \ref{s6.3.1} that orbits of $\Gamma_0$ can use 
their own 
histories
for binding, defining a sequence of shrinking parameter sets $\Delta_n
:=\{a \in \Delta_0:z_i(a)$ obeys (IA2) and (IA4) for all 
$z_0 \in \Gamma_0$ and $i \leq
n\}$. At time $3N_0$, $\Gamma_{3\theta N_0}$ is introduced
for $a \in \Delta_{3N_0}$, and deletions are made for orbits originating 
from $\Gamma_{3\theta N_0}$ up to time $3^2N_0$. 

For purposes of estimating the measure of parameters deleted in the first
$3N_0$ iterates, we consider one $\hat z_0 \in \Gamma_0$ at a time 
and estimate
the set of parameters discarded on account of $\hat z_0$ alone.
The argument for $3^iN_0<n \leq 3^{i+1}N_0, \ i \geq 1$,
is identical except it has to be made for a larger
set of critical points. For simplicity, let us first {\it decouple}
the situation, i.e. pretend(!) while considering $\hat z_0$ that no deletions are made for any other $z_0 \in \Gamma_0$. 

Let $\hat z_0$ be fixed in the next two paragraphs, and let $\hat \gamma_0$ be 
the curve $a \mapsto \hat z_0(a), \ a \in \Delta_0$. 
We ``iterate" $\hat \gamma_0$ until it gets into ${\cal C}^{(0)}$.
More precisely, let $i_0(\cdot)$ be the first time a point enters ${\cal C}^{(0)}$.
Then $i_0$ is a function on $\Delta_0$, and with $f_{a^*}$ being
a Misiurewicz map and $\Delta_0$ chosen as in the first paragraph, 
it is easy to  modify $i_0$ slightly so that either $\{i_0=n\} = \emptyset$
or $\hat \gamma_n|\{i_0=n\}$
contains $I_{\pm\mu^*}$, one of the outermost 
$I_\mu$.
(This is to ensure that part of $\hat \gamma_n|\{i_0=n\}$ will be retained.)
We partition $\hat \gamma_n|\{i_0=n\}$ into $I_{\mu j}$,
letting each element of this partition play the role of the start-up segment
$\omega_N$ in the model problem in Sect. \ref{s6.1}.
These segments are iterated independently, partitioned at free returns, 
and the deletion process  begins.
Care is taken to retain segments of full length
each time something is discarded.
For (IA2) we follow the procedure described in 
Sect. \ref{s6.1}. For (IA4), it is clear from the $1$-dimensional proof 
in \cite{BC2} that segments of full length are retained at each stage. 
The technical justifications for treating $\gamma_i$ as phase curves 
are given in Sect. \ref{s6.4}.  

The procedure in the last paragraph defines for each $n$
a set $\Delta_{\hat z_0, n}:=\{a \in \Delta_0: a$ is retained through
step $n\}$ and a partition ${\cal Q}_{\hat z_0, n}$ of 
$\Delta_{\hat z_0, n}$. Aside from $\{i_0>n\}$, the elements of 
${\cal Q}_{\hat z_0, n}$ are intervals $J$ such that $\hat \gamma_n|J$ 
in its last free return prior to time $n$ is a whole $I_{\mu j}$. 
From our estimates in Sect. \ref{s6.1}, Proposition \ref{nproposition6.2} 
(distortion estimate) and Proposition \ref{nproposition6.3} (large deviation), 
it follows that there exist $\alpha_1>0$ and $K>1$ such that
for all $k, \ n$, 
$$
|(\Delta_0 \setminus \Delta_{\hat z_0, n}) \cap \{i_0=k\}| \ 
\leq \ Ke^{-\alpha_1 k}
\ |\{i_0=k\}|.
$$
In particular, 
\begin{equation}
|\Delta_0 \setminus \Delta_{\hat z_0, n}| 
\ \leq \ Ke^{-\alpha_1 n_0} \ |\Delta_0|.
\label{nequation6.1}
\end{equation}

The discussion in the last two paragraphs applies to every 
$z_0 \in \Gamma_0$ -- under the same
erroneous assumption that deletions due to distinct critical orbits 
do not interfere with each other. We now remove this assumption:

Suppose all is well through time $n-1$. Let $J \in {\cal Q}_{\hat z_0, n-1}$
be such that all or part of $\hat \gamma_n|J$ makes a free return.
(There is no problem if $\hat \gamma_n|J$ is in the middle of
a free period or a bound period.) In order to continue, 
we need to know that the necessary binding structure exists for all
$a \in J$. Observe that $J$ is not necessarily contained in
 $\Delta_{n-1}:=\cap_{z_0 \in \Gamma_0} \Delta_{z_0, n-1}$. 
Indeed, it may have been
deleted in its entirety before time $n$ without $\hat z_0$ knowing 
about it. If that is the case,  we should not have been looking at it 
in the first place (and hence no parameter is deleted on account of 
$\hat z_0$ at this step -- or thereafter). 
If $J \cap \Delta_{n-1} \neq \emptyset$, we claim that $J \subset \tilde
\Delta_{K\alpha n}$, so that the necessary binding structure
for the part of $\hat \gamma_n|J$ to be retained is in place. 
The claim above is verified as follows: since $|\hat \gamma_n|<1$, 
we have, by (IA4) and Proposition 6.1, $|J|<\lambda_1^{-1}e^{-cn}$. 
With $\lambda_1^{-1}e^{-cn}<<\rho^{2K\alpha n}$, we are guaranteed 
by Lemma \ref{nlemma6.4} that $J \subset \tilde \Delta_{K\alpha n}$.
(Note that $J$ may not be contained in $\tilde \Delta_n$.)
Thus the estimate (\ref{nequation6.1}) remains valid even as we take 
into consideration
deletions due to other $z_0 \in \Gamma_0$.

The total measure deleted, therefore, is estimated by
$$
\sum_{i=0}^\infty |\Delta_{3^iN_0} \setminus \Delta_{3^{i+1}N_0}|
\ \leq \ {\rm card}(\Gamma_0)Ke^{-\alpha_1 n_0}|\Delta_0|
\ + \ \sum_{i=1}^\infty {\rm card}(\Gamma_{3^i\theta N_0})
Ke^{-\alpha_1 3^iN_0}|\Delta_0|.
$$

To estimate card$(\Gamma_{\theta N})$, let $I_1, \cdots, I_r$ be the monotone
intervals of $f$, and let
$K_0 = \max_i \{$ number  of $I_j$   counted 
with multiplicity $:I_j \cap f(I_i) \neq \emptyset \}$.

\begin{lemma} 
\label{nlemma6.10}
$$
{\rm card}(\Gamma_{\theta N}) < K_0^{\theta N}
$$
\end{lemma}

\noindent {\bf Proof} \ Partition $\partial R_k$ into segments 
by orbits of critical points
of generation $\leq k$. Then each segment has at most one free component, 
and each free component meets $\leq K_0$ of 
the monotone intervals, giving rise
to $\leq K_0$ new critical points. For more details, see Sect. \ref{s8.1}. 
\hfill $\square$

\bigskip
We conclude that the fraction of $\Delta_0$ deleted tends to $0$
as $n_0 \to \infty$ and $b \to 0$.

\newpage
%%%%%%%%%%%%%%%%%%%%Section 7
{\bf In the remainder of this paper, $T$ is assumed to be  $T_{a,b}$  where 
$(a,b)$ is a pair of ``good" parameters, i.e. $(a,b) \in \Delta$ where 
$\Delta$ is as in Theorem \ref{theorem1}.}

\section{Nonuniform Hyperbolic Behavior}
\label{s7}

Recall that $\Gamma$ is the set to which $\Gamma_{\theta N}$ converges
as $N \to \infty$. One of the properties guaranteed by parameter selection
is that orbits starting from $\Gamma$ have some hyperbolic behavior
(Theorem \ref{theorem1}(2)(ii)).
The purpose of this section is to show that this behavior
is passed on to a large set of points on the attractor and in the basin,
proving Theorem \ref{theorem2} except for the assertion in (1)(iii), the proof
of which we postpone to Sect. \ref{s9.add}.

%%%%%%%%%%%%%%%%%%%%%%%%%%%Section 7.1
\subsection{Control and hyperbolicity of non-critical orbits}
\label{s7.1}

We recapitulate the ideas developed in Sections \ref{s3}--\ref{s5} with a view 
toward
proving hyperbolicity for an arbitrary (non-critical) orbit. 
Given arbitrary $z_0 \in R_0$, we let 
$$
0 \ \leq \  n_1 \ <\ n_1+p_1 \ \leq \ n_2\ <\ n_2+p_2 \ \leq \ n_3
\ < \ \cdots
$$ 
be such that $z_{n_j} \in {\cal C}^{(0)}$ and is bound
to a suitable point in $\Gamma$,
$p_j$ is the ensuing bound period, 
and $n_{j+1}$ is the first return after $n_j+p_j$.
Then:

\begin{itemize}
\item[(1)] During its free periods, i.e. between times $n_j+p_j$ and $n_{j+1}$, the orbit is  outside of ${\cal C}^{(0)}$, where $DT^i$ is essentially 
uniformly hyperbolic (Lemma \ref{nlemma2.8}). 

\item[(2)] During its bound periods, i.e. between times $n_j$ and $n_j+p_j$,
$DT^i(z_{n_j})$ copies the derivative of its guiding orbit 
from $\Gamma$ (see (IA6)),
which has been guaranteed through parameter selection 
to have some form of hyperbolicity 
((IA4)).

\item[(3)] The concatenation of hyperbolic segments, however, need not result in a hyperbolic orbit, for the direction expanded at the end of 
one segment may be near the contractive direction of the next.
Indeed, this happens at times $n_j$, when there is a ``confusion" of 
stable and unstable directions, leading to a loss of hyperbolicity 
(see Sect. \ref{s3.1}).

\item[(4)] The properties that guarantee that hyperbolicity is preserved
through these concatenations are precisely the {\it h-relatedness} and
{\it correct splitting} properties at free returns. At time $n_j$,
the correct splitting of an expanded vector limits the magnitude of 
the loss (Lemma \ref{nlemma2.12} and Sect. \ref{s3.3.2}), while 
the h-relatedness
of $z_{n_j}$ to some $\hat z \in \Gamma$ guarantees that the ensuing bound 
period is long enough for this loss to be compensated (see (IA5)).
\end{itemize}

In particular,
if $z_0$ has a unit tangent vector $w_0$ such that
$(z_0, w_0)$ is {\it controlled} by $\Gamma$ for all $n \geq 0$
in the sense of Definition \ref{ndefinition4.2}, then
\begin{equation}
\limsup_{n \to \infty} \frac{1}{n} \log \|DT^n(z_0)w_0\| \geq c' >0
\label{nequation7.1}
\end{equation}
where $e^{c'}$ is a lower bound
of the growth rates of $b$-horizontal 
vectors outside of ${\cal C}^{(0)}$ and net derivative gains during
bound periods. Assuming that the rate of growth outside of ${\cal C}^{(0)}$
is $> e^{\frac{c}{3}}$ where $c$ is as in Theorem 1, we may take 
$c'=\frac{c}{3}$.
We remark that in general, the growth of $\|DT^n(z_0)w_0\|$ 
is not regular: without any assumptions on how close to $\Gamma$ 
the free returns are allowed to be,
i.e. without a condition in the spirit of (IA2), the loss of hyperbolicity
at time $n_j$ can be arbitrarily large; for example, the $\liminf$ in 
(\ref{nequation7.1})
can be negative.

\medskip
Recall that to establish control of $(z_0, w_0)$, it suffices to look at
free returns (Lemmas \ref{nlemma4.2} and \ref{nlemma4.5}). We record below a 
condition
at free returns that enables us to extend control 
through another bound-free cycle. Lemma \ref{nlemma7.1} plays a crucial role
in all the results in this section. First, we identify certain locations 
that are potentially problematic. For $k \geq 0$, let
$$
Z^{(k)}:=\{z \in {\cal C}^{(k)}: \ d_{\cal C}(z)<b^{\frac{k}{20}} \}.
$$

%%%%Lemma 7.1
\begin{lemma} 
Let $z_0$ be an arbitrary point in $R_0$, $w_0$ an arbitrary unit vector, 
and suppose that $(z_0, w_0)$ is controlled by $\Gamma$ up to time $k-1$.
Let $z_k$ be a free return. If $z_k \in {\cal C}^{(i)} \setminus
Z^{(i)}$ for some $ i< \frac{5}{4}k$, then $w_k$ splits correctly. 
\label{nlemma7.1}
\end{lemma}

\noindent {\bf Proof:} \ The proof of this lemma is virtually identical
to that of Proposition \ref{nproposition5.2}. Let $j= \min\{i, k\}$, so that $z_{k-j}$
makes sense. (The reason we allow $i$ to exceed $k$ has to do with the way
this lemma is used.) Claims \ref{nclaim5.1}-\ref{nclaim5.3} 
in Proposition \ref{nproposition5.2} continue to be
valid because the only requirement on $(z_0, w_0)$ is that the pair be 
controlled. Note that we have already controlled $\partial R_0$ and its 
tangent vectors for all times.
The proof here differs from that in Section \ref{s5} 
only at the end, where under present conditions we have
$$
b^{\frac{j'}{4}} \ \leq \ b^{\frac{j}{12}}\ \leq \ b^{\frac{1}{12}\frac{4}{5} i } \ 
<< \ b^{\frac{i}{20}} \ \leq \ d_{\cal C}(z_k).
$$ \hfill $\square$

%%%%%%%%%%%%%%%%%%%%%%%Section 7.2
\subsection{Typical derivative behavior in the basin}
\label{s7.2}
Let $m$ denote the 2-dimensional Lebesgue measure.
%%%%Proposition 7.1
\begin{proposition} 
\label{nproposition7.1} 
Assuming the additional regularity condition (**) in Sect. \ref{s1.2}, we have
$$
m\ \{z_0 \in R_0: z_k \in Z^{(k)}  \  {\rm infinitely \ often }  \} \ = 0.
$$
\end{proposition}

To prove this result, we need more refined estimates on the
width of $Q^{(k)}$ than that given in Lemma \ref{nlemma4.1}. 

%%%%%%%%%%%%%%%%%%%Lemma 7.2
\begin{lemma} 
\label{nlemma7.2}
There exists $K>0$ such that
if $Q^{(k)}$ is a component of ${\cal C}^{(k)}$, and $d_v$ is the vertical distance between the two
horizontal boundaries of $Q^{(k)}$ measured anywhere
along the length of $Q^{(k)}$, then
$$
(K^{-1}b)^{k+1} \ < \ d_v \ < \ (Kb)^{\frac{99}{100}k}.
$$
\end{lemma} 
\noindent {\bf Proof:} \ 
First we prove the lower bound, which relies heavily on the condition (**). 
Let $\omega_k$ be a vertical line segment joining two points in $\partial Q^{(k)}$. For $i<k$, let $\omega_i=T^{-k+i}\omega_k$. 
If $\omega_0$ connects the
two components of $\partial R_0$, then $d_v >(K^{-1}b)^k \cdot K^{-1}b$ since
by (**), $\|DT v \| \geq K^{-1} |\det(DT)| \geq K^{-1} K_1^{-1} b$ for every 
unit vector $v$.
If not, we will need to rule out
the possibility that $\omega_0$ may be extremely short.
Let $z_0, z_0' \in \omega_0 \cap \partial R_0$, and let $\gamma_0$ be the
shorter of the two segments of $\partial R_0$ between $z_0$ and $z_0'$. 
We consider $\gamma_i:=T^i\gamma_0$, and remember that points on
$\partial R_0$ together with their tangent vectors are controlled
(Proposition \ref{nproposition5.1}). Since $z_k$ and $z'_k$ are both free, 
and they
do not lie on a $C^2(b)$-curve, we conclude that a critical point
is created on $\gamma_i$ for some $i<k$. Let $i$ be the first time 
this happens. If $|z_i-z'_i|>\delta$,
then $|\omega_k|>\delta (K^{-1}b)^k$. If not, then both $z_i$ and 
$z'_i$ are in ${\cal C}^{(0)}$. Since both of their bound periods have 
expired by 
time $k$, 
it follows from (IA5) that $d_{\cal C}(z_i)$ and 
$d_{\cal C}(z'_i)$ are $>e^{-K(k-i)}$. 
We claim that $d_{\cal C}(z_i)+d_{\cal C}(z'_i)$ is approximately
the horizontal distance between these two points (see Lemma \ref{nlemma8.1}
 for more details). This gives
$|\omega_k|>2(e^{-K}K^{-1}b)^k$. 

For the upper estimate, we pick an arbitrary $z_k \in \partial Q^{(k)}$,
and borrow the argument in the proof of Claim \ref{nclaim5.1}  
with $j=k$, pivoting the line $L$ at $L \cap \{x=\frac{1}{100}k\}$
(instead of $L \cap \{x=k-\frac{1}{3}j\}$) as we rotate clockwise.
This gives $i_0$ with $0 \leq i_0 \leq \frac{1}{100}k$ such that
$\|DT^i(z_{i_0})\|>\|DT\|^{-100i}$. Iterating forward once if necessary
(and possibly losing a factor of $K^{-1}$
in the last estimate), we may assume that $z_{i_0} \not \in
{\cal C}^{(0)}$, so that it
lies on an integral curve $\gamma_0$ of $e_{k-i_0}$ which 
joins the two components of $\partial R_0$. Note that $\gamma_0$
meets $\partial R_0$ only at its end points.
Iterating forward, this curve brings in two segments of $\partial R_{k-i_0}$.
They must lie on the two horizontal boundaries of 
$Q^{(k-i_0)}(z_k)$ because $\gamma_{k-i_0}$ passes through $z_k$ and
intersects no other point of $\partial R_{k-i_0}$. This proves
that $d_v$ measured at $z_k$ has length at most that of $\gamma_{k-i_0}$, 
which by Lemma \ref{nlemma2.3} is 
$<(\|DT\|^{200}b)^{k-i_0}<(Kb)^{\frac{99}{100}k}$.
\hfill $\square$

\bigskip

\noindent {\bf Proof of Proposition \ref{nproposition7.1}}:  
By the Borel-Cantelli 
Lemma, it suffices
to show that $\sum_k m(T^{-k}Z^{(k)}) < \infty$. We estimate 
$m(T^{-k}Z^{(k)})$ by
\begin{eqnarray*}
m(T^{-k}Z^{(k)}) & = & \sum m(T^{-k}(Q^{(k)} \cap Z^{(k)})) \\ 
 & \leq & {\rm max}   
\frac{m(T^{-k}(Q^{(k)} \cap Z^{(k)}))}{m(T^{-k}Q^{(k)})}
\sum m(T^{-k}Q^{(k)})
\end{eqnarray*}
where the summations and maximum are taken over all components 
$Q^{(k)}$ of ${\cal C}^{(k)}$. Note also that $\sum m(T^{-k}Q^{(k)})<1$. 
Using Lemma \ref{nlemma7.2} and the regularity of $\det(DT)$ in (**), we obtain
\begin{eqnarray*}
\frac{m(T^{-k}(Q^{(k)} \cap Z^{(k)}))}{m(T^{-k}Q^{(k)})} & \leq & 
K^{2k} \cdot \frac{m(Q^{(k)} \cap Z^{(k)})}{m(Q^{(k)})} \\ 
& \leq & K^{2k} \cdot \frac{(Kb)^{\frac{99}{100}k} \cdot b^{\frac{1}{20}k}}
{(\frac{b}{K})^{k+1} \cdot \rho^k}
\leq K^{4k} \frac{1}{b} \cdot \frac{b^{\frac{1}{25}k}}
{\rho^k}
\end{eqnarray*}
which decreases geometrically in $k$ as desired.
\hfill $\square$

\bigskip

\noindent {\bf Proof of Theorem \ref{theorem2}(2):} \  Let $\xi_0 
\in R_0$. From the discussion in Sect. \ref{s7.1}, it follows that
$$
\limsup_{n \to \infty} \frac{1}{n} \log{\|DT^n(\xi_0)\|} \geq \frac{c}{3}
$$
holds if we are able to produce $k_0 > 0$ and a vector $w_0$ 
such that if $z_0=\xi_{k_0}$, then $(z_0, w_0)$ is 
controlled by $\Gamma$ for all $n \geq 0$. 
In light of Proposition \ref{nproposition7.1}, it suffices to
consider the following two cases.

\smallskip

{\it Case 1.} \  $\xi_k \not \in Z^{(k)}$ for all $k\geq 0$. 
We take $k_0=0$ and let $w_0={\tiny \left(\!\!\begin{array}{c} 1 \\ 0 \end{array}\!\!\right)}$ if $\xi_0 \not \in {\cal C}^{(0)}$, 
$w_0 ={\tiny \left(\!\!\begin{array}{c} 0 \\ 1 \end{array}\!\!\right)}$
if $\xi_0 \in {\cal C}^{(0)}$. We assume $(z_0, w_0)$ is controlled
up to time $k-1$, and let $z_k$ be a free return. 
The hypothesis of Lemma \ref{nlemma7.1} is verified at time $k$ as follows:
Let $j$ be the largest integer such that $z_k \in {\cal C}^{(j)}$. 
Then if $j \geq k$, $i=k$ meets the requirements of Lemma \ref{nlemma7.1}
since $\xi_k \not \in Z^{(k)}$; 
and if $j<k$, then $z_k$ must be in $\hat Q^{(j+1)}
\setminus Q^{(j+1)}$ for some $Q^{(j+1)}$ since it is in $R_k$,
and so we may take $i=j+1$.

\smallskip
{\it Case 2.} \ $\xi_{k_0} \in Z^{(k_0)}$ for some $k_0$ and $\xi_k \not \in Z^{(k)}$ for all $k>k_0$. Here we let $z_0=\xi_{k_0}$ and 
$w_0={\tiny \left(\!\!\begin{array}{c} 0 \\ 1 \end{array}\!\!\right)}$. 
There is a critical point $\hat z$ in $Q^{(k_0)}(z_0)$ to which $z_0$ is
bound for  $k_1$ iterates. Since $\|DT\|^{k_1}b^{\frac{k_0}{20}}
>e^{-\beta k_1}$, we have $k_1 \sim k_0 \theta^{-1} >> k_0$. During this period, we may regard  $(z_0, w_0)$ as controlled by $\Gamma$. 
For $k \geq k_1$, the situation is identical to that in Case 1 
except that $z_k \in R_{k+k_0}$ and we can only guarantee 
$z_k \not \in Z^{(k+k_0)}$. To verify the hypothesis of Lemma 7.1
for $z_k$, we proceed as above, distinguishing between the cases
$j \geq k+k_0$ and $j<k+k_0$ and noting that for $k \geq k_1$, $k+k_0<(1+K \theta)k$.
\hfill $\square$

\bigskip

\noindent {\bf Remark.} \ The results in this paper that use (**) 
in Sect. \ref{s1.2}
remain valid if (**) 
is replaced by

\medskip

${\rm (**)}^{\prime}$ \  
There exist $\eta \geq 1$ and $K_1, K_2 >0$ such
that for all $z \in R_0$,

\centerline{$K_1^{-1} b^{\eta} \leq |\det(DT)| \leq K_2 b^{\eta}$.}

\medskip

\noindent To Prove this, it suffices to check that Proposition 
\ref{nproposition7.1}
is valid under ${\rm (**)}^{\prime}$. Observe that the results in Sect. 
\ref{s2.1} are abstract, so that if $\|DT^i(z_0)\| \geq \kappa^i$ for all
$i \leq n$, then $\|DT^i e_n\| \leq (K b^{\eta} \kappa^{-2})^i$ for
all $i \leq n$. Using this and $\|DT v \| \geq K^{-1} b^{\eta}$ for
all $\| v \| = 1$, one checks easily that under ${\rm (**)}^{\prime}$, the 
conclusion of Lemma \ref{nlemma7.2} is valid if $b$ is replaced by 
$b^{\eta}$. Moreover, the number $\frac{99}{100}$ can be replaced
by $1-\varepsilon_0$ for any prespecified $\varepsilon_0 > 0$. Choosing 
$\varepsilon_0$ such that $\varepsilon_0 \eta < \frac{1}{20}$, we 
check that the proof of Proposition \ref{nproposition7.1} goes through
as is.

%%%%%%%%%%%%%%%Section7.3
\subsection{Uniform hyperbolicity away from $\cal C$}
\label{s7.3}

Recall that
$$
\Omega_{\varepsilon} = \{ z_0 \in \Omega: \ \ d_{\cal C}(z_n) \geq \varepsilon  
\ \ {\rm for \ all}  \ n \in {\mathbb Z} \}.
$$
The purpose of this subsection is to prove that $\Omega_{\varepsilon}$ is a uniformly hyperbolic 
invariant set \footnote{Technically, $z^i \to z$ does not imply
$d_{\cal C}(z^i) \to d_{\cal C}(z)$ when 
$z^i \not \in {\cal C}^{(k)}$ and $z \in {\cal C}^{(k)}$, but let us assume
$\Omega_{\varepsilon}$ is closed by taking its closure if necessary.}  
for every
$\varepsilon > 0$. This result together with the fact that the
strength of hyperbolicity deteriorates as $\varepsilon \to 0$
justifies our identification of $\cal C$ as the critical set and confirms that $d_{\cal C}(\cdot)$ is a valid notion of ``distance" to 
the critical set.
The approximation of $\Omega$ by $\Omega_\varepsilon$  is a concrete example of the use of uniformly
hyperbolic invariant sets to approximate systems that have (weak)
hyperbolic properties. See \cite{K} and \cite{P} for  
results in the same spirit.

Proofs of uniform hyperbolicity often rely on {\it a priori}
knowledge of invariant cones. In our setting, these cones are easily 
identified for $\Omega_{\varepsilon}$ with $\varepsilon>\sqrt b$; 
see Sect. \ref{s2.5}. 
As $\varepsilon \to 0$, the situation becomes considerably more
delicate: the stable and unstable directions at points in
$\Omega_\varepsilon$ become increasingly
confused, both ranging over nearly all possible directions within very
small neighborhoods. Our line of proof, which does not rely on 
{\it a priori} knowledge of cones,
can be formulated  as follows:

\smallskip
Let $g:X \to X$ be a self-map of a compact metric space, and let
$M:X \to GL(2, {\mathbb R})$ be a continuous map. For $x \in X$ and $n \geq 0$,
we define $M^{(n)}(x)=M(g^{n-1}x) \cdots M(gx) M(x)$ and $M^{(-n)}(x)=
M(g^{-n}x)^{-1} \cdots M(g^{-1}x)^{-1}$. It is clear what it means for the {\it cocycle} $(g, M^{(n)})$
to be {\it uniformly hyperbolic} (think of $g$ as a diffeomorphism
and $M(x)=Dg(x)$). Since the condition of interest to us is
projective in nature, we will state our result assuming that $M$ 
takes its values in $SL(2, {\mathbb R})$.

%%%%%%%%%%%%%Lemma 7.3
\begin{lemma} 
Let $(g, M^{(n)})$ be as above. If there exist $\lambda>1$
and $N \in {\mathbb Z}^+$ such that at each $x \in X$, there exists a unit 
vector
$v=v(x)$ such that
$$
\|M^{(n)}(x)v\| \leq \lambda^{-n} \ \ \ \ {\rm for \  all } \ \ n \geq N,
$$
then $(g, M^{(n)})$ is uniformly hyperbolic.
\label{nlemma7.3}
\end{lemma}
 
\noindent {\bf Proof:} \ Let $E^s(x)$ be the subspace spanned by
$v(x)$, and observe that $M(x)E^s(x)=E^s(gx)$: if not, then there
are two linearly independent vectors, $v_1 \in M(x)E^s(x)$
and $v_2=v(gx)$ such that both $\|M^{(n)}(gx)v_1\|$ and
$\|M^{(n)}(gx)v_2\|$ decrease exponentially as $n \to \infty$,
contradicting $M \in SL(2, {\mathbb R})$. The continuity of
$x \mapsto E^s(x)$ is proved similarly.

Using the uniform contraction of
$M^{(N)}$ on vectors in $E^s$ and the fact that $|\det(M)|$ $=1$, we choose 
$\delta_0>0$ such that for all $x \in X$ and
$w \not = 0 \in {\mathbb R}^2$, if $\angle(w, v(x))<\delta_0$, then
$\angle(M^{(N)}(x)w, v(g^Nx))>\frac{1}{2} \lambda^{2N} \angle(w, v(x))$.
Let $C^s(x)=\{w: \angle(w, v(x))< \delta_0\}$ and  
$C^u(x)={\mathbb R}^2 \setminus C^s(x)$. We claim that  
$E^u(x):= \cap_{n=1}^\infty M^{(nN)}(g^{-nN}x) C^u(g^{-nN}x)$ is a
1-dimensional subspace. This is true because from the angle separation between 
vectors in $E^s$ and $E^u$, it follows that for all $w \in E^u$, 
$\|M^{(-nN)}w\|$ decreases exponentially. 
The $M$-invariance of $E^u$ is 
checked easily.
\hfill $\square$

%%%%%%%%%%%%%%%%%%%%%Proposition 7.2
\begin{proposition} 
For every $\varepsilon>0, \ 
\Omega_\varepsilon$ is uniformly hyperbolic with
$$
\|DT^iu\| \ \geq \ K_\varepsilon^{-1}e^{c'i}
$$
for all $u \in E^u$. Here $K_\varepsilon$ is a constant depending on
$\varepsilon$, and $c'$ can be taken to be $\approx \frac{c}{3}$.
\label{nproposition7.2}
\end{proposition}

\noindent {\bf Proof:} \ We fix $\varepsilon$ and let $k_\varepsilon$ be the smallest integer $k$ such that $\varepsilon>b^{\frac{k}{20}}$.

%%%%%%%%%%Claim 7.1
\begin{claim}
\label{nclaim7.1}
For every $\xi_0 \in \Omega_\varepsilon$, there exists
$k(\xi_0) \leq 2k_\varepsilon$ and a unit vector $w_0$ such that if $z_0=\xi_{k(\xi_0)}$, then for all $i >0$,
$$
\|DT^i(z_0)w_0\| \ \geq \ e^{\frac{c}{3}i}  
b^{\frac{k_\varepsilon}{20}} 
K^{-\frac{k_\varepsilon}{10}}.
$$
\end{claim}

\noindent {\it Proof of Claim \ref{nclaim7.1}:} \ We consider separately the 
following cases:

{\it Case 1.} \ $\xi_i \not \in {\cal C}^{(0)}$
for all $i \leq k_{\varepsilon}$. In this case we let $k(\xi_0)=0$ and 
$w_0 ={\tiny \left(\!\!\begin{array}{c} 1 \\ 0 \end{array}\!\!\right)}$.

{\it Case 2.} \ $\xi_{i_0} \in {\cal C}^{(0)}$ for some 
$i_0 \leq k_{\varepsilon}$
and $\xi_{i_0+k} \not \in Z^{(k)}$ for all $k \geq 0$. We let
$k(\xi_0)=i_0$ and $w_0 ={\tiny \left(\!\!\begin{array}{c} 0 \\ 1 \end{array}\!\!\right)}$. 

{\it Case 3.} \ $\xi_{i_0} \in {\cal C}^{(0)}$ for some 
$i_0 \leq k_{\varepsilon}$
and $\xi_{i_0+k} \in Z^{(k)}$ for some $k \geq 0$.
We let $k$ be the last time this happens, and choose $k(\xi_0)=i_0+k$, 
$w_0 ={\tiny \left(\!\!\begin{array}{c} 0 \\ 1 \end{array}\!\!\right)}$. 
Note that $k(\xi_0) \leq 2 k_{\varepsilon}$.

\smallskip

In each of the three cases, we first show that $(z_0, w_0)$ is controlled by
$\Gamma$ for all $n \geq 0$. This is done by verifying inductively at free
returns the hypothesis of Lemma \ref{nlemma7.1}. The arguments are essentially 
the same as those for Theorem \ref{theorem2}(2).

From the control of $(z_0, w_0)$, it follows that at free returns, 
$\|w_i\| > e^{\frac{c}{3}i}$. Next we consider the drop in $\|w^*_i\|$ one step
later. This is given by $d_{\cal C}(z_i)$, which by the definition of 
$\Omega_{\varepsilon}$ is $\geq b^{\frac{k_{\varepsilon}}{20}}$. 
Further drops at bound returns are exponentially small.
For comparisons between $w^*_i$- and $w_i$- vectors, since the fold period 
$\ell$  initiated at time $i$ is $\leq \frac{k_{\varepsilon}}{10}$, we have, 
for $j < \ell$, 
$\|w_{i+j}\| \geq K^{-\frac{k_{\varepsilon}}{10}}\|w_{i+\ell}\|=
K^{-\frac{k_{\varepsilon}}{10}}\|w_{i+\ell}^*\|$.
\hfill $\diamondsuit$

\medskip

Let $z_0$ be as above.  
From Claim \ref{nclaim7.1},  
the fields of most 
contracted
directions of sufficiently high orders are defined at $z_0$, and their uniform 
contractive estimates are passed on to $e_\infty:= \lim_n e_n$ 
(see Corollary \ref{ncorollary2.1}).
Let $v(z_0)=e_\infty (z_0)$. For other $\xi_0 \in \Omega_\varepsilon$,
let $v(\xi_0)=DT^{-k(\xi_0)}(\xi_{k(\xi_0)})v(\xi_{k(\xi_0)})$.
Using the fact that $k(\xi_0)<2k_\varepsilon$ and letting 
$M(z)=\frac{1}{|\det DT(z)|^{1/2}} DT(z)$, 
we see that the conditions of Lemma \ref{nlemma7.3} are satisfied.
Uniform hyperbolicity follows. 

It remains to prove that a lower bound for $\|DT^i | E^u \|$ is as claimed.
In the argument above we have produced for each $\xi_0 \in \Omega_{\varepsilon}$
a vector $u_0$ uniformly bounded away from $E^s(\xi_0)$ such that $\|u_i\|
\geq K_{\varepsilon}^{-1} e^{\frac{c}{3}i}$. Since 
$\angle (u_n, E^u(\xi_n)) \to
0$ uniformly, we have $\|u_n\| \sim \|DT^n | E^u(\xi_0) \|$. The assertion 
in Theorem \ref{theorem2}(i) on periodic points is proved similarly.
\hfill $\square$

\bigskip
\noindent {\bf Proof of Theorem \ref{theorem2}(1)(ii):} \ We now prove
that the deterioration of hyperbolicity on $\Omega_\varepsilon$
as $\varepsilon \to 0$ is not only a possibility but a fact.
To do this, it suffices to produce a point $z \in \Omega_\varepsilon$
with the property that $\angle(E^u(z), E^s(z))<K\varepsilon$.
We can choose this point to be on the unstable manifold $W^u(\hat z)$
of any $\hat z  \in \Omega_\delta$. For $\xi_0 \in 
W^u(\hat z)$, let $\tau_0$ be its unit tangent vector to $W^u(\hat z)$,

%%%%%%%%%%%%%%%%%%%%%%Claim 7.2
\begin{claim} 
For all $\xi_0 \in 
W^u(\hat z), \ (\xi_0, \tau_0)$ is controlled by $\Gamma $ for all $n \geq 0$. 
\label{nclaim7.2}
\end{claim}

\noindent {\it Proof of Claim \ref{nclaim7.2}:} \ It suffices to prove the 
result for 
$\xi_0 \in W^u_{loc}(\hat z)$.
Suppose that $(\xi_0, \tau_0)$ is 
controlled up to time $k-1$, $\xi_k$ is a free return, and $\xi_k \in 
{\cal C}^{(j-1)} \setminus {\cal C}^{(j)}$ for some $j$. Since $\xi_{k-j} \in 
\Omega$, it follows that $\xi_k \in R_j$, so that 
$\xi_k \in {\hat Q}^{(j)} \setminus Q^{(j)}$ for some $Q^{(j)}$.
If $j \leq k$, then Lemma \ref{nlemma7.1} 
applies directly. If not, we let $z_0=\xi_{k-j}$ and apply
Lemma \ref{nlemma7.1} to the orbit of $(z_0, \tau_0(z_0))$.  \hfill $\diamondsuit$

\medskip

Let $\gamma=W^u_{\delta/2}(\hat z)$. We will show that there exists 
$z \in (T^n\gamma \cap \Omega_{\varepsilon})$ for some $n > 0$ such that 
$d_{\cal C}(z) < 2 \varepsilon$.  
As $\gamma$ is iterated,
it gets long and eventually meets the region $\{d_{\cal C}(\cdot)<
\varepsilon \}$. Let $n_0$ be the first time this happens, and let
$\omega_0 \subset T^{n_0} \gamma$ correspond to some $I_{\mu j}$ in the 
region  $\{\varepsilon \leq d_{\cal C}(\cdot)
\leq 2\varepsilon\}$. (See the beginning of Sect. \ref{s6.1} for notation.)
Note that $\omega_0$ is free. We set binding for $\omega_0$ and iterate until 
it becomes free again at time $n_1$. 
We then subdivide the image into segments corresponding to $I_{\mu j}$ 
(by which we include  pieces outside of ${\cal C}^{(0)}$),
and let $\omega_1$ be the longest of the divided subsegments. 
We iterate $\omega_1$ until it becomes
free again at time $n_2$. Then divide and 
choose $\omega_2$ to be the longest of the
subsegments etc. Let $z \in \cap_{i \geq 0} T^{-(n_i-n_0)}\omega_i$.
Using Corollary \ref{ncorollary4.3}, 
we verify that 
$\omega_i \cap \{d_{\cal C}(\cdot)
<\varepsilon \} = \emptyset$ for all $i \geq 0$, so that $z \in
\Omega_\varepsilon$.

It remains to estimate $\angle(E^u(z), E^s(z))$. First, since
$\tau(z)$ splits correctly, we have 
$\angle(E^u(z), \tau(\phi(z)))<\varepsilon_0
d_{\cal C}(z)<2 \varepsilon_0 \varepsilon$. Note that 
$\tau(\phi(z))=e_\infty(\phi(z))$ and $E^s(z)=e_\infty(z)$. 
We leave it as an easy exercise to show that 
$\|DT^n(z) \tau_0(z_0)\| \geq 1$
for all $n>0$ (use Claim \ref{nclaim7.2} and Corollary \ref{ncorollary4.3}),
so that at both
$z$ and $\phi(z)$, $\angle (e_n, e_\infty)={\cal O}(b^n)$.
Let $n$ be such that $\lambda^n \sim \varepsilon$ where $\lambda$ is as
in Lemma \ref{nlemma2.2}. Then $\angle(e_n(z), e_n(\phi(z)))<K\varepsilon$,
and ${\cal O}(b^n)<<\varepsilon$, proving $\angle(\tau(\phi(z)), E^s(z))
< K^{\prime} \varepsilon$. This completes the proof.
\hfill $\square$

\vspace{.3in}

%%%%%%%%%%%%%%%%%%%%Section 10
\section{Statistical Properties of SRB Measures}
\label{s10}
We follow \cite{Y3} and \cite{Y4}, which put forward a scheme for obtaining 
statistical information for general dynamical systems with some hyperbolic 
properties. In this approach, one constructs reference sets and studies
regular returns to these sets. Sufficient conditions 
in terms of return times are then given for various statistical properties.

In Sect. \ref{s10.1}, we indicate how this setup is arranged for the class of
attractors in question. For technical details on this construction, 
we refer the reader 
to \cite{BY2}, where a similar construction is carried out for the H\'enon 
maps. 
SRB measures and their statistical properties are
discussed in Sects. \ref{s10.2} and \ref{s10.4}. A feature of the present 
setting
is that depending on the transitivity properties of $T$, our attractor 
may admit multiple SRB measures.

Obviously, the method of \cite{Y3} and \cite{Y4}
gives information only on orbits that 
pass through the reference sets constructed. 
To complete the picture, we prove in Sect. 
\ref{s10.3} that all SRB measures 
are captured by our reference sets, and Lebesgue-almost every 
initial condition in the basin is accounted for.

\subsection{Positive-measure horseshoes with infinitely many branches and 
variable return times}
\label{s10.1}

In \cite{Y3}, a unified way of looking at nonuniformly hyperbolic systems is
proposed. This dynamical picture requires that one constructs a reference set 
and a return map with Markov properties. The purpose of this subsection is to 
recall this construction in the context of the maps under consideration, and
to give a summary of the facts needed in the discussion to follow.

\subsubsection{Construction of reference set}
\label{s10.1.1}

Let $\{x_1, \cdots, x_r\}$ 
be the set of critical points of $f$.
Our reference set $\Lambda$ is the disjoint union of $2r$ Cantor sets
$\Lambda^\pm_1, \cdots, \Lambda^\pm_r$ where $\Lambda^+_i$ and 
$\Lambda^-_i$ are located in the component of ${\cal C}^{(0)}$ containing
$(x_i, 0)$, one on each side of $(x_i,0)$.
We define $\Lambda^+_i$ (respectively $\Lambda^-_i$) by
specifying two transversal families of curves $\Gamma^{+,s}_i$ and 
$\Gamma^{+,u}_i$ and letting
$$
\Lambda^+_i=\{z \in \gamma^u \cap \gamma^s:
\gamma^u \in \Gamma^{+,u}_i, \ \gamma^s \in \Gamma^{+,s}_i\}. 
$$

The family $\Gamma^{+,s}_i$ (no relation to the critical set $\Gamma_i$ in 
Sections \ref{s3}--\ref{s6}) is defined as follows.
Let $\cal P$ be the partition in Sect. \ref{s6.1} centered at 
$(x_i, b) \in \partial R_0$. (To simplify notation, 
$\partial R_0$ in this section refers to the top boundary of $R_0$.)
Let $\omega_0 \subset \partial R_0$ be the outermost 
$I_{\mu j}$ on the right, and let $\omega_\infty=\{z_0 \in \omega_0: 
d_{\cal C}(z_n)>\delta e^{-\alpha n}$ for all $n \geq 0\}$. 
Letting $m_\gamma (\cdot)$ denote the measure on a curve $\gamma$
induced by arc length, it
 is proved in Sect. \ref{s6.1} that $m_{\omega_0}(\omega_\infty)>0$.
For every $z_0 \in \omega_\infty$, since 
$\|DT^i(z_0)\tau_0\| \geq \delta
e^{\frac{cn}{3}}$ for all $n \geq 0$ (use (IA5) and the definition of
$\omega_\infty$), there is a stable curve of every order passing through it. 
These curves converge to  a stable curve $\gamma^s(z_0)$
of infinite order (Sect. \ref{s2.1}). Moreover,  $\gamma^s(z_0)$ 
has slope $>K^{-1} \delta$ and  connects the two boundaries 
of $R_0$. We define $\Gamma^{+,s}_i:=\{\gamma^s(z_0): z_0 \in \omega_\infty\}$.

To define $\Gamma^{+,u}_i$, we first let
${\tilde \Gamma^{+,u}_i}$ be the set of all free segments $\gamma$
of $\partial R_n$, all $n\geq 0$, such that $\gamma$ is three times as long
as $\omega_0$ and has its midpoint vertically aligned with that of
 $\omega_0$. Let $\Gamma^{+,u}_i$ be the set of curves 
that are pointwise limits of sequences in $\tilde \Gamma^{+,u}_i$.
We remark that since the curves in $\tilde \Gamma^{+,u}_i$ are
$C^2(b)$, their slopes as functions in $x$ form an equicontinuous family.
This implies that the curves in $\Gamma^{+,u}_i$ are at least $C^{1+Lip}$,
and that the tangent vectors of curves in 
$\tilde \Gamma^{+,u}_i$ converge uniformly to the tangent vectors of
curves in $\Gamma^{+,u}_i$.

Recalling that $\Lambda^+_i$ and $\Lambda^-_i$ are the Cantor sets that
straddle $x_i$, we may, for convenience, choose $\Gamma^{-,s}_i$ and $\Gamma^{+,s}_i$
in such a way that their elements are paired, i.e. 
the $T$-image of each element in $\Gamma^{-,s}_i$ lies on a stable curve containing the $T$-image of an element of $\Gamma^{+,s}_i$, and {\it vice versa}.

This completes the construction of $\Lambda = \cup_{i=1}^r \Lambda^\pm_i$.
A similar construction is carried out for the
H\'enon maps in \cite{BY2}, Sects. 3.1-3.4. 

\subsubsection{Structure of return map}
\label{s10.1.2}

Next we define a return map $T^R:\Lambda \to \Lambda$ with the following 
properties: 
Topologically, $T^R:\Lambda \to \Lambda$ has the structure of an
infinite horseshoe. For simplicity of notation, we  write 
$\Lambda_i=\Lambda^+_i$ or $\Lambda^-_i$.
A set $X \subset \Lambda_i$ 
is called an {\it $s$-subset} of $\Lambda_i$
if there exists a subcollection of $\Gamma \subset \Gamma^{s}_i$ 
such that $X=\{z \in \gamma^s \cap \gamma^u: \gamma^s \in \Gamma, \ 
\gamma^u \in \Gamma^{u}_i\}$; {\it $u$-subsets} are defined similarly. 
If $X$ is an $s$-subset of $\Lambda_i$, we say {\it $X=\Lambda_i \
mod \ 0$} if $m_{\partial R_0}(\Lambda_i \triangle X)=0$.

\begin{lemma} \ There is a map $T^R:\Lambda \to \Lambda$ with the following 
properties: every $\Lambda_i$ has a collection of pairwise disjoint  
$s$-subsets $\{\Lambda_{i,j}\}_{j=1,2, \cdots }$  
with $\Lambda_i=\cup_j
\Lambda_{i,j} \ mod \ 0$ such that for each $j$, 

- $T^R|\Lambda_{i,j} = T^{n_{i,j}}|\Lambda_{i,j}$ for some $n_{i,j}
\in {\mathbb Z}^+$; 

- $T^R(\Lambda_{i,j})$ is a $u$-subset of $\Lambda_k$ for some $k=k(i,j)$.
\label{nlemma10.1}
\end{lemma}

We stress that the partition of $\Lambda_i$ into $\{\Lambda_{i, j} \}$ is an 
infinite one, and that the return times $n_{i, j}$ are not bounded. The 
{\bf return time function} $R: \Lambda \to  {\mathbb Z}^+$ is defined to be 
$R|\Lambda_{i, j} = n_{i, j}$. As we will see, the tail of this function, that
is, the distribution of its large values, plays a crucial role in determining
the statistical properties of the system. 
Note that $T^R$ is not
necessarily the first return map; we have settled for possibly larger
return times in favor of a Markov structure. Lemma \ref{nlemma10.1}
corresponds to 
Proposition A(1) in \cite{BY2}; its proof is given in Sects. 3.4 and 3.5
of \cite{BY2}.

\subsubsection{Two important analytic estimates}
\label{s10.1.3}

Technical estimates corresponding to (P1)-(P5) in \cite{Y3} or
Proposition B of \cite{BY2}
are needed. Referring the reader to Section 5 of \cite{BY2} for their precise
statements and proofs, we state below two of the most relevant facts.

\begin{lemma} {\bf (Distortion estimate for controlled segments)} \ 
\label{nlemma10.2}
There exists $K>0$ such
that the following holds: Let $\gamma_0$ be a curve and $\tau_0$ 
its unit tangent vectors. We assume that 

(i) \ for all $z_0 \in \gamma_0$, $(z_0, \tau_0)$ is controlled up to 
time $n-1$;

(ii) \ $\gamma_i$ is bound or free simultaneously for each $i$,
and $\gamma_i$ is contained in three  

\ \ \ \ contiguous $I_{\mu j}$ at all free returns;

(iii) \ $\gamma_n$ is a free return.

\noindent Then for all $z_0, z_0' \in \gamma_0$, 
$$
\frac{1}{K} \leq \frac{\|\tau_n(z_0)\|}{\|\tau_n(z_0')\|}
\leq K.
$$
\end{lemma}

The proof is similar to that of Proposition \ref{nproposition6.2}
(it is,  in fact,
a little simpler) and will be omitted. 
In the construction of $T^R:\Lambda \to \Lambda$, it is
important to arrange that $\gamma_0$, the shortest subsegment of $\partial R_0$ that spans $\Lambda_{i,j}$ in the $u$-direction, satisfies (ii)
above up to time $n_{i,j}$.

\medskip

Let $\cup \Gamma_i^s:= \cup \{z \in \gamma^s:\gamma^s \in \Gamma^s_i\}$.
If $\gamma$ and $\gamma'$ are curves transversal to  the elements of $\Gamma^s_i$ and intersecting them, we define 
$\psi: \gamma \cap (\cup \Gamma_i^s)
\to \gamma'$ by sliding along the curves in $\Gamma^s_i$, and 
say $\Gamma^s_i$ is {\it absolutely continuous} if for every pair of
$C^2(b)$-curves $\gamma$ and $\gamma'$ as above,
$\psi$ carries sets of $m_\gamma$-measure zero to sets of 
$m_{\gamma'}$-measure zero. Recall that if $\gamma$ is the subsegment of
$\partial R_0$ in $\Gamma^u_i$, then $m_\gamma(\gamma \cap (\cup \Gamma_i^s))>0$; in particular, the definition above is not vacuous.

\begin{lemma} {\bf (Absolute continuity of $\Gamma^s_i$)} \ 
\label{nlemma10.3}
$\Gamma^s_i$ is absolutely continuous with  
$$
\frac{1}{K} \ < \ \frac{d}{dm_{\gamma'}}\psi_*(m_\gamma|\cup \Gamma_i^s)
\ < \ K \ \ \ \ \ \ \ {\rm on} \ \gamma' \cap (\cup \Gamma^s_i).
$$
\end{lemma}

Except for one minor technical difference, the proof of Lemma \ref{nlemma10.3}
is identical to that of Sublemma 10 in Section 5 of \cite{BY2}: 
in the latter, the transversals
are taken to be curves in $\tilde \Gamma^u_i$, whereas we need them to be 
arbitrary $C^2(b)$ curves here. Clearly, it suffices to show that 
Sublemma 10 of \cite{BY2} is valid with $\gamma \in \tilde \Gamma^u_i$ 
and $\gamma'$ arbitrary $C^2(b)$, and for that we need distortion estimates for the $\tau_i$-vectors on certain subsegments of $\gamma'$ ($\omega'$
in the proof of Sublemma 10). We have them because
these subsegments are connected to subsegments of 
$\gamma$ by (temporary) stable curves, and the corresponding $\tau_i$-vectors are comparable.
 
\subsubsection{Tail of return times}
\label{s10.1.4}
Finally we state an estimate on which the statistical properties of $T$
depend crucially. Its proof is identical to that of Proposition A(4) 
in \cite{BY2}.

\begin{lemma} 
\label{nlemma10.4}
 \ There exists $K$ and $\theta_0<1$ such that for every 
$\Lambda_i$, 
$$
m_{\partial R_0}\{z \in \partial R_0 \cap \Lambda_i: R(z)>n\} < K\theta_0^n.
$$
\end{lemma}

%%%%%%%%%%%%%%%%%%%%%%%%%%%%
\subsection{SRB measures}
\label{s10.2}
\subsubsection{Construction of SRB measures}
\label{s10.2.1}
We describe below a recipe for constructing SRB measures using the 
reference sets $\{\Lambda^\pm_i\}$.
For the definition of SRB measures, see Sect. \ref{s1.4}. For more details
on the technical justification of the steps below, see \cite{Y3} or \cite{BY2},
Sect. 6.2. The construction consists of three steps.

\medskip

\noindent {\it Step 1. \ Construction of a $T^R$-invariant measure $\nu$ on 
$\cup \Lambda^\pm_k$ with absolutely continuous conditional measures on the 
leaves of $\Gamma^u:=\cup_k
\Gamma^{\pm,u}_k$.} \ 
We fix some $\Lambda_i=\Lambda^+_i$ or $\Lambda^-_i$, and let 
$m_0=m_{\partial R_0} \mid (\Lambda_i \cap \partial R_0)$. Let $\nu$ be an
accumulation point of the sequence of measures 
$$ 
\frac{1}{n} \sum_{j=0}^{n-1} (T^R)^j_*m_0, \ \ \ n=1,2, \cdots.
$$
Then $\nu$ is a $T^R$-invariant measure. By Lemma \ref{nlemma10.2}, the 
conditional measures of
$(T^R)^j_*m_0$ on the curves of $\tilde \Gamma^u:=\cup_k
\tilde \Gamma^{\pm,u}_k$ have uniformly bounded densities. From 
Lemma \ref{nlemma10.2} and the Markov property of $T^R$ 
(see Sect. \ref{s10.1.2}), it follows that for $\gamma \in \tilde \Gamma^u$, 
the conditional densities of $(T^R)^j_*m_0$ on $\gamma$ when restricted to 
$\gamma \cap (\cup \Gamma^s)$ are uniformly bounded away from 0. 
These properties are 
passed on to the conditional measures of $\nu$ on the leaves of $\Gamma^u$. 
(The curves in $\Gamma^u$ are pairwise disjoint except possibly
for a countable number of pairs; this is nothing more 
than a technical nuisance.) 

\medskip

\noindent {\it Step 2. \ Construction of a $T$-invariant probability measure $\mu$
given $\nu$.} \ It follows from the bounded densities of $\nu$, Lemma
\ref{nlemma10.3} and Lemma \ref{nlemma10.4} that 
$\int_\Lambda R d\nu_0 < \infty$. Let 
$$
\mu \ = \ \frac{1}{\int R d\nu_0} \ \sum_{j=0}^\infty T^j_*(\nu_0 \mid \{R>j\}).
$$
It is straightforward to check that $\mu$ is a $T$-invariant 
probability measure.

\medskip

\noindent {\it Step 3. \ Proof of SRB property.} Let $\mu$ be as in Step 2. First we check
that $T$ has a positive Lyapunov exponent $\mu$-a.e. 
At $z_0 \in (\cup \Gamma^s) \cap (\cup \tilde \Gamma^u)$, let
 $\tilde \tau(z_0)$ be a unit tangent vector to $\tilde \Gamma^u(z_0)$, 
the $\tilde \Gamma^u$-curve through $z_0$.
Just as on $\omega_\infty$, we have
$\|DT^n(z_0)\tilde \tau\| \geq \delta e^{\frac{cn}{3}}$ for all $n \geq 0$.
This uniform growth is passed on to the tangent vectors $\tau$
to $\Gamma^u$-curves at every $z \in \Lambda=\cup_k \Lambda^\pm_k$.
The existence of a positive Lyapunov exponent $\mu$-a.e. follows 
from the fact that the orbit of $\mu$-almost every point passes through 
$\Lambda$. General nonuniform hyperbolic theory (see e.g. \cite{P} or \cite{R4}) then tells us that stable and unstable manifolds exist $\mu$-a.e. 

To prove that $\mu$ is an SRB measure, we need to show that
its conditional measures  on unstable manifolds are absolutely continuous.
Since $\mu$ is the sum of forward images of $\nu$, it suffices
to prove this for $\nu$. We know from Step 1 that $\nu$ has
absolutely continuous conditional measures on the leaves of $\Gamma^u$.
Thus it remains to prove

\begin{claim} 
\label{nclaim10.1}
For $\nu$-a.e. $z_0$, $\Gamma^u(z_0)$ is a local unstable manifold, i.e. 
$$
\limsup_{n \to \infty} \frac{1}{n}\log  \sup_{\xi_0 \in \Gamma^u(z_0)} |\xi_{-n}-z_{-n}| \ < \ 0.
$$
\end{claim}

\noindent {\it Proof of Claim \ref{nclaim10.1}:} \  From the construction of 
$\nu$,
it follows that that for $\nu$-a.e. $z_0 \in \Lambda$, there is a 
sequence of $\tilde \Gamma^u$-curves $\{\tilde \gamma_i\}$ such that
$\tilde \gamma_i \to \Gamma^u(z_0)$.  Let $n_i$ be such that 
$T^{-n_i}\tilde \gamma_i \subset \partial R_0$. Since $\tilde \gamma_i$
is free, we have that for all tangent vectors $\tilde \tau$ of $\tilde \gamma$, $\|DT^{-n}\tilde \tau\| <e^{-c''n}$ for some $c''>0$ and
$0<n \leq n_i$ (Proposition \ref{nproposition5.1} and Lemma \ref{nlemma4.8}). 
These uniform estimates
for backward iterates of $T$ are passed on to all tangent vectors of
$\Gamma^u(z_0)$, proving that it is a local unstable manifold of $z_0$.
\hfill $\diamondsuit$

\newpage
\subsubsection{Ergodic decomposition of SRB measures}
\label{s10.2.2}
We begin by considering the ergodic decompositions of the $T^R$-invariant 
measures constructed in Step 1 in Sect. \ref{s10.2.1}.

\begin{definition}
\label{ndefinition10.1} 
Let $g: X \to X$ be a continuous map of a compact 
metric space, and let $\nu$ be a $g$-invariant Borel probability
measure on $X$. We say $z \in X$ is {\bf generic} or {\bf future-generic}
with respect to $\nu$ if for every continuous function $\varphi: X
\to {\mathbb R}$,
$$
\frac{1}{n}\sum_{i=0}^{n-1} \varphi (z_i) \to \int \varphi d\nu.
$$
\end{definition}

Let ${\cal M}(T^R)$ be the set of all normalized invariant measures
constructed in Step 1 of Sect. \ref{s10.2.1}. Let $\nu \in {\cal M}(T^R)$, and
suppose that $\nu(\Lambda_i) >0$ for some $\Lambda_i = \Lambda_i^+$
or $\Lambda_i^-$. From the positivity of the conditional densities of $\nu$ on
$\Lambda_i \cap \gamma$, Lemma \ref{nlemma10.3}, and a standard argument due 
to Hopf, we know that there is an ergodic component $\nu^e$ of
$\nu$ such that

(i) $\nu$-a.e. $z \in \Lambda_i$ is generic with respect to $\nu^e$;

(ii) for every $C^2(b)$-curve $\gamma$, $m_{\gamma}$-a.e. 
\ $z \in \gamma \cap (\cup \Gamma_i^s)$ is generic with respect to $\nu^e$.

\noindent We abbreviate this by saying $\nu^e$ ``occupies" $\Lambda_i$.

Let ${\cal M}_e(T^R)$ denote the set of normalized ergodic components of 
measures in ${\cal M}(T^R)$. Then each $\Lambda_i^+$ (resp. $\Lambda^-_i$) 
is occupied by an elememt of ${\cal M}_e(T^R)$. Because the
stable curves of $\Lambda_i^+$ and $\Lambda_i^-$ are joined, 
$\Lambda_i^-$ and $\Lambda_i^+$ are in fact occupied by the same 
element of ${\cal M}_e(T^R)$. Thus the cardinality of ${\cal M}_e(T^R)$
is $\leq r$.

To further study the structure of ${\cal M}_e(T^R)$ we borrow some
ideas from finite state Markov chains. Let $\Lambda_i^{\pm}: = \Lambda_i^+
\cup \Lambda_i^-$. We think of each the sets 
$\Lambda_i^{\pm}, i = 1, \cdots, r$, as a state, and write $``i \to j"$
if $T^R(\Lambda_i^{\pm}) \cap \Lambda_j^{\pm} \not = \emptyset$. We say $i$ is
{\it transient} if there exists $j$ such that there is a chain 
$i \to \cdots \to j$
but no chain with $j \to \cdots \to i$. Non-transient states are called
{\it recurrent}. The following are consequences of simple facts about directed 
graphs.

\begin{itemize}
\item[(a)] The set of recurrent states is partitioned into equivalence classes
where $i \sim j$ if there is a chain $i \to \cdots \to j$. On the union of the 
$\Lambda_i^{\pm}$ corresponding to the states in each equivalence class is
supported exactly one element of ${\cal M}_e(T^R)$, which occupies each of 
these $\Lambda_i^{\pm}$. 

\item[(b)] If $i$ is transient, then clearly $\nu(\Lambda_i^{\pm}) = 0$ for
every $\nu \in {\cal M}_e(T^R)$. The following claim is a consequence of
the structure of $T^R$ (Lemma \ref{nlemma10.1}) and the fact that for every
transient state $j$, there exists a recurrent $k$ such that 
$j \to \cdots k$.

\begin{claim}
$\Lambda^{\pm}_i$ is the mod 0 union of a collection of pairwise disjoint 
$s$-subsets $\{ \hat \Lambda_{i, \ell} \}_{\ell = 1, 2, \cdots}$ with
the property that for each $\ell$, there exists $n_{\ell} > 0$ such that 
$(T^R)^{n_{\ell}} {\hat \Lambda}_{i, \ell}$ is a $u$-subset of 
some recurrent state.
\label{chain}
\end{claim}
\end{itemize}

The discussion in Sects. \ref{s10.2.1} and \ref{s10.2.2} are summarized as
follows:

\begin{proposition}
\label{nproposition10.1}
Let $r$ be the number of critical points of $f$.
Then there exist ergodic SRB measures
$$
\mu_1,\ \mu_2, \ \cdots, \ \mu_{r'}, \ \ \ \ 1 \leq r' \leq r, 
$$
such that for every $C^2(b)$-curve $\gamma$, 
$m_\gamma$-a.e. $z \in \gamma \cap (\cup \Gamma^s)$
is generic with respect to some $\mu_i$.
\end{proposition}
\noindent {\bf Proof:} Let 
${\cal M}_e(T^R) = \{\nu_1, \nu_2, \cdots, \nu_{r^{\prime \prime}} \}$. Then 
the $\mu_i$ in this proposition are saturations of the $\nu_j \in
{\cal M}_e(T^R)$ in the sense of Step 2 in Sect. \ref{s10.2.1}. Clearly 
$r^{\prime} \leq r^{\prime \prime} \leq r$; it may happen that $r^{\prime}
\leq r^{\prime \prime}$ because the saturations of distinct $T^R$-invariant 
measures may merge. The genericity assertion is proved as follows. If $k$ is
a recurrent state, then it is occupied by some $\nu_j$, and hence 
$m_{\gamma}$-a.e. $z \in \gamma\cap(\cup\Gamma_k^s)$ is generic with
respect to some $\mu_i$. Via Claim \ref{chain}, the same conclusion holds if 
$k$ is a transient state.

\hfill $\square$

%%%%%%%%%%%%%%%%%%%%%%%%%%%%%%%%%%%%%%%%%%%%%%%%%%%%%

\subsection{Accounting for almost every initial condition}
\label{s10.3}
It is a fact from general nonuniform hyperbolic theory that
if an SRB measure has nonzero Lyapunov exponents,
then the set of points generic with respect to it has positive
Lebesgue measure. This is a consequence of the absolute continuity of
stable foliations \cite{PS}. In general, the set of generic points 
may not have full Lebesgue measure in any neighborhood
of the attractor. 

Let $m$ denote the Lebesgue measure on $R_0$.

\begin{proposition} 
\label{nproposition10.2}
Let $\{\mu_i\}$ be the ergodic SRB measures in 
Proposition \ref{nproposition10.1}. Then for $m$-a.e. $z_0 \in R_0$, 
there exists 
$\mu_i$ with respect to which $z_0$ is generic, and $z_n \in \cup \Gamma^s$
for some $n>0$.
\end{proposition}

It follows that $\{ \mu_i \}$ is the set of all the ergodic SRB measures 
that $T$ admits.
Propositions \ref{nproposition10.1} and \ref{nproposition10.2}
together comprise the proof of Theorem \ref{theorem6}.

\bigskip
\noindent {\bf Proof:} \ 
Let $B$ be the set of points {\it not} generic with respect to any of
the $\mu_i$. We remark that $B$ is a Borel measurable set, for genericity
with respect to a given measure is determined by a countable number of
test functions. Let $Z^{(k)}$ be as in Sect. \ref{s7.3}.
Let $Y_0=\{z_0 \in R_0:z_k \not \in Z^{(k)}$ for any $k \geq 0 \}$, and
for $i \geq 1$, let 
$$
Y_i=\{z_0 \in R_0: z_i \in Z^{(i)} \ {\rm and} \  z_k \not \in
Z^{(k)} \ {\rm for \  all} \ k>i \}.
$$

Suppose $m(B \cap Y_i)>0$ for some $i>0$. Then $m(B \cap T^i Y_i)>0$,
and there is a vertical line
$\gamma$ with $m_\gamma(B \cap T^i Y_i)>0$. Let $\varepsilon>0$ be 
a small number. By the Lebesgue density theorem, there exists 
a short segment $\gamma_0 \subset \gamma$ with the property that
$m_\gamma(B \cap T^i Y_i \cap \gamma_0)>(1-\varepsilon)m_\gamma(\gamma_0)$.
We will show in the next paragraphs that points generic with respect
to some $\mu_i$ make up a definite fraction of  $\gamma_0$,
contradicting our choice of $\gamma_0$ if $\varepsilon$ is sufficiently
small. (The argument we present also works if 
$m(B \cap Y_0 \cap {\cal C}^{(0)})>0$. For
the case $m(B \cap Y_0 \cap (R_0 \setminus {\cal C}^{(0)}))>0$, 
use horizontal instead of vertical lines.)

Let $\tau_0$ denote the tangent vectors to $\gamma_0$, and let $\gamma_j=
T^j\gamma_0$. We regard all of $\gamma_0$ (which can be taken to be arbitrarily short) as bound to its nearest critical point, and 
let $n_1$ be the first time when part of $\gamma_j$ makes a free return 
to ${\cal C}^{(0)}$. As before, let $\Lambda_k=\Lambda^+_k$ or $\Lambda^-_k$. Let $D(\Lambda_k)$ denote the smallest rectangular
region bounded by $\Gamma^u$ and $\Gamma^s$-curves that contains $\Lambda_k$. If $\gamma_{n_1}$ crosses 
some $D(\Lambda_k)$ with two segments of at least comparable lengths extending beyond the two sides of $D(\Lambda_k)$, we consider the segment $\gamma_{n_1} \cap D(\Lambda_k)$ as having reached its final destination
and take it out of circulation. We then divide what remains of $\gamma_{n_1}$ into $I_{\mu j}$ and delete those subsegments that 
do not contain a point of $T^{n_1}(B \cap T^i Y_i)$.

Observe that 
for $z_0 \in \gamma_0 \cap T^iY_i$, $(z_0, \tau_0)$ is controlled through
time $n_1-1$, and by Lemma \ref{nlemma7.1}, $\tau_{n_1}$ splits correctly 
(see the
proof of Theorem \ref{theorem2}(2)). This is true not only 
for $z_0 \in \gamma_0 \cap T^iY_i$ but also for $z'_0 \in \gamma_0$ such that 
$z'_{n_1}$ is in the same $I_{\mu j}$ as $z_{n_1}$.
We iterate independently each one of the $I_{\mu j}$-segments 
that are kept. At the next free return 
we repeat the same procedure, namely we take out subsegments 
that cross some $D(\Lambda_k)$,
divide the rest into $I_{\mu j}$, delete those that do not contain
a point in the image of $B \cap T^i Y_i$, and observe that for the
remaining segments control is extended to the next free return.

Let $\gamma_0^d=\{z_0 \in \gamma_0:z_j$ is deleted
at a free return for some $j>0\}$, and let 
$\hat \gamma_0=\{z_0 \in \gamma_0: z_j$ reaches $D(\Lambda_k)$ for some $j$ and $k$ in the
required manner$\}$. 
We note that $m_{\gamma_0}(\gamma_0 \setminus (\gamma_0^d \cup \hat \gamma_0))=0$. This follows from a sublemma which is the first step in
the proof of Lemma \ref{nlemma10.4} (see \cite{BY2}, 
Sublemma 4 and its corollary).

Since $(\gamma_0 \cap B \cap T^i Y_i) \cap \gamma_0^d = \emptyset$, 
we have $(\gamma_0 \cap B \cap T^i Y_i) \subset \hat \gamma_0$
mod $0$ and that $\hat \gamma_0$ is the
disjoint union of a countable number of subsegments $\{\omega\}$
with the following properties:

\smallskip
-- each $\omega$ is mapped under some $T^{n(\omega)}$ onto a $C^2(b)$-curve that connects two $\Gamma^s$-sides 

\ \ \ \ of some $D(\Lambda_k)$;

--  $(z_0, \tau_0)$ is controlled up to time $n(\omega)$ for every
$z_0 \in \omega$.

\smallskip
\noindent From Lemmas \ref{nlemma10.2}, \ref{nlemma10.3} and Proposition
\ref{nproposition10.1}, it follows that
there exists $c_1>0$ independent of the choice of $\gamma_0$ such 
that for each $\omega$,

\bigskip
\centerline{$m_{\gamma}\{z_0 \in \omega:z_{n(\omega)} \in \cup \Gamma^s$ and 
is  generic w.r.t. 
some $\mu_k \}
> c_1 m_{\gamma}(\omega)$.}

\bigskip
\noindent This implies that 
$m_{\gamma}\{z_0 \in \gamma_0 :z_0$ is  generic w.r.t. some $\mu_k \}> c_1 m_{\gamma}(\hat \gamma_0)>c_1(1-\varepsilon)$ $m_\gamma(\gamma_0)$,
contradicting our choice of $\gamma_0$ if $c_1(1-\varepsilon)
>\varepsilon$. \hfill $\square$

\bigskip

\smallskip

For $i=1,2, \cdots, r'$, the set $B_i:=\{z \in R_0: z$ is generic with respect to $\mu_i\}$ 
can be thought of as the {\bf measure-theoretic basin} of $\mu_i$. 
When there exist multiple SRB measures, the $B_i$'s
can be quite delicately interwined (although they are not
``riddled"). We leave it as an exercise for the reader to construct an example of a $1$-dimensional map $f$ which when perturbed according to
the rules in Sect. \ref{s1.1} gives rise to a positive measure set of maps
$T$ with the following properties: 

(i) \ $T$ admits $n$ ergodic SRB measures for any $n \geq 2$; 

(ii) \ there is a Cantor set of stable curves
that do not meet any $B_i$;

(iii) \ every open set that meets any one of the curves in (ii) 
intersects every $B_i$ in  

\ \ \ \ a positive Lebesgue measure set.

%%%%%%%%%%%%%%%%%%%%%%%%%%%%%%%%%%%%%%%%%%%%%%%%%%%%%%%%

\subsection{Correlation decay and Central Limit Theorem}
\label{s10.4}
We indicate how Theorem \ref{theorem7} is proved.
The setup $T^R:\Lambda \to \Lambda$ is  designed so that
the statistical properties in question are easily read off from the tail
properties of the return time function $R$.
To use the results in \cite{Y3} or \cite{Y4} directly, however,
we need to consider returns to a single recurrent state.
Let $\tilde \mu$ be one of the $\mu_j$ in Proposition \ref{nproposition10.1}, 
and let
$\tilde \Lambda$ be one of the $\Lambda_i$ such that $\tilde \mu(\Lambda_i)>0$. For $z \in \tilde \Lambda$, we define a return time 
$\tilde R(z)$ of $z$ to $\tilde \Lambda$ as follows:
$$
\tilde R(z) \ = \ t_0+t_1+ \cdots +t_n,
$$
where $t_0=R(z), \ t_1=R(T^R(z)), \ \cdots, \ t_n=R((T^R)^n z)$
and $(T^R)^{n+1}z$ is the first return to $\tilde \Lambda$ under $T^R$. 
The results in [Y3] or [Y4] allow us to read off
information on the statistical properties of $(T, \tilde \mu)$
via the asymptotics of $m_{\partial R_0}\{z \in \partial R_0 \cap \tilde \Lambda: \tilde R(z)>n\}$.

\begin{lemma} 
\label{nlemma10.5}
There exists $K>0$ and $\tilde \theta_0<1$ such that
for every $n>0$,
$$
m_{\partial R_0}\{z \in \partial R_0 \cap \tilde \Lambda: \tilde R(z)>n\} < K\tilde \theta_0^n.
$$
\end{lemma}

This lemma, which we leave as an exercise, is an easy consequence 
of Lemma \ref{nlemma10.4}. The results in \cite{Y3} and \cite{Y4} state that
if the quantity estimated in Lemma \ref{nlemma10.5} 
is of order
${\cal O}(\frac{1}{n^{2+\varepsilon}})$ for some $\varepsilon>0$, then
the Central Limit Theorem holds in the context of Theorem \ref{theorem7}.
This condition is evidently satisfied here.
They also tell us that if this quantity is exponentially small, then every mixing component of $\tilde \mu$ has exponential decay of 
correlations as asserted.

%\vspace{.3in}
\newpage
\section{Global Geometry}
\label{s8} 
\subsection{Motivation}
\label{s8.1} 
Nonuniformly hyperbolic attractors have very complicated local structures. The purpose of this section is to develop an understanding of the 
{\it coarse geometry} of the attractor $\Omega$ for the maps in question,
that is to say, to describe
in a finite way the approximate shape and complexity of $\Omega$.
 
To illustrate the idea of coarse
geometry, consider the standard solenoid constructed from $z \mapsto z^2$. A good approximation of the attractor is given by the $k$th forward
image of $S^1 \times D_2$, which is a tubular neighborhood of 
a simple closed curve 
winding around the solid torus $2^k$ times.
For another example, consider piecewise monotonic maps in $1$-dimension.
Iterates of these maps continue to be piecewise monotonic and can be
understood in terms of their monotone pieces.
 
Returning to the maps under consideration, the standard solenoid example suggests that $R_k$ may be a good approximation of $\Omega$.
In analogy with $1$-dimension, one may also guess that
$R_k$ is a tubular neighborhood
of a simple closed curve whose $x$-coordinates vary in a
piecewise monotonic fashion. 
The latter is false, as is evident from the following sequence of pictures:
Depicted in (a) is a section of $R_k$ lying between two $C^2(b)$-curves;
(b) is the image of (a). As (b) is iterated, the horizontal distance
between the tips of the two parabolas increases as shown in (c), 
until at some point they fall on opposite sides of a component 
of the critical set, resulting in (d).
Since this happens to every ``turn" that is created, the geometry of
$R_k$ for large $k$ is quite complicated.
 
\begin{picture}(12, 8)
\put(1.5,0){
\psfig{figure=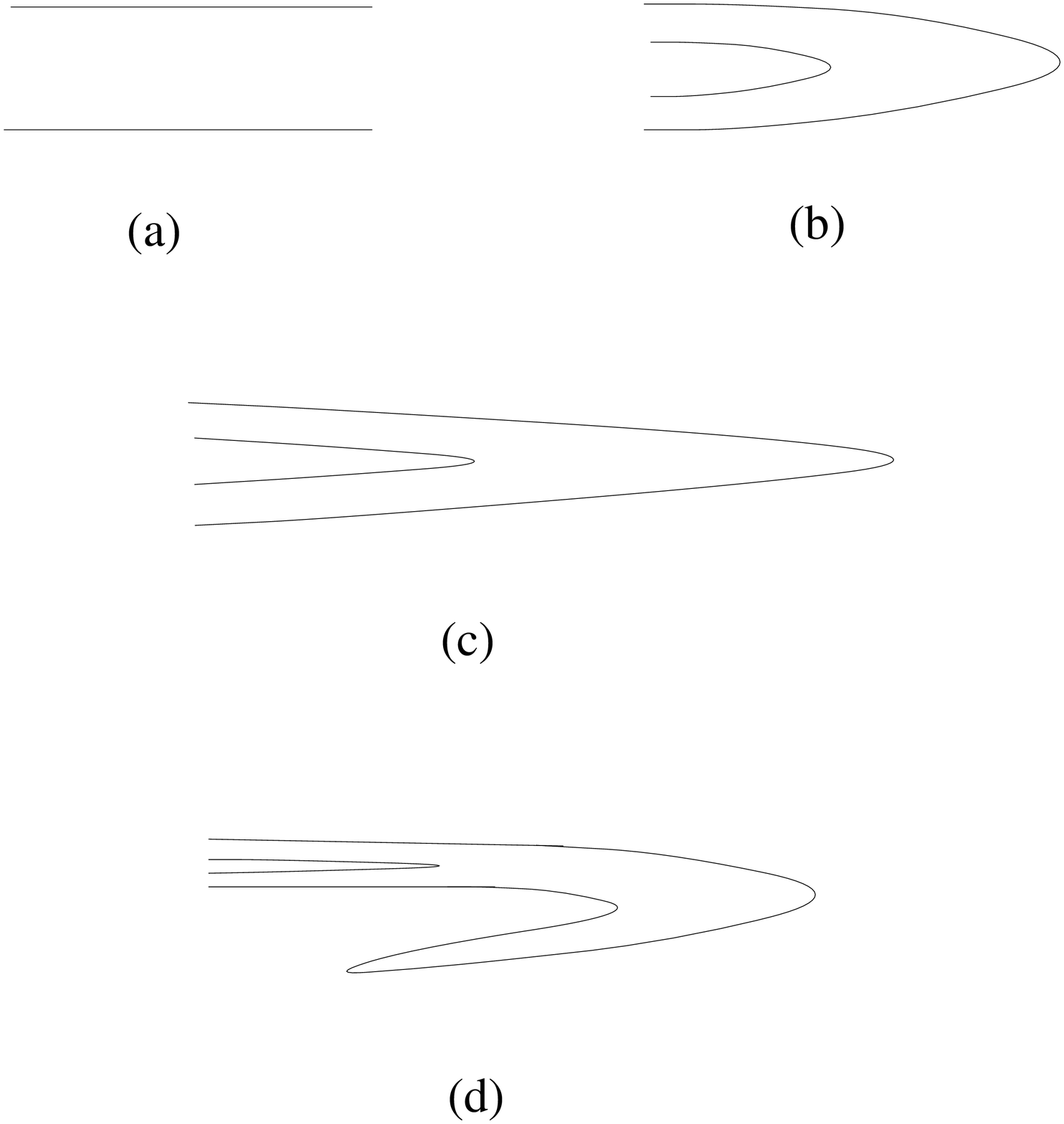,height=7cm,width=10cm}
}
\end{picture}

\vspace{0.7cm}

\centerline{Figure 3 \ \ The geometry of $R_k$}

\medskip
 
The purpose of this section is to introduce the idea of 
{\bf monotone branches}
as basic building blocks for understanding the global structure of $\Omega$.
To each map $T$ we will associate a {\bf combinatorial tree} whose edges correspond 
to monotone branches, and we will show that $\Omega$ has arbitrarily fine  neighborhoods made up of unions of finitely many monotone branches. Moreover, the way these branches fit together will tell us exactly how, in finite approximation, $T$ differs  from a $1$-dimensional map.

%%%%%%%%%%%%%%%%%%%%%%%%%%%%%%Sect. 8.2
\subsection{Monotone branches}
\label{s8.2} 
Let $\Gamma$ be the set of critical points as in Sections \ref{s3}--\ref{s6}.
For $z_0 \in R_0$, let $O_+(z_0)=\{ z_1, z_2, z_3, \cdots \}$
denote the positive orbit of $z_0$, and write $O_+(\Gamma)=
\cup_{z_0 \in\Gamma} O_+(z_0)$.
%%%%Definition 8.1 
\begin{definition} 
Let $\gamma$ be a connected subsegment of
$\partial R_k$. We say $\gamma$ is a {\bf (maximal) monotone segment}
if

(i) the two end points of $\gamma$ are in $O_+(\Gamma)$;

(ii) $\gamma$ does not intersect $O_+(\Gamma)$ in its interior.
\label{ndefinition8.1}
\end{definition}
When we say $\xi_i$ is an end point of a monotone segment,
it will be understood that $\xi_0$ is a critical point.
We record below some simple facts about monotone segments.

%%%%Lemma 8.1
\begin{lemma}
Let $\gamma \subset \partial R_k$ be a monotone segment. Then:

\noindent (a) All points near the two ends of $\gamma$ are in their fold periods; the part of $\gamma$ not in a 

fold period
(respectively bound period), if nonempty, is connected.

\noindent 
(b) If part of $\gamma$ is free, then its geometry is as follows:
$\gamma$ consists of a relatively long 

$C^2(b)$-curve connecting two sets of relatively small diameters at the
two ends; more 

precisely, there exists $p$ such that the $C^2(b)$-curve has length $>e^{-\beta p}$ while the 

diameters of the two small sets are $<b^{\frac{p}{2}}$; also, 
the curvature of $\partial R_k$ at the end 

point $\xi_i$ of $\gamma$ is $>b^{-i}$.

\noindent (c) If $\gamma$ meets $\Gamma$ in $r$ points, $r \geq 0$, then
$T(\gamma)$ is the union of $r+1$ monotone segments 

joined together
at the $T$-images of these points.

\label{nlemma8.1}
\end{lemma}
\noindent {\bf Proof:} \ (c) follows from the definition of a monotone
segment. (a) follows from the way monotone segments are created and 
from the monotonicity of bound and fold periods
(see the proof of Lemma \ref{nlemma4.10}).
The first assertion in (b) follows from estimates on the relative
 sizes of the parts of $\gamma$ that are in bound versus fold periods; 
the second follows from the
curvature formula in the proof of 
Lemma \ref{nlemma2.4}. \hfill $\square$

\bigskip
We now begin to study the geometry of certain $2$-dimensional objects.

%%%%Definition 8.2
\begin{definition} 
A simply connected region $S \subset R_k$ is called a {\bf
monotone branch} if it is bounded by
 two monotone segments $\gamma, \ \gamma' \subset \partial R_k$ 
and two ends $E_\xi$ and $E_\zeta$ with the following properties:

\noindent (i) If the end points of $\gamma$ are $\xi_i$ and $\zeta_j$,
then the end points of $\gamma'$ are $\xi_i'$ and $\zeta_j'$ where

$\xi_0$ and $\xi_0'$ lie on the upper and lower boundaries of
the same component $Q^{(k-i)}$ of 

${\cal C}^{(k-i)}$, and $\zeta_0$ and $\zeta_0'$ are related in the same way. 

\noindent (ii) $E_\xi = T^i\{z \in Q^{(k-i)}(\xi_0):|z-\xi_0|<b^{\frac{k-i}{4}}\}$; its 
{\bf time of creation} is said to be $k-i$;

$E_\zeta$ and its time of creation are defined analogously.

\noindent (iii) We define the {\bf age} of $E_\xi$ to be $i$
and require that $i<\theta^{-1} (k-i+1)$; there is an 

analogous limit on the age of $E_\zeta$.
\label{ndefinition8.2}
\end{definition}

The definitions of $E_\xi$ and $E_\zeta$ are quite arbitrary,
subject only to the following considerations:
We want $E_\xi$  to be large enough to contain all the critical orbits
that originate from $Q^{(k-i)}(\xi_0)$. On the other hand, we want it 
to remain relatively small during the life span of the monotone branch,
so that the phenomenon depicted in Figure 3 does not occur.
We assume $\theta$ is chosen such that for $i < \theta^{-1}(k-i+1)$, 
$\|DT\|^i b^{\frac{k-i}{4}} < b^{\frac{k-i}{8}} << e^{-\alpha i}$,
which is $< d_{\cal C}(z_i)$ for $z_0 \in \Gamma$ by (IA2) in Section \ref{s3};
that is to say, if $S$ is a monotone branch of $R_k$, then its ends are at 
least a certain distance away from ${\cal C}^{(k)}$. 
It is not always easy to visually identify monotone branches, particularly
when their boundary segments are in fold periods. When part of $\gamma$ 
is free, it follows from Lemma \ref{nlemma8.1}(b) that $S$
consists of a (relatively long) horizontal strip with two small
blobs at the two ends.

\bigskip
\noindent {\bf Tree structure of a class of monotone branches} \ 
\smallskip

Monotone branches can be constructed as follows. 
First we declare that $R_0$ is a monotone branch 
(even though it has no ends). Then if $x_i<x_{i+1}$ 
are adjacent critical points of the $1$-dimensional map $f$, 
the $T$-image of
$\{z=(x,y):x_i-b^{\frac{1}{4}}<x$ $<x_{i+1}+b^{\frac{1}{4}}\}$ is a
monotone branch of $R_1$. 
In general, let $S$ be a monotone branch of $R_k$. 
If one of the ends of $S$ is at its maximum allowed age, 
then $S$ is ``discontinued", meaning we do not iterate it further.
If not, $T(S)$ is the union of a finite number of monotone branches
of $R_{k+1}$. More precisely, 
if $S \cap {\cal C}^{(k)}=\emptyset$, then $T(S)$ is a monotone branch. 
If $S \cap Q^{(k)} \neq \emptyset$, 
then $S \supset Q^{(k)}$ (in fact, $S$  extends beyond $Q^{(k)}$ by 
$>e^{-\alpha k}$ in both directions).
 If $S$ contains $r$ components of ${\cal C}^{(k)}$, 
then $T(S)$ is the union of $r+1$ monotone branches 
split roughly along the $T$-images of the middle of each of the $Q^{(k)}$
contained in $S$ (cf. Lemma \ref{nlemma8.1}(c)).

Let ${\cal T}=\cup_k {\cal T}_k$ denote the set of all monotone branches inductively constructed this way, with ${\cal T}_k$ consisting of 
branches of $R_k$. More precisely, ${\cal T}_0=\{R_0\}$, and 
${\cal T}_{k+1}$ is obtained
from ${\cal T}_k$ via the procedure described above.
{\bf We will be working exclusively with monotone branches in ${\cal T}$}, which is a proper subset of the set of all monotone branches in 
Definition \ref{ndefinition8.2}. 
The set ${\cal T}$ has a natural tree structure: we call the branches 
obtained by mapping forward and subdividing a given branch
its {\bf descendants}. Note that every branch in ${\cal T}_k$ has a unique 
{\bf ancestor} in ${\cal T}_i$ for every $i<k$, but not all branches in 
${\cal T}$ have offsprings: the ones with no offsprings are exactly those
 one of whose ends has reached its maximum allowed age.

We have elected to discontinue a branch before its geometry ``deteriorates".
An immediate question that arises is what happens to the part of the attractor 
contained in a discontinued branch.
We will show in the next subsection 
that branches farther down the tree ${\cal T}$ can be used to take
 its place.
We will, in fact, prove the following stronger version of 
Theorem \ref{theorem3}.

\bigskip

\noindent {\bf Theorem $\rm{5}^{\prime}$ } \ \  
{\it One can construct special neighborhoods  
$\tilde R_n$ as in Theorem \ref{theorem3} using only monotone branches from 
${\cal T}_k,
\ n \leq k <(1+K\theta)n$.} 

\bigskip

\begin{picture}(8, 5)
\put(3, 0.5){
\psfig{figure=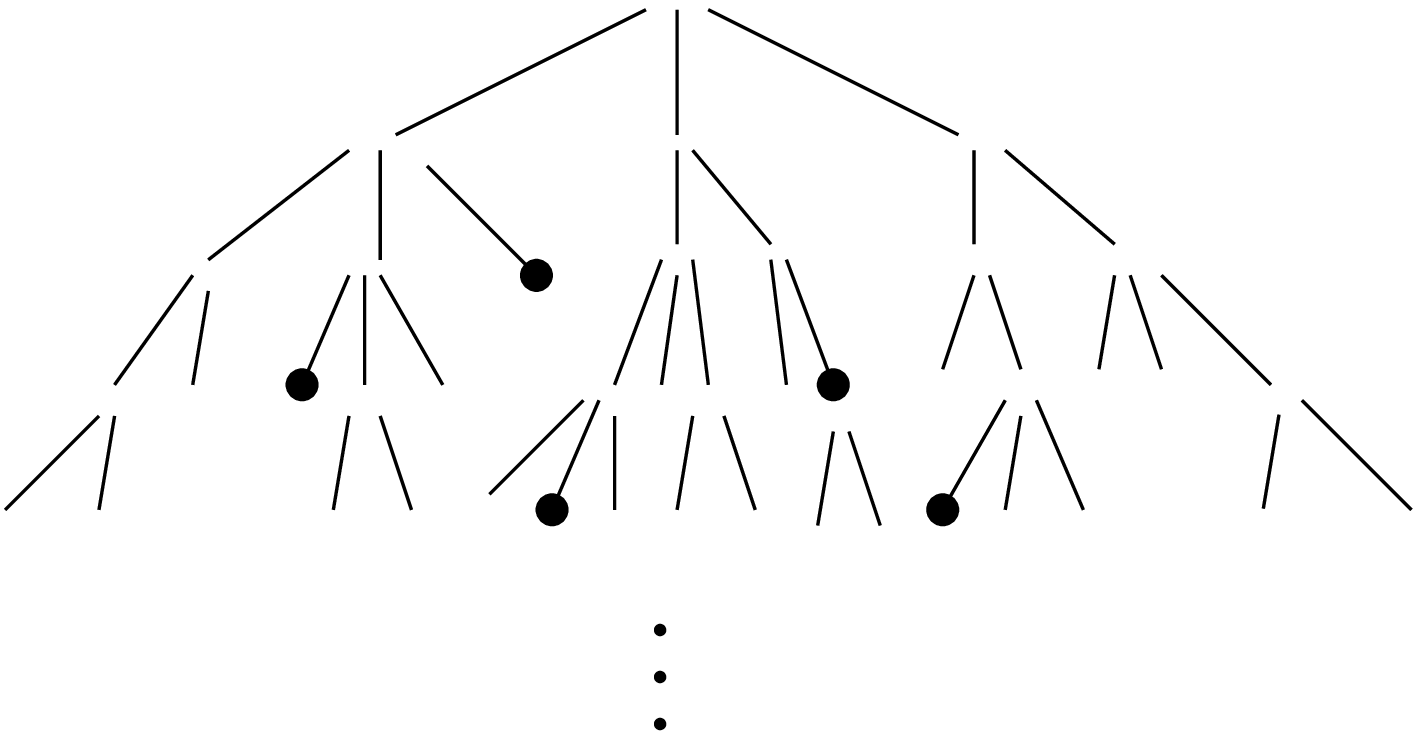,height=4cm,width=8cm}
}
\end{picture}

\medskip

\centerline{Figure 4 \ \ Tree of monotone branches: branches ending in $\bullet$ 
are discontinued}

\medskip

%%%%%%%%%%%%%%%%%%%%Section 8.3 

\subsection{Replacement of branches}
\label{s8.3} 

Let $S \in {\cal T}_k$ be a branch whose ends are denoted by $E_\xi$ and
$E_\zeta$. In the discussion to follow, we assume that $E_\xi$ is fairly advanced in age, meaning $(k-i) \sim \theta i$ where $i$ is the age of
$E_\xi$ and $k-i$ is its time of creation. As we search for replacements
for $S$, the picture we hope to have is the following.
There is a finite collection of branches 
$\{B\} \subset \cup_{k<j \leq (1+K\theta)k} {\cal T}_j$ such that 

\smallskip
(i) \ the ends of $B$ are contained in those of $S$; and

(ii) \ if $S \in {\cal S}$ where ${\cal S} \subset {\cal T}$ is a cover
of $\Omega$, then replacing $S$ by $\{B\}$ does not 

\ \ \ \ \ leave any part of $\Omega$ exposed.

\smallskip
\noindent Let $Q^{(k-i)}$ be the component of ${\cal C}^{(k-i)}$
containing $T^{-i}E_\xi$. 
We hope to show that $T^{-i}S \subset Q^{(k-i)}$, so that
the picture described above pulled back to $Q^{(k-i)}$ is as shown in
Figure 5. 

\begin{picture}(13,4.5) 
\put(0,0){
\psfig{figure=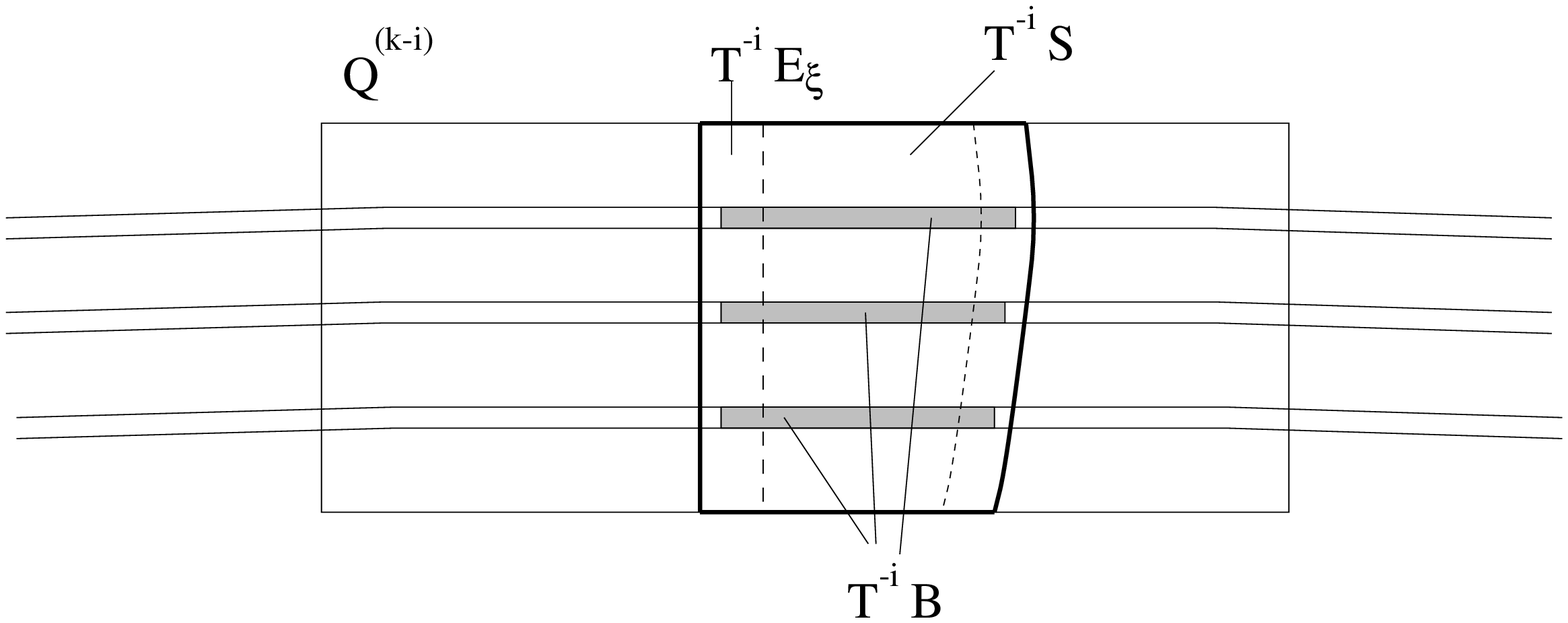,height = 4cm, width = 13cm}
}
\end{picture}

\medskip

\centerline{Figure 5 \ \ Replacing $S$ by $\{ B \}$}

We begin to systematically justify this picture. For $j = 0, 1, \cdots, i-1$, 
let $S_j \in {\cal T}_{k-i+j}$ 
be the ancestor of $S$, so that $S_0$ is the
monotone branch of $R_{k-i}$ containing $Q^{(k-i)}$.
Let $E_0$ denote the end of $S_0$ contained in $Q^{(k-i)}$, 
and let $E_j=T^jE_0$. Let the other end of $S_j$ be called $E'_j$. 
Let $t>k-i$, and let $P \in {\cal T}_t$ be such that $P \cap Q^{(k-i)}$
is a horizontal strip bounded by two $C^2(b)$-curves stretching
 all the way across $Q^{(k-i)}$. 
We think of $P$ as a {\it pre-branch} with respect to $S_0$
in the sense that $P \subset S_0$ and it is not yet born 
when $S_0$ is created.
If $P$ is not discontinued,
then we let $P_1$ be the (unique) child of $P$ with one end in
$E_1$, and assuming $P_1$ is not discontinued, we let $P_2$
be the child of $P_1$ with one end in $E_2$.
Similarly, we define $P_3, P_4, \cdots$ up to  $P_i$ if it makes sense.

\begin{lemma}
\label{nlemma8.2}
There exists $K_1$ depending on $\rho$ such that 

\noindent (i) \ for all $j$ with $K_1(k-i) <j \leq i$, $T^{-j}S_j \subset Q^{(k-i)}$;  

\noindent (ii) if $P_j$ is defined for all $j \leq K_1(k-i)$, then 
it is defined for all $j \leq i$; moreover,

for each $j \geq K_1(k-i)$, $P_j \subset S_j$, and the two ends of 
$P_j$ are contained in the two 

ends of $S_j$.
\end{lemma}

We isolate the following sublemma, the ideas in which are also used elsewhere.
See Sect. \ref{s6.1} for notation.

\begin{sublemma} 
\label{nsublemma8.1}
Let one of the horizontal boundaries of $Q^{(s)}$, any $s$,
 be identified with $[-\rho^s, \rho^s]$, with the critical point corresponding to $0$. Then for every $I_{\mu_0 j_0} \subset [-\rho^s, \rho^s]$, there exists $n<K|\mu_0|$ such that $T^n I_{\mu_0 j_0}$ 
traverses completely a component of ${\cal C}^{(0)}$.
\end{sublemma}

\noindent {\bf Proof:} \ Let $\omega_0=I_{\mu_0 j_0}$, and let
$r_0$ be the first time when part of $\omega_0$ makes a free return with
$T^{r_0}\omega_0$ containing an $I_{\mu j}$ of full length.
By Corollary \ref{ncorollary4.3}, either $T^{r_0}\omega_0$ contains one of the
outermost
$I_{\mu j}$ (which we will call $\tilde I$) or it contains
some $I_{\mu_1 j_1}$ with $|\mu_1|<K\beta |\mu_0|$.
In the latter case, we let $\omega_1=I_{\mu_1 j_1}$ and continue 
to iterate until $r_1$ iterates later  when
part of $T^{r_1}\omega_1$ is free and contains either $\tilde I$ or some 
$I_{\mu_2 j_2}$ with $|\mu_2|<K\beta |\mu_1|$.
After a finite number of iterates, we have $T^{r_q}w_q
\supset \tilde I$. 

From Corollary \ref{ncorollary4.3}, we see
that at the end of its bound period, $T^p \tilde I$ has length
$>>\delta$. Inductively define $\tilde I_{p+j}= T(\tilde I_{p+j-1}) \setminus {\cal C}^{(0)}$ for $j=1,2, \cdots$. Then $\tilde I_{p+j}$ 
is a connected $C^2(b)$-curve which grows essentially exponentially 
-- until it crosses completely a component of ${\cal C}^{(0)}$. 
Since $r_i \sim |\mu_i|$ up to the point when $T^{r_q}w_q
\supset \tilde I$, and the growth
is exponential thereafter, we conclude that the end game is reached
in a total of $<K|\mu_0|$ iterates. \hfill $\square$

\bigskip

\noindent {\bf Proof of Lemma \ref{nlemma8.2}:} \ 

\begin{claim} 
\label{nclaim8.1}
There exists $K_1$ (depending on $\rho$) such that $T^{-K_1(k-i)}S_{K_1(k-i)} \subset Q^{(k-i)}$.
\end{claim}

\noindent {\it Proof of Claim \ref{nclaim8.1}:} \ We 
 identify the upper horizontal boundary of $Q^{(k-i)}$
with the interval $[-\rho^{k-i}, \rho^{k-i}]$, 
with the critical point corresponding to $0$, and let
 $n_1$ be the smallest $n$ such that $T^n[0, \frac{1}{2}\rho^{k-i}]$
intersects the horizontal boundary of some $Q^{(k-i+n)}$.
From Sublemma \ref{nsublemma8.1}, $n_1<K_1(k-i)$ for some $K_1=K(\rho)$.
The claim is proved once we show that $T^{-(n_1+1)}S_{n_1+1}
\subset Q^{(k-i)}$.
Let $[0, \ell]$ be the shortest interval such that $T^{n_1}[0,\ell]$
contains the entire horizontal boundary of a $Q^{(k-i+n_1)}$.
Since this boundary is free, $\ell < \frac{1}{2}\rho^{k-i} +
e^{-c'n_1}\rho^{k-i+n_1}$, which is $\approx \frac{1}{2}\rho^{k-i}$.
Let $\hat S_{n_1}$ be the section of $R_{k-i+n_1}$
from $E_{n_1}$ to the middle of $Q^{(k-i+n_1)}$.
Since $b^{\frac{k-i}{4}}K^{K_1(k-i)} << \rho^{k-i+K_1(k-i)}$, 
we have that $T^{n_1}Q^{(k-i)} \supset \hat S_{n_1}$.
It remains to show $S_{n_1+1}=T(\hat S_{n_1})$,
for which we need only to check that 
$T(\hat S_{n_1})$ is a monotone branch. To do that, it suffices
to show that for $j< n_1$, $T^j[0,\ell]$ does not contain
the horizontal boundary of any $Q^{(k-i+j)}$. Suppose it does for some $j$. By 
our choice of $n_1$, this can happen only if $\ell >\frac{1}{2}\rho^{k-i}$
and $\mid T^j[\frac{1}{2}\rho^{k-i},\ell] \mid \geq \rho^{k-i+j}$,
which is impossible, for 
$\mid T^j[\frac{1}{2}\rho^{k-i},\ell] \mid <e^{-c'(n_1-j)} \rho^{k-i+n_1}$.
\hfill $\diamondsuit$

\medskip
Suppose we are guaranteed that $P_{n_1}$ exists. We show next that
$P_{n_1+1}$ exists and has the properties in Lemma \ref{nlemma8.2}(ii).
Let $\gamma$ be the part of a horizontal boundary of $P$ that lies
below $[0,\ell]$.
From the estimates above, we know that $T^{n_1}\gamma$ is $C^0$ very near
$T^{n_1}[0,\ell]$.  Let $\hat P_{n_1}$ be the section of 
$T^{n_1}(P \cap Q^{(k-i)})$ that runs from $E_{n_1}$ to the middle of some $Q^{(t+n_1)} \subset Q^{(k-i+n_1)}$.
We claim that $P_{n_1+1}=T(\hat P_{n_1})$. Clearly, $P_{n_1+1} \subset
S_{n_1+1}$.
To see that $P_{n_1+1}$ is a monotone branch, it suffices to observe
that for $j<n_1$, $T^{-n_1+j}\hat P_{n_1} \cap {\cal C}^{(t+j)}
= \emptyset$, which is an immediate consequence of the fact that
$T^{-n_1+j}\hat S_{n_1}
\cap {\cal C}^{(k-i+j)} = \emptyset$.

We are now ready to show that $P_j$ exists for all $j \leq i$.
Suppose that $P_{j-1}$ exists. The only reason why $P_j$ may not exist
is that one of its ends has reached its maximum allowed age.
Of the two ends of $P_{n_1+1}$, the one contained in $E_{n_1+1}$ is
clearly created earlier, which means that of the two ends of $P_{j-1}$,
the one contained in $E_{j-1}$ is created earlier. It suffices therefore
 to check that this end survives
the step from $P_{j-1}$ to $P_j$. It does, because it is created later than $E_{j-1}$ and has the same age as $E_{j-1}$, and, by definition, $E_{j-1}$
has not reached its maximum allowed age.
 
From here on we argue inductively that 
the relations in Lemma \ref{nlemma8.2}(ii)  between $P_j$ and $S_j$ hold
from $j=n_1+2$ to $j=i$. 
Assume this is true for $j-1$, and that $S_{j-1}$ has more than one child. Then $S_j=T(\hat S_{j-1})$ where $\hat S_{j-1}$ is the section of 
$S_{j-1}$ from $E_{j-1}$ to the middle of some $Q^{(k-i+j-1)}$. 
Since by inductive assumption $P_{j-1}$ has its ends contained in those 
of $S_{j-1}$, we are assured that it traverses some
$Q^{(t+j-1)} \subset Q^{(k-i+j-1)}$. Letting $\hat P_{j-1}$
be the section of $P_{j-1}$ from its end in $E_{j-1}$ to the middle of
$Q^{(t+j-1)}$, we see that $P_j=T(\hat P_{j-1})$ has the desired properties. 

This completes the proof of Lemma \ref{nlemma8.2}.
\hfill $\square$
 
\bigskip
\noindent {\bf Proof of Theorem $\rm{5}^{\prime}$:} \ Let ${\cal S}_0=\{R_0\}$,
and assume that for each $n \leq m$, a collection of monotone branches
${\cal S}_n$ is selected so that $\tilde R_n:= \cup_{S \in {\cal S}_n} S$ 
is a neighborhood of the attractor, and each $S \in {\cal S}_n$ has the 
following properties:
\medskip

(i) $S \in {\cal T}_k$ for some $n \leq k \leq (1+3\theta)n$;

(ii) if an end of $S$ is of age $i$, i.e. it is
created at time $k-i$, then
$2\theta i \leq k-i+1$.

\medskip

\noindent Note that (ii) is a more stringent requirement than the
definition of monotone branches.

The collection ${\cal S}_{m+1}$ is defined as follows.
For each $S \in {\cal S}_m$, if the ends of $S$ have not reached their
maximum ages as allowed by (ii) above, then we put the children of $S$
in ${\cal S}_{m+1}$. If one of its ends has reached this age, then
we choose a collection of branches $\{P\}$ to be specified in the next
paragraph, construct from each $P$ a monotone branch $P_i$ as in Lemma \ref{nlemma8.2},
replace $S$ by $\{P_i\}$ and put the children of $P_i$ in ${\cal S}_{m+1}$.

Suppose for definiteness that $S \in {\cal T}_k$, and its end
$E$ has reached age $i$ where
\begin{equation}
\label{aformula8.1}
2\theta i=k-i+1.
\end{equation}
Let $Q^{(k-i)}$ be the component of ${\cal C}^{(k-i)}$ containing $T^{-i}E$.
Let $\{P\}$ be the subcollection of ${\cal S}_{k-i+1}$ with the property
that $P \cap Q^{(k-i)} \neq \emptyset$. Observe immediately that by
our inductive hypotheses, $P$ is a monotone branch of $R_{\tilde k}$
for some $\tilde k$ with
\begin{equation}
\label{aformula8.2}
k-i+1 \ \leq \ \tilde k \ \leq \ (1+3\theta)(k-i+1).
\end{equation}
Since $e^{-\alpha(1+3\theta)(k-i+1)}>>\rho^{k-i}$, it follows that
$P$ intersects $Q^{(k-i)}$ in a horizontal strip bounded by $C^2(b)$
curves. Note also that since the union of the elements of ${\cal S}_{k-i+1}$
covers $\Omega$, we have $\cup P \supset (Q^{(k-i)} \cap \Omega)$.

To justify the validity of this replacement procedure, we need to show that

\smallskip
\noindent (a) \ for each $P$ as above, 
$P_{K_1(k-i)}$ is well defined where $K_1$
is as in Lemma \ref{nlemma8.2};

\noindent (b) $P_i$ is a monotone branch of 
$R_j$ for some $j \leq (1+3\theta)m$.

\smallskip
Suppose that an end of $P$, which is a branch of $R_{\tilde k}$,
is of age $\tilde i$. Then 
\begin{equation}
\label{aformula8.3}
2\theta \tilde i \ \leq \tilde k- \tilde i +1.
\end{equation}
To prove (a), if suffices to verify that this end lasts another
$K_1(k-i)$ iterates, i.e. 
$$
\theta \ [ \tilde i +K_1(k-i) ] \ \leq \ \tilde k- \tilde i +1.
$$
This is true because  $\theta \tilde i \leq \frac{1}{2}
(\tilde k- \tilde i +1)$ by (\ref{aformula8.3}), and  
\begin{eqnarray*}
K_1 \theta (k-i) & \leq & K_1 \theta (\tilde k+1) \ = \ K_1 \theta [(\tilde k- \tilde i +1)+ \tilde i] \\
& \leq & K_1 \theta (\tilde k- \tilde i +1) (1+\frac{1}{2\theta})
\ << \ \frac{1}{2} (\tilde k- \tilde i +1).
\end{eqnarray*}
The first inequality above is by (\ref{aformula8.2}) and the second by 
(\ref{aformula8.3}).

To prove (b), we need to check that the age of the end of $P_i$ that is
contained in $E$, namely $\tilde k +i$, is $\leq (1+3\theta)m$.
Observe first that $i \leq m$. This is because the replacement
procedure described in Lemma \ref{nlemma8.2}
does not change the ages of the respective ends of the monotone
branch in question. (The age of an end is equal to
the ``age" of the critical orbits it contains.)
Thus it remains to check that 
$$
\tilde k \ \leq \ (1+3\theta)(k-i+1) \ = \ (1+3\theta) 2\theta i
\ < \ 3\theta i \ \leq \ 3\theta m,
$$
the first inequality above coming from (\ref{aformula8.2}) and the equality 
from 
(\ref{aformula8.1}). 
This completes the proof of Theorem $\rm{5}^\prime$. \hfill $\square$

\bigskip
 
We mention two bonuses of this construction.
 
First, it can be seen inductively that for every $S \in {\cal S}_n$, 
if $S$ is a branch of $R_k$, then the two monotone segments of $\partial R_k$ that bound $S$ must necessarily be from different components of 
$\partial R_0$. This is used in Sect. \ref{s9.5}. 
 
Second, we claim that if $\deg(f) \neq 0$, then all of our monotone branches
$S \in {\cal S}_m$ intersect the attractor $\Omega$ in an essential way.
Let us call a monotone branch $S$ {\bf essential} if every curve connecting the two monotone segments $\gamma$ and $\gamma'$ in $\partial S$ meets $\Omega$. Observe first that $R_0 \in {\cal S}_0$ is essential
if $\deg(f) \neq 0$. If not, then there exists a curve
 $\omega$ connecting the two components of $\partial R_0$ that does not meet 
$\Omega$. Since $\Omega= \cap_k R_k$, this implies that for some $k$, 
$R_k \cap \omega = \emptyset$, which is absurd since $R_k$ is not contractible. Assuming that $S \in {\cal S}_m$ is essential, 
then clearly all the monotone branches that comprise
 $T(S)$ are essential if no end
replacements are needed in the next step. If an end replacement is
required, then since the new branches are the images of parts of earlier 
essential branches,
they are again essential.

%%%%%%%%%%%%%%%%%%%%%%%%%Section 8.4
\subsection{The coarse geometry of $\Omega$}
\label{s8.4} 
We explain in the following sequence of pictures
exactly how, in finite approximation,
the geometry of $\Omega$ differs from that of a small tubular
neighborhood of a single curve. These pictures are justified by
Lemma \ref{nlemma8.2}. Referring back to Figure 3(c),
we may think of the region between the parabolas as made up
to two ends belonging to adjacent branches. We know from Lemma \ref{nlemma8.2}
 that long before 
the tips of these parabolas ``separate", that is,  before the ends 
in question reach their maximum allowed age, there are {\it pre-branches} 
inside running parallel to these parabolas.
In Figure 6 below, the pre-branches are shown in grey, and 
the zig-zagging cut-lines represent pre-images of the
critical set. These cut-lines will become ``turns" before the ends in question reach their maximum allowed age.
 
\begin{picture}(13, 4)
\put(1,0){
\psfig{figure=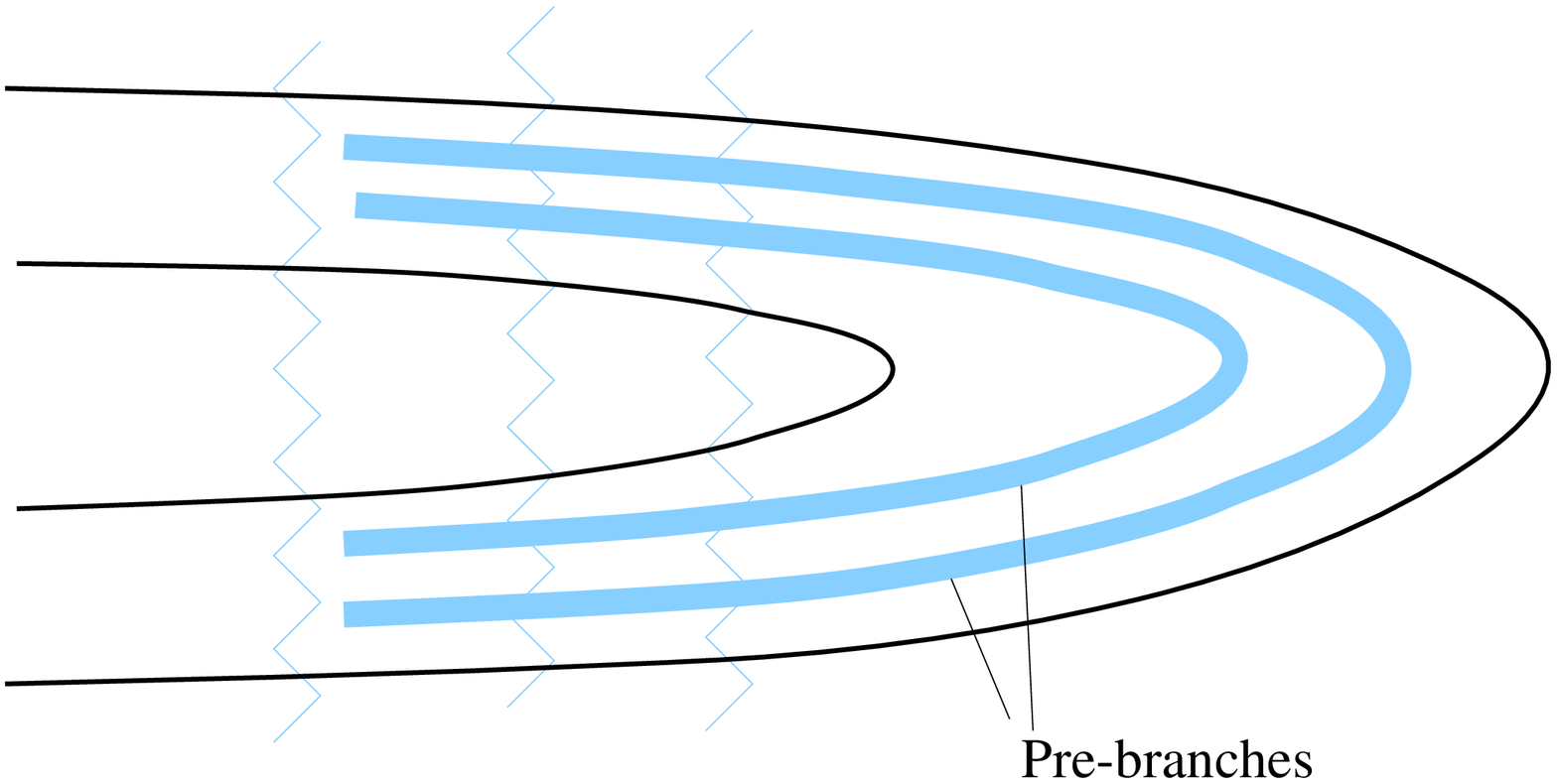,height = 3.5cm, width = 10cm}
}
\end{picture}

\medskip

\centerline{Figure 6 \ \ Pre-branches waiting to be released} 

\bigskip

As this age is reached, the pre-branches are released. Figure 7(a) shows four 
newly released 
montone branches
grafted onto a branch created earlier. Once released,
the new branches evolve independently, resulting possibly
in the configuration in Figure 7(b) (cf. Figure 3(d)).
 
\begin{picture}(13, 3.5)
\put(1,0){
\psfig{figure=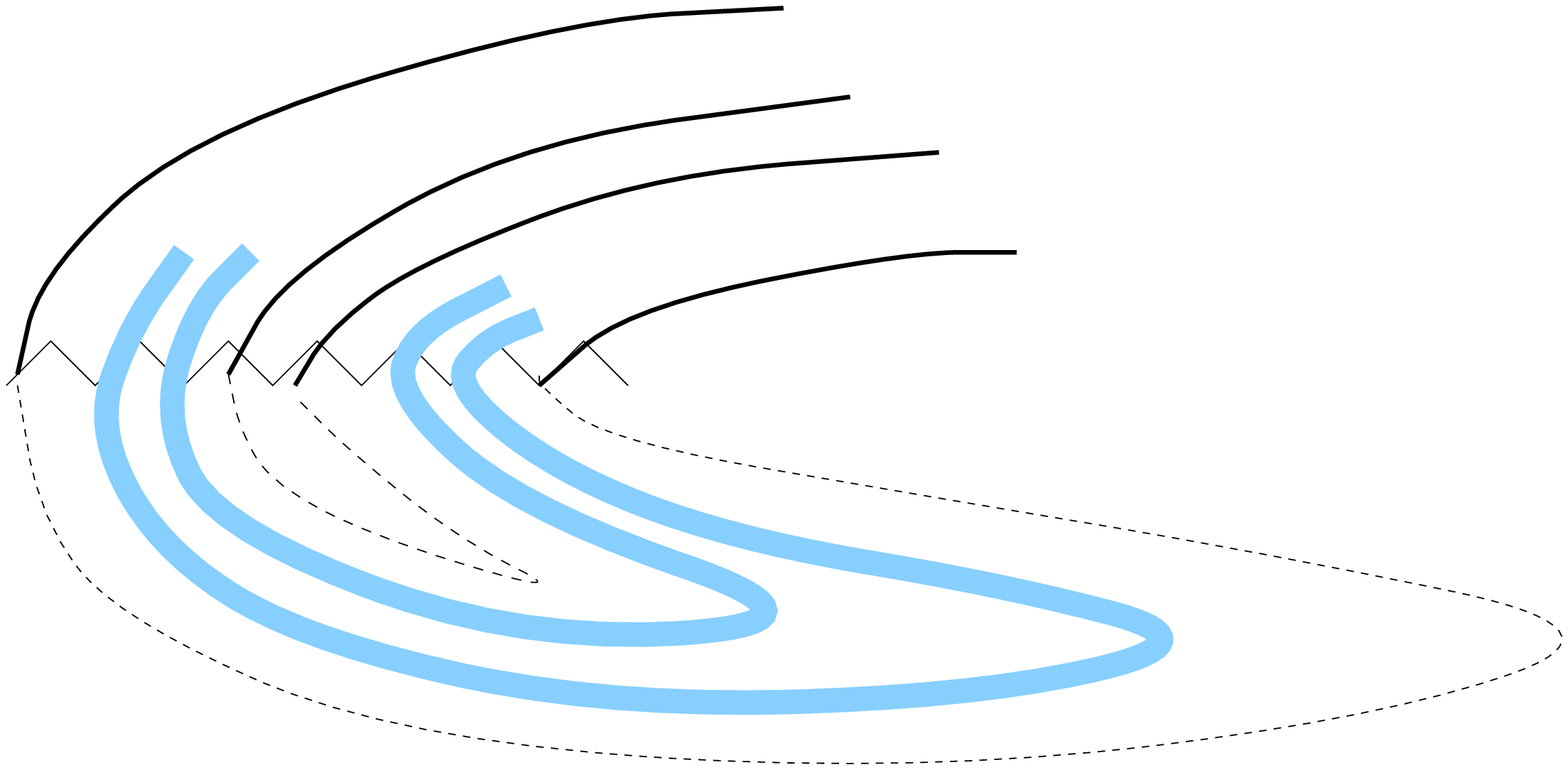,height = 3cm, width = 10cm}
}
\end{picture}

\vspace{0.5cm}

\centerline{\bf (a)}

\begin{picture}(13, 3)
\put(1,0){
\psfig{figure=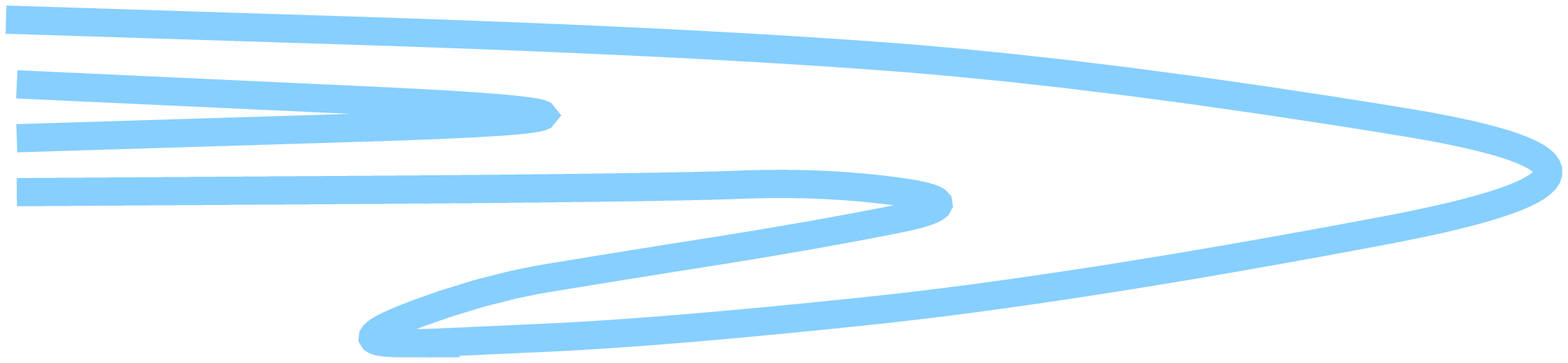,height = 1.5cm, width = 10cm}
}
\end{picture}

\vspace{0.5cm}

\centerline{\bf (b)}

\medskip

\centerline{ Figure 7 \ \ Newly released monotone branches evolving independently}

\bigskip
 
The boundaries of every turn (or pair of ends) created every step 
of the way will in time separate, releasing new branches grafted 
onto ones born earlier.
As the new branches evolve, they create new turns, 
which again will last for only a finite duration of time.
In terms of global geometry, this, in a sense, is the {\it only} way
in which $T$ differs from a $1$-dimensional map.
Tip replacements are scheduled to take place roughly once every
$\sim \log\frac{1}{b}$ iterates, so that in the limit as $b$ tends to $0$,
no replacement is needed -- as it should be for $1$-dimensional maps.

\vspace{.3in}

%%%%%%%%%%%%%%%%%%%%Section 9 
\section{Symbolic Dynamics and Topological Entropy}
\label{s9} 
The goals of this section are  (1) to introduce a natural and 
unambiguous coding of all points on the attractor $\Omega$ for the maps
 in question, and (2) to use this coding to obtain results on
topological entropy and equilibrium states.
%%%%%%%%%%%%%%%%Section 9.1 
\subsection{Coding of points on the attractor}
\label{s9.1} 

Abusing notation slightly, let $x_1<x_2< \cdots <x_r<x_{r+1}=x_1$ 
be the critical
points of $f$ in the order in which they appear on the circle,
and let ${\cal C}_i:={\cal C} \cap {\cal C}_i^{(0)}$ 
where ${\cal C}_i^{(0)}$ is the component of
${\cal C}^{(0)}$ containing $x_i$. We remark that ${\cal C}_i$ 
may be a fractal set,
and that for an arbitrary $z \in R_0$ near ${\cal C}_i$, it does not 
always make sense to think of $z$ as being located on the left 
or on the right of ${\cal C}_i$. The goal of this subsection is to show
that points on $\Omega$ are special, in that for them this 
left/right notion is always well defined.
 
Recall that if $Q^{(k)}$ is a component of
${\cal C}^{(k)}$, then ${\hat Q}^{(k)}$ is the component of
$R_k \cap {\cal C}^{(k-1)}$ containing $Q^{(k)}$. In particular,
${\hat Q}^{(k)} \setminus Q^{(k)}$ has a left and a right component.

%%%%Lemma 9.1 
\begin{lemma} 
\label{nlemma9.1}
The critical set ${\cal C}$ partitions 
$\Omega \setminus {\cal C}$
into disjoint sets $A_1, \cdots, A_r$ as follows:

- For $z=(x,y) \not \in {\cal C}^{(0)}$, $z \in A_i$ if and only if
$x_i<x<x_{i+1}$.

- For $z \in {\cal C}_i^{(0)} \setminus {\cal C}_i$, let $Q^{(k)}$ be such that $z \in \hat Q^{(k)} \setminus
Q^{(k)}$. Then $z \in A_i$ if it lies in 

\ \ the right component of ${\hat Q}^{(k)} \setminus Q^{(k)}$; 
$z \in A_{i-1}$ if it lies in the left component of 

\ \ ${\hat Q}^{(k)} \setminus Q^{(k)}$.
\end{lemma}

\noindent {\bf Proof:} 
This lemma is an immediate consequence of our description of critical
regions (Theorem \ref{theorem1}(1)). The sets $\{A_i\}$ are  defined
by the conditions above. What sets points in $\Omega$
apart from arbitrary points in $R_0$ is that $z \in \Omega$ implies $z \in
R_k$ for all $k$, so that for $z \in {\cal C}^{(0)}$, there are only two
possibilities: either $z \in \cap_{k \geq 0} {\cal C}^{(k)}$, in which case
it is a critical point,
or there is a largest $k$ such that $z \in {\cal C}^{(k-1)}$.
In the latter case, it follows
from the geometric relation between ${\cal C}^{(k)}$ and 
${\cal C}^{(k-1)}$ that
$z \in \hat Q^{(k)} \setminus Q^{(k)}$ for some $Q^{(k)}$. \hfill $\square$

\bigskip
Lemma \ref{nlemma9.1} gives a well defined {\bf address} $a(z)$ for all 
$z \in \Omega \setminus {\cal C}$.
We write $a(z)=i$ if $z \in A_i$. Points in ${\cal C}$ have two addresses;
for example, for $z \in {\cal C}_i$, $a(z)=$ both $i-1$ and $i$. 
This in turn allows us
to attach
to each $z_0 \in \Omega$ with $z_i \not \in {\cal C}$ for all $i$
an {\bf itinerary} $\iota (z_0)=(\cdots, a_{-1}, a_0, a_1,
\cdots)$ where $a_i =a(z_i)$.
Orbits that pass through ${\cal C}$ have exactly two itineraries
as $T^i {\cal C} \cap {\cal C}=\emptyset$ for all $i$.
 
We would like to show that the symbol 
sequence
$\iota(z_0)$ uniquely determines $z_0$. This may fail in a trivial way:
Let $I_i=[x_i, x_{i+1}]$. Then our coding is clearly
not unique if for some $i$, $f(I_i)$ wraps all
the way around the circle, meeting some $I_j$ more than once.
For simplicity of exposition {\it we will assume this does not happen}.
If it does, it suffices to consider the partition on $\Omega$ whose elements 
correspond to  the connected components of $I_i \cap f^{-1}I_j$.
%%%%%%%%%%%%%%%%%%%%%%%%%%%%%%%%Section 9.2
\subsection{Coding of monotone branches}
\label{s9.2} 

\noindent {\bf Coding of monotone segments of $\partial R_k$.} \ 
Observe that points in $\partial R_k$ also have well-defined $a$-addresses
in the spirit of Lemma \ref{nlemma9.1}: if 
$z \in \partial R_k \cap {\cal C}^{(k)}$,
then its location
with respect to $\Gamma_k$ is obvious (except when
$z \in \Gamma_k$). This allows us to assign in a unique way a
$k$-block $[a_{-k}, \cdots, a_{-1}]$ to each monotone segment $\gamma$ of
$\partial R_k$. We write $\iota(\gamma)=[a_{-k}, \cdots, a_{-1}]$.

\medskip
\noindent {\bf Coding of monotone branches of $R_k$.} \ Each $S \in 
{\cal T}_k, \ k>0$, is associated with a block $\iota (S)=
[a_{-k}, \cdots, a_{-1}]$ defined inductively as follows:
Let $S \in {\cal T}_{k-1}$ be such that $\iota (S)=
[a_{-(k-1)}, \cdots, a_{-1}]$. If $S \cap {\cal C}^{(k-1)} =\emptyset$,
then it lies between two components of ${\cal C}$, say
${\cal C}_i$ and ${\cal C}_{i+1}$, and
$\iota (T(S)):=[a'_{-k}, \cdots, a'_{-1}]$ where $a'_{-1}=i$ and  $a'_{-j}=a_{-j+1}$ for $j>1$. If $S \cap {\cal C}^{(k-1)} \neq \emptyset$,
then $S=\hat S_1 \cup \cdots \cup \hat S_n$ where $\hat S_1$ is the
section of $S$ from one end to the middle of the first $Q^{(k-1)}$ 
that it meets, $\hat S_2$ is the section from the middle of this
$Q^{(k-1)}$ to the middle of the next component of ${\cal C}^{(k-1)}$ 
etc., and the $a'_{-1}$-entry of $\iota(T(\hat S_j))$ is defined
according to the location of $\hat S_j$. Note that
this coding of branches in ${\cal T}$ is injective, i.e. 
$S \neq S'$ implies $\iota(S) \neq \iota(S')$, and that if $\gamma$
and $\gamma'$ are monotone segments that bound $S$, then 
$\iota(\gamma)=\iota (\gamma') =\iota (S)$. Note also that the replacement
procedure in Sect. \ref{s8.3} corresponds to replacing 
$[a_{-k}, \cdots, a_{-1}]$
by blocks of the form $[*, \cdots, *, a_{-k}, \cdots, a_{-1}]$. 

\medskip
\noindent {\bf Coding of arbitrary points in $R_0$.} \ 
For points in certain locations of $R_0$, there is no meaningful way
of assigning to it an address as we did in Sect. \ref{s9.1}. Instead, for
each $k \geq 0$, we define the $\tilde a^{(k)}$-address(es) of 
$z \in R_k$ as follows:
$\tilde a^{(k)}(z)$ has the obvious definition if $z \not \in {\cal C}^{(k)}$;
if $z=(x,y) \in Q^{(k)}$ for some $Q^{(k)} \subset {\cal C}^{(0)}_i$, we let  
$\tilde a^{(k)}(z)
=i$ if $x> \hat x -b^{\frac{k}{4}}$ where
$\hat z=(\hat x, \hat y)$ is one of the critical points in $\partial Q^{(k)}$;
$\tilde a^{(k)}(z)=i-1$ if $x< \hat x +b^{\frac{k}{4}}$. Clearly
$\tilde a^{(k)}$-addresses are not unique: an open set of points in the middle part of each $Q^{(k)} \subset
{\cal C}^{(0)}_i$ have as their $\tilde a^{(k)}$-addresses both $i-1$ and $i$. 

\medskip
We further introduce the following notation:
$$
\pi_\Omega ([a_{n}, a_{n+1}, \cdots, a_m])=\{z_0 \in \Omega: a(z_{i})=a_{i},
\ n \leq i \leq m\};
$$
$$
\pi_{R_0} ([a_{-k}, a_{-k+1}, \cdots, a_{-1}])=\{z_0 \in R_k: 
\tilde a^{(k-i)}(z_{-i})=a_{-i}, \ 1 \leq i \leq k\};
$$
``$\tilde a^{(k-i)}(z_{-i})=a_{-i}$" above means
$a_{-i}$ is an admissible $\tilde a^{(k-i)}$-address of $z_{-i}$.

%%%%Lemma 9.2 
\begin{lemma} 
\label{nlemma9.2} \ 
(i) Every $S \in {\cal T}_k$, $k \geq 1$, is $= \pi_{R_0}(\iota(S))$ and 
contains a neighborhood 

\ \ of
$\pi_{\Omega}(\iota(S))$. 

\noindent (ii) Given $z_0 \in \Omega$ and $n \in {\mathbb Z}^+$, 
there exists $k$ with
$n \leq k \leq n(1+3\theta)$ and $S=
S(z_0, n)$

\ \ $\in {\cal T}_k$ such
that $z_0 \in \pi_\Omega(\iota(S))$.
\end{lemma}

\noindent {\bf Proof:} That $S =  \pi_{R_0}(\iota(S))$ follows inductively 
from the definitions of these two objects. That $S$ contains a neighborhood of 
$\pi_{\Omega}(\iota(S))$ is also obvious inductively. For (ii), we know from 
Theorem ${\rm 5}^{\prime}$ that there exists $S \in {\cal T}_n$ with 
$z_0 \in S$.
The only way one can have $z_0 \not \in \pi_{\Omega}(\iota(S))$ is that at the 
time $S$ is created, say at time $k-i$,  
$T^{-i} S$ meets the mid $b^{\frac{k-i}{4}}$-section $E$ of some 
$Q^{(k-i)}$ and extends to the left of $E$, while $z_{-i} \in E$ and lies to 
the ``right" of $\Gamma \cap Q^{(k-i)}$ in the sense of Lemma \ref{nlemma9.1}.
Let $S_0$ be the ancestor of $S$ in ${\cal T}_{k-i}$, and let $S_1$ be the
descendant of $S_0$ that contains the right half of $Q^{(k-i)}$. 
Our replacement procedure guarantees that there exists 
$S^{\prime} \in {\cal T}$ that is either a
descendant of $S_1$ or a replacement for a descendent of $S_1$
which contains $z_0$.
\hfill $\square$

\bigskip
Let
$$
\Sigma \ := \ \{ {\bf a} =(a_i)_{i=-\infty}^\infty :
\iota(z_0)={\bf a} \ {\rm for \ some} \ z_0 \in \Omega \},
$$
 
\noindent and let $(\sigma{\bf a})_i = ({\bf a})_{i+1}$ denote the shift 
operator.
It is easy to check that $\Sigma$ is a closed subset of $\Pi_{-\infty}^\infty
 \{1,2, \cdots, r\}$ with $\sigma^{-1}\Sigma \subset \Sigma$. 
Extending our definition of $\pi_{\Omega}$
to infinite sequences and writing $\pi = \pi_{\Omega}$, we have that 
$\pi({\bf a})$ is
the set of all points $z_0 \in \Omega$ with $\iota(z_0)={\bf a}$.
The following proposition, whose proof occupies all of the next subsection, 
completes the proof of Theorem \ref{theorem4}.
%%%%Proposition 9.1 
\begin{proposition} For every ${\bf a} \in \Sigma$,  $\pi({\bf a})$
consists of exactly one point, and $\pi:\Sigma \to \Omega$ is a
continuous mapping.
\label{nproposition9.1}
\end{proposition}

Let $B(z_0, \varepsilon)$ denote the ball of radius $\varepsilon$ centered
at $z_0$, and let us say $S \in {\cal T}_k$ is {\it compatible} with 
${\bf a} = (a_i)$ if $\iota(S) = [a_{-k}, \cdots, a_{-1}]$. 
Proposition \ref{nproposition9.1} follows immediately from  
Lemma \ref{nlemma9.2}(i) and Proposition ${\rm 10.1}^{\prime}$ below.

\bigskip

\noindent {\bf Proposition $\rm{\bf 10.1}^{\prime}$} \ \ {\it Given 
${\bf a} \in \Sigma$, $z_0 \in \pi({\bf a})$, and 
$\varepsilon >0$, there exists $S \in 
{\cal T}_{n+m}$ compatible with $\sigma^n {\bf a}$ such that 
$T^{-n}S \subset B(z_0, \varepsilon)$.}

%%%%%%%%%%%%%%%%%%%%%%%%%%%%%Section 9.3
\subsection{Uniqueness of point in $\Omega$ corresponding to each itinerary}
\label{s9.3} 

We begin with a situation that resembles that in $1$-dimension.

\begin{lemma}  Let ${\bf a}, z_0$ and
$\varepsilon$ be as in Proposition ${\rm 10.1}^{\prime}$.
Suppose that for some $k$, the component of
$R_k \cap B(z_0, \varepsilon)$ containing $z_0$, which we denote
by $H$, is bounded by two $C^2(b)$ subsegments $\gamma$ and $\gamma'$ of 
$\partial
 R_k$ cutting across $B(z_0, \varepsilon)$ as shown
with
$$
{\rm Hausdorff \ distance}\ (\gamma, \gamma')<\varepsilon^{10}.
$$
Then there exists $S \in {\cal T}_{n+m}$ compatible with $\sigma^n {\bf a}$ 
such that
$T^{-n}S \subset H$.
\label{nlemma9.3}
\end{lemma}

\begin{picture}(8, 5)
\put(3,0){
\psfig{figure=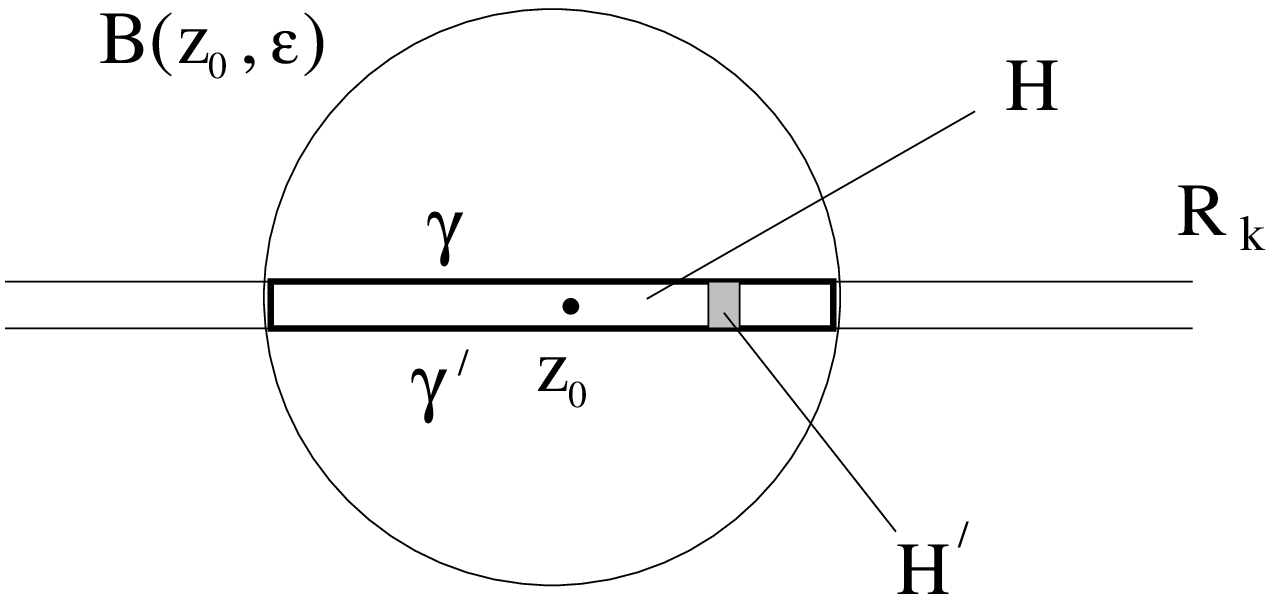,height=4cm,width=8cm}
}
\end{picture}

\bigskip

\centerline{Figure 8 \ \ The situation considered in Lemma \ref{nlemma9.3}}

\bigskip

\noindent {\bf Proof:} Our plan of proof is as follows. 
Since $z_0 \in \Omega$, we have, for $i = 1, 2, \cdots$, a monotone branch
$S_i \in {\cal T}_{k+n_i}$, $n_i > i$, such that $z_i \in S_i$. Then
$T^{-i}S_i \subset R_k$ for all $i$. We need to show that for $i$ 
sufficiently large, $T^{-i}S_i \subset H$.
To prevent $T^{-i}S_i$ 
from exiting $B(z_0, \varepsilon)$
via the right boundary of $H$, we will show that
for some section $H' \subset H$ as shown and $k'>0$,
$T^{k'}(H')$ is a component of ${\cal C}^{(k+k')}$, so that the left 
and
right boundaries of $H'$ have incompatible $\tilde a^{(k+k')}$-addresses. 
Assuming $i>k'$, it will follow (using Lemma \ref{nlemma9.2}) that $T^{-i}S_i$
cannot meet both the left and right boundaries of $H'$.
Being connected and contained in $R_k$, $T^{-i}S_i$
must meet
both boundaries of $H'$ in order to exit 
$B(z_0, \varepsilon)$ from the right.
The left boundary of $H$ is blocked off similarly.

The proof that $T^{k'}H$ crosses a component of ${\cal C}^{(k+k')}$
for some $k'$ is similar to that of Sublemma \ref{nsublemma8.1}, but there
are two differences: initially at least, we do not know the lengths of
$T^j\gamma$ relative to their distances to the critical set, and  
we must control the shearing between $\gamma$ and $\gamma'$ as we iterate.  
Details of the proof follow.

 \medskip
Consider first the case where $z_0 \not \in {\cal C}^{(0)}$.
Let $\gamma_0$ be a subsegment of $\gamma$ of length
$\frac{\varepsilon}{2}$ located half-way between $z_0$ and the right boundary 
of $H$. We first describe how to locate $\gamma_0 \cap H^{\prime}$. Let $n_1$ be
the first time when $T^i(\gamma_0)$ meets ${\cal C}^{(0)}$. 
If $T^{n_1} \gamma_0$ contains an $I_{\mu j}$ of full length, then we let 
$\gamma_1 \subset T^{n_1} \gamma_0$ correspond to the longest $I_{\mu j}$
or segment outside of ${\cal C}^{(0)}$, whichever is longer. If not, 
we let $\gamma_1 = T^{n_1} \gamma_0$. In both cases, we let $n_2$ be the 
first time when part of $T^{n_2 - n_1} \gamma_1$ makes a free return. Choose
$\gamma_2 \subset T^{n_2 -n_1} \gamma_1$ as before, let $n_3$ be the first time when 
part of $T^{n_3 - n_2} \gamma_2$ makes a free return, and so on. Using the fact 
that $\partial R_n$ is controlled (Proposition \ref{nproposition5.1}), we see
that the $\gamma_i$ increase in length, so that there exists 
some $i_0$ such that $\gamma_{i_0}$ contains an $I_{\mu j}$. From then on, the
argument in Sublemma \ref{nsublemma8.1} produces an $i_1$ such that 
$T^{n_{i_1} - n_{i_1-1}} \gamma_{n_{i_1-1}}$ traverses a component of 
${\cal C}^{(0)}$. 

We now proceed to construct $H^{\prime}$.
Letting $\tau_0$ denote unit tangent vectors to $\gamma$, we have that $\|DT^i(\xi_0)\tau_0\| \geq c>0$ for all $\xi_0 \in \gamma_0$ and 
$i \leq n_1$. Through each $\xi_0 \in \gamma_0$, therefore, is 
a stable curve of order $n_1$ connecting $\xi_0$ to a point in $\gamma'$
less than $\varepsilon^9$ away (see Sect. \ref{s2.2} and Lemma \ref{nlemma2.9}).
Let $H_0$ be the region between $\gamma$ and $\gamma'$
made up of the union of these stable curves. 

Since we do not know how close $T^{n_1}\gamma_0$ gets to the critical set, 
we cannot continue to claim the expanding property of $\tau_0$
beyond time $n_1$.
Instead, we observe that for $\xi_0 \in \gamma_1$, 
$\|DT^j(\xi_0){\tiny (\!\!\begin{array}{l} 0 \\ 1 \end{array}\!\!)}\| \geq 1$ 
for $j \leq n_2-n_1$ so that through
each $\xi_0 \in \gamma_1$, there is a stable curve of order $n_2-n_1$. 
Assuming that these stable curves meet $T^{n_1} \gamma'_0$, we define $H_1$ to be the region between $T^{n_1} \gamma_0$ and $T^{n_1} \gamma'_0$ spanned by these curves, and check that $H_1$ can be chosen 
to be a subregion of $T^{n_1}H_0$.
 
To justify the last sentence, observe first that if 
$\gamma_1 \subset I_{\mu_1 j}$, then $e^{-\mu_1} > 
\varepsilon$. This is true regardless of whether $\gamma_1=T^{n_1}\gamma_0$.
Second, since the contractive field $e_{n_2-n_1}$ near $\gamma_1$
makes angles $\sim e^{-\mu_1}$ with $T^{n_1}\gamma_0$ and with $T^{n_1}\gamma'_0$, every point in $\gamma_1$ is connected by a 
stable curve to a point in $T^{n_1}\gamma'_0$ not more 
than a distance of $(b^{n_1}
\varepsilon^9)/e^{-\mu_1} < b^{n_1}\varepsilon^8 <<e^{-\mu_1}$
away. This allows us to define $H_1$.
Finally, we may need to trim the edges of $H_1$ by a length $\sim b^{n_1}
\varepsilon^8$
in order to fit it inside $T^{n_1}H_0$. This is easily done since
$|\gamma_1|>\min(\varepsilon, \frac{1}{\mu_1^2}e^{-\mu_1})$.

At time $n_2$, we again do not know how close $T^{n_2 - n_1} \gamma_1$ is to 
the critical set, and so we use $\|DT^j {\tiny (\!\! \begin{array}{l}
0 \\ 1 \end{array} \!\!)} \| \geq 1$ for $j \leq n_3 - n_2$ to construct
new stable curves which are then used to construct $H_2$. Observe that 
compared to time $n_1$, the situation has improved: 
$|\gamma_2| \geq |\gamma_1|$, and the segments $\gamma_2$ and 
$\gamma_2^{\prime}$ are closer than before. We construct $H_3$, $H_4$, $\cdots$,
until time $n_i$, when $T^{n_{i_1} - n_{i_1-1}} H_{i_1-1} \supset Q$, a
component of ${\cal C}^{(k+n_{i_1})}$. Letting $k^{\prime} = n_{i_1}$ and
$H^{\prime} = T^{-n_{i_1}}(Q)$, the proof for the case 
$z_0 \not \in {\cal C}^{(0)}$ is complete.
 
For $z_0 \in {\cal C}^{(0)}\setminus {\cal C}$, let $j$ be such that
 $z_0 \in \hat Q^{(j)} \setminus Q^{(j)}$. If $k$ in the statement of
the lemma is $\geq j$, repeat the argument above with $n_1=0$. 
If not, replace $k$ by $j$ and $\varepsilon$ by $\min(\varepsilon,
\frac{1}{2}\rho^j)$ and let $n_1=0$. The case of $z_0 \in {\cal C}$
is dealt with similarly.
\hfill $\square$
 
\bigskip
Recall that for all $z_0 \in \Omega$, at every return to ${\cal C}^{(0)}$,
$z_i$ is h-related, and bound and fold periods are well defined.
(See Section \ref{s3} for definitions.)
 
\bigskip
 
\noindent {\bf Proof of Proposition $\rm{10.1}^{\prime}:$} \ 
Let ${\bf a} \in \Sigma$ and $z_0 \in \pi ({\bf a})$ be given. We wish to 
arrange for the scenario in 
Lemma \ref{nlemma9.3} at $z_0$, but it is not possible to do it directly 
when $z_0$ 
is near a ``turn". Intuitively, in order for $z_0$ to be near a ``turn", 
$z_{-i}$ must be near the critical set for some $i>0$. 
This motivates the following considerations.
 
\medskip \noindent
{\it Case 1.} \ There exists arbitrarily large $i$ such that
$d_{\cal C}(z_{-i})<\rho^k$ for $k \approx $ $K_0 (\log \|DT\|) \theta i$
where $K_0$ is to be specified shortly. \
Let $\varepsilon>0$ be given, and let $i$ and $k$ have the relationship
above with $\|DT\|^{-i}<\varepsilon^{10}$.
Let $j$ be such that $z_{-i} \in \hat Q^{(j)} \setminus Q^{(j)}$.
Then $j \geq k$. We wish to apply Lemma \ref{nlemma9.3} to $z'_0=z_{-i}$
with $\varepsilon'=\|DT\|^{-i}\varepsilon$ and $H$ bounded by 
$\partial R_j$. This result transported back to $z_0$ proves
the proposition.
To satisfy the hypotheses of Lemma \ref{nlemma9.3} at $z'_0$, it suffices to 
check
that the Hausdorff distance between the two horizontal boundaries of
$\hat Q^{(j)}$ is $<(\|DT\|^{-i}\varepsilon)^{10}$.
This is true provided $K_0$ is chosen to satisfy the inequality
$$
b^{\frac{k}{4}} = (b^{\frac{1}{4}K_0 \theta \log\|DT\|})^i
= (\|DT\|^{\log b^{\frac{1}{4}K_0  \theta}})^i < \|DT\|^{-11i}
=(\|DT\|^{-i} \varepsilon)^{10}.
$$

\smallskip \noindent 
{\it Case 2.} Not Case 1. \ Note that this means that $z_{-i}$ approaches
${\cal C}$ extremely slowly (if at all) as $i \to \infty$.
First we observe that with $d_{\cal C}(z_{-i})>>b^{\frac{i}{2}}$,
$z_0$ is out of all fold periods from the past.
To arrange for the scenario of Lemma \ref{nlemma9.3} at $z_0$, we will show:

\smallskip
\ \ (i) there exist $\kappa={\cal O}(1)$ and arbitrarily large $i$ such that 
$\|DT^j(z_{-i}){\tiny (\!\!\begin{array}{l} 0 \\ 1 \end{array}\!\!)}\|
 \geq \kappa^j$ 

\ \ \ \ \ for all $j \leq i$;

\ (ii) the stable curves near $z_{-i}$ when mapped forwards bring with them to $z_0$ a    

\ \ \ \ \ pair of curves from $\partial R_n$ with $z_0$ sandwiched in between;

(iii) these curves are $C^2(b)$, they have a minumum length 
$\varepsilon_1$ independent
of $i$ and

\ \ \ \ \  their Hausdorff distance can be made as small as need be
by choosing $i$ large.

\smallskip

We prove (i). Leaving the $\inf_i d_{\cal C}(z_{-i})>0$ case as an exercise,
we consider $i$ with $d_{\cal C}(z_{-i}) \leq d_{\cal C}(z_{-j})$
for all $0<j\leq i$. 
Suppose $d_{\cal C}(z_{-i}) \approx e^{-\mu}$, so that the ensuing
bound period is $>K^{-1}\mu$. 
Let $w_j=DT^j(z_{-i}){\tiny (\!\!\begin{array}{l} 0 \\ 1 \end{array}\!\!)}$,
and let $z_{-i+n}$ be the next free return.
Then $\|w_j\| \geq 1$ for $j \leq n$. We argue that $w_n$ splits correctly:
If $z_{-i+n} \in {\cal C}^{(n)}$, then $d_{\cal C}(z_{-i+n}) \geq
d_{\cal C}(z_{-i}) \approx e^{-\mu} >> b^{\frac{1}{20}K^{-1}\mu} 
\geq b^{\frac{1}{20}n}$; if $z_{-i + n} \not \in {\cal C}^{(n)}$, then it is
$\in {\hat Q}^{(j)} - Q^{(j)}$ for some $j < n$. In both cases, 
Lemma \ref{nlemma7.1} applies, and we have
$\|w^*_{n+1}\| \geq
e^{\frac{cn}{3}}e^{-\mu} \geq e^{(\frac{c}{3}-K)n}$. 
Since the situation at subsequent
free returns is clearly improved ($d_{\cal C}(\cdot) \geq d_{\cal C}(z_{-i})$
and the derivative
has built up), 
we have $\|w_j\| \geq e^{(\frac{c}{3}-K)j}$ for all $j \leq i$.
 
To prove (ii), suppose $z_{-i} \in {\hat Q}^{(k)}\setminus Q^{(k)}$
for some $k$.
We consider the stable curve of order $i$ through $z_{-i}$ and let
 $\zeta_0$ be its intersection with the upper boundary of $\hat Q^{(k)}$.
A subsegment $\gamma_0$ of this upper boundary centered at
$\zeta_0$ is constructed by iterating forward $i$ times and trimming whenever necessary
so that $T^j\gamma_0$ stays inside three consecutive $I_{\mu \ell}$
for all $j \leq i$. Clearly, stable curves of order $i$ can be constructed through
all points in $\gamma_0$, and these curves ``tie together" the two
subsegments of $\partial {\hat Q}^{(k)}$.
 
We leave it as an exercise to show the existence of 
$\varepsilon_1$ (which depends only on the slow rate of approach to 
${\cal C}$ in backward time). The curves brought in are sebsegments of 
$\partial R_{k+i}$ and they are out of all fold periods.
This completes the proof of Proposition $\rm{10.1}^{\prime}$. 
\hfill $\square$

%%%%%%%%%%%%%%%%%%%%%%%%%%%%%%%%%%%%%%%%%%%%%%%%%%%%%%%%%%%%%%%%%%
\subsection{Proof of Theorem 2(1)(iii)}
\label{s9.add}
We explain how $\Omega = 
\overline{\cup_{\varepsilon > 0} \Omega_{\varepsilon}}$
follows readily from the ideas in the last two subsections and the surjectivity
condition (*) in Sect. \ref{s1.2}.

In view of Proposition ${\rm 10.1}^{\prime}$, it suffices to show that every
$S \in {\cal T}$ contains a point in $\Omega_{\varepsilon}$ for some 
$\varepsilon > 0$. Recall the way monotone branches in 
${\cal T}$ are constructed. Given $S \in {\cal T}$, let $\ell > 0$ be the 
smallest integer such that $T^{-\ell}S \not \in {\cal T}$. Then
$T^{-\ell}S$
contains half of some $Q^{(k)}$. Let $H$ be the middle half of 
$T^{-\ell}S \cap Q^{(k)}$, with length $\frac{1}{4} \rho^k$. An argument 
similar to that 
in Lemma \ref{nlemma9.3} but carried on indefinitely in time gives a sequence 
of domains $H \supset H_1^{\prime} \supset H_2^{\prime} \supset \cdots$ and a
curve $\omega_0 \subset \cap_{n \geq 1} H_n^{\prime}$ with the following properties:

\medskip

--- $\omega_0$ connects the top 
and bottom boundaries of $Q^{(k)} \cap T^{-\ell}S$;
 
\smallskip

--- there exists $\varepsilon > 0$ such that  
$\forall z \in \omega_0, \ d_{\cal C}(z_n) \geq \varepsilon \ \forall n \geq 0$.

\medskip

To finish, it suffices to produce $\hat z_0 \in \omega_0$ such that 
$\hat z_{-i} \not \in {\cal C}^{(0)} \ \forall i > 0$. Let $D_i$ be the
component of $R_0\setminus {\cal C}^{(0)}$ between the $i$-th and
$(i+1)$-st components of ${\cal C}^{(0)}$, and let $\hat D_i$ be the 
union of $D_i$ with the two components of ${\cal C}^{(0)}$ adjacent to it. 
Then we may assume from condition (*) that for every $i$, there exists
$j$ such that $T(D_j) \cap {\hat D}_i$ contains a horizontal strip
traversing the full length of ${\hat D}_i$. Suppose $\omega_0
\subset {\hat D}_i$,  and let $j$ be as above. Then there is a subsegment 
$\omega_1 \subset \omega_0$ such that $T^{-1} \omega_1 \subset D_j$ and 
connects the top and bottom boundaries of $D_j$. Similarly, we produce for 
$n = 2, 3, \cdots$ segments $\omega_n \subset \omega_{n-1}$ 
such that $T^{-n} \omega_n$ is contained in some $D_{j(n)}$ and connects the two 
horizontal boundaries 
of $D_{j(n)}$. Let $\hat z_0 \in \cap_{n \geq 0} \omega_n$. 
\hfill $\square$ 

\subsection{Existence of Equilibrium states}
\label{s9.4}
This is a corollary to the symbolic dynamics we have developed.
Let $\varphi: R_0 \to {\mathbb R}$ be a continuous function,
and let $P(T; \varphi)$ denote the {\bf topological pressure} of $T$
for the potential $\varphi$. (See e.g. \cite{Wa}, Chapter 9, for
definitions and basic facts.)
A well known variational principle says that
$$
P(T; \varphi) \ = \ sup \ \ P_\nu(T;\varphi)
$$
where the supremum is taken over all $T$-invariant Borel probability
measures $\nu$ and
$$
P_\nu(T;\varphi):=h_\nu(T)+\int \varphi d\nu,
$$
where $h_{\nu}(T)$ denotes the metric entropy of $T$ with respect to $\nu$.
An invariant measure for which this supremum 
is attained is
called an {\bf equilibrium state} for $(T; \varphi)$.
 
Let $\sigma:\Sigma \to \Sigma$ and $\pi:\Sigma \to \Omega$ be as in
Theorem \ref{theorem4}.
 
\bigskip
 
\noindent {\bf Proof of Corollary \ref{ncoro2}:} \ 
Let $\varphi: R_0 \to {\mathbb R}$ be given.
We need to prove that there exists $\nu$ such that
$P_\nu(T;\varphi)=P(T; \varphi)$.
Let $\tilde \varphi$ be the function on $\Sigma$ defined by $\tilde
\varphi = \varphi \circ \pi$. Then 
$P(T; \varphi)=P(T|\Omega; \varphi|\Omega) \leq 
P(\sigma;
\tilde \varphi)$. Since $\sigma:\Sigma \to \Sigma$ has a natural finite 
generator
without boundary, $(\sigma, \tilde \varphi)$ has an
equilibrium state which we call $\tilde \nu$. Let $\nu=\pi_* \tilde \nu$.
 It suffices to show that
$P_\nu(T|\Omega; \varphi|\Omega)=P_{\tilde \nu}(\sigma; \tilde \varphi)$.
This follows from the fact that $\pi$ is one-to-one over $\Omega \setminus
\cup T^i{\cal C}$, and $\mu(\pi^{-1}(\cup T^i {\cal C}))=0$ for 
any $\sigma$-invariant probability
measure $\mu$
because $\sigma^i(\pi^{-1}{\cal C}) \cap \pi^{-1}{\cal C} = \emptyset$ for all 
$i \in {\mathbb Z}$. \hfill $\square$

\bigskip
Since the {\bf topological entropy} of $T$, written $h_{\rm top}(T)$, is
equal to $P(T;0)$, the discussion above gives immediately
 
\begin{corollary}
\label{ncorollary9.1} 
\begin{itemize}
\item[(i)] $T$ has an invariant measure of maximal entropy.
\item[(ii)] Let $N_n$ be the number of distinct blocks of symbols of
length $n$ that appear in $\Sigma$. Then
$$
\lim_{n \to \infty} \ \frac{1}{n} \log N_n \ = \ h_{\rm top}(T).
$$
\end{itemize}
\end{corollary}
%%%%%%%%%%%%%%%%%%%%%%%%%%%%%%%%%%%%%%%%%%%%%%%%%%%%%%%%
 
\subsection{Topological entropy}
\label{s9.5}
Topological entropy is, in general, defined in terms of open covers
of arbitrarily small diameters, $\varepsilon$-separated or spanning sets. None 
of the standard definitions is easy to compute with. 
Corollary \ref{ncorollary9.1} 
gives a concrete way to think about this invariant for the class
of dynamical systems under consideration.
Three other characterizations and estimates of geometric interest are
discussed here.
 
\medskip
Recall the notion of $\tilde a^{(k)}$-addresses for
$z \in R_k$ (see Sect. \ref{s9.2}). For $z_0 \in R_0$, we define
its (future) $\tilde a$-itinerary to be $(a_i)_0^{\infty}$ if
for each $i, \ \tilde a^{(i)}(z_i)=a_i$. These itineraries are clearly not 
unique. Let
 
\bigskip
\centerline{$\tilde N_n \ = \ $ the number of $n$-blocks
appearing in the $\tilde a$-itineraries of points in $R_0$,}
\bigskip
\noindent overcounting whenever ambiguities arise, that is, if an orbit
has $j$ different admissible $\tilde a$-itineraries of length $n$,
they will be counted as $j$ distinct blocks in $\tilde N_n$.
Obviously, $N_n \leq \tilde N_n$.

\begin{lemma}
\label{nlemma9.4}
$$
\limsup_{n \to \infty} \ \frac{1}{n} \log \tilde N_n \ \leq \
h_{\rm top}(T).
$$
\end{lemma}
 
\noindent {\bf Proof:} \ We fix some arbitrarily small $\varepsilon>0$, and 
choose $n_0$ 
so that
$$
\frac{1}{n_0} \log N_{n_0} < h_{\rm top}(T)+ \varepsilon
\ \ \ {\rm and} \ \ \  \frac{1}{n_0}\log (2n_0) <\varepsilon.
$$
Let $n_1>n_0$ be large enough that $b^{\frac{n_1}{10}}\|DT\|^{n_0}
<e^{-\beta n_0}$, so that no orbit segment in $R_0$ of length $\leq n_0$
can pass through the region $D:=\{\xi_0 \in {\cal C}^{(n_1)}: |\xi_0-\hat z_0|<
b^{\frac{n_1}{10}}$ for some $\hat z_0 \in {\cal C} \cap Q^{(n_1)}(\xi_0)\}$
more than once.
For each $z_0$, let $S_{z_0}=T^{-n_0} S(z_{n_0}, 2n_1)$ where
$S(z_{n_0}, 2n_1)$ is as 
in Lemma \ref{nlemma9.2}(ii). By part (i) of the same lemma, $S_{z_0}$ is 
a neighborhood of $z_0$. Let $n_2>n_1$ be such that
$R_{n_2} \subset \cup_{z_0 \in \Omega} S_{z_0}$. Define
 
\bigskip
\centerline{$\tilde N(n_2, n_2+n_0) \ = \ $ the number of
distinct blocks of $[a_{n_2}, \cdots, a_{n_2+n_0-1}]$ }
\centerline{that appear in the $\tilde a$-itineraries
of all points in $R_0$.}
 
\begin{claim}
\label{nclaim9.1}
$\tilde N(n_2, n_2+n_0) \leq 2n_0 N_{n_0}$.
\end{claim}
 
\noindent {\it Proof of Claim \ref{nclaim9.1}:} Let $\xi_0 \in R_0$, and let 
$(a_i)$ be
any one of its $\tilde a$-itineraries. Let $\xi_{n_2} \in
S_{z_0}$ for some $z_0 \in \Omega$, and let $\iota(z_0)=(b_i)$.
We compare the two blocks $[a_{n_2}, \cdots, a_{n_2+n_0-1}]$
and $[b_0, \cdots, b_{n_0-1}]$. The $i$th entry of the first block
is an $\tilde a^{(n_2+i)}$-address of $\xi_{n_2+i}$. Since $\xi_{n_2+n_0} \in 
S = S(z_{n_0}, 2n_1)$, it follows from Lemma \ref{nlemma9.2} that the $i$-th 
entry of the second block is an $\tilde a^{(n(S) - n_0 + i)}$-address of 
$\xi_{n_2 + i}$ where $n(S)$ is such that $S \in {\cal T}_{n(S)}$.
Since the indices in both of these $\tilde a$-addresses exceed $n_1$,
they may differ only if $\xi_{n_2+i} \in D$.
This can happen at most once in the time period in question.
In other words,  $[a_{n_2}, \cdots, a_{n_2+n_0-1}]$ and $[b_0, \cdots, b_{n_0-1}
]$ can differ in at most one entry, and the difference
is either $+1$ or $-1$. Since $[b_0, \cdots, b_{n_0-1}]$
is one of the sequences counted in $N_{n_0}$, the claim is proved.

\hfill $\diamondsuit$
 
\bigskip
Similar reasoning shows that $\tilde
N(n_2+k n_0, n_2+(k+1)n_0) \leq 2n_0 N_{n_0}$ for all $k \geq 0$,
giving
$$
\tilde N_{n_2+k n_0} \ \leq \ K^{n_2} \cdot (2n_0 N_{n_0})^k.
$$
This combined with the properties we imposed on $n_0$ at the beginning
of the proof
gives the desired inequality.
\hfill $\square$
 
\bigskip
 
To complete the proof of Theorem \ref{theorem5}(i), recall that $P_n$ is the
number of fixed points of $T^n$ in $\Omega$.
 
\begin{lemma}
\label{nlemma9.5}
$$
\lim_{n \to \infty} \ \frac{1}{n} \log P_n \ = \ h_{\rm top}(T).
$$
\end{lemma}
 
\noindent {\bf Proof:} \ Since no point in ${\cal C}$ is periodic, 
there is a one-to-one
correspondence between the fixed points of $T^n$ and the periodic symbol 
sequences
of period $n$ in $\Sigma$,
proving ``$\leq$" in the lemma. That
$$
\liminf_{n \to \infty} \ \frac{1}{n} \log P_n \ > \ h_{\rm top}(T)
-\varepsilon
$$
for every $\varepsilon>0$ follows from a general theorem of Katok
for all $C^2$ surface diffeomorphisms \cite{K}. \hfill $\square$
 
\bigskip
 
Perhaps the most concrete geometric quantity of all is
the rate of growth of the number of monotone segments of a curve such as
$\partial R_0$. Our next lemma compares this
growth rate to the topological entropy of $T$.
Let $\partial R_0^+$ and $\partial R_0^-$ denote the two components of
$\partial R_0$, and define
 
\bigskip
\centerline{$M_n^{\pm}\ = \ $ the number of monotone segments in $\partial
R_n^{\pm}$}
 
\bigskip
\noindent where ``monotone segments" are as defined in Sect. \ref{s8.1}.
 
\bigskip

\noindent {\bf Proof of Theorem \ref{theorem5}(ii):} \
First we prove $M_n^{\pm} \leq \tilde N_n$. This 
follows from the fact that
for every monotone segment $\gamma$ in
$\partial R_n^{\pm}$,  $\iota (\gamma)$ is counted in 
$\tilde N_n$, and the mapping $\gamma \mapsto \iota (\gamma)$ is injective.
 
To prove the second inequality, we associate to each $n$-block
$[a_{-n}, \cdots, a_{-1}]$ that appears in $\Sigma$ first a point
$z_0 \in \Omega$ with $a(z_{-i})=a_{-i}$ and then a monotone branch 
$S = S(z_0, n)$
as in Lemma \ref{nlemma9.2}.
Then $S \in {\cal T}_k$ for some $k$ with
$n \leq k \leq n(1+\varepsilon_0)$,  
$\varepsilon_0 = 3(\log \frac{1}{b})^{-1}$.
We remarked at the end of Sect. \ref{s8.3} that every $S \in {\cal T}$ has a 
boundary component $\gamma^+$ in $\partial R_k^+$ and one 
in $\partial R_k^-$. 
We have thus defined, for each fixed $n$, a mapping from the set of 
$n$-blocks in $\Sigma$
to the set of monotone segments of $\partial R_k^+, \
n \leq k \leq n(1+\varepsilon_0)$.
This mapping is clearly injective since 
$\iota(\gamma^+) = \iota(S) = [*, \cdots, *, a_{-n}, \cdots, a_{-1}]$, proving
$$
N_n \ \leq \ \sum_{n \leq k \leq n(1+\varepsilon_0)} M_k^+.
$$
From this one deduces easily that
$$
\lim \frac{1}{n} \log N_n \ \leq \
(1+\varepsilon_0) \ \liminf  \frac{1}{(1+\varepsilon_0)n} \log M^+_{(1+
\varepsilon_0)n}.
$$
\hfill $\square$

\footnotesize
\newpage
\appendix
\setcounter{section}{0}
\renewcommand{\thesection}{Appendix \Alph{section}}
%%%%%%%%%%%%%%%%%%%%%%%%%%%Appendix A
\section{Examples}
\label{appendix-A}

\renewcommand{\thesection}{\Alph{section}}
\subsection{Attractors arising from interval maps including 
\\ the H\'enon attractors}
\label{appendix-A.1}

\noindent {\bf Reduction of Theorem \ref{theorem8} to Theorems 
\ref{theorem1}--\ref{theorem5}:} \ 
Let $I_0$ be a closed interval such that $f(I) \subset$ int$(I_0)
\subset I_0 \subset$ int $(I)$, and let $J_1$ and $J_2$ be the two 
components of $I \setminus I_0$. Choosing $b_0<<|J_1|, |J_2|$, one obtains easily from the formulas 
for $T_{a,b}$ in Sect. 1.1 that there exist $K>0$ and $\hat \Delta:=
[a_0, a_1] \times (0,b_0]$ such that for all $(a,b) \in \hat \Delta$, 
$T_{a,b}$ maps $R:= I \times [-Kb, Kb]$ strictly into 
$I_0 \times [-Kb, Kb]$. 

Our plan is to replace $\partial I \times [-Kb, Kb]$ by two curves
$\omega_1$ and $\omega_2$ so that each $\omega_i \subset J_i 
\times [-Kb, Kb]$, joins the top and bottom boundaries of $R$, 
and lies on the stable curve of a periodic orbit. We may assume that 
these periodic orbits stay outside of ${\cal C}^{(0)}$.
Replacing $R$ by $R_0$, the subregion of $R$ bounded by $\omega_1$
and $\omega_2$, the situation is now virtually indistinguishable from 
that of the annlus maps treated in Theorems \ref{theorem1}--\ref{theorem5}:
the top and bottom
boundaries of $R_0$ play the role of $\partial R_0$ 
in the previous situation, and the left and right boundaries 
shrink exponentially as we iterate. (There are small differences, such as
the existence of monotone branches with one end bounded by images of
$\omega_i$. These differences are inessential.)

To produce $\omega_1$ and $\omega_2$, we claim that pre-periodic points 
of $f$ are dense in $I$. This claim is justified as follows.
 First, Misiurewicz maps have no homtervals, 
so that there is a coding of the orbits of $f$ by a subshift 
$\sigma: \Sigma \to \Sigma$ with the property that each element of $\Sigma$ 
corresponds to the itinerary of exactly one point in $I$.
Second, $\Sigma$ is the closure of $\cup_n \Sigma_n$ where 
$\{\Sigma_n\}$ is an 
increasing sequence of subshifts of finite type, and third, 
pre-periodic points are dense in shifts of finite type. 

To finish, we fix pre-periodic points $p_1$ and $p_2$ of $f$
near the middle of $J_1$ and $J_2$. Shrinking
$\hat \Delta$ if necessary, we may assume that for $T_{a,b}$ with 
$(a,b) \in \hat \Delta$, the periodic orbits related to 
$p_1$ and $p_2$ persist and the stable curves through 
the continuation of $p_i$ have the desired properties. 
This is possible because the slopes of these
stable curves are bounded away from zero (see Lemma \ref{nlemma2.9}(a)).

\hfill $\square$

\bigskip

\noindent {\bf Proof of Corollary 3:} \ For the quadratic family,
the transversality condition in
Step II in Sect. 1.1 hold at all Misiurewicz points [T]. 
The nondegeneracy condition in Step IV is obviously satisfied. 
(To ensure that $f(I) \subset$ int$(I)$ for some $I$ in the case 
$a^*=2$, consider $a$ slightly less than $2$.)
\hfill $\square$ 

%%%%%%%%%%%%%%%%%%%%%%%%%%%%%%%%%%
\subsection{Homoclinic bifurcations}
\label{appendix-A.2}
We verify here the conditions in Sect. \ref{s1.1} and condition (**) in Sect. 
\ref{s1.2} for homoclinic bifurcations in 2-dimensions, setting the stage to 
apply Theorems \ref{theorem1}--\ref{theorem5}. See Sect. \ref{s1.5} for a
more detailed description of the bifurcation in question.

Following \cite{PT}, pages 47-51, we assume that linearizing coordinates 
have been
chosen in which $g_{\mu}, \ \mu \in [0, \mu^*]$, has the following properties:
\begin{itemize}
\item[(i)] On $\{ |\xi|, |\eta| < 2 \}$, $g_{\mu}$ is the linear map
$$
g_{\mu}(\xi, \eta) = (\sigma_{\mu} \xi, \lambda_{\mu} \eta)
$$
where $0 < \lambda_{\mu} < 1 < \sigma_{\mu}$, $\lambda_{\mu} \sigma_{\mu}<1$, 
and $\lambda_{\mu}, \sigma_{\mu}$ depend continuously on $\mu$.

\item[(ii)] There exists $N \in \mathbb Z^+$ such that $g_0^N$ maps the point
$(1, 0)$ to $(0, 1)$, carrying the unstable curve at $(1, 0)$ to a curve making
a quadratic tangency with the stable curve at $(0, 1)$. Near $(1, 0)$, 
$g_{\mu}^N$ has the form 
\begin{equation}
g_{\mu}^N(\xi, \eta) = ( \alpha (\xi-1)^2 + \beta \eta + \gamma \mu 
+ H_1(\mu, \xi, \eta), \ 1+
H_2(\mu, \xi, \eta))
\end{equation}
where $\alpha, \beta, \gamma \not = 0$ are constants. Furthermore, we have that
at $(\mu, \xi, \eta)
=(0,1,0)$, $H_1 = H_2 = 0$,
$\partial_{\xi} H_1 = \partial_{\eta} H_1 = \partial_{\mu}
H_1 = 0$ and $\partial_{\xi \xi} H_1 = \partial_{\xi \mu} H_1 = 
\partial_{\mu \mu} H_1 
=0$.
\end{itemize}

\begin{picture}(6, 5.5)
\put(4,0){
\psfig{figure=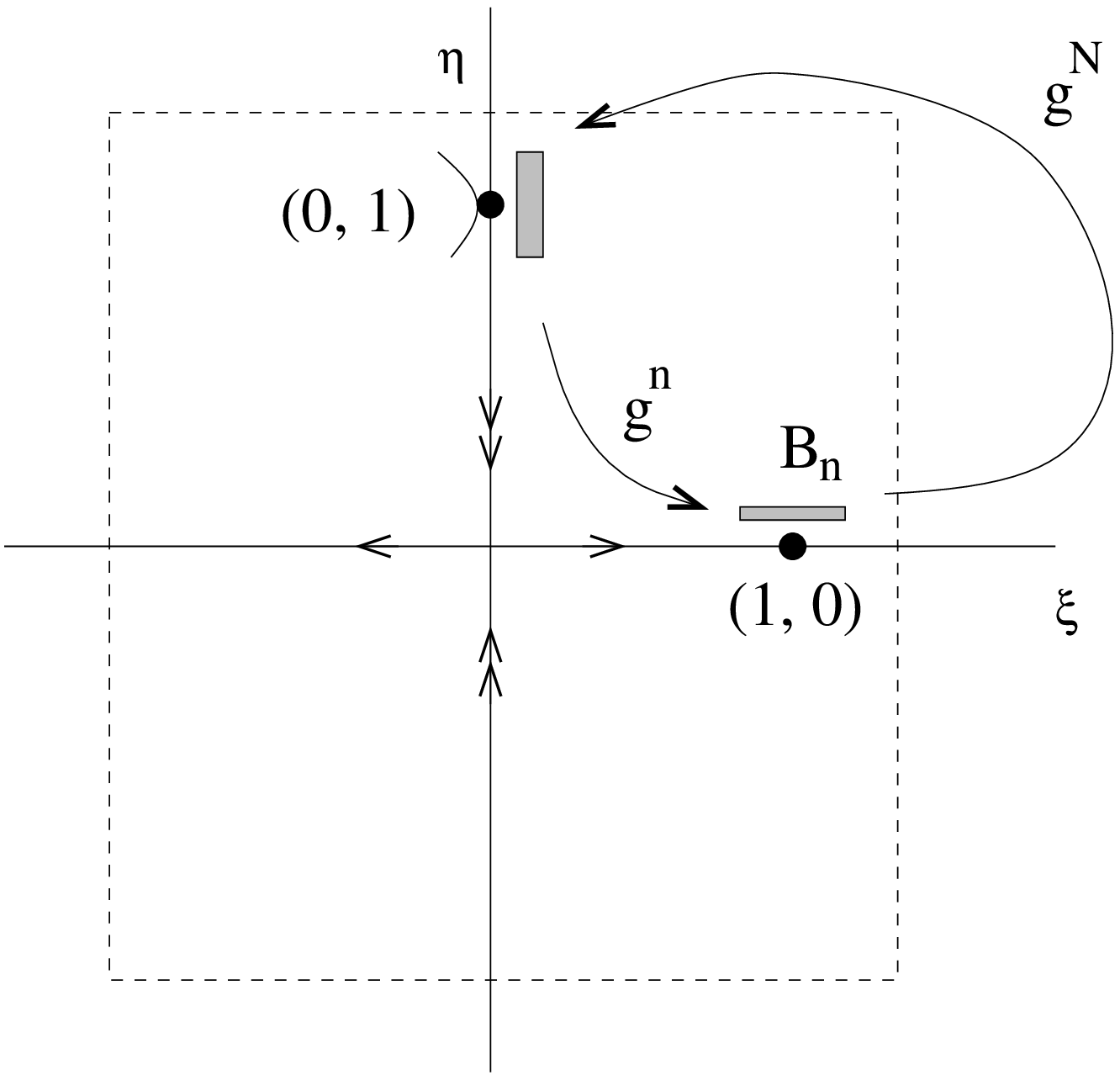,height = 5cm, width = 5cm}
}
\end{picture}

\medskip

\centerline{Figure 9 \ \  Attractors arising from homoclinic bifurcations}

\bigskip

It is not hard to see that for each fixed $n$, $n$ large, there exist a
box $B_n$ (with ${\rm diam}(B_n) \to 0$ as $n \to \infty$) and a range of parameters $\mu$ (also depending on  $n$) such that 
$(g^n_{\mu} \circ g^N_{\mu})(B_n) \subset B_n$. 
The attractors of interest to 
us have $(n+N)$ components permuted cyclically by $g_\mu$, with one of 
these components residing in $B_n$. 

To maneuver $g^n \circ g^N$ into the setting in Sect. \ref{s1.1}, we apply
the coordinate transformation $\Phi = \Phi_2 \circ \Phi_1$ where 
$$
\Phi_1(\xi, \eta) = (\xi-1, \eta -\lambda^n), \ \ \ \
\Phi_2(\xi, \eta) = (- \frac{\sigma^n}{a} \xi, -\frac{\sigma^{2n}}{a} \eta ).
$$
The purpose of $\Phi_1$ is to shift the center of $B_n$ to the origin. The
map $\Phi_2$ magnifies the attractor to unit length; its scaling in the 
$\eta$-direction is chosen with the standard quadratic family in mind. 
A straightforward computation yields
$$
T:=\ \ \Phi \circ g^n \circ g^N \circ \Phi^{-1}: \ \ 
\left(\begin{array}{c} x \\ y \end{array}\right) 
\mapsto \left(\begin{array}{c}
\frac{1}{a}[\sigma^n - \sigma^{2n}(\lambda^n + \mu)] - ax^2 + y - 
\frac{\sigma^{2n}}{a} H_1(\mu, \Phi^{-1}(x, y)) \\ 
-\frac{\sigma^{2n}}{a} \lambda^n H_2(\mu, \Phi^{-1}(x, y))
\end{array}\right).
$$
Letting $a = \Psi(\mu):= \sigma^n - \sigma^{2n}(\lambda^n + \mu)$ and
$\tilde H_i(a, x, y) := H_i(\mu, \Phi^{-1}(x, y)), \ i = 1, 2$, we have
$$
T: \ \ \left( \begin{array}{c} x \\ y \end{array} \right)
\mapsto \left( \begin{array}{c}
1- ax^2 + y - \frac{\sigma^{2n}}{a} {\tilde H}_1(a, x, y) \\
-\frac{\sigma^{2n}}{a} \lambda^n {\tilde H}_2(a, x, y) \end{array} \right).
$$
Since $\mu = \sigma^{-n} - a \sigma^{-2n} - \lambda^n$,
 the range of $a$ of interest to us, namely $a \in [1.5, 2)$ 
(see Appendix \ref{appendix-A.1}), corresponds to a subset of
 $(0, \mu^*]$ for $n$ large.

What we have so far is a 1-parameter family $\{ T_a \}$, which we regard
as defined on $U:=\{|x|,|y|<2\}$. The role of $b \to 0$ here is played by 
$n \to \infty$. Our next task is to choose $b$ (as a 
function of $n$) in such a way that $T_{a, b}$ has the form
$$
T_{a, b} : \ \ \left( \begin{array}{c} x \\ y \end{array} \right)
\mapsto \left( \begin{array}{c} 1- ax^2 + y + b u \\
                                bv \end{array} \right)
$$
where $u = u(a, x, y)$ and $v = v(a, x, y)$ have uniformly 
bounded $C^3$-norms. This will put us in the setting of Theorem 8
(see the proof of Corollary 3).

We begin by examing the $C^3$-norms of $\sigma^{2n} {\tilde H}_1$ and $\sigma^{2n} \lambda^n {\tilde H}_2$. Using the facts that the 
leading terms in $H_1$ are $\eta (\xi-1 + \eta + \mu)$,
and that $|\xi| < 3 \sigma^{-n}$ and $|\eta| < 3 \sigma^{-2n}$ for
$(\xi, \eta) \in \Phi^{-1}(U)$, we have $\|{\tilde H}_1\|_{C^0} =
{\cal O}(\sigma^{-3n})$. Similarly, $\|{\tilde H}_2 \|_{C^0} = 
{\cal O}(\sigma^{-n})$. Let $\partial^i, i = 1, 2, 3$, denote any one of the 
$i$-th partial derivatives. Using again the special form of $H_1$ and the nature of the coordinate transformations $\Phi$ and $\Psi$, we have 
$\| \partial^i {\tilde H}_1 \|  = {\cal O}(\sigma^{-3n})$ and
$\| \partial^i {\tilde H}_2 \| = {\cal O}(\sigma^{-n})$. Together this gives
$$
\| \sigma^{2n} {\tilde H}_1 \|_{C^3} < K \sigma^{-n}, \ \ \
\| \sigma^{2n} \lambda^n {\tilde H}_2 \|_{C^3} < K (\sigma \lambda)^n.
$$ 
The following choices of $b$ therefore give the desired result:

\medskip

If $\sigma^2 \lambda \leq 1$, let $b = \sigma^{-n}$.

If $\sigma^2 \lambda \geq 1$, let $b = (\sigma \lambda)^n$.

\medskip

This completes the verification of the conditions in Sect. \ref{s1.1} for 
the family $\{T_{a, b}\}$. 
We finish with the observation that
all the results in Section 1 that assume (**) are valid in the present 
setting: In the case $\sigma^2 \lambda \leq 1$,  
$|\det(DT)| \sim b$, so (**) is satisfied. When $\sigma^2
\lambda \geq 1$, 
$|\det(DT)| \sim (\sigma \lambda)^n = b^{\eta}$ where  
$\sigma^{-1} = (\sigma \lambda)^{\eta}$. This is condition 
${\rm (**)}^{\prime}$, a variant of (**) discussed in Sect. \ref{s7.2}

\vskip .5in

\renewcommand{\thesection}{Appendix \Alph{section}}
%%%%%%%%%%%%%%%%%%%%%%%%%%%%Appendix B
\section{Computational Proofs}
\label{appendix-B}

\renewcommand{\thesection}{\Alph{section}}

%%%%%%%%%%%%%%%%%%%%%%%%%%%%%%Appendix B.1
\subsection{Linear algebra (Sect. \ref{s2.1})}
\label{app-B.1}

%%%%%Sublemma B.1
\begin{sublemma}
Let $e$ be a unit vector in the most contracted direction of  
$$
M = \left(  \begin{array}{cc}
               A & C \\
               B & D \end{array} \right)
$$
with $ \| M e \| = \lambda^{min}$.
Then
\begin{eqnarray}
e & = & \pm \frac{1}{\rho}( C^2 + D^2 - (\lambda^{min})^2,
\  -(AC+BD) ) \ ,
\label{nB.1(1)} \\
M e & = & \pm \frac{1}{\rho}(-A(\lambda^{min})^2 + D \det(M), 
\ -B (\lambda^{min})^2 - C \det(M))
\label{nB.1(2)}
\end{eqnarray}
and
$$
(\lambda^{min})^2 = \frac{1}{2}(A^2+B^2+C^2+D^2 - \sqrt{(A^2+B^2+C^2+D^2)^2
- 4 (\det(M))^2})
$$
where $\rho$ is the normalizing constant in (\ref{nB.1(1)}).
\label{nsublemmaB.1}
\end{sublemma}

The proof is left as an easy exercise.

\bigskip

\noindent {\bf Proof of Lemma \ref{nlemma2.1}:} Let $O_1$ and $O_2$ be 
orthogonal 
matrices such that
$$
O_2 M^{(i-1)} O_1 = \left(  \begin{array}{cc}
               \lambda^{min}_{i-1} & 0 \\
               0 & \lambda^{max}_{i-1} 
             \end{array} \right).
$$
Then the tangent of the angle between $e_{i-1}$ and $e_i$ is given by
the slope of the most contracted direction of the matrix
$$
M_i O^{-1}_2  \left(  \begin{array}{cc}
               \lambda^{min}_{i-1} & 0 \\
               0 & \lambda^{max}_{i-1} 
             \end{array} \right)  := \left(  \begin{array}{cc}
               A & C \\
               B & D
             \end{array} \right) \left(  \begin{array}{cc}
               \lambda^{min}_{i-1} & 0 \\
               0 & \lambda^{max}_{i-1} 
             \end{array} \right) = \left(  \begin{array}{cc}
               \lambda^{min}_{i-1}A & \lambda^{max}_{i-1}C \\
               \lambda^{min}_{i-1}B & \lambda^{max}_{i-1}D
             \end{array} \right).
$$
From Sublemma \ref{nsublemmaB.1}, we see that the slope in question is equal to
$$
\frac{(AC+BD)\lambda^{min}_{i-1} \lambda^{max}_{i-1}}{(C^2+D^2) (\lambda^{max}_{i-1})^2 -(\lambda_{i}^{min})^2} \  .
$$
This is $ \leq \left( \frac{Kb}{\kappa^2} \right)^{i-1}$ because 
$\lambda^{min}_{i-1} \lambda^{max}_{i-1}=|\det(M^{(i-1)})
|
<b^{i-1}$, $\lambda^{min}_i < (\frac{b}{\kappa})^i$ and $(C^2+D^2) (\lambda^{max}_{i-1})^2
>K^{-1}\kappa^{2(i-1)}$, the last inequality being a 
consequence of the fact that
$\| M^{(i)} \|>\kappa^i$ and $(A^2+B^2) (\lambda^{min}_{i-1})^2
<K(\frac{b}{\kappa})^{2(i-1)}$. \hfill $\square$

\bigskip

Before giving the proof of Corollary \ref{ncorollary2.2} we state another 
lemma
the proof of which is also a straightforward computation.

%%%%%%Sublemma B.2
\begin{sublemma}
\label{nsublemmaB.2}
Let 
$$
M_i = \left(  \begin{array}{cc}
               A & C \\
               B & D
             \end{array} \right),  \ \ \ \ \ \ \ \
M^{(j)} = \left( \begin{array}{cc}
                A_j & C_j \\
                B_j & D_j
             \end{array} \right), \ \ \ \ j=i-1, i.
$$
Then
\begin{eqnarray}
\| e_i \times e_{i-1}\| & = &
\frac{1}{\rho^{(i)} \rho^{(i-1)}}
\mid  \det(M^{(i-1)})[(AC+BD)(C_{i-1}^2+D_{i-1}^2) \nonumber \\
& & + (A^2 + B^2 - C^2 - D^2)C_{i-1}D_{i-1}] + \Delta_i \mid
\label{nB.1(3)}
\end{eqnarray}
where $\rho^{(i-1)}$ and $\rho^{(i)}$ are the normalizing constants
for $e_{i-1}$ and $e_i$ as in Sublemma \ref{nsublemmaB.1}, and
$$
\Delta_i = - (\lambda_i^{min})^2(A_{i-1}C_{i-1}+B_{i-1}D_{i-1})
+ (\lambda_{i-1}^{min})^2
(A_iC_i + B_i D_i).
$$
\end{sublemma}

Observe that each the terms in the numerator of (\ref{nB.1(3)}) has a factor
$|\det(M^{(i-1)})|$, $\lambda^{min}_{i-1}$ or $\lambda^{min}_i$, all of which 
are $\leq (\frac{b}{\kappa})^{i-1}$. Observe also that if both $e_{i-1}$ and $e_i$ are nearly parallel to the $x$-axis, then $\rho^{(i)}, \rho^{(i-1)}$ are $> K^{-1}\kappa^{2i}$ (see the proof of Lemma \ref{nlemma2.1}).

\medskip

\noindent {\bf Proof of Corollary \ref{ncorollary2.2}:} \ We begin with some useful derivative estimates. First, we claim that
\begin{equation}
\|\partial^1 M^{(i)}\| < K^i.
\label{nB.1(4)}
\end{equation}
This is because $\partial^1 M^{(i)}$ is the sum of $i$ terms of the form
$M_i  \cdots M_{j+1} (\partial^1 M_j) M_{j-1} \cdots M_1$
and the norm of this product is $<K_0^{2i}$. A similar argument gives
\begin{equation}
|\partial^1 \det{M^{(i)}}| \leq (Kb)^i.
\label{nB.1(5)}
\end{equation}
Since $\lambda^{max}_i =\|M^{(i)}\|$, it follows from (\ref{nB.1(4)}) 
that $|\partial^1 \lambda^{max}_i|<K^i$; and since
$\lambda^{min}_i=|\det M^{(i)}|/\lambda^{max}_i$, we have
$|\partial^1 \lambda^{min}_i|<(\frac{Kb}{\kappa^2})^i$. 

\smallskip
Pre-composing with a suitable orthogonal matrix as in the proof of 
Lemma \ref{nlemma2.1}, 
we may assume that $\rho^{(i)}, \rho^{(i-1)}$ are $> K^{-1}\kappa^{2i}$. 
The estimate for $\partial^j \theta_1$ is obtained by differentiating
(\ref{nB.1(1)}). To prove (\ref{nformula2.1}), we differentiate 
(\ref{nB.1(3)}), and observe using the 
inequalities above that after differentiation, the numerator is the sum of a finite number
of terms each one of which is bounded above by $(\frac{Kb}{\kappa^2})^{i-1}$. 

\smallskip

To prove (\ref{nformula2.2}), we write
\begin{eqnarray*}
M^{(i)} e_n &=& M^{(i)} e_i +
M^{(i)} (e_n - e_i) \\
 &=&  M^{(i)} e_i +
\sum_{k = i}^{n-1}  M^{(i)}(e_{k+1} - e_k)
\end{eqnarray*}
and take partial derivative one term at a time.
First we have
$$
\partial^1 M^{(i)}(e_{k+1} - e_k)
= \partial^1 M^{(i)} \cdot (e_{k+1} - e_k) +
M^{(i)} \cdot \partial^1 (e_{k+1} - e_k).
$$
The norm of the first term on the right side is bounded by 
$(\frac{K b}{\kappa^2})^k$
because $\|\partial^1 M^{(i)}\| \leq K^i$ and $\| e_{k+1} - e_k \|
< (\frac{K b}{\kappa^2})^k$. The norm of the second term is bounded by 
$(\frac{K b}{\kappa^2})^k$ according to (\ref{nformula2.1}).
It remains to show 
$\|\partial^1 M^{(i)}e_i\| < (\frac{K b}{\kappa^{2}})^i$.
This follows by differentiating (\ref{nB.1(2)}) 
and using the inequalities above. The proof for $j =2$ is similar. 
\hfill $\square$

\bigskip
%%%%Sublemma B.3
\begin{sublemma}
Let $M_i$ and $M_i^{\prime}$ be as in Lemma \ref{nlemma2.2}, 
let $m < \frac{n}{2}$, and write
$$
M_{i, m} = M_{i+m} M_{i-1+m} \cdots M_m \ , \ \ \ 
M_{i, m}^{\prime} =  M_{i+m}^{\prime} M_{i-1+m}^{\prime} \cdots M_m^{\prime}.
$$
Then 
\begin{equation}
\| M_{i, m} - M_{i, m}^{\prime} \| < \frac{1}{4} (K \lambda)^m
\label{nB.1(6)}
\end{equation}
for all $i$, $0 \leq i \leq m$. \label{nsublemmaB.3}
\end{sublemma}
\noindent {\bf Proof:} 
Set $\rho_{k} = \| M_{k, m} - M_{k, m}^{\prime} \|$. Then
\begin{eqnarray*}
M_{k+1, m} - M_{k+1, m}^{\prime} &=& M_{k +1 + m} M_{k, m} - 
M_{k +1+ m}^{\prime} M_{k, m}^{\prime} \\
&=&  M_{k+1+m} (M_{k, m} - M_{k, m}^{\prime})
+ (M_{k+1+m} - M_{k+1+m}^{\prime}) M_{k, m}^{\prime}.
\end{eqnarray*}
Since $\| M_{k, m}^{\prime} \| < K_0^{k}$ and
$\|M_{k+1+m} - M_{k+1+m}^{\prime}\| < \lambda^{k + m}$, we have
$$
\rho_{k+1} \leq K \rho_{k} + K^{k} \lambda^{m+k},
$$
which implies (\ref{nB.1(6)}).
\hfill $\square$

\medskip

\noindent {\bf Proof of Lemma \ref{nlemma2.2}:} \ (\cite{BC2}, p. 108): 
\ We prove the assertion for  
all the indices 
that are  
powers of two and leave the rest as an exercise.
To prove (b), write $m_j = 2^j$, and let
$$
u_j = \frac{w_{m_j}}{\| w_{m_j} \|},
\ \ \ u_j^{\prime} =
\frac{w_{m_j}^{\prime}}{\| w_{m_j}^{\prime} \|}
$$
where $w_{m_j} = M^{(m_i)} w$ and $w_{m_j}^{\prime} = M^{\prime (m_i)} w$.
We will show inductively that
\begin{equation}
\| u_j \times u_j^{\prime} \| < \lambda^{\frac{m_j}{4}}.
\label{nB.1(7)}
\end{equation}
Assume that (\ref{nB.1(7)}) is true up to
index $j$. Let
$$
A = M_{m_{j+1} - m_j, m_j} \ \ \ and \ \ \ 
A^{\prime} = M^{\prime}_{m_{j+1} - m_j, m_j}.
$$
Since
$\| w_{m_j} \| < K^{m_j}$ and
$\| w_{m_{j+1}} \| > \kappa^{m_{j+1}}$,
we have
\begin{eqnarray}
\| Au_j \| &=&
\frac{\| w_{m_{j+1}} \|}{\| w_{m_j} \|}
> \left( \frac{\kappa^2}{K} \right)^{m_j}, 
\label{nB.1(8)} \\
\| Au^{\prime}_j \| & \geq &
\| Au_j \| - \| A \| \| u_j - u^{\prime}_j \|  
 \geq \left(\frac{\kappa^2}{K}\right)^{m_j} - K^{m_j} \lambda^{\frac{m_j}{4}}
\geq \frac{3}{4}\left(\frac{\kappa^2}{K}\right)^{m_j}.
\label{nB.1(9)}
\end{eqnarray}
Writing
$\| A^{\prime}u^{\prime}_j \|
= \| A^{\prime} u^{\prime}_j - A^{\prime} \hat{u}_j + A^{\prime} \hat{u}_j
- A \hat{u}_j + A \hat{u_j} \|$
where $\hat{u}_j = u_j$ if the angle between $u_j$ and $u_j^{\prime}$ is
smaller than $\frac{\pi}{2}$, $\hat{u}_j = -u_j$ otherwise, we obtain 
$\| A^{\prime}u^{\prime}_j \|
\geq \| Au_j \| - \| A \|
\| u_j \times u_j^{\prime} \| - \| A - A^{\prime} \|$.
Using Sublemma \ref{nsublemmaB.3} to bound $\| A - A^{\prime} \|$,
we again have
\begin{equation}
\parallel A^{\prime}u^{\prime}_j \parallel
\geq \frac{3}{4}\left(\frac{\kappa^2}{K}\right)^{m_j}.
\label{nB.1(10)}
\end{equation}
We are now ready to prove (\ref{nB.1(7)}) for index $j+1$:
\begin{eqnarray*}
\| u_{j+1} \times u_{j+1}^{\prime} \|
&=& \frac{\| Au_j \times 
A^{\prime}u_j^{\prime} \|}{\| A u_j \| \cdot \| A^{\prime}_j u^{\prime}_j \|} 
 = \frac{\| Au_j \times (A - A + A^{\prime})u_j^{\prime} \|}{\| A u_j \| \cdot \| A^{\prime
}_j u^{\prime}_j \|} \\ 
& \leq &
\frac{\| Au_j \times Au_j^{\prime} \|}{\| A u_j \| \cdot \|
A^{\prime
}_j u^{\prime}_j \|} + 
\frac{\| Au_j \times (A - A^{\prime}) u_j^{\prime} \|}{\| A u_j \| \cdot \|
A^{\prime
}_j u^{\prime}_j \|} .
\end{eqnarray*} 
The first term is fine since 
$\| A u_j \times A u_j^{\prime} \| = |\det(A)| \ \| u_j \times u_j^{\prime} \|$ 
and $|\det(A)| < b^{m_j}$. To estimate the second term, we use
Sublemma \ref{nsublemmaB.3} and
(\ref{nB.1(8)})-(\ref{nB.1(10)}).

To prove (a), we again let 
$i = 2^k$. Then for 
$0 < j \leq k$, we have
$$
\| w_{m_{j+1}}^{\prime} \| =  
\| w_{m_j}^{\prime} \|
\| A^{\prime} u^{\prime}_j \| \\
= \| w_{m_j}^{\prime} \|
\| A^{\prime} u^{\prime}_j -  A u^{\prime}_j
+  A u^{\prime}_j - A \hat{u}_j + A \hat{u}_j \|,
$$
so that
\begin{eqnarray*}
\frac{\| w_{m_{j+1}}^{\prime} \|}{\|
w_{m_j}^{\prime} \|} & \geq & \| Au_j \|
-\| A^{\prime} -A \| \| u^{\prime}_j \|
-\| A \| \| u^{\prime}_j - \hat{u}_j \|
\\
& = & \frac{\| w_{m_{j+1}} \|}{\|
w_{m_j} \|}
\left(1- \frac{\| w_{m_j} \|}{\|
w_{m_{j+1}} \|}(
\| A^{\prime} -A \| +
\| A \| \| u^{\prime}_j - \hat{u}_j \| ) \right).
\end{eqnarray*}
Using Sublemma \ref{nsublemmaB.3} to bound 
$\| A - A^{\prime} \|$ and part (b) of this lemma to bound
$\| \hat{u}_j - u^{\prime}_j \|$, we obtain 
$$
\frac{\| w_{m_{j+1}}^{\prime} \|}{\|
w_{m_j}^{\prime} \|} \geq
\frac{\| w_{m_{j+1}} \|}{\|
w_{m_j} \|}(1- 4^{-m_j}),
$$
which implies (a). \hfill $\square$

%%%%%%%%%%%%%%%%%%%%%%%%%%%Appendix B.2
\subsection{Stable curves (Sect. \ref{s2.2})}
\label{app-B.2}
On a ball of radius $\frac{\lambda}{2K_0}$ centered at $z_0$,
we have $\|DT\| \geq \frac{\kappa}{2}$
so that $e_1$, the field of most contracted directions
of $DT$, is well defined. Let $\gamma_1$ be the integral curve to $e_1$ of 
length
 $\sim \lambda$
passing through $z_0$.

To construct $\gamma_2$, let $B_1$ be the $\frac{\lambda^2}{2K_0}$-neighborhood
of $\gamma_1$.
For $\xi \in B_1$, let $\xi'$ be a point in $\gamma_1$
with $|\xi-\xi'|<\frac{\lambda^2}{2K_0}$.
Then $|T\xi-Tz_0| \leq |T\xi-T\xi'| +|T\xi'-Tz_0|
\leq \frac{\lambda^2}{2} + \frac{Kb}{\kappa^2} \lambda <\lambda^2$,
so by Lemma \ref{nlemma2.2},
$\|DT^2 \xi\| \geq \frac{\kappa^2}{2}$. This ensures that $e_2$,
the field of most
contracted directions for $DT^2$, is defined on all of
$B_1$. Let $\gamma_2$ be
the integral curve through $z_0$ in $B_1$. We leave it as an exercise
to show that the Hausdorff distance between $\gamma_1$ and $\gamma_2$
is ${\cal O}(\frac{b}{\kappa^2} \lambda)<< \lambda^2$, so that $\gamma_2$ has 
essentially the same length as $\gamma_1$.
This uses the fact that $e_1$ has Lipschitz constant $K$ 
(Corollary \ref{ncorollary2.2})
and that the angle between $e_1$ and $e_2$ is $<\frac{Kb}{\kappa^2}$
(Corollary \ref{ncorollary2.1}).
 
Next we let $B_2$ be the $\frac{\lambda^3}{2K_0}$-neighborhood
of $\gamma_2$ and repeat the argument above to get $e_3$ and $\gamma_3$.
Using the Lipschitzness of $e_2$
and the fact that $\parallel e_3 \times e_2 \parallel
\leq (\frac{Kb}{\kappa^2})^2$,
we conclude again that $\gamma_3$ has essentially the same length
as $\gamma_2$.
This process is continued for $n$ steps. 

%%%%%%%%%%%%%%%%%%%%%%%%%%%Appendix B.3
\subsection{Curvature estimates (Sect. \ref{s2.3})}
\label{app-B.3}
Recall that 
$$
k_i(s) = \frac{\| \gamma^{\prime}_i(s) 
\times \gamma^{\prime \prime}_i(s) \|}{\| \gamma^{\prime}_i(s) \|^3}.
$$
Write
$$
DT = DT(\gamma_i(s)) = \left( \begin{array}{cc}
            A & C \\
            B & D 
            \end{array} \right) 
$$
and 
$$
X = \left( \begin{array}{cc}
            <\nabla A, \gamma_{i-1}^{\prime}>  & <\nabla C, \gamma_{i-1}^{\prime}> \\
            <\nabla B, \gamma_{i-1}^{\prime}>  & <\nabla D, \gamma_{i-1}^{\prime
}> 
            \end{array} \right)
$$
where $< \ , \ >$ is the usual inner product. Since 
$\gamma_i^{\prime} = DT \cdot \gamma_{i-1}^{\prime}$ and
$\gamma_i^{\prime \prime} = 
DT \cdot \gamma^{\prime \prime}_{i-1} 
+  X \cdot \gamma_{i-1}^{\prime}$,
we have
\begin{equation}
k_i  =  \frac{1}{\| \gamma_i^{\prime} \|^3} \| DT \cdot \gamma_{i-1}^{\prime}
\times (DT \cdot \gamma^{\prime \prime}_{i-1})
+ X \cdot \gamma_{i-1}^{\prime})\| 
 \leq  \frac{1}{\| \gamma_i^{\prime} \|^3} (I + II)
\label{nB.3(1)}
\end{equation}
where
$$
I = 
|\det(DT)| \cdot \|\gamma_{i-1}^{\prime} \times \gamma^{\prime \prime}_{i-1}\|, 
\ \ \ \ II = \| DT \cdot \gamma_{i-1}^{\prime}
\times
 X \cdot \gamma_{i-1}^{\prime}\|.
$$
Term II is degree three homogeneous in $\gamma_{i-1}^{\prime}$.
Moreover, the second component of each vector involved
in the cross product has a factor $b$.  
Thus there exist $K>0$ such 
that 
\begin{equation}
k_i \leq (b \cdot k_{i-1} + K \cdot b)\cdot 
\frac{\|\gamma_{i-1}^{\prime}\|^3}{\|\gamma_{i}^{\prime}\|^3}.
\label{nB.3(2)}
\end{equation}

Lemma \ref{nlemma2.4} follows by recursively applying
inequality (\ref{nB.3(2)}). 

%%%%%%%%%%%%%%%%%%%%%%%%%%%Appendix B.4
\subsection{One-dimensional dynamics (Sect. \ref{s2.4})}
\label{app-B.4}

Let $\delta_0:=$inf$\{d(f^n \hat x,C):\hat x \in C, n>0\}$. We begin with three easy observations:

(i) There exists $k_0>0$ such that for all $\delta<\frac{1}{2}\delta_0$,
if $x$ is such that $f^nx \in C_\delta$, then $|(f^n)'x|\geq k_0$.
This is true because there is an interval $(x_1,x_2)$ containing $x$
on which $f^n$ is monotone and $f^n(x_1,x_2) \supset (\hat x-2\delta,
\hat x+2\delta)$ for some $\hat x \in C$. It then follows from the negative
Schwarzian property that restricted to $f^{-n}(\hat x-\delta, \hat x+\delta) \cap (x_1,x_2)$, $|(f^n)'| \geq $ some $k_0>0$ independent of $x$.

(ii) There exists $\lambda_0>1$ such that for all sufficiently small
$\delta$,  if $d(x,C)<\delta$, then
there exists $p=p(x)$ such that $f^ix \not \in C_\delta$
for all $i<p$ and $|(f^p)'x| \geq \lambda_0^p$. This is an easy computation
using the fact that the forward critical orbits of $f$ are contained in a
uniformly expanding invariant set. Let $\hat p(\delta)=$inf$\{p(x):d(x,C)
<\delta\}$.

(iii) For all sufficiently small $\delta$, there exist 
$N_1(\delta) \in {\mathbb Z}$
and $\lambda_1(\delta)>1$ such that if $x, \cdots, f^nx \not \in C_\delta$
for some $n>N_1$, then $|(f^n)'x| \geq \lambda_1^n$. This is proved in 
\cite{M1}.

\smallskip
We now prove the assertion in Lemma \ref{nlemma2.5}. Fix $\delta_1$ sufficiently small for (i)--(iii) above, and with
the property that $\lambda_0^{\hat p(\delta_1)}>>k_0^{-1}$.  Consider $\delta<\delta_1$ and an orbit segment $x, \cdots, f^nx$ with
$f^ix \not \in C_\delta$ for $i<n$ and $f^n x \in C_\delta$. To estimate $(f^n)'x$, we let
$n_j$ be the $j$th time $f^ix \in C_{\delta_1}$, and let $p_j=p(f^{n_j}x)$.
Then $|(f^{p_j})'(f^{n_j}x)| \geq \lambda_0^{p_j}$, and between the times
$n_j+p_j$ and $n_{j+1}$, the derivative is bounded below by 
$\lambda_1(\delta_1)^{n_{j+1}-(n_j+p_j)}$ if $n_{j+1}-(n_j+p_j)>N_1(\delta_1)$, 
by $k_0$ otherwise. The same estimate holds for the initial stretch up to
time $n_1$. Noting that the factor $k_0$ can be absorbed into $\lambda_0^{p_j}$, we see that $|(f^n)'x| \geq e^{\hat c_1 n}$
where $e^{\hat c_1}$ can be taken to be slightly smaller than
$(\min(\lambda_0, \lambda_1(\delta_1))^{\frac{\hat p(\delta_1)}
{\hat p(\delta_1)+N_1(\delta_1)}}$. Also, $\hat c_0$ can be taken to be $k_0 \lambda^{-N_1
(\delta_1)}$. This completes the proof of part (ii) of Lemma \ref{nlemma2.5}.

To prove (i), let $n_q<n$ be the last time $f^ix \in C_{\delta_1}$,
and observe that $|(f^{n-n_q})'(f^{n_q}x)| \geq K^{-1}k_0 \delta$
if $n-n_q<N_1(\delta_1)$, $ \geq K^{-1}\delta \lambda_1(\delta_1)^{n-n_q}$
otherwise. 
\hfill $\square$
%%%%%%%%%%%%%%%%%%%%%%%%%%%%%%%%%%Appendix B.5
\subsection{Critical points inside ${\cal C}^{(0)}$ (Sect. \ref{s2.6})}
\label{app-B.5}
{\bf Proof of Lemma \ref{nlemma2.9}:} \ Write
\begin{equation}
\frac{dq_1(s)}{ds} = \partial_x q_1(x, y) \frac{dx(s)}{ds}
+ \partial_y q_1(x, y) \frac{dy(s)}{ds}.
\label{nB.5(1)}
\end{equation}
Since $\gamma$ is $b$-horizontal, we have $\frac{dx(s)}{ds} \approx 1$ and
$\mid \frac{dy(s)}{ds} \mid <  {\cal O}(b) \cdot \mid \frac{dx(s)}{ds} \mid$. 
By (\ref{nB.1(1)})
\begin{equation}
q_1(s) = \frac{AC + BD}{C^2 + D^2 - (\lambda^{min})^2}, 
\label{nB.5(2)}
\end{equation}
so
\begin{eqnarray*}
\partial_x q_1(x, y) & = & \frac{A_x C + A C_x + {\cal O}(b)}{C^2 + D^2 - 
(\lambda^{min})^2}
- 2 \frac{(AC + BD)(C C_x + D D_x+ \lambda^{min} \lambda^{min}_x)}{(C^2 + D^2 -
(\lambda^{min})^2)^2} \\
& := & I + II
\end{eqnarray*}
where
$$
A = F_x + b u_x, \ \ \ C = F_y + bu_y,
$$
$$
B = b v_x, \ \ \ D = b v_y.
$$
We will show that $|I| \geq K^{-1}$ and $|II| = {\cal O}(\delta)$. To estimate
$I$, observe that the denominator is $> K^{-1}$, and that for $(x, y) \in 
{\cal C}^{(0)}$, $|AC_x| = {\cal O}(\delta)$,
while $|A_xC| = |F_{xx}F_y|(1+{\cal O}(b)) \geq K^{-1}$ since $|F_y| > K^{-1}$ 
(non-degeneracy condition). 
Term II follows from the fact that its denominator is 
$\geq K^{-1}$, and $AC + BD = {\cal O}(\delta)$. \hfill $\square$

\bigskip

\noindent {\bf Proof of Lemma \ref{nlemma2.10}:} \
Using the results in Sect. \ref{s2.1} and Lemma \ref{nlemma2.9}, we have that
at $\gamma(s)$ with 
$|s|<(Kb)^{\frac{m}{2}}$,
$e_{3m}$ is defined  with
$|q_{3m}-q_m|<(Kb)^m$ (Lemma \ref{nlemma2.2})
and $|\frac{d}{ds}q_m| \geq K^{-1}$
(Corollary \ref{ncorollary2.2}
and Lemma \ref{nlemma2.9}).
Let $\tau(s)$ denote the slope of $\gamma'(s)$, and
assume for definiteness that $\frac{d}{ds}q_{3m}>0$.
Then
\begin{eqnarray*}
q_{3m}((Kb)^{\frac{m}{2}}) - \tau((Kb)^{\frac{m}{2}}) & = &
(q_{3m}((Kb)^{\frac{m}{2}}) - q_m((Kb)^{\frac{m}{2}}))
 + (q_m((Kb)^{\frac{m}{2}}) - q_m(0)) \\
& & + (q_m(0)- \tau(0))+(\tau(0)-\tau((Kb)^{\frac{m}{2}})) \\
&\geq &  -(Kb)^m + K^{-1}(Kb)^{\frac{m}{2}} + 0 - K_1b(Kb)^{\frac{m}{2}} \\
&\geq & \frac{K^{-1}}{2}(Kb)^{\frac{m}{2}}.
\end{eqnarray*}
 
Similarly, $q_{3m}(-(Kb)^{\frac{m}{2}}) - \tau(-(Kb)^{\frac{m}{2}})<0$,
giving a unique critical point of
order $3m$ in between. 

\hfill $\square$

\bigskip

\noindent {\bf Proof of Lemma \ref{nlemma2.11}:} \ 
Let $\tau(s)$ be the slope of the tangent vector to $\gamma$ at 
$\gamma(s)$, and let $q_m(s)$ be the slope of $q_m$ at $\gamma(s)$. 
Let $\hat{\tau}(s)$ and $\hat{q}_m(s)$ denote the corresponding quantities
at $\hat \gamma(s)$. First we claim that
\begin{equation}
\mid \tau(0) - \hat{\tau}(0) \mid \leq 2 \sqrt{\varepsilon}.
\label{nB.5(3)}
\end{equation}
An easy calculation (which we omit) shows that if this was not the case,
then $\gamma$ and $\hat \gamma$ would meet at 
$\gamma(s)$ for some $|s|<\sqrt{\varepsilon}$. 

Let $\hat{m}$ be the largest integer $j \leq m$ such that 
$4K_1\sqrt{\varepsilon} < \|DT\|^{-13j}$.
Then by Lemma \ref{nlemma2.2},  $\|DT^i(\gamma(s))\| > \frac{1}{2}$ for 
$0 < i < \hat{m}$ and $s \in [- 4K_1\sqrt{\varepsilon}, 
4K_1\sqrt{\varepsilon}]$. This guarantees that $q_{\hat m}$ is defined 
everywhere on 
$\gamma$ and on $\hat \gamma$. 
Let $\hat \sigma(s):=\hat{q}_{\hat{m}}(s) - \hat \tau(s)$.
We have
\begin{eqnarray*}
|\hat \sigma(0)| & \leq & |\hat q_{\hat{m}}(0) -q_{\hat{m}}(0)|
+ |q_{\hat{m}}(0) - q_{m}(0)| + |q_m(0) - \tau(0)|
+ |\tau(0) - \hat \tau(0)| 
\\
& < & K\varepsilon+ (Kb)^{\hat m}+0+2\sqrt{\varepsilon} \ 
< \ 3\sqrt{\varepsilon}.
\end{eqnarray*}

To prove the existence of a critical point of order $\hat m$ on $\hat
\gamma$, we will compare the signs of $\hat \sigma$ at the two
end points of $\hat \gamma$. First, 
$$
\hat \sigma(4K_1 \sqrt{\varepsilon})=\hat{q}_{\hat{m}}(0)
+ \frac{d}{ds}q_{\hat m}(s_1) \cdot 4K_1 \sqrt{\varepsilon}
-\hat \tau(0) - \frac{d}{ds}\hat \tau(s_2) \cdot 4K_1 \sqrt{\varepsilon}
$$
for some $s_1, s_2 \in [0,4K_1 \sqrt{\varepsilon}]$. This is
$$=\hat \sigma(0)+ (\frac{d}{ds}q_{\hat m}(s_1)+{\cal O}(b)) 
\cdot 4K_1 \sqrt{\varepsilon}.
$$
Since the second term has absolute value 
$>(K_1^{-1}-{\cal O}(b))\cdot 4K_1\sqrt{\varepsilon}>|\hat \sigma(0)|$, 
it follows that $\hat \sigma(4K_1 \sqrt{\varepsilon})$ 
has the same sign as $\frac{d}{ds}q_1$.
An analogous computation shows that $\hat \sigma(-4K_1 \sqrt{\varepsilon})$ has the opposite sign as $\frac{d}{ds}q_1$.\hfill $\square$

%%%%%%%%%%%%%%%%%%%%%%%%%%%Appendix B.6
\subsection{Growth of $w_i$ and $w^*_i$ \ (Sect. \ref{s4.2})}
\label{app-B.6}
%%%%%%%%%%%%%%%%%%%%%%%%%%Sublemma B.4
\begin{sublemma}
 Let $z_0$ be h-related to $\hat z_0 \in
\Gamma_{\theta N}$ with bound period $p<\frac{2}{3}N$, and let 
$w_0={\tiny \left(\!\!
\begin{array}{c} 0 \\ 1 \end{array}\!\!\right)}$.
Then for $i \leq p$, $\|w^*_i\| >K^{-1}e^{c''i}$ for some $c'' \approx c$.
\label{sublemmaB.4}
\end{sublemma}
\noindent {\bf Proof:} Let $\hat w^*_i$ be as defined in (IA6).
Then (IA4) and (IA6) together imply that
$\|\hat w^*_i\| >\frac{c_0}{2}e^{ci}$. The only difference between $\hat w^*_i$
and $w^*_i$ is that contractive fields of order $\ell(\hat z_i)$
are used for splitting for the former and $\ell (z_i)$ the latter at returns to
${\cal C}^{(0)}$. By Lemma \ref{nlemma4.2}, $\ell (z_i)=
\ell(\hat z_i) \pm 1$, so that recombination times may differ by one.
This is clearly of no consequence. Assuming these times are
synchronized, we observe next that $w^*_i$ has the same direction as 
$\hat w^*_i$.
This can be seen inductively (using the nested property of fold periods).
Finally, a vector split using a field of order
$\ell$ or $\ell+1$ may differ in length by a factor of
$1 \pm {\cal O}(b^\ell)$.  Thus $\|w^*_i\| \geq (1 -{\cal O}(b))^i
\| \hat w^*_i\|$.
\hfill $\square$

\bigskip

\noindent {\bf Proof of Lemma \ref{nlemma4.6}:} \ We may assume $z_i$ is in a 
fold period,
otherwise there is nothing to prove. Let $i_1 < i \leq i_2$ be
the longest fold period containing $i$. By Lemma \ref{nlemma4.4},
which applies also to controlled orbits satisfying
$d_{\cal C}(z_j)>e^{-\alpha j}$, we have $i_2-i_1 \leq \varepsilon i$.
Let
$w_{i_1}=Ae+B {\tiny \left(\!\! \begin{array}{c} 0 \\ 1 \end{array}\!\!\right)}$
be the usual splitting. Then
$$
\|w_i^*\| \leq K^{i-i_1}|B| \leq K^{i-i_1}\|w_{i_2}^*\|
=K^{i-i_1}\|w_{i_2}\| \leq K^{i-i_1}(K^{i_2-i}\|w_i\|)
\leq K^{\varepsilon i} \|w_i\|.
$$
The first ``$\leq$" uses the fact that $\frac{\|w^*_{j+1}\|} {\|w^*_j\|} \leq$ 
some $K$, the second uses Sublemma \ref{sublemmaB.4}, 
and the third $\|DT\| \leq K$.
The reverse estimate follows from $\|w_i\| \leq K^{i-i_1}\|w_{i_1}\| \leq
K^{i-i_1} d_{\cal C}(z_{i_1})^{-1} \|w^*_i\|$ and
$d_{\cal C}(z_{i_1})>e^{-\alpha i}$.
\hfill $\square$
 
\bigskip
 
\noindent {\bf Proof of Lemma \ref{nlemma4.7}:} \ We give a proof in the case 
where $j$ exists;
the other case is simpler.
Let $k \leq i_1< i_1+p_1 \leq i_2 < i_2+p_2 \leq \cdots
\leq i_r=j < n$ be defined as follows:
we let $i_1$ be the first return to
${\cal C}^{(0)}$ at or after time $k$, $p_1$ the 
bound period of $z_{i_1}$, $i_2$
the first
return after $i_1+p_1$, and so on until $i_r=j$.
Writing $k=i_0+p_0$,  we have that $\frac{\|w^*_n\|}{\|w^*_k\|}$ is a 
product of
factors of the following three types:
$$
I:=\frac{\|w^*_{i_{s+1}}\|}{\|w^*_{i_s+p_s}\|}, \ \ \ \ \
II:=\frac{\|w^*_{i_s+p_s}\|}{\|w^*_{i_s}\|} \ \ \ \ {\rm and} \ \ \ \
III:=\frac{\|w^*_n\|}{\|w^*_j\|}.
$$
 
First we prove the lemma assuming that no fold 
periods initiated before time $k$
expires between times $k$ and $n$. By Lemma \ref{nlemma2.8}, $I \geq
c_0e^{c_1(i_{s+1}-(i_s+p_s))}$. Since $w^*_{i_s}$ splits correctly,
we have, by (IA5), $II
\geq K^{-1}e^{\frac{c}{3}p_s}$. Moreover, we may assume that
$c_0$ and $K$ above
can be absorbed into the exponential estimate for the bound period
$[i_s, i_s+p_s]$.
For $III$, let $\ell$ be the fold period initiated at time $j$.
If $\ell>n-j$, then $III \geq K^{-1}d_{\cal C}(z_j)
e^{c''(n-j)}$ by Sublemma \ref{sublemmaB.4}. If not, we split $w^*_j$ into
$w^*_j=Ae_{n-j}+B 
{\tiny \left(\!\!\begin{array}{c} 0 \\ 1 \end{array}\!\!\right)}$,
noting that $e_{n-j}$ is defined at $z_j$ by Sublemma \ref{sublemmaB.4}
and Lemma \ref{nlemma4.6}.
Then $III \geq K^{-1}d_{\cal C}(z_j)e^{c''(n-j)} - (Kb)^{n-j}$.
The last term is negligible because $d_{\cal C}(z_j)
\sim b^{\frac{\ell}{2}}>>(Kb)^{n-j}$.
Altogether, this gives $\frac{\|w^*_n\|}{\|w^*_k\|}
\geq K^{-1}d_{\cal C}(z_j)e^{c'(n-k)}$ for some $c'>0$ as claimed.
 
In the rest of the proof, we view contributions from fold periods initiated 
before time $k$ as perturbations of the estimates above, and verify that
they are inconsequential.
For $I$, we claim that for each $t$ in question,
$$
\frac{\|w^*_{t}\|}{\|w^*_{t-1}\|} = (1 \pm {\cal O}(\sqrt b))
\frac{\|DT(z_{t-1})w^*_{t-1}\|}{\|w^*_{t-1}\|},
$$
so that $I$ has the same estimate as before with possibly a slightly
smaller $c_1$. This claim follows from the fact that when a fold period
initiated $\ell$ steps earlier expires at time $t$, the vector to rejoin
the main term has magnitude 
$\|DT(z_{t-1})w^*_{t-1}\| {\cal O}(b^{\frac{\ell}{2}})$.
(See Sect. \ref{s2.7})
 
Next we turn to $III$, which is similar to and a little more complicated than 
$II$. Given $z_t$ and a vector $u$, we let $u, \ T^1_*(z_t)u, \ T^2_*(z_t)u,\
\cdots$ denote the vectors given by the splitting algorithm for the
orbit segment beginning at $z_t$ with initial vector $u$ -- neglecting
recombinations from fold periods initiated before time $t$.
Then
$$
w^*_n=T^{n-j}_*(z_j)w^*_j + \sum_{t=j+1}^n T^{n-t}_*(z_t)E_t
$$
where $E_t$ is the sum of the vectors to be rejoined at time $t$.
For fixed $t$, let $\ell$ be the shortest  fold period initiated before $k$ to 
expire at time $t$. From Sect. \ref{s2.7}, we have $\|E_t\| \leq
(Kb)^{\frac{\ell}{2}}\|w^*_t\|$.
Also, since this fold period contains the one initiated
at $j$, we have, by Sect. \ref{s4.1},  $K \alpha \ell> (n-j)$. Together this gives
$$
\|T^{n-t}_*(z_t)E_t\| \leq K^{n-t}(Kb)^{\frac{\ell}{2}}\|w^*_t\|
\leq (Kb^{\frac{1}{K\alpha}})^{n-j}\|w^*_t\|.
$$
Assuming inductively that the assertion in the lemma has been proved 
for shorter
time intervals, we have $\|w^*_t\| \leq
Kd_{\cal C}(z_{j_t})^{-1}\|w^*_n\|$ where $j_t$ is a return between
times $t$ and $n$. Thus
$$
\sum_{t=j+1}^n \|T^{n-t}_*(z_t)E_t\|
< (n-j) (Kb^{\frac{1}{K\alpha}})^{n-j} e^{\alpha (n-j)}  \|w^*_n\|
<< \|w^*_n\|,
$$
which together with our earlier estimate on $\|T^{n-j}_*(z_j)w^*_j\|$
gives the disired result.
\hfill $\square$

\bigskip
\noindent {\bf Proof of Lemma \ref{nlemma4.8}:} \  The case where $z_k$ is not 
in a fold
period is contained in Lemma \ref{nlemma4.7}. Let $j$ be the starting point of
the largest
fold period covering $z_k$. Observe that its length $\ell$ is
$< K \theta (n-j)< 2K \theta (n-k)$
because the bound period initiated at time $j$
has expired by time $n$.
Then by Lemma \ref{nlemma4.7},
$$
\|w_n\|>K^{-1}e^{c'(n-j)}\|w_j\| \geq K^{-1}e^{c'(n-j)}K^{-\ell}\|w_k\|.
$$
\hfill $\square$

\bigskip

\noindent {\bf Proof of Lemma \ref{nlemma4.5}:} \
The proof proceeds inductively. Consider a bound return $z_i$,
and assume that the $w_j^*$-vectors split correctly at {\it all} returns
prior to time $i$. Let $\angle(\cdot, \cdot)$ denote the
angle between two vectors, and let
$u = {\tiny \left(\!\!\begin{array}{c} 0 \\ 1 \end{array}\!\!\right)}$.
 
{\bf Case 1.} $z_i$ is in a fold period. Let $j<i$ be the largest
integer such that the fold period initiated at time $j$ remains in effect
at time $i$, and let $\hat z_0=\phi(z_j)$. Then proving $w_i^*$
splits correctly is equivalent to proving
$$
\angle(DT^{i-j}(z_j)u, \ \tau(\phi(z_i)) \ < \varepsilon_0 d_{\cal C}(z_i).
$$
We compare this inequality to
$$
\angle(DT^{i-j}(\hat z_0)u, \ \tau(\phi(\hat z_{i-j})) \ < \varepsilon_0 d_{\cal
 C}(\hat z_{i-j}),
$$
which we know to be true by (IA3). Suppose $\hat z_{i-j} \in {\cal C}^{(k)}$. 
Then
 
- $\angle(DT^{i-j}(z_j)u, \ DT^{i-j}(\hat z_0)u)
<<e^{-\beta(i-j)}<<d_{\cal C}(z_i)$ by (IA6);
 
- $\angle(\tau(\phi(z_i), \ \tau(\phi(\hat z_{i-j}))<b^{\frac{k-1}{4}}
<<d_{\cal C}(z_i)$ by Lemma \ref{nlemma4.1};
 
- $|d_{\cal C}(z_i)- d_{\cal C}(\hat z_{i-j})|
<e^{-\beta(i-j)}+b^{\frac{k-1}{4}}<<d_{\cal C}(z_i)$.

\medskip
 
{\bf Case 2.} $z_i$ is not in any fold period. In this case let $j<i$
be the last free return, so that the bound period initiated
at $j$ remains in effect at $i$ and $w_i^*=DT^{i-j}(z_j)w_j^*$.
We split
$$
w^*_j(z_0) = A e_{i-j} + Bu;
$$
$e_{i-j}(z_j)$ is defined (even though $i-j>\ell(z_j)$) by
(IA6) and Lemma \ref{nlemma4.6}. We argue as above that $DT^{i-j}(z_j)u$ splits
correctly at $z_i$.
It remains to check that adding $A \cdot DT^{i-j}(z_j)e_{i-j}$
will only change the angle of $B \cdot DT^{i-j}(z_j) u$ by
$<<e^{-\alpha(i-j)} <d_{\cal C}(z_i)$. This is true because
$$
\|A \cdot DT^{i-j}(z_j)e_{i-j}\| < \frac{|B|}{\angle(e, w_j^*)} b^{i-j}
<|B| b^{\frac{i-j}{2}}
< b^{\frac{i-j}{2}}\|B \cdot DT^{i-j}(z_j)u\|.
$$ \hfill $\square$

%%%%%%%%%%%%%%%%%%%%%%%%%%%%Appendix B.7
\subsection{Distortion during bound periods \ (Sect. \ref{s4.3})}
\label{app-B.7}
%%%%%%%%%%%%%Sublemma B.5 
\begin{sublemma}
$$
\sum_{i=1}^{\mu} K \frac{\Delta_i}{d_{\cal C}(z_i)} << 1.
$$
\label{nsublemmaB.5}
\end{sublemma} 
\noindent {\bf Proof:}  
Since $|\xi_s - z_s| < e^{-\beta s}$ for all $s<\mu$, we have 
$\Delta_i < 2 e^{-\beta i}$. 
Let $h_0 = \frac{1}{2\|DT\|}|\log{\delta}|$. Then 
$$
\sum_{i=h_0+1}^{\mu} K \frac{\Delta_i}{d_{\cal C}(z_i)} 
< K \sum_{i=h_0+1}^{\infty} K e^{-(\beta - \alpha) i}
< K \frac{e^{-(\beta - \alpha) h_0}}{1- e^{-(\beta - \alpha)}}
<<1
$$
provided $\delta$ is sufficiently small. Also,
$$
\sum_{i=1}^{h_0}  K \frac{\Delta_i}{d_{\cal C}(z_i)} < 
(\sum_{i=1}^{h_0} K (e^{\alpha}\|DT\|)^i) \delta << 1
$$
by our choice of $h_0$.
\hfill $\square$

\bigskip
 
\noindent {\bf Proof of Lemma \ref{nlemma4.9}}: (cf. \cite{BC2}, Lemma 7.8)
\ Assuming the lemma
for all $i < \mu$, we give the proof of (\ref{formula3.2})
for step $\mu$; the bound in (\ref{formula3.3}) is proved similarly. 
For notational simplicity, we drop the hat in ${\hat w}_{\mu}^*(\xi_0)$.

\medskip
 
\noindent {\bf Case 1} \ No fold period expires at $z_{\mu}$ and $\mu-1$ is not
a return
time. In this case $w^*_{\mu}(\cdot) = DT(\cdot) w^*_{\mu-1}(\cdot)$. Writing
$C = DT(z_{\mu-1})$, $C^{\prime} = DT(\xi_{\mu-1})$, 
$$
u = \frac{w^*_{\mu-1}(z_0)}{\|w^*_{\mu-1}(z_0)\|} \ \ \  and \ \ \ 
u^{\prime} = \frac{w^*_{\mu-1}(\xi_0)}{\|w^*_{\mu-1}(\xi_0)\|},
$$
we have
\begin{eqnarray*}
\frac{M_{\mu}^{\prime}}{M_{\mu}} & = &
\frac{M_{\mu-1}^{\prime}}{M_{\mu-1}} \cdot
\frac{\|C^{\prime} u^{\prime}\|}{\|C u\|} 
\leq \frac{M_{\mu-1}^{\prime}}{M_{\mu-1}}
\left( 1 + \frac{\| C^{\prime}u^{\prime}-Cu \|}{\| Cu \|}  \right) \\
& \leq & \frac{M_{\mu-1}^{\prime}}{M_{\mu-1}} \left( 1 +
\frac{\| C^{\prime}-C\|}{\| Cu \|} +
\frac{\| C(u-u^{\prime}) \|}{\| Cu \|} \right).
\end{eqnarray*}
Since $\| Cu \| > K^{-1} \delta$,
$\| C - C^{\prime} \| <
K |\xi_{i-1} - z_{i-1} |$ and $\|u - u^{\prime}\| \sim
| \theta_{\mu-1}^{\prime} - \theta_{\mu-1} | <
K b^{\frac{1}{2}} \Delta_{\mu-2}$,
we have
$$
\frac{M_{\mu}^{\prime}}{M_{\mu}} \leq  \frac{M_{\mu-1}^{\prime}}{M_{\mu-1}}
\cdot \left( 1 + K \frac{\Delta_{\mu-1}}{d_{\cal C}
(z_{\mu-1})} \right).
$$
 
\medskip
 
\noindent {\bf Case 2} \ $\mu-1$ is a return time. Then
$$
w^*_{\mu-1}(z_0) = A(z_{\mu-1}) \cdot e(z_{\mu-1})
+ B(z_{\mu-1}) \cdot w_0.
$$
Let
$$
A_0 = \frac{A(z_{\mu-1})}{\|w^*_{\mu-1}(z_0)\|}; \ \ \ \
B_0 = \frac{B(z_{\mu-1})}{\|w^*_{\mu-1}(z_0)\|}.
$$
Then since
$w^*_{\mu}(z_0) = B(z_{\mu-1}) \cdot DT(z_{\mu-1}) w_0$,
we have
$$
\frac{M_{\mu}^{\prime}}{M_{\mu}} = \frac{M_{\mu-1}^{\prime}}{M_{\mu-1}} \cdot
\frac{|B_0^{\prime}|}{|B_0|} \cdot \frac{\| C^{\prime} w_0 \|}{\|C w_0\|}.
$$
Also with $|B_0| \sim d_{\cal C}(z_{\mu-1})$ and
$|B^{\prime}_0 - B_0| \leq |\theta_{\mu-1}^{\prime} - \theta_{\mu-1}|
+ \| e - e^{\prime} \|$,
we get
\begin{equation}
\left| \frac{B_0^{\prime}}{B_0} - 1 \right| <
K \frac{\Delta_{\mu-1}}{d_{\cal C}(z_{\mu-1})}.
\label{nB.7(1)}
\end{equation}
For the last ratio,
$$
\frac{\|C^{\prime}w_0\|}{\|Cw_0\|} \leq
1 + K | \xi_{\mu-1} - z_{\mu-1} |.
$$
This finishes the computation for the magnitude. We record also the estimate
$$
| A_0 - A_0^{\prime} | < K \Delta_{\mu-1}
$$
for use in Case 3. 

\medskip
\noindent {\bf Case 3} \ There exists a return time $j$ whose fold period 
expires
at time $\mu$. In this case
$$
w^*_{\mu}(z_0) = B(z_j) \cdot DT^{\mu-j}(z_j) w_0
+ A(z_j) \cdot DT^{\mu-j}(z_j) e(z_j).
$$
Let
\begin{eqnarray*}
B_0 = \frac{B(z_j)}{\|w^*_j(z_0)\|} \ ,  & \ \ \ \ & A_0 =
\frac{A(z_j)}{\|w^*_j(z_0)\|} \ , \\
C = DT^{\mu-j}(z_j)w_0 \ , & \ \ \ \ & Y = DT^{\mu-j}(z_j) 
e(z_j) \ .
\end{eqnarray*}
As before, all the corresponding quantities for $\xi_0$ carry a prime. 
Then
$$
\frac{M_{\mu}^{\prime}}{M_{\mu}}  =  \frac{M_j^{\prime}}{M_j} \cdot
\frac{\|B_0^{\prime} C^{\prime} + A_0^{\prime} Y^{\prime}\|}{\|B_0 C +
A_0 Y \|}
\leq  \frac{M_j^{\prime}}{M_j} \cdot  \frac{\|C^{\prime}\|}{\|C \|} \cdot
\frac{|B_0^{\prime}|}{|B_0|} \cdot \left( 1 + \frac{\left\| \frac{C^{\prime}}{\|C^{\prime}\|}
- \frac{C}{\|C\|} + 
\frac{A_0^{\prime}Y^{\prime}}{B_0^{\prime}\|C^{\prime}\|}-
\frac{A_0Y}{B_0\|C\|} \right\|}{\left\| \frac{C}{\|C\|} + 
\frac{A_0 Y}{B_0\|C\|}\right\|} \right). 
$$
Since 
$\frac{|A_0|}{|B_0|} \sim \frac{1}{d_{\cal C}(z_j)}$ and $\frac{\|Y\|}{\|C\|} 
\leq
d_{\cal C}^2(z_j)$, it follows that 
$\frac{\|A_0Y\|}{\|B_0C\|} << 1$, giving
$$
\frac{M_{\mu}^{\prime}}{M_{\mu}} \leq \frac{M_j^{\prime}}{M_j} \cdot
\frac{\|C^{\prime}\|}{\|C \|} \cdot \frac{|B_0^{\prime}|}{|B_0|} \cdot
\left( 1+ 2 \left\|\frac{C^{\prime}}{\|C^{\prime}\|}
- \frac{C}{\|C\|}\right\| + 2
\left\|\frac{A_0^{\prime}Y^{\prime}}{B_0^{\prime}\|C^{\prime}\|}+
\frac{A_0Y}{B_0
\|C\|}\right\| \right).
$$
Since both $\{z_s \}_{s=j}^{\mu}$
and
$\{ \xi_s \}_{s=j}^{\mu}$ are bound to a critical segment
$\{\eta_s\}_{s=0}^{\mu-j}$, $\eta_0 \in \Gamma_{\theta N}$, we have
$$
\frac{\|C^{\prime}\|}{\|C\|}  \leq  1 + K \sum_{s=1}^{\mu-j-1}
\frac{\hat{\Delta}_s}{d_{\cal C}(\eta_s)} 
\leq  1 + K \sum_{s=1}^{\mu-j-1} \frac{{\Delta}_{s+j}}{d_{\cal C}(z_{s+j})} 
 =  1 + K \sum_{i=j+1}^{\mu-1} \frac{\Delta_i}{d_{\cal C}(z_i)}
$$
where
\begin{equation}
\hat{\Delta}_s = \sum_{j=1}^s (Kb)^{\frac{j}{4}} |z_{s-j} - \xi_{s-j}|.
\label{nB.7(2)}
\end{equation}
The factor $\frac{|B_0^{\prime}|}{|B_0|}$ is estimated in (\ref{nB.7(1)}). 
This term has no cumulative effect because it is a one-time
addition 
to the exponent in the distortion formula for any given return. Next 
$$
 \left\|\frac{C^{\prime}}{\|C^{\prime}\|}
- \frac{C}{\|C\|}\right\| < {\hat \theta}
$$
where ${\hat \theta}$ is the angle between $C$ and $C^{\prime}$, which is 
smaller than ${\hat \Delta}_{\mu-j-1}$. Now 
$$
\left\|
\frac{A_0^{\prime}Y^{\prime}}{B_0^{\prime}\|C^{\prime}\|}-
\frac{A_0Y}{B_0\|C\|}
\right\| \leq \frac{|A_0|}{|B_0|} \cdot \frac{\|Y^{\prime}-Y\|}{\|C\|}
+ \left\| \frac{A_0}{B_0\|C\|} - \frac{A_0^{\prime}}{B_0^{\prime}\|C^{\prime}\|}
\right\| \|Y^{\prime}\|.
$$
For the first term we have
$$
\frac{|A_0|}{|B_0|} \sim \frac{1}{d_{\cal C}(z_j)}, \ \ \ \ \ \ \
\|Y^{\prime} - Y\| \leq (Kb)^{\mu-j}|\xi_j - z_j|,
$$
and $\|C\| > 1$, where the 
estimate on $\|Y^{\prime}-Y\|$ is from  
(\ref{nformula2.2}) in Corollary \ref{ncorollary2.2}. For the second term, 
\begin{eqnarray*}
\left\| \frac{A_0}{B_0\|C\|} - \frac{A_0^{\prime}}{B_0^{\prime}\|C^{\prime}\|}
\right\| \|Y^{\prime}\| & \leq &
(Kb)^{\mu-j} \frac{|A_0^{\prime}|}{|B_0^{\prime}|} \cdot  \frac{1}{\|C\|} \cdot
\left(
\left| \frac{A_0}{A_0^{\prime}}  \cdot  \frac{B_0^{\prime}}{B_0} -1
\right| +
\left| 1-\frac{\|C\|}{\|C^{\prime}\|} \right| \right) \\
& \leq & \frac{(Kb)^{\mu-j}}{d_{\cal C}(z_j)} \left( \frac{|A_0|}{|A_0^{\prime}|} 
\left|
\frac{B_0^{\prime}}{B_0} -1 \right|
+ \left| \frac{A_0}{A_0^{\prime}} -1 \right| + \left|
1-\frac{\|C\|}{\|C^{\prime}\|}\right|
\right).
\end{eqnarray*}
We again estimate term by term: For the first term,
$$
\frac{(Kb)^{\mu-j}}{d_{\cal C}(z_j)} \cdot \frac{|A_0|}{|A_0^{\prime}|} \cdot
\left| \frac{B_0^{\prime}}{B_0} -1 \right|
\leq
\frac{(Kb)^{\mu-j}}{d_{\cal C}(z_j)} \cdot \frac{\Delta_j}{d_{\cal C}(z_j)} <
\frac{(Kb)^{\frac{\mu-j}{2}}}{d_{\cal C}(z_j)} \cdot \Delta_j
$$
because $b^{\frac{\mu-j}{2}} < d_{\cal C}(z_j)$ by the definition of fold
period.
For the second term,
$$
\frac{(Kb)^{\mu-j}}{d_{\cal C}(z_j)} \cdot \left| \frac{A_0 -
A_0^{\prime}}{A_0^{\prime}}\right|  \leq  \frac{(Kb)^{\mu-j}}{d_{\cal C}(z_j)}
\Delta_j.
$$
Finally, for the third term, we have
$$
\frac{(Kb)^{\mu-j}}{d_{\cal C}(z_j)} \left| 1- \frac{\|C^{\prime}\|}{\|C\|} 
\right|
 \leq 
\frac{(Kb)^{\mu-j}}{d_{\cal C}(z_j)} \sum_{s=1}^{\mu-j-1}
\frac{\hat{\Delta}_s}{d_{\cal C}(z_{s+j})}
$$
where ${\hat \Delta}_s$ is as in (\ref{nB.7(2)}). 
We also have $b^{\frac{\mu-j}{2}} < d_{\cal C}(z_{s+j})$, 
for no fold period starting 
at time $s+j$ 
extends beyond index $\mu$.
Also, $b^{\frac{\mu-j}{2}} \cdot \hat{\Delta}_s \leq \Delta_j$
for all $s$, $0 < s < \mu-j$. Therefore the third term is again 
bounded by
$K\frac{\Delta_j}{d_{\cal C}(z_j)}$.
 
Observe further that if we replace $z_0$ by another point
$\xi_0^{\prime}$ which is bounded to $z_0$, the same argument above 
continues to work
with $\Delta_i(\xi_0, z_0)$ replaced by $\Delta_i(\xi_0, \xi_0^{\prime})$.
This completes the proof.

\hfill $\square$
 
%%%%%%%%%%%%%%%%%%%%%%%%%%%%%Appendix B.8
\subsection{Quadratic behavior \ (Sect. \ref{s4.3})}
\label{app-B.8}
 
Let $\xi_0(s)$ and $z_0$ be as in Lemma \ref{nlemma4.11}.
We begin with the following
 
\bigskip
\noindent {\bf A priori estimate on $\xi_\mu(s)-z_\mu$ :} \ (cf. \cite{BC2},
p.144-147) \ Let $t_0(s)$ be a unit vector to $\gamma$ at $\xi_0(s)$, and let
$t_\mu=DT^\mu t_0$. We split $t_0$ using $e_\mu$ to get
$$
t_0 = A_0 e_\mu + B_0 {\tiny \left(\!\!\begin{array}{c} 0 \\ 1 
\end{array}\!\!\right)},
\ \ \ \ {\rm so \ that} \ \ \ \ t_\mu = A_0 DT^\mu e_\mu + B_0 w_\mu.
$$
Writing
$$
w_{\mu} = w_{\mu}(0) +(w_{\mu} - w_{\mu}(0))
= w_{\mu}(0) + (w_{\mu}^* - w_{\mu}^*(0)) + (E_{\mu} - E_{\mu}(0))
$$
where
$$
E_{\mu} = \sum_{j \in S_{\mu}} A_j DT^{\mu-j} e_{\ell_j}
$$
and $S_{\mu}$ is the collection of $j$ such that the fold period
begun at time $j$ extends beyond time $\mu$, we have
\begin{equation}
\xi_\mu(s)-z_\mu\ = \ \int_0^s t_\mu(u)du \ =\ w_{\mu}(0) \int_0^{s} B_0(u) du \
 +\  I\  +\  II
\  +\ III
\label{nB.8(1)}
\end{equation}
where
$$
I=\int_0^{s} A_0 DT^\mu e_\mu, \ \ \ \
II= \int_0^{s} B_0(w^*_\mu-w^*_\mu(0)), \ \ \ \
III=\int_0^{s} B_0(E_\mu-E_\mu(0)).
$$
Since $A_0 \approx 1$, $\|I\| \leq (Kb)^\mu s$. We claim that
$$
\|II\|, \ \|III\| \ \leq K e^{2\alpha \mu} \| w^*_{\mu}(0) \|
 \int_{0}^{s} u \ (\sup_{i \leq \mu} |z_i - \xi_i(u)|) \ du.
$$
The norm of $II$ is estimated using the distortion estimate in 
Lemma \ref{nlemma4.9}.
To estimate $\|III\|$, we have, for each $j \in S_\mu$,
\begin{eqnarray*}
\| A_j DT^{\mu-j}(\xi_j) e(\xi_j) - A_j(0) DT^{\mu-j}(z_j)
e(z_j) \|  \ \ \ \ \ \ \ \ \ \ \ \ \ \ \ \ \ \ \ \  \\
\leq (Kb)^{\mu-j} |A_j -  A_j(0)|  + |A_j(0)| \  \|  DT^{\mu-j}(\xi_j) e(\xi_j)
- DT^{\mu-j}(z_j)e(z_j) \|.
\end{eqnarray*}
From the distortion estimate in appendix \ref{app-B.7},
$$
|A_0 -  A_0(0)| < K \| w^*_j(0) \|
e^{2\alpha j} \sup_{i \leq j} |z_i - \xi_i| .
$$
For the second term we have $|A_j(0)| \leq \| w^*_j(0) \| e^{\alpha j}$
because $w^*_j(0)$ splits correctly at time $j$, and
$\| w^*_j(0) \| \leq e^{\alpha \mu}\| w^*_\mu(0) \|$ by Lemma \ref{nlemma4.7}.
Finally, $\| DT^{\mu-j}(\xi_j) e(\xi_j) - DT^{\mu-j}(z_j)e(z_j) \| \leq (Kb)^{
\mu-j}|\xi_j-z_j|$ by Corollary \ref{ncorollary2.2}, 
and $B_0(s) \approx 2K_1 s$.

\bigskip
 
\noindent {\bf Proof of Lemma \ref{nlemma4.11}:} \ We will
show that for the $\mu$ and $s$ in question, the first term in
(\ref{nB.8(1)}) is the dominating one. Let
$$
U_\mu \ := \ Ke^{4\alpha \mu} \sup_{j \leq \mu} \|w^*_j(0)\|
$$
where $K$ is the constant in the bound for $\|II\|$ and $\|III\|$ above.
We choose $\mu_0$ large enough that $e^{7\alpha \mu_0}e^{-\beta \mu_0}<<1$,
and assume $\delta$ is small enough that $U_{\mu_0}\delta^2<<1$.
We will show inductively first the weaker statement
 
(i) \ $|\xi_j(s)-z_j|<U_j s^2$
 
\noindent and then the stronger statement
 
(ii) \ $|\xi_j(s)-z_j|=K_1(1 \pm \varepsilon_1)\|w_j(0)\| s^2$.
 
Assume this has been done for all $j<\mu$. To prove (i) for $j=\mu$,
we need to first verify that $U_\mu s^2<<1$. In the case where $\mu>\mu_0$,
we use
$$
\sup_{j \leq \mu} \|w^*_j(0)\| < e^{\alpha \mu}\|w^*_{\mu}(0)\| \leq e^{2\alpha
\mu}\|w_\mu(0)\| \leq Ke^{2\alpha \mu}\|w_{\mu-1}(0)\|
$$
(see Lemmas \ref{nlemma4.7} and \ref{nlemma4.6}) 
combined with (ii) for step $\mu -1$ to get
$$
U_\mu s^2
\leq Ke^{4\alpha \mu} \ (Ke^{2\alpha \mu}\|w_{\mu-1}(0)\|) s^2
\approx K^2 e^{6\alpha \mu} K_1^{-1} |\xi_{\mu-1}(s)-z_{\mu-1}|
$$
$$
< K^2 e^{6\alpha \mu} K_1^{-1}e^{-\beta (\mu-1)} <<1.
$$
Noting that $\int_0^{s} B_0 \approx K_1 s^2$, we see from our {\it a priori} 
estimate that
$$
|\xi_\mu(s)-z_\mu| \ \leq \ \ \|w_{\mu}(0)\| K_1 s^2 \ +\  (Kb)^\mu s
\  +\ 2U_\mu \int_0^{s} Ku \ [\sup_{i < \mu} |z_i - \xi_i(u)|] \ du.
$$
With the quantity inside square brackets being $<U_\mu u^2$ by (i)
from the previous step, this is
$$
<(U_\mu s^2)e^{-\alpha \mu} + (Kb)^{\frac{\mu}{2}}s^2 +K(U_\mu s^2)^2
\ <\ U_\mu s^2.
$$
The proof of (ii) for step $\mu$ now follows immediately.
\hfill $\square$
 
%%%%%%%%%appendix B.9
\subsection{Proof of Lemma \ref{nlemma6.2} \ (Sect. \ref{s6.2})}
\label{app-B.9}
 
We begin with a scenario for which one sees easily that the assertion
in this lemma holds: Suppose for $1 \leq j \leq i-s$, $\|DT^j(z_s)\|
>\kappa^j$ for some $\kappa>>b^{\frac{1}{2}}$, and that $z_s$ is
bounded away from ${\cal C}^{(0)}$. Then $e_{i-s}(z_s)$ is well defined
and has slope $>K^{-1}$. Suppose, in addition, that $z_s$ is out of
all fold periods,
so that $w_s$ is a $b$-horizontal vector. Then
$$
\|DT^{i-s}(z_s)\|  \ \|w_s\| \leq K \|DT^{i-s}(z_s) w_s\|
=  \|w_i\|.
$$
This together with $\|w_s\|>c^{c''s}$ (which follows from $\|w^*_s\|>
e^{cs}$) gives the desired estimate.
 
Now, intuitively, the behavior of  $\|DT^j(z_s)\|$ is a little different
 just before or after a return to ${\cal C}^{(0)}$. This motivates the
following definition: If $t$ is a return time to ${\cal C}^{(0)}$
for $z_0$, let $\ell_t$ denote
its fold period and let $I_t:=(t-5\ell_t, t+\ell_t)$.

%%%%%%%%%%%Claim B.1 
\begin{claim} By modifying $I_t$ slightly to $\tilde I_t=(t-(5\pm
\varepsilon)\ell_t, \ t+(1\pm \varepsilon)\ell_t)$,
we may assume they have a nested structure.
\label{nclaimB.1}
\end{claim}
 
\noindent {\it Proof of Claim \ref{nclaimB.1}:} \ We consider $t=0,1,2, \cdots$
in this order, and determine, if $t$ is a return time, what  $\tilde I_t$ will
 be. The right end point of $\tilde I_t$ is determined by the following
algorithm: Go to $t+\ell_t$, and look for the largest $t'$
{\it inside the bound period initiated at time $t$} with the property
that $t'-5\ell_{t'}<t+\ell_t$. If no such $t'$ exists, then $t+\ell_t$
is the right end point of $\tilde I_t$. If $t'$ exists, then the new candidate 
end point is $t'+\ell_{t'}$, and the search continues.
For the same reasons as in Sect. \ref{s4.1}, the increments in length are
exponentially small and the process terminates.
 
As for the left end point of $\tilde I_t$, it is possible that
$t-5\ell_t \in \tilde I_{t'}$ for some $t'$ the bound period initiated
at which time does not extend to time $t$.
This means that $\ell_{t'}<<\ell_t$, and since we assume a nested structure
has been arranged for $\tilde I_{t'}$ for all $t'<t$, we simply extend
the left end of $\tilde I_t$ to include the largest $\tilde I_{t'}$ that
it meets. \hfill $\diamondsuit$
 
\medskip
  
Let us assume this nested structure and write $I_t$ instead of $\tilde I_t$
from here on.

%%%%%%%%%%%%Claim B.2
\begin{claim}
\label{nclaimB.2}
For $s \not \in \cup I_t$, we have, for all $j$ with
$1 \leq j <i-s$,
$$
\|w_{s+j}\| \geq b^{\frac{j}{9}}\|w_s\|.
$$
\end{claim}
 
\noindent {\it Proof of Claim \ref{nclaimB.2}:}
We fix $j$ and let $r$ be such that $z_r$ makes the deepest
return between times $s$ and $s+j$. Let $j'$ be the smallest integer
$\geq j$ such that $z_{s+j'}$ is outside of all fold periods.
Then from Sect. \ref{s4.2} , it follows that
\begin{equation}
\|w_{s+j}\| \geq K^{-K\theta(j'-j)}\|w_{s+j'}\| \geq
K^{-K\theta (j'-j)}d_{\cal C}(z_r) \|w_s\| \approx
K^{-K\theta (j'-j)}b^{\frac{\ell_r}{2}} \|w_s\|.
\label{nequation-B.9(1)}
\end{equation}
 
{\it Case 1.} $s+j \not \in I_r$. \ In this case, $6 \ell_r<j$
since $I_r$ is sandwiched between $s$ and $s+j$,
and $j'-j\leq \ell_r$ because $r$ is the deepest return.
The rightmost quantity in (\ref{nequation-B.9(1)}) is therefore
$>K^{-\ell_r}b^{\frac{\ell_r}{2}} \|w_s\|
>b^{\frac{j}{9}}\|w_s\|$.
 
{\it Case 2.} \  
$s+j \in I_r$. \ The argument is as above, except
we only have  $5 \ell_r<j$.
 
\noindent This completes the proof of the claim. \hfill $\diamondsuit$
 
\medskip
 
As noted in the first paragraph, Claim \ref{nclaimB.2} implies the assertion in
Lemma \ref{nlemma6.2} for $s \not \in \cup I_t$ provided $z_s$ is bounded away 
from
${\cal C}^{(0)}$. This last proviso is easily removed by considering
$z_{s+1}$ if necessary.
 
It remains to prove the lemma for $s \in \cup I_t$.
Let $I_r$ be the maximal $I_t$-interval containing $s$. Observe
that $6\ell_r<K\alpha \theta s$ (recall that $z_0$ obeys (IA2))
and $\|w_i\|>e^{c''i}$ for some $c''>0$.
If $i \in I_r$, then $\|DT^{i-s}(z_s)\|<K^{6\ell_r}<<
e^{\frac{1}{2} c''i}<e^{-\frac{1}{2}c''s}e^{c''i}
<e^{-\frac{1}{2}c''s}\|w_i\|$.
If $i \not \in I_r$, let $s'=r+\ell_r$. Then $s' \not \in \cup I_t$, and
$$
\|DT^{i-s}(z_s)\| \leq \|DT^{s'-s}(z_s)\| \cdot \|DT^{i-s'}(z_{s'})\|
\leq K^{6\ell_r} \cdot Ke^{-c's'}\|w_i\|.
$$

%%%%%%%%%%%%%%%%%%%%%%%%%%%%Appendix B.10
\subsection{Initial data for critical curves \ (Sect. \ref{s6.3})}
\label{app-B.10}
 
\noindent {\bf Proof of Lemma \ref{nlemma6.4}:} \ Let
$J_i:=[\hat a-\rho^{2i}, \hat a+\rho^{2i}]$.
Assume for all $i<n$ that the
following has been proved:
 
(i) $J_i \subset \tilde \Delta_i$;
 
(ii) $\Gamma_i^i(\hat a)$ has a smooth continuation on $J_i$
and ${\cal C}^{(i)}$ deforms continuously;
 
(iii) for all $z \in \Gamma_{i,i}$, $\|\frac{dz}{da}\| \leq K^i$.
 
\noindent We now prove (i)--(iii) for $i=n$.
 
First we verify that for all $a \in J_n$ and $z_0 \in \Gamma_{n-1,n-1}$,
(IA2) and (IA4) hold up to time $n$. This is true for $a=\hat a$.
For $a \in J_n$,
$|z_0(a)-z_0(\hat a)|<\rho^{2n}K^{n-1}$,
so that $|z_j(a)-z_j(\hat a)|<\rho^{2n}K^{2n}$ for all $j \leq n$.
We may assume that $\rho K$ is $<<1$. It then follows from
 the discussion at the beginning of Sect. \ref{s6.3.1} that $\Gamma_{n,n}(a)$
is well defined,
proving (i).

To prove (ii), we fix an arbitrary $\tilde a \in J_n$,
a component $Q^{(n-1)}$ of ${\cal C}^{(n-1)}$,
and show that
every segment of $\partial R_n(\tilde a) \cap Q^{(n-1)}(\tilde a)$ has
a continuation to a segment of $\partial R_n(a) \cap Q^{(n-1)}(a)$.
Let $\tilde \omega$ be a segment of this kind, and let $\omega(a):=
T_a^n (2T^{-n}_{\tilde a} \tilde \omega)$ where
$2T^{-n}_{\tilde a} \tilde \omega$ refers to the segment in $\partial R_0$
with the same midpoint as $T^{-n}_{\tilde a} \tilde \omega$ and two times
as long.
Observe that as we vary our parameter from $\tilde a$ to $a$,
the segment $\omega(a)$  cannot intersect the horizontal boundaries
of $Q^{(n-1)}(a)$. Thus the only way
$\omega (a)$ can fail to traverse fully $Q^{(n-1)}(a)$ is that it has
moved sufficiently far from $\omega ({\tilde a})$
in the {\it horizontal} direction.
We know this cannot happen because
$|T_{\tilde a}^n-T_a^n| \leq \rho^{2n} K^n$
which is $<<\rho^n$. This proves (ii).

It remains to prove (iii). Consider
$\bar z(a)=(\bar x(a), \bar y(a)) \in \Gamma_{n,n}(a)$, and let
$y = \psi(x, a)$ denote the
$C^2(b)$-curve in $\partial R_n(a)$ containing $\bar z(a)$. Then
$$
q_n(\bar x(a), \psi(\bar x(a), a), a) = \partial_x \psi(\bar x(a), a)
$$
where $q_n(x, y, a)$ is the slope of the contractive vector of order
$n$ at $z = (x, y)$. Taking derivative with respect to $a$ on both sides
of the last equation, we have
$$
\partial_x q_n \cdot \frac{d \bar x}{da} + \partial_y q_n \cdot
(\partial_x \psi \cdot \frac{d \bar x}{da} + \partial_a \psi)
+ \partial_a q_n =
\partial_{xx} \psi \cdot \frac{d \bar x}{da} + \partial_{ax} \psi.
$$
This implies
\begin{equation}
\frac{d \bar x}{da} =
\frac{\partial_{xa} \psi -  \partial_y q_n \cdot \partial_a \psi
- \partial_a q_n}{\partial_x q_n + \partial_y q_n \cdot \partial_x
\psi - \partial_{xx} \psi}.
\label{B.14(1)}
\end{equation}
Since $\partial_x \psi, \ \partial_{xx} \psi = {\cal O}(b)$,
$|\partial_x q_n| > K_1$ and $|\partial_y q_n|
< K$ (see Corollary \ref{ncorollary2.2} and Lemma \ref{nlemma2.9}),
the denominator on the right-hand side is
bounded away from zero.
In the numerator, we have $|\partial_y q_n|, | \partial_a q_n|<K$,
and we need to estimate $\partial_a \phi(x, a)$ and 
$\partial_{ax} \phi(x, a)$.
 
For this purpose we write the horizontal curve $y = \psi(x, a)$
in parametric form $x = X(t, a), y = Y(t, a)$ where $t$ is the x-coordinate
of $T_a^{-n}(x,y)$, i.e., $(t, \pm b) \in \partial R_0$ and
$$
(X(t, a), Y(t, a)) = T_a^n(t, \pm b).
$$
Let $t = t(x, a)$ be defined by
$\psi(x, a) = Y(t(x, a), a)$.
Then
$$
\partial_a \psi = \partial_t Y(t, a) \cdot \partial_a t(x, a) +
\partial_a Y(t, a).
$$
Clearly, $|\partial_t Y(t, a)|<K^nb$ and $|\partial_a Y(t, a)|<K^n$.
One way to bound $\partial_a t(x, a)$ is to write it as
$$
\partial_a t(x, a) = - \frac{\partial_a X(t, a)}{\partial_t X(t, a)}.
$$
Since $|\partial_t X(t, a)| > 1$ (recall that
$T_a^n|\partial R_0$ is controlled), this term is also $<K^n$.
Similar considerations yield $|\partial_{ax} \psi(x, a)| < K^n$.
We have proved $\frac{d\bar x}{da} < K^n$.
The corresponding estimate for $\frac{d \bar y}{da}$ follows immediately since
$$
\frac{d \bar y}{da} = \partial_x \psi \frac{d\bar x}{da} + \partial_a \psi.
$$
 
\bigskip
We record an estimate needed in the proof of Lemmas \ref{nlemma6.5} and
\ref{nlemma6.6}. 
Taking derivatives with respect to $a$ one more time on both sides of
(\ref{B.14(1)}) and estimating corresponding terms (using again Corollary \ref{ncorollary2.2} and Lemma \ref{nlemma2.9}), we obtain
$|\frac{d^2 \bar x}{da^2}| < K^n$. This estimate requires that $T_{a,b}$
be $C^3$.
 
Recall the following lemma due to Hadamard:
 
\begin{lemma}{\rm (Hadamard)}
\label{nlemmaB.1} \
Let $g \in C^2(0, L)$ be such that $|g| \leq M_0$ and
$|g^{\prime \prime}| < M_2$.
If $4M_0 < L^2$, then
$$
|g^{\prime}| \leq \sqrt{M_0}(1 + M_2).
$$
\end{lemma}
 
\bigskip
\noindent {\bf Proof of Lemma \ref{nlemma6.5}:} \ 
Let $z^{(n)}=(x^{(n)}, y^{(n)})$.
For our pusposes, let $g(a)=x^{(n)}(a)-x^{(n-1)}(a)$ and
$L=2\rho^{2n}$. Then $M_0 = b^{\frac{n}{4}}$ and
$M_2=K^n$. Thus $|\frac{d x^{(n)}}{da}|<b^{\frac{n}{8}}K^n
<b^{\frac{n}{9}}$. A similar estimate holds for $y^{(n)}$. \hfill
$\square$
 
\bigskip
\noindent {\bf Proof of Lemma \ref{nlemma6.6}:} \ 
Let $z^n(a) = (x^n(a), y^n(a))$,
$z^m(a) = (x^m(a), y^m(a))$, and let $y = \psi(x, a)$ be the $C^2(b)$-curve
segment in $\partial R_n$ containing both $z^n(a)$ and $z^m(a)$.
Arguments similar to those used to prove
$|\partial_{xa} \psi | < K^n$ can also be used to prove that
the $C^3$-norm of $\psi$ is $< K^n$.
 
Let $P_n$ and $Q_n$ be the numerator and denominator on the right hand side
of (\ref{B.14(1)}), and similarly for $P_m$ and $Q_m$. Then
$$
\frac{dx^n}{da} - \frac{dx^m}{da}
=\frac{P_nQ_m-P_mQ_n}{Q_nQ_m} = \frac{(P_n-P_m)Q_n + (Q_m-Q_n)P_n}{Q_nQ_m}.
$$
As observed in the proof of Lemma \ref{nlemma6.4}, $|Q_m|, |Q_n|>K^{-1}$.
$|Q_m|, |P_n| < K^n$.
It remains therefore to estimate
$|Q_m - Q_n|$ and $|P_m - P_n|$.
Let $q_n$ and $q_m$ denote the slopes of $e_n$ and $e_m$ respectively.
Fixing $a$ and omitting it in the arguments of the functions
below, we have
\begin{eqnarray*}
|Q_m-Q_n| & \leq & |\partial_x q_m(z^m) - \partial_x q_n(z^n)| \ + \
|\partial_y q_m (z^m)\cdot \partial_x \psi (z^m)
- \partial_y q_n (z^n) \cdot \partial_x \psi (z^n)| \\
& & + \ |\partial_{xx} \psi (z^m) - \partial_{xx} \psi (z^n)|.
\end{eqnarray*}
The second difference, for example, is
\begin{eqnarray*}
\ \ \ \ \ & \leq & \ |\partial_y q_n(z^n)| \ |\partial_x \psi(z^m)-\partial_x \psi(z^n)|
\ +\ |\partial_x \psi(z^m)| \ |\partial_y q_m(z^m)-\partial_y q_n(z^m)| \\
& & + \ |\partial_x \psi(z^m)| \ |\partial_y q_n(z^m)-\partial_y q_n(z^n)|.
\end{eqnarray*}
This is $<(Kb)^{\frac{n}{4}}$ since $\|\psi\|_{C^3} <K^n$,
$|\partial_{xx}q|, |\partial_x \partial_y q|
<K$ (Corollary \ref{ncorollary2.2}), 
$|q_m(z^n)-q_n(z^n)|<(Kb)^n$ 
(Lemma \ref{nlemma2.1}) and
$|z^m-z^n|<(Kb)^{\frac{n}{4}}$ (Lemma \ref{nlemma2.10}).
The other terms in $|Q_m-Q_n|$ and $|P_m-P_n|$
are estimated similarly. \hfill $\square$

%%%%%%%%%%%%%%%%%%%%%%%%%%%%%%%Appendix B.11
\subsection{Dynamics of critical curves \ (Sect. \ref{s6.4})}
\label{app-B.11}
 
\noindent {\bf Proof of Lemma \ref{nlemma6.8}:} \ Let  $\hat z_0$ be an arbitrary
 critical point. First we observe that as functions
of $a$, $z_i(a)$ and $\hat z_0(a)$ move at very different speeds:
$\|\frac{d}{da} z_i(a)\|\sim \|w_i(a)\| >e^{ci}$ by Proposition 
\ref{nproposition6.1},
whereas from Sect. \ref{s6.3} we have $\|\frac{d}{da} \hat z_0(a)\|
<K$.
 
Next we consider $z_i(a) \in Q^{(k-1)}(a) \setminus Q^{(k)}(a)$
for some $k<<i$, so that $\phi_a(z_i(a)) \in \partial Q^{(k-1)}(a)$,
 and study the relative movements of $z_i$, $\phi(z_i)$ and the relevant
critical regions as $a$ varies. For definiteness, let us assume that
$z_i$ is in the right component of $(Q^{(k-1)} \cap R_k) \setminus Q^{(k)}$
(which we call $A$), and that it moves left as $a$ increases.  (See Fig. 1
in Sect. \ref{s1.2}.) In horizontal distance,
it follows from the first paragraph that
{\it relative to $\phi(z_i)$}, $z_i$ is moving left at a speed
$>K^{-1}e^{ci}-K$, which we assume to be $>>1$. We do not have analytic
estimates on the relative vertical movements of $\phi(z_i)$ and $z_i$,
but note that since $z_i \not \in \partial R_k$, it must enter
$A$ through its right vertical boundary and exit through the left.
As $z_i$ meets these vertical boundaries, it crosses them instantaneously
due again to the horizontal speed differential between $z_i$ and
the critical points which determine these regions.
 
What we have shown is that the function $a \mapsto \phi_a(z_i(a))$ is
continuous  except at a  discrete set of points corresponding to when
$z_i(a)$ crosses a vertical boundary of some $Q^{(k)}$.
If $a$ and $a'$ are the entry and exit parameters for $Q^{(k-1)}
\setminus Q^{(k)}$ as above, we have $|a-a'|<\rho^{k-1}(K^{-1}e^{ci}-K)^{-1}$,
and consequently $|\phi_a(z_i)-\phi_{a'}(z_i)|<K' \rho^{k-1}e^{-ci}$.
As $z_i$ crosses the vertical boundary into $Q^{(k)}$,
 a jump in $\phi_a(\cdot)$ occurs
due to our rule for selecting binding points;
this jump is  $<b^{\frac{k-1}{4}}$.
 
As we continue to move toward the cricial set,
either $\gamma_i$ ends or we enter
the ``last" $Q^{(k)}$ available at step $i$, with $k \sim \theta i$.
Let $\bar a$ be the parameter that corresponds to the end point
of $\gamma_i$ or where $d_{{\cal C}(\bar a)}(z_i(\bar a)) =
e^{-\frac{\alpha i}{2}}$, whichever is reached first,
and let $\bar z= \phi_{\bar a}(z_i(\bar a))$.
We will use $\bar z$ as our ``binding point" for $\gamma_i$.
The ``error" in this choice for $z_i(a)$, i.e.
$|\ |z_i(a) - \bar z|_h - d_{{\cal C}(a)}(z_i(a)) \ |$,
is less than the total variation of $a \mapsto
\phi_a(z_i(a))$ between $a$ and $\bar a$.
We have proved that this is $<Ke^{-ci}d_{{\cal C}(a)}(z_i(a))
+Kb^{\frac{k-1}{4}}$ where
$\rho^k \sim d_{{\cal C}(a)}(z_i(a))$. \hfill $\square$
 
\bigskip
\noindent {\bf Proof of Lemma \ref{nlemma6.9}:} \ Let
$\tilde p=$ min $\{p_a(z_i(a)):z_i(a) \in I_{\mu j}\}$.
Then by Corollary \ref{ncorollary4.2}(a) and the last lemma,
 $\tilde p <K|\mu|$.
Assertion (a) in
Lemma \ref{nlemma6.9} is obvious for $j \leq \ell$ where $\ell$ is the common
fold period.
For $\ell < j \leq \tilde  p$, we have:
$$
|z_{i+j}(a)-z_{i+j}(a')| \ \leq \ {\rm length}(\omega_j) \ = \
\int_{\omega_0} \frac{\|\tau_{i+j}\|}{\|\tau_i\|}
 \ \sim \ \int_{\omega_0} \frac{\|w_{i+j}\|}{\|w_i\|}.
$$
Since $z_i$ is outside of fold periods and $w_i$ splits correctly,
we have, for $j>\ell$, $\|w_{i+j}\|/\|w_i\| \sim e^{-\mu}\|w_j(z_i)\|$.
Furthermore, if $\hat z_0=\phi(z_i)$ and $\tilde p=p(z_i(\tilde a))$, then
$$
e^{-\mu}\|w_j(z_i)\| \sim e^{-\mu}\|w_j(\hat z_0)\| \sim
e^{-\mu}\|w_j(\hat z_0(\tilde a))\|,
$$
the first $\sim$ coming from Lemma \ref{nlemma4.9},
and the second from
the fact that $|\phi_a(z_i(a))-\phi_{\tilde a}(z_i(\tilde a))|
<e^{-ci}<<K^j$ for $j<K|\mu| <K\alpha i$. Thus using the distance formula
(\ref{nequation4.3(1)}) in Lemma \ref{nlemma4.11} for $T_{\tilde a}$, we have
$$
\int_{\omega_0} \frac{\|\tau_{i+j}\|}{\|\tau_i\|}
\ \sim \ e^{-2\mu}\frac{1}{\mu^2} \| w_j(\hat z_0(\tilde a))\|
\ \sim \ \frac{1}{\mu^2} |z_{i+j}(\tilde a)-\hat z_j(\tilde a)|
\ < \ \frac{1}{\mu^2} e^{-\beta j}.
$$
 
This completes the proof of (a); (c) follows from
$\|w_{\tilde p}(z_i(a))\| \sim \|w_{\tilde p}(z_i(\tilde a))\|$
and Proposition \ref{nproposition6.1}.
 
It remains to prove (b). From (a) we have that $z_{i+\tilde p}$ is out of
all fold periods whenever $z_{i+\tilde p}(\tilde a)$ is.
To show that the slopes of $\tau_{i+\tilde p}$ are $<K(\delta)$,
we use Lemma \ref{nlemma6.3}: the $w_{i+\tilde p}$-vectors
are b-horizontal, so it suffices to show
that $\|w_s\| \leq K\|w_{i+\tilde p}\|$ for all $s< i+ \tilde p$.
For $s \geq i$, this is true by comparison with $a=\tilde a$;
for $s<i$, $\|w_s\| \leq \|w_i\|$ because $z_i$ is a free return.
Finally, the small
slope of $\omega_{\tilde p}$ allows us to reverse the inequalities
displayed above to conclude that
$|\omega_{\tilde p}| \geq \frac{1}{\mu^2} e^{-\beta \tilde p}$.
\hfill $\square$
 
%%%%%%%%%%%%%%%%%%%%%%%%%%%%%%%Appendix B.12
\subsection{Distortion estimate for critical curves \ (Sect. \ref{s6.4})}
\label{app-B.12}
Let $J$ be a parameter interval satisfying all the
assumptions made in Proposition 
\ref{nproposition6.2}, and let $a, a^{\prime} \in J$.
Assume that
$z_i(a)$ and $z_i(a^{\prime})$ are free returns, and that they lie in the 
same $I_{\mu j}$ with $\mu < \alpha i$. 
Write $\xi_0(a) = z_i(a)$ and 
$w_k(\xi_0(a)) = DT_a^k(\xi_0(a))
{\tiny \left(\!\!\begin{array}{c} 0 \\ 1 \end{array}\!\!\right)}$. Let
${\tilde p} = p(z_i(\tilde a))$ be the bound period and ${\hat z}_0(\tilde a)
= \phi(z_i(\tilde a))$ the binding point in the proof of Lemma \ref{nlemma6.9}.
For $k < {\tilde p}$, let $\{ w^*_k(\xi_0(a)) \}$ be given by the splitting
algorithm taken with respect to the orbit segment 
$\{ {\hat z}_k(\tilde a) \}_{k=0}^{\tilde p}$, and write $w^*_k(\xi_0(a)) = 
M_k e^{i  \theta_k(\xi_0(a))}$. The corresponding quantities for 
$\xi_0(a^{\prime}) = z_i(a^{\prime})$ are defined analogously.

\begin{sublemma}
\label{nsublemmaB.6}
For $k < {\tilde p}$,  
$$
\frac{M_k(\xi_0(a^{\prime}))}{M_k(\xi_0(a))}, \ \ \  
\frac{M_k(\xi_0(a))}{M_k(\xi_0(a^{\prime}))} \leq
\exp \{ K \sum_{j=1}^{k-1}
\frac{\Delta_j(a, a^{\prime})}{d_{\cal C}(\hat{z}_j(a))} \}
$$
and
$$
|\theta_k(\xi_0(a)) - \theta_k(\xi_0(a^{\prime}))| 
< (Kb)^{\frac{1}{2}} \Delta_{k-1}(a, a^{\prime})
$$
where
$$
\Delta_j(a, a^{\prime}) =
\sum_{s=1}^j (Kb)^{\frac{s}{4}}(| \xi_{j-s}(a) - \xi_{j-s}(a^{\prime}) | +
 | a - a^{\prime} |).
$$
\end{sublemma}
\noindent {\bf Proof:} The computation is similar to that in Appendix 
\ref{app-B.7}, modulo the following adaptations to accommodate for the
fact that
different parameter values are involved in the present situation:
\begin{itemize}
\item[(i)] Replace 
$|DT(\xi) - DT(\xi^{\prime})|  < K |\xi - \xi^{\prime}|$ 
by
$$
|DT_a(\xi(a)) - DT_{a^{\prime}}(\xi(a^{\prime}))|
< K (|\xi(a) - \xi(a^{\prime})| + |a - a^{\prime}| ).
$$
 
\item[(ii)] Replace 
$|e - e^{\prime}|< K|\xi - \xi^{\prime}|$ 
by
$$
|e(a) - e(a^{\prime})| <
K ( |\xi(a) - \xi(a^{\prime})| + |a - a^{\prime}|).
$$
 
\item[(iii)] Replace 
$|Y - Y^{\prime}| < (Kb)^{\mu - j} |\xi - \xi^{\prime}|$ by
$$
|Y(a) - Y(a^{\prime})| <
(Kb)^{\mu - j}(|\xi(a) - \xi(a^{\prime})| + |a - a^{\prime}|).
$$
\hfill $\square$
\end{itemize}

\medskip

Next we prove a version of Sublemma \ref{nsublemmaB.6} with
${\tiny \left(\!\!\begin{array}{c} 0 \\ 1 \end{array}\!\!\right)}$ 
replaced by 
$u_i(a) := \frac{w_i(z_0)(a)}{\|w_i(z_0)(a)\|}$.
\begin{sublemma}
\label{nsublemmaB.7}
$$
\frac{\|DT^{\tilde p}_a(\xi_0(a)) u_i(a)\|}{\| DT^{\tilde p}_{a^{\prime}}
(\xi_0(a^{\prime})) 
u_i(a^{\prime}) \| }
< \exp \{ K \frac{|\xi_0(a) - \xi_0(a^{\prime})| }{ e^{-\mu}} \}
$$
\end{sublemma}
\noindent {\bf Proof:} 
The proof uses the fact that both
$u_i(a)$ and $u_i(a^{\prime})$ split correctly. 
Writing
$$
u_i(a) = A(a)e(a) 
+ B(a){\tiny \left(\!\!\begin{array}{c} 0 \\ 1 \end{array}\!\!\right)},
$$
we have
$$
DT^{\tilde p}_a(\xi_0(a)) u_i(a) = A(a) DT^{\tilde p}_a(\xi_0(a)) e(a) + 
B(a) w_{\tilde p}(\xi_0(a)).
$$
The proof is similar to that of Case 3 of Lemma \ref{nlemma4.9}, and 
Sublemma \ref{nsublemmaB.6} is used to compare 
$w_p(\xi_0(a))$ and $w_p(\xi_0(a^{\prime}))$.  
\hfill $\square$

\medskip
 
\noindent {\bf Proof of Proposition \ref{nproposition6.2}:} \ 
In view of Proposition \ref{nproposition6.1}, it suffices 
to show that there exists a constant  $K > 0$ such that
$$
\frac{1}{K} < \frac{\mid w_n(z_0(a)) \mid }{\mid w_n(z_0(a^{\prime}))\mid }
< K.
$$
Divide the time interval $(1, n)$ into
bound and free period according to Lemma \ref{nlemma6.9}. As usual we denote 
free return times as
$t_k$, $1 \leq k < q$, and the bound period at $t_k$ as $p_{t_k}$.
Write
$$
\log{\frac{\| w_n(z_0(a)) \| }{\| w_n(z_0(a^{\prime}))\| }}
= \sum_{k<q} S_k^{\prime} + \sum_{k<q} S_k^{\prime \prime}
$$
where
$$
S^{\prime}_k =
\log{\frac{\| DT^{p_k}_{a}(z_{t_k}(a)) u_{t_k}(a) \|}{\| DT^{p_k}_{a^{
\prime
}}(z_{t_k}(a^{\prime})) u_{t_k}(a^{\prime}) \| }}, \ \ \ \
S^{\prime \prime}_k =
\log{\frac{\| DT^{t_{k+1}-p_k}_{a}(z_{t_k+p_k}(a)) u_{t_k+p_k}(a) \|}{\| DT
^{t_{k+1}- p_k}_{a^{\prime}}(z_{t_k+p_k}(a^{\prime})) u_{t_k+p_k}(a^{\prime}) 
\| }}.
$$
 
First we prove that $\sum_{k<q} S_k^{\prime \prime} < K$. 
Since $\gamma_j \cap {\cal C}^{(0)} = \emptyset$ for 
$t_k + p_k \leq j \leq t_{k+1}$,
it is straightforward to see using Sublemma \ref{nsublemmaB.6} that
$$
S^{\prime \prime}_k <
\frac{K}{\delta} \sum_{j=t_k+p_k}^{t_{k+1}} (|z_j(a) - z_j(a^{\prime})| 
+ | a - a^{\prime} |).
$$
The effect of $|a - a^{\prime}|$ can be ignored
since $|a - a^{\prime}| < e^{-cn}$. 
By Lemma \ref{nlemma6.7}, the slopes of $\gamma_j$ are uniformly bounded 
and the length of $\gamma_j$ grows exponentially, so
$$
\sum_{j=t_k+p_k}^{t_{k+1}} |z_j(a) - z_j(a^{\prime})|   <
K |\gamma_{t_{k+1}}|.
$$
Again by Lemma \ref{nlemma6.9}(b),
$|\gamma_{t_{k+1}}| > K |\gamma_{t_k}|$. 
Therefore $\sum_{k<q} S_k^{\prime \prime} < K$.

To estimate  $\sum_{k<q} S_k^{\prime}$ we apply Sublemma \ref{nsublemmaB.7}.
The effect of 
the term $|a - a^{\prime}|$
can again be ignored, so that
$$
\sum_{k<q} S_k^{\prime} \leq K \sum_{k=1}^{q-1} \frac{\gamma_{t_k}}{e^{-\mu_k}}
$$
where $\gamma_{t_k} \in I_{\mu_k j_k}$. To estimate this sum, let
$m(\mu) = max\{ t_k: \mu_k = \mu \}$ for each $\mu$. Using the fact that
$\mid \gamma_{t_{k+1}} \mid \geq K|\gamma_{t_k}|$, we conclude
that 
$$
\sum_{k<q} S_k^{\prime}
< K \sum_{k < q} \frac{\mid \gamma_{t_k} \mid }{e^{-\mu_k}}
< K \sum_{\mu} \frac{\mid \gamma_{m(\mu)} \mid }{e^{-\mu}}
< K \sum_{\mu} \frac{1}{\mu^2}.
$$
This completes the proof. \hfill $\square$

\end{document}